\documentclass[12pt]{article}

\usepackage[english]{babel}
\usepackage{amsmath}
\usepackage{amssymb}
\usepackage{amsthm}
\usepackage{thmtools}
\usepackage{amsfonts}
\usepackage{float}
\usepackage{aliascnt}
\usepackage{fancyhdr} 
\usepackage{bm}
\usepackage{cleveref}
\usepackage{verbatim}
\crefname{enumi}{}{}
\usepackage{graphicx}
\usepackage{xcolor}
\usepackage{enumitem}
\usepackage{changepage} 
\usepackage{bbm}
\usepackage{tabularx,colortbl}
\usepackage{stmaryrd}
\usepackage[refpage]{nomencl}
\usepackage[margin=1in]{geometry}
\usepackage[utf8]{inputenc}
\usepackage{centernot}
\usepackage{marvosym} 
\usepackage{stmaryrd}

\makeatletter
\def\thmt@refnamewithcomma #1#2#3,#4,#5\@nil{%
  \@xa\def\csname\thmt@envname #1utorefname\endcsname{#3}%
  \ifcsname #2refname\endcsname
    \csname #2refname\expandafter\endcsname\expandafter{\thmt@envname}{#3}{#4}%
  \fi
}

\makeatother

\def\independenT#1#2{\mathrel{\rlap{$#1#2$}\mkern2mu{#1#2}}}
\usepackage{etoolbox}
\author{Martin Emil Jakobsen}
\title{TITLE}

\newcommand\independent{\protect\mathpalette{\protect\independenT}{\perp}}
\def\independenT#1#2{\mathrel{\rlap{$#1#2$}\mkern2mu{#1#2}}}
\renewcommand{\t}{\intercal}
\newcommand{\my}{\mu}
\newcommand{\ep}{\varepsilon}
\renewcommand{\epsilon}{\varepsilon}

		\renewcommand{\i}{\infty}
		
		\renewcommand{\phi}{\varphi}
\newcommand{\dcov}{\mathrm{dCov}}

\newcommand{\convp}{\stackrel{P}{\longrightarrow}}
\newcommand{\convas}{\stackrel{a.s.}{\longrightarrow}}
\newcommand{\convd}{\stackrel{\mathcal{D}}{\longrightarrow}}

\newcommand{\eqd}{\stackrel{\mathcal{D}}{=}}
\newcommand{\inner}[2]{\langle #1 , #2 \rangle}
\newcommand{\sigalg}{$\sigma$-algebra }

\newcommand{\pspace}{$(\Omega,\mathcal{F},P)$ }

\newcommand{\define}[4]{\expandafter#1\csname#3#4\endcsname{#2{#4}}}
\newcommand{\R}{\mathbb{R}}
\newcommand{\F}{\mathbb{F}}
\newcommand{\N}{\mathbb{N}}
\forcsvlist{\define{\newcommand}{\mathcal}{c}}{A,B,C,D,E,F,G,H,I,J,K,L,P,M,N,U,S,T,W,V,X,Y,Z}
\forcsvlist{\define{\newcommand}{\mathbb}{b}}{A,B,C,D,E,F,G,H,I,J,K,L,P,Q,M,N,Z,T,X,Y}

\newcommand{\lp}{\left( }
\newcommand{\rp}{\right) }
\newcommand{\lb}{\left\lbrace}
\newcommand{\rb}{\right\rbrace}

\newcommand{\ra}{\rangle}
\newcommand{\la}{\langle}

\newcommand{\lv}{\left\vert}
\newcommand{\rv}{\right\vert}

\declaretheoremstyle[
  spaceabove=\topsep, 
  spacebelow=\topsep,
  headfont=\normalfont\bfseries,
  notefont=\bfseries, notebraces={(}{)},
  bodyfont=\normalfont\itshape,
 postheadspace=\newline,
 qed=$\circ$
]{bombom}
\declaretheoremstyle[
spaceabove=\topsep, 
spacebelow=\topsep,
headfont=\normalfont\bfseries,
notefont=\normalfont\bfseries, notebraces={- }{},
bodyfont=\normalfont,
postheadspace=\newline,
qed=$\circ$
]{exstyle}
\declaretheorem[style=bombom,title=Theorem,numberwithin=section,refname={theorem}]{theorem}
\declaretheorem[style=bombom,name=Assumption,sibling=theorem,refname={assumption}]{assumption}
\declaretheorem[style=bombom,name=Lemma,sibling=theorem,refname={lemma}]{lemma}
\declaretheorem[style=bombom,name=Definition,sibling=theorem,refname={definition}]{definition}
\declaretheorem[style=bombom,name=Corollary,sibling=theorem,refname={corollary}]{corollary}
\declaretheorem[style=exstyle,name=Problem,numbered=no,refname={problem}]{problem}
\declaretheorem[style=exstyle,name=Remark,sibling=theorem,refname={remark}]{remark}

\declaretheoremstyle[
  spaceabove=\topsep, 
  spacebelow=30pt,
  headfont=\normalfont\bfseries,
  notefont=\mdseries, notebraces={(}{)},
  bodyfont=\normalfont,
 postheadspace=\newline,
 numbered=no,
  qed=\qedsymbol,
  name=Proof
]{mythmstyle}
\declaretheorem[style=mythmstyle]{p}
\makenomenclature

\begin{document}
\bibliographystyle{alpha}
\newpage
	\begin{titlepage}
	\centering
	{\scshape\LARGE University of Copenhagen\par}
	\vspace{1cm}
	{\scshape\Large Master's Thesis in Actuarial Mathematics\par}
	\vspace{1.5cm}
	{\huge\bfseries Distance Covariance in Metric Spaces\\} {\LARGE Non-Parametric Independence Testing in Metric Spaces\par}
	\vspace{2cm}
	{\Large Martin Emil Jakobsen\par}
	\vspace{3cm}
	\large Supervised by \par
	Professor Thomas Valentin Mikosch \\
	\vfill
	
	{\small Thesis for the Master Degree in Actuarial Mathematics.\\
	\vspace{-0.1cm}
	 Department of Mathematical Sciences, University of Copenhagen \\}
 	\vspace{0.3cm}
	 {\small Speciale for cand.act graden i forsikringsmatematik. \\
 	\vspace{-0.2cm}
 Institut for matematiske fag, Københavns Universitet}
	\vspace{1cm}
	
	{\large January 24, 2017\par}
\end{titlepage}
\newpage\null\thispagestyle{empty}\newpage
\clearpage
\thispagestyle{empty}
\vspace*{\fill}
\begin{center}
	{\Large\textbf{Acknowledgements} } \\
	\begin{minipage}{.7\textwidth}
		\vspace{0.4cm}
		I would like to thank my thesis supervisor Professor Thomas Valentin Mikosch, for great guidance and rewarding discussions throughout the writing of this thesis. Furthermore, I would also like to extend my sincere gratitude to Mads Bonde Raad for the countless and fruitful  discussions about various mathematical problems and concepts. I would also like to thank Russell Lyons for taking the time to both confirm problems, and in the case of \cref{lemma_strong_negative_type_beta_tensor_injective} (3.8 in \cite{lyons2013distance}) providing a smart workaround idea that yielded the new and correct proof. Finally, I would like to thank my parents for always being  supportive during my studies.
	\end{minipage}
\end{center}
\vfill 
\clearpage
\newpage

\clearpage
\thispagestyle{empty}
	\vspace*{\fill}
	\begin{center}
		{\Large\textbf{Abstract} }\\
	\begin{minipage}{.7\textwidth}
		\vspace{0.4cm}
		The aim of this thesis is to find a solution to the non-parametric independence problem in separable metric spaces. Suppose we are given finite collection of samples from an i.i.d. sequence of paired random elements, where each marginal has values in some separable metric space. The non-parametric independence problem raises the question on how one can use these samples to reasonably draw inference on whether the marginal random elements are independent or not. We will try to answer this question by utilizing the so-called distance covariance functional in metric spaces developed by Russell Lyons. We show that, if the marginal spaces are so-called metric spaces of strong negative type (e.g.\, seperable Hilbert spaces), then the distance covariance functional becomes a direct indicator of independence. That is, one can directly determine whether the marginals are independent or not based solely on the value of this functional. As the functional formally takes the simultaneous distribution as argument, its value is not known in the posed non-parametric independence problem. Hence, we construct estimators of the distance covariance functional, and show that they exhibit asymptotic properties which can be used to construct asymptotically consistent statistical tests of independence. Finally, as the rejection thresholds of these statistical tests are non-traceable we argue that they can be reasonably bootstrapped. 
	\end{minipage}
	\end{center}
	\vfill 
	\clearpage
\newpage
\newpage

\clearpage
\thispagestyle{empty}
\vspace*{\fill}
\begin{center}
	{\Large\textbf{Resumé} }\\
	\begin{minipage}{.7\textwidth}
		\vspace{0.4cm}
		Det primære formål med dette speciale er at finde en løsning til det såkaldte ikke-parametriske uafhængighedsproblem i separable metriske rum. Antag, at vi er givet en endelig samling af stikprøver fra en uafhængig og identisk fordelt følge af parvise stokastiske elementer med marginaler, der antager værdier i et separabelt metrisk rum. Det ikke-parametriske uafhængighedsproblem stiller nu spørgsmålet om, hvordan disse stikprøver kan bruges, på fornuftig vis, til at drage inferens omkring, hvorvidt de marginale stokastiske elementer er uafhængige eller ej. Vi vil besvare dette spørgsmål i en tilfredsstillende grad ved at anvende det såkaldte distance covariance funktionale udviklet af Russell Lyons. Dette gøres ved at vise, at hvis de marginale metriske rum er af såkaldt stærk negativ type (f.eks. separable Hilbert rum), så er distance covariance funktionalet en såkaldt direkte uafhængigheds indikator. Dette betyder, at vi direkte kan bestemme om marginalerne er uafhængige ved at aflæse værdien af dette funktionale. Da distance covaraince funktionalet formelt tager den simultane fordeling som argument, kan vi i den givne problemstilling ikke aflæse værdien af funktionalet. Derfor konstruerer vi estimatorer for distance covariance funktionalet og viser at de besidder asymptotiske egenskaber, der muliggør konstruktionen af asymptotisk konsistente statistiske tests for uafhængighed. Da disse tests har forkastelses-niveauer, der ikke direkte kan identificeres, redegøres der for, at man på fornuftig vis kan bootstrappe dem i stedet for. 
	\end{minipage}
\end{center}
\vfill 
\clearpage
\newpage\null\thispagestyle{empty}\newpage

\clearpage
\thispagestyle{empty}
\tableofcontents
\newpage
\setcounter{page}{1}
\section{Introduction}
Consider the following set-up applicable throughout the thesis. Let  $(Z_n)_{n\in \N}=((X_n,Y_n))_{n\in \N}$ be an independent and identically distributed sequence of paired random elements, defined on a probability space $(\Omega,\F,P)$. It is assumed that, each pair of random elements $Z_n=(X_n,Y_n)$ takes values in some product space $\cX\times \cY$. Furthermore, throughout the thesis we let $\theta$ denote the simultaneous  distribution and let $\mu$ and $\nu$ denote the marginal distributions on $\cX$ and $\cY$ respectively. That is, 
\begin{align*}
(X_n,Y_n)\sim \theta, \quad \quad X_n\sim \mu:=\pi_1(\theta), \quad \quad Y_n\sim \nu := \pi_2(\theta),
\end{align*}
where $\pi_1$ and $\pi_2$ are the coordinate projections onto the marginal spaces $\cX$ and $\cY$ respectively (see \cref{Appendix_Product_spaces_section} for further details on product spaces). The purpose of this thesis is to answer the following problem in a set-up as general as possible. 
\begin{problem}[The Non-Parametric Independence Problem]
	Suppose that we are given a finite collection of paired sample points $z_{1,n}=[(x_i,y_i)]_{1\leq i \leq n}$, where each pair $(x_i,y_i)$ is a realization of $(X_i,Y_i)$. Given this collection of samples, how can we without restricting $\theta$ to a specific parametric class of distributions, draw inference on whether to reject the null-hypothesis of independence
	\begin{align*}
	H_0: \theta = \mu \times \nu,
	\end{align*}
	in favor of the alternative hypothesis of dependence 
	\begin{align*}
	H_1:\theta \not = \mu\times \nu.
	\end{align*}
\end{problem}
A solution to the above problem was proposed by Gábor J. Székely, Maria L. Rizzo and Nail K. Bakirov,  in the widely cited article "Measuring and Testing Dependence by Correlation of Distances" from 2007, published in The Annals of Statistics; \cite{szekely2007measuring}. In this article, a solution  to the above problem is proposed, in the case that both $\cX$ and $\cY$ are finite-dimensional Euclidean spaces. This is done by introducing the so-called distance covariance measure between two random vectors $X\in \R^n$ and $Y\in \R^m$, with simultaneous distribution $\theta$ on $\R^n\times \R^m$. This distance covariance measure, is given by
\begin{align*}
\dcov(X,Y) = \sqrt{\frac{1}{c_nc_m}\int_{\R^{n+m}} \frac{|\hat{\theta}(t,s)-\hat{\mu}(t)\hat{\nu}(s)|^2}{\|t\|_{\R^n}^{n+1}\|s\|^{m+1}_{\R^m}} \, d\lambda^{n}\times \lambda^m(t,s)},
\end{align*}
a weighted $L^2$ difference between the characteristic functions of $\theta$ and $\mu\times \nu$. The distance covariance measure $\dcov(X,Y)$ is easily seen to be zero if and only if $X\independent Y$. They furthermore introduce a plug-in estimator of this distance covariance measure, based on empirical characteristic functions. Hereafter they showed, that the estimator possesses asymptotic properties that allow the construction of an asymptotically consistent test of independence. \\ \\ \\

In 2013, Russell Lyons published the article "Distance covariance in metric spaces" in The Annals of Probability; \cite{lyons2013distance}. This article proposes a solution to the above problem, under the weaker assumption that the marginal spaces $(\cX,d_{\cX})$ and $(\cY,d_\cY)$ are so-called metric spaces of strong negative type. This is done by introducing another so-called distance covariance measure (a generalization of $\dcov$; see \cref{theorem_equivalence_distance_covariance_in_metric_spaces_and_Euclidean}) between the random Borel elements $X\in \cX$ and $Y\in \cY$ with simultaneous distribution $\theta$. This distance covariance measure, is equivalently (see  \cref{eq.temp.987}) given by
\begin{align*}
dcov(X,Y):=& \,dcov(\theta)\\
=& \,E\Big[ \Big( d_\cX(X_1,X_2) -E(d_\cX(X_1,X_2)|X_1) - E(d_\cX(X_1,X_2)|X_2)+ Ed_\cX(X_1,X_2) \Big) \\
&\times \Big( d_\cY(\,Y_1\,,\,Y_2\,) -E(d_\cY(\,Y_1\,,\,Y_2\,)|Y_1) - E(d_\cY(\,Y_1\,,\,Y_2\,)|Y_2)+ Ed_\cY(\,Y_1\,,\,Y_2\,) \Big) \Big],
\end{align*}
where $(X_1,Y_1)$ and $(X_2,Y_2)$ are independent copies of $(X,Y)$. An important but non-trivial property of this distance covariance measure, is that it can be used as a direct indicator of independence. That is, $dcov(X,Y)=0$ if and only if $X\independent Y$, whenever $(\cX,d_{\cX})$ and $(\cY,d_\cY)$ are metric spaces of strong negative type (a superset of separable Hilbert spaces; see \cref{theorem_separable_hilbert_spaces_are_of_strong_negative_type}). Russell Lyons then introduces a plug-in estimator for the distance covariance measure and show that it possesses  asymptotic properties that can be used to construct an asymptotically consistent test of independence.\\ \\ 
In this thesis we will answer the non-parametric independence problem using the theory developed in \cite{lyons2013distance}. This thesis is therefore essentially best described as, a  very detailed exposition of discoveries made by Russell Lyons. The original article leaves a surprisingly large amount of details to the reader, therefore it has not been easy or without problems to make this thesis.  

Some of the mathematical concepts and constructions needed to understand and describe the theory of distance covariance in metric spaces, were at the beginning unknown to me. So in order to keep the thesis self-contained, appendices have been added to introduce these concepts in a degree which suffices for our needs.  

In writing this thesis I also stumbled upon several discrepancies in the original article, ranging from negligible to serious. Whenever the non-negligible discrepancies are met, I have explicitly added remarks  explaining the problems and how they are solved. I am grateful that Russell Lyons has taken the time to both confirm problems, and in the case of \cref{lemma_strong_negative_type_beta_tensor_injective} (lemma 3.8 in \cite{lyons2013distance}) providing a smart workaround idea that yielded the new and to some extend quite different proof.   \\ \\
\newpage
\noindent We will now provide a brief overview of the content of the following sections. \\ \\
\begin{tabularx}{1\linewidth}{>{\raggedleft}p{2cm}X}
	\textit{Section 2} & We construct the so-called  distance covariance measure $dcov$. This so-called measure $dcov$, is formally a real-valued functional with domain given by a space of sufficiently nice Borel probability measures on the product space $\cX\times \cY$ of metric spaces. From the definition of $dcov$ it is easily realized that $\theta=\mu\times \nu$ implies $dcov(\theta)=0$. However, the converse implication which would render the distance covariance measure a direct indicator of independence, is not true for general metric spaces $\cX$ and $\cY$.  \vspace{0.2cm}\\
	
	\textit{Section 3} & To answer the question regarding which metric spaces would yield the converse implication mentioned above, we define metric spaces of negative and strong negative type. Metric spaces of negative type, are metric spaces that can be isometrically embedded into Hilbert spaces. If both marginal metric spaces $\cX$ and $\cY$ are of negative type,  then we  show that the functional $dcov$ has an alternative representation in terms of the isometric embeddings. This alternative representation leads us to the definition of metric space of strong negative type. The essential property of these spaces are, if both $\cX$ and $\cY$  are metric spaces of strong negative type, then
	$$
	dcov(\theta)=0\iff \theta=\mu\times \nu.
	$$
	It is furthermore shown that, when disregarding the unimportant singleton spaces,  it is necessary for the marginal metric spaces to be of strong negative type in order to have the implication $dcov(\theta)=0\implies \my\times \nu$. This section is concluded with a theorem identifying all separable Hilbert spaces as metric spaces of strong negative type. \vspace{0.2cm}\\
	
	\textit{Section 4} & This section is dedicated to proving some properties and bounds on the functional $dcov$. We also establish the connection between the distance covariance in Euclidean spaces from \cite{szekely2007measuring} and the distance covariance measure in metric spaces. That is, if $\cX$ and $\cY$ are finite-dimensional Euclidean spaces, then $\dcov(X,Y)^2=dcov(X,Y)$, proving that the distance covariance measure in metric spaces indeed is a generalization of the former.   \vspace{0.2cm}\\
	
	\textit{Section 5} & \Cref{section_asyp_test_main} is divided into three subsections. 
	In \cref{section_estimators}, we introduce two different estimators for $dcov(\theta)$. It is seen that, $dcov$ is a so-called regular functional, and one may recall that such functionals are the building blocks of the so-called $U$- and $V$-statistic estimators. Our choice of estimators for $dcov(\theta)$ are therefore given by such estimators.
	In \cref{section_asymp_prop_of_estimators}, we show that these estimators are both strongly consistent and if scaled correctly also possess rather complicated asymptotic distributions. In \cref{section_tests}, we formally describe the statistical models for which the assymptotic properties from \cref{section_asymp_prop_of_estimators} yield asymptotically consistent tests of independence. These tests turns out to have non-traceable rejection thresholds, so we end this last section by describing how one may reasonably bootstrap the rejection thresholds. \vspace{0.1cm}\\
\end{tabularx}
\newpage
\section{Distance covariance in metric spaces} \label{section_distance_covariance}
As mentioned in the introduction the main objective of this thesis is to establish a measure of dependence that can  be used to create an asymptotically consistent statistical test for independence. In this section we will construct the so-called distance covariance measure $dcov$ of a probability measure $\theta$ on a product space $\cX\times \cY$, which can be used to directly establish whether or not the probability measure $\theta$ is in fact given by the product of its marginals $\mu\times \nu$. That is, we will construct a functional 
\begin{align*}
dcov: M^{1,1}_1(\cX\times \cY) \to \R,
\end{align*}
with the desired property, that whenever the marginal spaces $\cX$ and $\cY$ are sufficiently nice, $dcov(\theta)=0$ if and only if $\theta=\mu\times \nu$. Here $M^{1,1}_1(\cX\times \cY)$ is the space of all Borel probability measures on $\cX\times \cY$ with sufficient integrability. Exactly what this sufficient integrability entails, is the content of the first definition below. 

Note that $dcov$ is not a measure in the usual sense, nevertheless we will still refer to it as the distance covariance measure rather than functional.  Before proceeding, we make an initial restriction on what kind of marginal spaces $\cX$ and $\cY$ we will consider. This restriction is the content of the following universal assumption of this thesis:
\begin{assumption}
	Every metric space $(\cX,d_\cX)$ and $(\cY,d_\cY)$ considered in this thesis is assumed separable.
\end{assumption}
As we shall see later, this restriction on the marginal spaces is not sufficient for the distance covariance measure to have the desired property. This is indeed solved by assuming that the marginal metric spaces are of strong negative type, which is the focus of attention in \cref{section_negative_and_strong_negative_type}. Before continuing we present a short remark on the above assumption.
\begin{remark}
	In the article of Russell Lyons \cite{lyons2013distance}, it is nowhere stated that we move beyond the realm of general metric spaces. This is an obvious error in the article as one  has to require as a minimum, that the metric spaces considered have cardinality less than or equal to the continuum. 
	
	The reason for this, is that in order to define the distance covariance, we need that the metrics on our marginal spaces are jointly measurable, i.e. $d_\cX:\cX\times \cX\to \R$ needs to be $\cB(\cX)\otimes \cB(\cX)/\cB(\R)$-measurable. Due to Nedomas pathology (see prop. 21.8 \cite{schechter1996HandbookOfAnalysisAndItsFoundations} or example 6.4.3 \cite{MeasureTheory2007bogachev2}) we get that, every metric $d$ on space with cardinality strictly greater than the continuum $\mathfrak{c}$, is not jointly measurable (the diagonal is not measurable). An example of such a space could be $\cX=\{f|f:\R\to \R\}$ endowed with the discrete metric, since $\mathrm{card}(\{f|f:\R\to \R\})> \mathfrak{c}$.
	
	This problem is of course eliminated by the assumption of separability of $(\cX,d_\cX)$, which implies that $\mathrm{card}(\cX)\leq \mathfrak{c}$ but also implies that $\cB(\cX\times \cX)=\cB(\cX)\otimes \cB(\cX)$ (see \cref{theorem_product_sig_alg_is_tensor_when_separable}), rendering $d_\cX$ jointly measurable, since it is continuous.  	There are indeed other places in this thesis, that utilize the separability of the considered metric spaces. Some of these are \cref{lemma_tensor_product_map_of_isometric_embeddings_is_pettis_integrable} which uses  that $\cB(\cX\times \cY)=\cB(\cX)\otimes \cB(\cY)$, furthermore in \cref{theorem_dcov(XX)=0_iff_X_degenerate} and \cref{theorem_bar_h_is_degenerate_and_identity} where we explicitly use the separability.

	In personal communication with Russell Lyons he acknowledges the problems, and agrees with me that this discrepancy is best solved by only considering separable metric spaces. 
\end{remark}
Throughout the thesis we will  use a variety of Borel measureas on our metric spaces, so we start by defining some commonly used spaces of measures. \\ \\ 
 In order to do so we need to define moments of measures. \textit{For any $k>0$ we say that a finite signed Borel measure on $(\cX,d_\cX)$ has finite $k$'th moment if }
 \begin{align*}
 \int d_\cX(x,o)^k \, d|\mu|(x) < \i,
 \end{align*}
 for some $o\in \cX$, where $|\mu|=\mu^+ + \mu^-$ is the total variation of $\mu$ and $\mu=\mu^+-\mu^-$ is the Jordan-Hahn decomposition. We may also note that, if the above holds for some $o \in \cX$, then it holds for all $o\in \cX$. In order to see this, note that, if there is an $o_1\in \cX$ such that $\int d_\cX(x,o_1)^k \, d|\mu|(x)<\i$, then for any $o_2\in \cX$ the $c_r$-inequality allows for the following finite bound
 \begin{align*}
 \int d_\cX(x,o_2)^k \, d|\mu|(x) &\leq \int (d_\cX(x,o_1)+d_\cX(o_1,o_2))^k \, d|\mu|(x) \\
 &\leq c_k \int d_\cX(x,o_1)^k d|\mu|(x) + c_k \int d_\cX(o_1,o_2)^k d|\mu|(x) \\
 &= c_k \int d_\cX(x,o_1)^k d|\mu|(x) + c_k  d_\cX(o_1,o_2)^k |\mu|(\cX) \\
 &<\i,
 \end{align*}
where $c_k=1$ when $k\leq 1$ and $c_k=2^{k-1}$ when $k\geq 1$. We may also note that in the case that $\cX=\R^n$  and $\mu$ is a probability measure on $\R^n$, the above definition of moments coincides with the regular definition of moments of random vectors. That is, if $X\sim \mu$ then choose $o=0$ and note that
 \begin{align*}
 \int_{\R^n} d_{\cX}(x,0)^k \, d|\mu|(x) = \int_{\R^n} \|x\|^k \, d\mu(x) = E\|X\|^k.
 \end{align*}
\begin{definition}
	Let $(\cX,d_\cX)$ and $(\cY,d_\cY)$ be two metric space. We define the following spaces\\\\
\begin{tabularx}{1\linewidth}{>{\raggedleft}p{2.4cm}X}
	 $M(\cX)$ & The space of all finite signed measures on $(\cX,\cB(\cX))$. \vspace{0.1cm}\\
	 $M^k(\cX)$ & The space of all finite signed measures on $(\cX,\cB(\cX))$ with finite $k$'th moment. \vspace{0.1cm}\\
	 $M_1(\cX)$ & The space of all probability measures  on $(\cX,\cB(\cX))$. \vspace{0.1cm}\\
	 $M_0(\cX)$ & The space of all finite signed measures on $(\cX,\cB(\cX))$ that assigns the entire space $\cX$ to zero. That is, $\mu(\cX)=0$ for all $\mu\in M_0(\cX)$. \vspace{0.1cm}\\
	 $M(\cX\times \cY)$ & The space of all finite signed measures on $(\cX\times \cY,\cB(\cX)\otimes \cB(\cY))$. \vspace{0.1cm} \\
	 $M^{k,k}(\cX\times \cY)$ & The space of all finite signed measures  $\theta$ on $(\cX\times \cY,\cB(\cX)\otimes \cB(\cY))$ for which it holds that $\pi_1(|\theta|)\in \cM^k(\cX)$ and $\pi_2(|\theta|)\in M^k(\cY).$ \vspace{0.1cm} \\
	 $M_1(\cX\times \cY)$ & The space of all probability measures on $(\cX\times \cY,\cB(\cX)\otimes \cB(\cY))$. \vspace{0.0cm} \vspace{0.1cm} \\
	 $M_{1,nd}(\cX\times \cY)$ & The space of all probability measures on $(\cX\times \cY,\cB(\cX)\otimes \cB(\cY))$ that has non-degenerate marginal distributions. That is, the marginal distributions are not concentrated on a singleton or equivalently not Dirac measures. \vspace{0.0cm} \vspace{0.0cm} \\ \\
\end{tabularx}
Whenever we put both a subscript and superscript it denotes the intersection. For example, $M^1_1(\cX)=M^1(\cX)\cap M_1(\cX)$ is the space of probability measures with finite 1st moment.
\end{definition}
 We may furthermore note that $M_1(\cX)\subset M_1^1(\cX)\subset M(\cX)$. Whenever we consider a measure $\theta\in M^1(\cX\times \cY)$ then we indirectly assume that $\mu$ is the marginal measure on $(\cX,\cB(\cX))$ and $\nu$ is the marginal measure on $(\cY,\cB(\cY))$, i.e. $\mu=\pi_1(\theta)$ and $\nu=\pi_2(\theta)$. Also note that if $\theta\in M^{1,1}_1(\cX\times \cY)$ then $\mu\in M^1_1(\cX)$ and $\nu\in M^1_1(\cY)$, since $|\theta|=\theta$. \\ \\
 It turns out that some of these spaces are in fact $\R$-vector spaces, if we define some sensible addition and multiplication on them. This fact is not important for the definition of $dcov$, but it will later play a very important part in the further analysis of the  distance covariance measure.
\begin{lemma} \label{lemma_M(X)_and_M^1(X)_are_vector_spaces} \label{lemma_M(XtimesY)_and_M11(XtimesY)_are_vector_spaces}
	For any metric space $(\cX,d_\cX)$, it holds that $M(\cX)$ is an $\R$-vector space, and $M^1(\cX) $  is a  linear subspace of $M(\cX)$. If furthermore $(\cY,d_\cY)$ is yet another metric space, then $M(\cX\times \cY)$ is a $\R$-vector space and $M^{1,1}(\cX\times \cY)$ is a linear subspace of $M(\cX\times \cY)$.
\end{lemma}
\begin{p}
	On $M(\cX)$ - the space of finite signed Borel measures - we define scalar multiplication and addition of these measures by
	\begin{align*}
	(\mu_1+\mu_2)(A) = \mu_1(A)+\mu_2(A) \quad \quad \text{ and } \quad \quad (a\mu_1)(A)=a \cdot \mu_1(A).
	\end{align*}
	for any $\mu_1,\mu_2\in M(\cX)$, $A\in \cB(\cX)$ and $a_1,a_2\in \R$. 
	It is obvious that $M(\cX)$ is closed under any finite linear combination and satisfies every other axiom of vector spaces, meaning that $M(\cX)$ is a vector space. The question is now, if the subset of signed measures with finite first moments $M^1(\cX)\subset M(\cX)$, indeed is a linear subspace. Since the zero measure (maps every measurable set to zero) clearly has a finite first moment (rendering $M^1(\cX)$ non-empty), we note that it suffices to show that $ a_1\mu_1+a_1\mu_2$	is a measure with finite first moment for all $\mu_1,\mu_2\in M^1(\cX)$ and $a_1,a_2\in \R$, in order to prove that $M^1(\cX)$ is a linear subspace of $M(\cX)$. Hence we see that
	\begin{align*}
	\int d_\cX(x,o) \, d|a_1\mu_1+a_2\mu_2|(x) &\leq  \int d_\cX(x,o) \, d(|a_1||\mu_1|+|a_2||\mu_2|)(x)  \\
	&=|a_1|\int d_\cX(x,o) \, d|\mu_1|(x) +|a_2|\int d_\cX(x,o) \, d|\mu_2|(x)  \\
	&<\i,
	\end{align*}
	where we used that $|a_1\mu_1+a_2\mu_2|\leq |a_1||\mu_1|+|a_2||\mu_2|$ and that the Lebesgue integral is monotone in measure, when the integrand is non-negative. \\ \\	
	The inequality $|a_1\mu_1+a_2\mu_2|\leq |a_1||\mu_1|+|a_2||\mu_2|$ follows by standard arguments, but in order to keep the thesis self-contained we show it regardless.
	First note that if $\lambda_1$ and $\lambda_2$ are two positive measures and $\nu = \lambda_1-\lambda_2$ , then $\lambda_1 \geq \nu^+$ and $\lambda_2 \geq \nu^-$ (see p. 88 \cite{folland1999real}). In our setup, we have that
	$
	a_1\mu_1+a_2\mu_2 = ((a_1\mu_1)^+-(a_1\mu_1)^-) +((a_2\mu_2)^+-(a_2\mu_2)^-) = ((a_1\mu_1)^+ +(a_2\mu_2)^+)- ((a_1\mu_1)^-+(a_2\mu_2)^-),
$ 
	and as a consequence $	|	a_1\mu_1+a_2\mu_2| = (	a_1\mu_1+a_2\mu_2)^+ + (	a_1\mu_1+a_2\mu_2)^- \leq 	((a_1\mu_1)^+ +(a_2\mu_2)^+)+((a_1\mu_1)^-+(a_2\mu_2)^-) 
	= ((a_1\mu_1)^+ + (a_1\mu_1)^-) + ((a_2\mu_2)^++(a_2\mu_2)^-) = |a_1\mu_1| + |a_2\mu_2|$.
	Now note that the total variation measure $|a_1\mu_1|=(a_1\mu_1)^++(a_1\mu_1)^-$  equivalently  can be stated as the expression $	|a_1\mu_1|(A)= \sup \sum_{i=1}^\i |a_1\mu_1(A_i)| = |a_1| \sup \sum_{i=1}^\i |\mu_1(A_i)| = |a_1||\mu_1|(A)
$
	for any $A\in \cB(\cX)$, where the supremum is over all mutually disjoint sequences $(A_i)$ in $\cB(\cX)$ with $\cup A_i=A$ (see section A.1 \cite{sokol2014introduction} or p. 177 \cite{MeasureTheory2007bogachev}), proving the wanted inequality. \\ \\
	The proof for $M(\cX\times \cY)$ and $M^{1,1}(\cX\times \cY)$ follows by analogous arguments. That is, for $\theta_1,\theta_2\in M^{1,1}(\cX\times \cY)$ and $a,b\in \R$ we have that $a\theta_1 + b\theta_2\in M(\cX\times \cY)$ and that $|a\theta_1 + b\theta_2| \leq |a||\theta_1|+|b||\theta_2|$. Hence\begin{align*}
	\pi_1(|a\theta_1 + b\theta_2|) \leq \pi_1 (|a||\theta_1|+|b|\theta_2|) = |a|\pi_1(|\theta_1|)+|b|\pi_1(|\theta_2|) \in M^1(\cX)
	\end{align*}
	since $M^1(\cX)$ is a $\R$-vector space and $\pi_1(|\theta_1|),\pi_1(|\theta_2|)\in M^1(\cX)$ by definition of $M^{1,1}(\cX\times \cY)$. A similar derivation follows for the $\pi_2$ projection, so $a\theta_1+b\theta_2\in M^{1,1}(\cX\times \cY)$.
	
\end{p}
Now that we have defined the important spaces of measures we are almost ready to define the distance covariance measure. The distance covariance measure is defined in terms of integrals of certain mappings, hence we start by proving that these are sufficiently integrable. Before doing this, we need to establish some common ground, on how to define the integral of a mapping with respect to the product of two finite signed measures.\\ \\
Recall the integral of a measurable mapping $f$ with respect to a signed measure $\mu$ is defined as
\begin{align*}
\int f \, d\mu = \int f\, d\mu^+ - \int f\, d\mu^-,
\end{align*}
whenever $f\in \mathcal{L}(\mu):=\mathcal{L}^1(\mu^+)\cap \mathcal{L}^1(\mu^-)=\mathcal{L}^1(|\mu|)$, where  $\mu=\mu^+ -\mu^-$ is the Jordan-Hahn decomposition and $|\mu|=\mu^++\mu^-$ is the total variation of $\mu$.
We remind the reader that one constructs the product measure of two signed measures $\mu\in M(\cX)$ and $\nu\in M(\cY)$ by utilizing the Jordan-Hahn decomposition theorem. 

If we let $\mu = \mu^+ - \mu^-\in M(\cX)$ and $\nu= \nu^+ - \nu^-\in M(\cY)$,  then we define the product measure $\mu\times \nu\in M(\cX\times \cY)$ directly by its Jordan-Hahn decomposition. That is, as the difference between the mutually singular measures $\mu^+ \times \nu^+ + \mu^- \times \nu^-$ and  $\mu^+ \times \nu^- +\mu^- \times \nu^+$. To see that they are mutually singular simply realize that the first measure is concentrated on $(\cX^+ \times \cY^+)\cup (\cX^- \times \cY^-)$ and the other is concentrated on $(\cX^+ \times \cY^-)\cup (\cX^- \times \cY^+)$, where $\cX=\cX^+ \cup \cX^-$ and $\cY=\cY^+\cup \cY^-$ are the disjoint decompositions  of the spaces given by the Jordan-Hahn decomposition of the marginal measures. Hence $\mu\times \nu \in M(\cX\times \cY)$ is the finite signed measure given by
\begin{align*}
\mu\times \nu = \mu^+ \times \nu^+ + \mu^- \times \nu^- - \mu^+ \times \nu^- - \mu^- \times \nu^+.
\end{align*}
and the above decomposition is equal to its Jordan-Hahn decomposition.
 We say that a $\cB(\cX)\otimes \cB(\cY)/\cB(\R)$-measurable mapping $f:\cX\times \cY \to \R$ is integrable with respect to $\mu\times \nu$, written $f\in \cL^1(\mu\times \nu)$ , if $f$ is integrable with respect to all the above product measures. The integral of $f$ with respect to $\mu\times \nu$ is defined in the natural way as the sum and difference of the four regular product integrals. In terms of integrability conditions, it suffices to check that $\int |f| \, d|\mu|\times |\nu| <\i$, since 
	\begin{align*}
	\int |f| \, d\mu^\pm \times \nu^\pm   \leq \int |f| \, d|\mu| \times |\nu|,
	\end{align*}
	which is seen by using Tonelli's theorem and then successively making an upper bound for the inner and then the outer integral, by changing integration measures to $|\mu_1|$ and $|\mu_2|$. That is, if $|f|\in \mathcal{L}^1(|\mu|\times |\nu|)$ then $f\in \mathcal{L}^1(\mu\times \nu)$. In the case of an integrable mapping $f\in \mathcal{L}^1(\mu\times \nu)$, we have the following Fubini theorem for the product integral of signed measures
	\begin{align*}
	\int f \, d\mu \times d\nu = \int\int f \, d\mu \, d\nu = \int\int  f \, d\nu \, d\mu,
	\end{align*}
	which is seen be utilizing Fubini's theorem for the marginal integrals with respect to $\mu_1^\pm \times \mu_2^\pm$. For further details on the integration with respect to the product of finite signed measures, we refer the reader to section 3.3 of \cite{MeasureTheory2007bogachev}.
	 
\begin{lemma} \label{lemma_metric_is_integral_product_of_M_with_finte_first_moment}
	If $(\cX,d_\cX)$ is a metric space, then $d_\cX\in \cL^1(\mu_1\times \mu_2)$ for any two measures  $\mu_1,\mu_2\in M^1(\cX)$.
\end{lemma}
\begin{p}
	By the above remark if suffices to show that $\int d_\cX \, d  |\mu_1| \times |\mu_2| <\i$. This is easily seen by using the triangle inequality: $d_\cX(x_1,x_2) \leq d_\cX(x_1,o)+d_\cX(o,x_2)$ for any $o\in\cX$. Thus for some $o\in \cX$, we get
	\begin{align*}
	\int d_\cX \, d  |\mu_1| \times |\mu_2| &\leq \int d_\cX(x_1,o) \, d  |\mu_1| \times |\mu_2|(x_1,x_2)+\int d_\cX(o,x_2) \, d  |\mu_1| \times |\mu_2|(x_1,x_2) \\
	&= \mu_2(\cX)\int d_\cX(x_1,o) \, d  |\mu_1| (x_1)+\mu_1(\cX)\int d_\cX(x_2,o) \, d  |\mu_2|(x_2) \\
	&<\i,
	\end{align*}
 where we used that $\mu_i$ is finite and $d_\cX(\cdot , o)\in \mathcal{L}^1(|\mu_i|)$ for both $i\in \{1,2\}$, since they are finite signed measures with finite first moment.
\end{p}
This lemma allows for the definition of the mappings relevant for the distance covariance measure.
\begin{definition}
	Let $(\cX,d_\cX)$ be a metric space. For any measure $\mu\in M^1(\cX)$ we may define the $M^1(\cX)$-integrable and continuous mapping $a_\mu:\cX \to \R$ given by
	\begin{align} \label{defi_a_mu}
	a_\mu(x) = \int d_\cX(x,y) d\mu(y),
	\end{align}
	and the mapping $D:M^1(\cX) \to \R$  by
	\begin{align}
	D(\mu)= \int a_\mu (x) \, d\mu(x)=\int d_\cX(x_1,x_2) \, d\mu \times \mu(x_1,x_2). \label{defi_D_mu}
	\end{align}
	Lastly for any $\mu\in M^1(\cX)$ we may define the $\mu$-modified "distance" $d_\mu:\cX \times \cX \to \R$ by 
	\begin{align}
	d_\mu(x_1,x_2)= d_\cX(x_1,x_2) -a_\mu(x_1) - a_\mu(x_2) +D(\mu). \label{defi_d_mu}
	\end{align}
\end{definition}
The quotation sign in "distance", signifies that it is not a real metric. To see the continuity of $a_\mu$ for all $\mu\in M^1(\cX)$, we note that $d_\cX:\cX^2 \to \R$ is continuous and that by applying the reverse triangle inequality $
\left| a_\mu(x) - a_\mu(z) \right|  \leq \int |d_\cX(x,y)-d_\cX(z,y)| \, d\mu(y) \leq \int d_\cX(x,z) \, d\mu(y) = d_\cX(x,z) \mu(\cX),
$ which tends to zero as $x\to z$ for any $z\in\cX$, proving continuity of $a_\mu$. \\ \\
We may note that the $\mu$-modified distance $d_\mu$ lies within $\mathcal{L}^1(\mu_1\times \mu_2)$ for any two $\mu_1,\mu_2\in M^1(\cX)$, but as shown below it actually possess stronger integrability than $d_\cX$.
\begin{lemma} \label{lemma_d_mu_is_squre_integrabel_with_respect_to_product_measure_in_M11}
For any metric space $(\cX,d_\cX)$	 and any $\mu,\mu_1,\mu_2 \in M^1_1(\cX)$, we have that $d_{\mu}\in \mathcal{L}^2(\mu_1 \times \mu_2)$
\end{lemma}
\begin{p}
	First of all $d_\mu$ is the composition of $\cB(\cX)\otimes \cB(\cX)/\cB(\R)$-measurable mappings, hence it is itself jointly measurable. 
	By the triangle inequality $d(x,y)\leq d(x,z)+d(z,y)$ we get that 
	\begin{align*}
	D(\mu)&= \int d_\cX(x_1,x_2) \, d\mu\times \mu(x_1,x_2) \\
	&\leq  \mu(\cX)\int d_\cX(x_1,x) \, d\mu(x_1) + \mu(\cX)\int d_\cX(x_2,x) \, d\mu(x_2) \\
	&= 2a_\mu(x),
	\end{align*}
	for any $x\in \cX$, by Fubini's theorem and the fact that $\mu(\cX)=1$ since it is a probability measure. This also implies that $D(\mu) \leq a_\mu(x) +a_\mu(y)$
	for all $x,y\in\cX$. We also have that
	\begin{align} \label{eq_showing_that_|d(xy)-a_mu(x)|<=a_mu(y)}
	d_\cX(x,y) \leq a_\mu(x) + a_\mu(y), \quad \quad \text{and} \quad \quad a_\mu(x) \leq d_\cX(x,y) + a_\mu(y),
	\end{align}
	which is seen by integrating on both sides of the triangle inequality with respect to $\mu$ of different arguments. Hence with $A=\{(x,y)\in \cX\times \cX: d_\cX(x,y)+D(\mu)\geq a_\mu(x)+a_\mu(y)\}$, all the above inequalities yield that
	\begin{align*}
	|d_\mu(x,y)| &= 1_A(x,y)(d_\cX(x,y)+D(\mu)- a_\mu(x)-a_\mu(y))\\
	&+1_{A^c}(x,y)( a_\mu(x)+a_\mu(y)-d_\cX(x,y)-D(\mu)) \\
	&\leq 1_A(x,y)D(\mu) +1_{A^c}(x,y)(2a_\mu(y)-D(\mu)) \\
	&\leq 2a_\mu(y),
	\end{align*}
	and $|d_\mu(x,y)|\leq 2a_\mu(x)$ (by symmetry) for all $x,y\in \cX$.  Thus
	\begin{align*}
	\int d_\mu(x,y)^2 \, d \mu_1\times \mu_2(x,y) &\leq 4\int a_\mu(x)a_\mu(y)   \, d\mu_1\times \mu_2(x,y) \\
	&= 4\int d_\cX(x,y) \, d\mu_1\times \mu(x,y)  \int d_{\cX}(x,y) \, d\mu\times \mu_2(x,y) <\i,
	\end{align*}
	by the above inequalities, Fubini's theorem and \cref{lemma_metric_is_integral_product_of_M_with_finte_first_moment}.
\end{p}
Having established the square integrability of the mapping $d_\mu:\cX\times \cX\to \R$ we define the distance covariance functional in the following way
\begin{definition}[Distance Covariance] \label{defi_Distance_covariance}
	For any two metric spaces $(\cX,d_\cX)$ and $(\cY,d_\cY)$, we define the distance covariance measure as  $dcov:M^{1,1}_1(\cX\times \cY)\to\R$ given by
	\begin{align*}
	dcov(\theta ) =& \int d_\mu(x_1,x_2) d_\nu(y_1,y_2) \, d\theta\times \theta((x_1,y_1),(x_2,y_2)),
	\end{align*}
	for any $\theta\in M^{1,1}_1(\cX\times \cY)$ with marginal probability measures $\mu=\pi_1(\theta)\in M_1^1 (\cX)$ and $\nu=\pi_2(\theta)\in M_1^1(\cY)$.
\end{definition}

We stress that $dcov(\theta)$ is indeed well-defined for any $\theta \in M_1^{1,1}(\cX\times \cY)$ by the Cauchy-Schwarz inequality. Simply note that
\begin{align*}
\int d_\mu(x_1,x_2)^2 \, d\theta \times \theta ((x_1,y_1),(x_2,y_2)) = \int d_\mu(x_1,x_2)^2 d\mu\times \mu(x_1,x_2) < \i,
\end{align*}
by Fubini's theorem and \cref{lemma_d_mu_is_squre_integrabel_with_respect_to_product_measure_in_M11}. By analougus aruguments for $d_\nu$, we conclude that the mappings $((x_1,y_1),(x_2,y_2))\mapsto d_\mu(y_1,y_2),d_\nu(y_1,y_2) \in \mathcal{L}^2((\cX\times \cY)^2,\theta \times \theta)$. Hence Cauchy-Schwarz inequality yields that 
$((x_1,y_1),(x_2,y_2))\mapsto d_\mu(x_1,x_2) d_\nu(y_1,y_2) \in \mathcal{L}^1((\cX\times \cY)^2,\theta \times \theta)$, proving that $dcov(\theta)$ is well-defined. \\ \\
It is immediately seen from the above definition, that if $\theta=\mu\times \nu\in M^{1,1}_1(\cX\times \cY)$ then
\begin{align*}
dcov(\theta)=& \int d_\mu(x_1,x_2) d_\nu(y_1,y_2) \, d(\mu\times \nu) \times (\mu\times \nu)((x_1,y_1),(x_2,y_2)) \\
=&\int d_\mu(x_1,x_2) \, d\mu\times \mu(x_1,x_2) \int d_\nu(y_1,y_2) \, d\nu\times \nu(y_1,y_2) \\
=& (D(\mu)-D(\mu)-D(\mu)+D(\mu))  (D(\nu)-D(\nu)-D(\nu)+D(\nu)) =0,
\end{align*}
 by Fubini's theorem. But it is not readily apparent what is needed in order to get the converse statement: if $dcov(\theta)=0$, then $\theta=\mu\times \nu$. As mentioned previously this implication does not hold for general metric spaces. The next section is dedicated to understand when it does.    \\ \\ 
 We may also derive a representation of $dcov$ in terms of conditional expectations of random variables. Let $(\cX,d_\cY)$ and $(\cY,d_\cY)$ be metric spaces and let $X\in \cX$ and $Y\in\cY$ be random Borel elements defined on a common probability space $(\Omega,\F,P)$ with distribution $P_X=\mu \in M^{1}_1(\cX)$ and $P_Y=\nu\in M^1_1(\cY)$ respectively, i.e. $P_{X,Y}=\theta\in M^{1,1}_1(\cX\times \cY)$. Now let $(X_1,Y_1)$ and $(X_2,Y_2)$ be independent copies of $(X,Y)$ and note that $ \theta\times \theta = ((X_1,Y_1),(X_2,Y_2))(P)$. Thus we may realize that
 \begin{align}
 dcov(X,Y) :=& dcov(\theta) \notag \\
 =& E[d_{P_X}(X_1,X_2)d_{P_Y}(Y_1,Y_2)]\notag \\
 =& E\Big[ \Big( d_\cX(X_1,X_2) -a_{P_X}(X_1) - a_{P_X}(X_2)+ D(P_X) \Big) \label{eq.temp.987}\\
 &\times \Big(  d_\cY(\, Y_1\, , \, Y_2\, ) -a_{P_Y}(\,Y_1\, ) - a_{P_Y}(Y_2)+ D(P_Y) \Big)\Big] \notag\\
=& \,E\Big[ \Big( d_\cX(X_1,X_2) -E(d_\cX(X_1,X_2)|X_1) - E(d_\cX(X_1,X_2)|X_2)+ Ed_\cX(X_1,X_2) \Big) \notag \\
&\times \Big( d_\cY(\,Y_1\,,\,Y_2\,) -E(d_\cY(\,Y_1\,,\,Y_2\,)|Y_1) - E(d_\cY(\,Y_1\,,\,Y_2\,)|Y_2)+ Ed_\cY(\,Y_1\,,\,Y_2\,) \Big) \Big], \notag
 \end{align}
In the case when $d_\cX$ and $d_\cY$ are sufficiently integrable, these linear combinations are actually orthogonal projections onto certain subspaces of $L^2(\Omega,\F,P)$. 
In the appendix on $U$- and $V$-statistics (\cref{Appendix_Hoeffding_decomposition_section}) we create formulas for the orthogonal projection onto certain spaces, which we use to prove the Hoeffding decomposition of $U$-statistics. For our purpose right here it suffices to consider the space $H_{\{1,2\}}$ consisting of all mappings in $L^2(\Omega,\F,P)$ that have the form $g(X_1,X_2)$ for some measurable map $g$ and for which it holds that
\begin{align*}
E(g(X_1,X_2)|X_1)=0, \quad E(g(X_1,X_2)|X_2)=0, \quad Eg(X_1,X_2)=0,
\end{align*}
with a similar space $H_{\{1,2\}}'$ constructed for the $Y_1,Y_2$ variables. Now if $\theta\in M^{2,2}_1(\cX\times \cY)$ we have that $d_\cX(X_1,X_2),d_\cY(Y_1,Y_2)\in L^2(\Omega,\F,P)$ and their orthogonal projection onto the spaces $H_{\{1,2\}}$ and $H_{\{1,2\}}'$  are given by
 \begin{align*}
 P_{H_{\{1,2\}}}(d_\cX(X_1,X_2))=&\sum_{B\subset \{1,2\}}(-1)^{2-|B|}E(d_\cX(X_1,X_2)|X_i:i\in B) \\
 =&d_\cX(X_1,X_2) -E(d_\cX(X_1,X_2)|X_1) - E(d_\cX(X_1,X_2)|X_2)+ Ed_\cX(X_1,X_2),
 \end{align*}
 and
 \begin{align*}
 P_{H_{\{1,2\}}'}(d_\cY(X_1,X_2)) =d_\cY(Y_1,Y_2) -E(d_\cY(Y_1,Y_2)|Y_1) - E(d_\cY(Y_1,Y_2)|Y_2)+ Ed_\cY(Y_1,Y_2),
 \end{align*}
 respectively, where we used the projection formula from \cref{lemma_projection_onto_H_A_spaces}. Thus we may say that, if $\theta\in M^{2,2}_1(\cX\times \cY)$, then $dcov(X,Y)$ is given as the expectation of $d_\cX(X_1,X_2)$ projected onto $H_{\{1,2\}}$ multiplied by $d_\cY(Y_1,Y_2)$ projected onto $H'_{\{1,2\}}$. This observation is not used in the next chapters, but we feel it was a connection worth mentioning. 
 
\section{Metric spaces of negative and strong negative type}\label{section_negative_and_strong_negative_type}
In this section we will examine what restriction on our marginal metric spaces $(\cX,d_\cX)$ and $(\cY,d_\cY)$ allows for the distance covariance measure to be used as a direct indicator of independence. That is, for which marginal metric spaces $(\cX,d_\cX)$ and $(\cY,d_\cY)$ are we allowed to conclude that
\begin{align*}
dcov(\theta)=0 \iff \theta= \mu\times \nu,
\end{align*}
for any $\theta\in M^{1,1}_1(\cX\times \cY)$ with marginals $\mu\in M^1_1(\cX)$ and $\nu\in M^1_1(\cY)$. The answer to this problem is: if we only consider spaces for which the independence problem is indeed valid, then it is necessary and sufficient to assume that both $(\cX,d_\cX)$ and $(\cY,d_\cY)$ are what is called metric spaces of strong negative type. \\ \\
The procedure for showing this is as follows. We define a certain subset of all metric spaces, called metric spaces of negative type. It turns out that whenever both marginal metric spaces $(\cX,d_\cX)$ and $(\cY,d_\cY)$ are of negative type, one may derive a Hilbert space representation of the distance covariance measure given by
\begin{align*}
dcov(\theta)= 4\| \beta_{\phi\otimes \psi}(\theta-\mu\times \nu)\|_{\cH_1\otimes \cH_2},
\end{align*}
for any $\theta\in M^{1,1}_1(\cX\times \cY)$, and some mapping $\beta_{\phi\otimes \psi}:M^{1,1}(\cX\times \cY)\to \cH_1\otimes \cH_2$ with values in the tensor product of two Hilbert spaces $\cH_1$ and $\cH_2$. This Hilbert space representation will also be essential in \cref{section_properties_of_dcov}, where we prove that the distance covariance measure in metric spaces (as we defined it) coincides with the distance covariance measure from \cite{szekely2007measuring}, when the marginal spaces are assumed to be finite-dimensional Euclidean spaces. 

If we restrict ourselves to an even smaller set of metric spaces, so-called metric spaces of strong negative type, then we get that the mapping $\beta_{\phi\otimes \psi}:M^{1,1}(\cX\times \cY)\to\cH_1\otimes \cH_2$ is an injective linear transformation. Since $\theta-\mu\times\nu \in M^{1,1}(\cX\times \cY)$, this of course entails that  $\beta_{\phi\otimes \psi}(\theta-\mu\times \nu)=0\iff \theta=\mu\times \nu$. Hence we get that, if both marginal  spaces are of strong negative type, then the distance covariance measure can be used as a direct indicator of independence. On the other hand, if we disregard the unimportant (from an independence perspective) singleton-spaces, then it is also necessary that both marginal spaces are  of strong negative type, in order for the distance covariance measure to be a direct indicator of independence. \\ \\
The above procedure is split into the next two subsections. In \cref{section_negative_type}, we define metric spaces of negative type and prove the alternative Hilbert space representation of the distance covariance measure. In \cref{section_strong_negative_type}, we define metric spaces of strong negative type and prove the above mentioned implication of injectivity, which yields the wanted property of the distance covariance measure.  In \cref{section_strong_negative_type}, we also identify all separable Hilbert spaces as metric spaces of strong negative type, and we furthermore prove that it is necessary for the marginals spaces to be of strong negative type in order for the distance covariance measure to be used as a direct indicator of independence. \newpage

\subsection{Metric spaces of negative type} \label{section_negative_type}
Metric spaces of negative type are defined in \cite{lyons2013distance} in the same way we do below, but without mention of negative definite kernels. The concept of negative definite kernels is presented in "Harmonic Analysis on Semigroups" by Christian Berg et al. \cite{berg1984harmonic} and here a connection between negative definite kernels and positive definite kernels (defined below) is established. This connection turns out to be  useful, since the concept of positive definite kernels is central in the theory of reproducing kernel Hilbert spaces. By utilizing this connection between negative definite kernels and reproducing kernel Hilbert spaces,  we can establish a more structured (yet longer) presentation of the properties of metric spaces of negative type, than the one presented in \cite{lyons2013distance}. 
\begin{definition}[Metric spaces of negative type]
	A metric space $(\cX,d_\cX)$ is said to be of negative type if the metric mapping $d_\cX:\cX\times \cX\to \R$ is a negative definite kernel. That is, if
	\begin{align*}
	\sum_{i=1}^n \sum_{j=1}^n \alpha_i \alpha_j d_\cX(x_i,x_j) \leq 0,
	\end{align*}
	for any $n\in\N$, $x\in \cX^n$ and $\alpha\in \R^n$ with $\sum_{i=1}^n \alpha_i =0$. 
\end{definition}

One of the main objectives of this section is to provide an equivalence between metric spaces of negative type and the isometric embeddability of $(\cX,d_\cX^{1/2})$ into a Hilbert space. This turns out to be a very helpful equivalence, because it is indeed these isometric embeddings that allow for the derivation of the alternative Hilbert space representation of the distance covariance measure.  \\ \\
But before we continue with proving the existence of previously mentioned isometric embeddings, we establish a connection between inner product spaces and negative definite kernels - a connection found in \cite{wells2012embeddings}. 
\begin{theorem}
	Let $\cX$ be a normed $\R$-vector space and let $d_\cX:\cX\times \cX\to \R$ denote the naturally induced metric. Then $\cX$ is an inner product space if and only if  $d_\cX^2:\cX\times \cX\to \R$ is a negative definite kernel. That is, if and only if the semi-metric space $(\cX,d_\cX^2)$ is of negative type.
\end{theorem}
\begin{p}
	First note that, if $d_\cX$ is a metric, then $d_\cX^2$ is in general not  a metric but rather a semi-metric, which by definition means that it satisfies all conditions of metrics except the triangle inequality. Let $\la \cdot,\cdot\ra :\cX\times \cX\to \R$ and $\|\cdot\|:\cX\to \R$ denote the inner product and its naturally induced norm on $\cX$. \\ \\
	If $\cX$ is an inner product space, we see that $\la x_i-x_j,x_i-x_j\ra = \|x_i\|^2+\|x_j\|^2-2\la x_i,x_j \ra$. Hence 
	\begin{align*}
	\sum_{i=1}^n \sum_{j=1}^n \alpha_i \alpha_j d_\cX(x_i,x_j)^2 &=\sum_{i=1}^n \sum_{j=1}^n \alpha_i \alpha_j \la x_i-x_j , x_i-x_j \ra \\
	&= \sum_{i=1}^n \alpha_i \|x_i\|^2 \sum_{j=1}^n \alpha_j + \sum_{j=1}^n \alpha_j \|x_j\|^2 \sum_{i=1}^n \alpha_i - 2 \Big\la \sum_{i=1}^n \alpha_i x_i,\sum_{j=1}^n \alpha_j,x_j \Big\ra \\
	&= -2\Big\| \sum_{i=1}^n \alpha_ix_i \Big\|^2 \leq 0,
	\end{align*}
	for all $x\in \cX^n$ and $\alpha\in \R^n$ with $\sum_{i=1}^n \alpha_i =0$, proving that $d_\cX^2$ is a negative definite kernel. Conversely assume that $d_\cX^2$ is a negative definite kernel. We know that $\cX$ is an inner product space  if and only if the norm satisfies the parallelogram law
	\begin{align*}
	\| x+y\|^2 + \|x-y\|^2 =2\|x\|^2 +2\|y\|^2,
	\end{align*}
	for all $x,y\in \cX$ (cf. theorem 6.9 \cite{applied_analysis_hunter2001}). Thus fix $x,y\in\cX$ and let $ a\in(0,1/2)$, $n=4,\,x_1=x,\,x_2=y,\,x_3=-y,\,x_4=0$ and $\alpha_1=1-2a,\,\alpha_2=\alpha_3=a,\,\alpha_4=-1$. After reduction we get that
	\begin{align*}
	\frac{1}{2}\sum_{i=1}^4 \sum_{j=1}^4 \alpha_i\alpha_j d_\cX(x_i,x_j)^2	=& a(1-2a)\|x+y\|^2 +a(1-2a)\|x-y\|^2 \\
	& -(1-2a)\|x\|^2-2a(1-2a)\|y\|^2,
	\end{align*}
	is less than or equal to zero. Rearranging this equation and dividing with $(1-2a)>0$ on both sides, we get that
   \begin{align*}
   \|x+y\|^2 + \|x-y\|^2 \leq a^{-1}\|x\|^2 +2\|y\|^2,
   \end{align*}
   and, if we let $a\uparrow 1/2$, then we get that $\|x+y\|^2 + \|x-y\|^2 \leq 2\|x\|^2 +2\|y\|^2$. Now this holds for any $x,y\in \cX$, hence we may also let $x'=x+y\in \cX$ and $y'=x-y\in \cX$ and note that
  \begin{align*}
  &\|x'+y'\|^2 + \|x'-y'\|^2 \leq 2\|x'\|^2 +2\|y'\|^2 \\
  \iff &\|2x\|^2 + \|2y\|^2 \leq 2\|x+y\|^2 +2\|x-y\|^2 \\
  \iff & 2\|x\|^2 + 2\|y\|^2 \leq \|x+y\|^2 +\|x-y\|^2,
  \end{align*}
  proving the reverse inequality. We conclude that $ \| x+y\|^2 + \|x-y\|^2 =2\|x\|^2 +2\|y\|^2$ for any $x,y\in \cX$, proving that $\cX$ is an inner product space.
\end{p}
 Now we return to our main topic of this section - to prove that metric spaces are of negative type if and only if they can be isometrically embedded into a Hilbert space. As mentioned above there is a connection between metric spaces of negative type and positive definite kernels. Thus we start by defining positive definite kernels and thereafter establish this connection.
\begin{definition}
	Let $\cX$ be a non-empty set. A symmetric mapping $k:\cX \times \cX\to \R$ is called a positive definite kernel if 
	\begin{align*}
	\sum_{i=1}^n\sum_{j=1}^n \alpha_i\alpha_j k(x_i,x_j) \geq 0,
	\end{align*}
	for any $n\in\N$, $\alpha\in\R^n$ and $x\in\cX^n$	 .
\end{definition}
Note that the definition of positive definite kernels can be equivalently represented in terms of positive semi-definite matrices. Let $k:\cX\times \cX \to \R$ be a symmetric mapping and let $K^n(x)$ denote the symmetric $n\times n$ dimensional real matrix with entries given by $K^n(x)_{i,j}=k(x_i,x_j)$ for any $x\in \cX^n$. Then
\begin{align*}
\alpha^\t K^n(x) \alpha =\sum_{i=1}^n\sum_{j=1}^n \alpha_i\alpha_j k(x_i,x_j),
\end{align*}
for any $\alpha\in \R^n$. Hence we get that $k:\cX\times \cX \to \R$ is a positive definite kernel if and only if $K^n(x)$ is a positive semi-definite matrix for all $n\in \N$ and $x\in \cX^n$. \\ \\
With the terminology in order, we are now ready to prove an equivalence between a metric space being of negative type and a certain mapping being a positive definite kernel.
\begin{theorem} \label{theorem_negative_type_iff_positive_kernel}
	A metric space $(\cX,d_\cX)$ is of negative type if and only if $d_0:\cX\times \cX\to \R$ given by
	\begin{align*}
	d_o(x,y) = d_\cX(x,o)+d_\cX(y,o)-d_\cX(x,y),
	\end{align*}
	for some $o\in\cX$, is a positive definite kernel.
\end{theorem}
\begin{p}
	Fix $o\in\cX$, $n\in \N$ and consider any $\alpha\in\R^n$ with $\sum_{i=1}^n\alpha_i=0$ and $x\in\cX^n$. Note that
	\begin{align*}
	\sum_{i=1}^n\sum_{j=1}^n \alpha_i\alpha_j d_\cX(x_i,x_j) &= -\sum_{i=1}^n\sum_{j=1}^n \alpha_i\alpha_j d_o(x_i,x_j) -\sum_{i=1}^n\sum_{j=1}^n \alpha_i\alpha_j [d_\cX(x_i,o)+d_\cX(x_j,o)] \\
	&=-\sum_{i=1}^n\sum_{j=1}^n \alpha_i\alpha_j d_o(x_i,x_j),
	\end{align*}
	using that the latter double sum can be written as two sums multiplied by $\sum_{i=1}^n\alpha_i=0$. This especially entails that, if $d_0$ is a positive definite kernel, then $(\cX,d_\cX)$ is a metric space of negative type. Now let $x_0=o$ and set $\alpha_0=-\sum_{i=1}^n\alpha_i$ such that $\sum_{i=0}^n\alpha_i =0$, note that
	\begin{align*}
	\sum_{i=0}^n\sum_{j=0}^n \alpha_i\alpha_j d_\cX(x_i,x_j) =& \sum_{i=1}^n \sum_{j=1}^n \alpha_i\alpha_j d_\cX(x_i,x_j) +  \sum_{i=1}^n \alpha_0\alpha_i d_\cX(x_0,x_i) +\sum_{j=1}^n \alpha_0\alpha_j d_\cX(x_0,x_j) \\
	&+ a_0^2d_\cX(x_0,x_0) \\
	=&\sum_{i=1}^n \sum_{j=1}^n \alpha_i\alpha_j d_\cX(x_i,x_j) -  \sum_{i=1}^n \sum_{j=1}^n\alpha_i\alpha_j [d_\cX(x_0,x_i) +d_\cX(x_0,x_j)] \\
	=&\sum_{i=1}^n \sum_{j=1}^n \alpha_i\alpha_j [d_\cX(x_i,x_j) - d_\cX(x_0,x_i) -d_\cX(x_0,x_j)] \\
	=& -\sum_{i=1}^n \sum_{j=1}^n \alpha_i\alpha_j d_o(x_i,x_j),
	\end{align*}
	proving that $d_0:\cX\times \cX\to \R$ is a positive definite kernel, when $(\cX,d_\cX)$ is a metric space of negative type.
\end{p}
As mentioned above, positive definite kernels are intertwined with the theory of reproducing kernel Hilbert spaces, so we summarize some useful relations from this theory.
\begin{lemma} \label{lemma_symmetric_positive_definite_kernel_implies_existence_of_hilbert_space_and_inner_prod_relation}
	If $\cX$ is a separable topological space and $k:\cX\times \cX\to \R$ is a continuous positive definite kernel, then there exist a separable $\R$-Hilbert space $\cH$ and a mapping $\phi:\cX\to \cH$ such that
	\begin{align*}
	k(x,y)=\la \phi(x), \phi(y) \ra_{\cH}.
	\end{align*}
	for any $x,y\in\cX$. We say that $k$ is the reproducing kernel of the reproducing kernel Hilbert space $\cH$.
\end{lemma}
\begin{p}
	Consider the real linear subspace \begin{align*}
	H_{\text{pre}}=\mathrm{span}(\{k(\cdot,y):\cX\to\R \, | \, y\in \cX\})\subset \{f\, |\, f:\cX\to\R\},
	\end{align*}
	 and note that any $f,g\in H_{\text{pre}}$ has some representation by  a finite linear combination $f=\sum_{i=1}^n a_i k(\cdot, x_i)$ and $g=\sum_{j=1}^m b_j k(\cdot, x_j')$ and for any such elements we define
	\begin{align*}
	\la f,g \ra = \sum_{i=1}^n \sum_{j=1}^m a_ib_j k(x_i,x_j').
	\end{align*}
	We may note that $\la \cdot , \cdot \ra : H_{\text{pre}}\times H_{\text{pre}} \to \R$ defined by the above relation is a well-defined mapping. To see this, note that 
	\begin{align*}
	\la f,g \ra = \sum_{i=1}^n \sum_{j=1}^m a_ib_j k(x_i,x_j') = \sum_{i=1}^n a_i g(x_i).
	\end{align*}
	Hence if for any two representations, $g=\sum_{j=1}^{m_1}b_{1,j}k(\cdot,x_{1,j}')=\sum_{j=1}^{m_2}b_{2,j}k(\cdot,x_{2,j}')$ we have that
	\begin{align*}
	\Big\la f , \sum_{j=1}^{m_1}b_{1,j}k(\cdot,x_{1,j}') \Big \ra  =\sum_{i=1}^n a_i g(x_i) =\Big\la f , \sum_{j=1}^{m_2}b_{2,j}k(\cdot,x_{2,j}') \Big \ra,
	\end{align*}
	so it evidently does not depend on the specific representation of $g$, and by symmetry of $k$ it is also independent of the representation of $f$. By the symmetry of $k$, we also have symmetry of $\la \cdot , \cdot \ra$ and for $f_1,f_2,g\in H_{\text{pre}}$ and $a_1,a_2\in\R$ and note that
	\begin{align*}
	\la a_1 f_1 + a_2f_2, g\ra = a_1\la f_1 , g \ra + a_2 \la f_2,g \ra,
	\end{align*}
	which is seen by taking a finite linear combination representation of $f_1,f_2,g$ and using the above definition followed by a separation of the terms. We also note that 
	\begin{align*}
	\la f,f \ra = \sum_{i=1}^n \sum_{j=1}^n a_ia_j k(x_i,x_j) \geq0,
	\end{align*}
	by the assumption that $k$ is a positive definite kernel. Hence $\la \cdot , \cdot \ra :H_{\text{pre}} \times H_{\text{pre}} \to \R$ is a semi-inner product and thus by the Cauchy-Bunyakowsky-Schwarz inequality (cf. proposition 1.4 \cite{conway1990ACourseInFunctionalAnalysis_reference_for_completion_inner_product_space}) we have that $	|\la f,g \ra |^2 \leq \la f, f\ra \la g,g \ra$.	As a consequence we get
	\begin{align*}
	|f(x)|^2  = \Big| \sum_{i=1}^n a_i k(x,x_i)\Big|^2 = |\la f, k(\cdot , x) \ra|^2 \leq \la f, f\ra \la k(\cdot,x ),k(\cdot ,x) \ra  =\la f, f\ra k(x,x),
	\end{align*}
	for any $x\in \cX$, hence $\la f, f\ra=0 \implies f\equiv0$. Conversely by the same considerations as above we have that $\la f,f \ra = \sum_{i=1}^na_if(x_i)$, proving that $f\equiv 0\implies \la f,f \ra=0$.  We conclude that $\la \cdot , \cdot \ra : H_{\text{pre}}\times H_{\text{pre}} \to \R$ is a inner product on $H_\text{pre}$, that is $(H_{\text{pre}},\la\cdot,\cdot\ra)$ is an inner product space. 
	
	Let $\cH$ denote the completion of $H_{\text{pre}}$ with respect to the inner product induced metric. By  proposition 1.9 \cite{conway1990ACourseInFunctionalAnalysis_reference_for_completion_inner_product_space}, we can extend the inner product $\la \cdot , \cdot \ra$ on $H_{\text{pre}}$ to an inner product $\la \cdot, \cdot \ra_{\cH}$ on $\cH$ in the following way. First one constructs the completion $\cH$ and a linear isometric embedding $\iota:H_{\text{pre}}\to \cH$ as done in  \cref{Appendix_Tensor_product_of_Hilbert_spaces}. The proposition from \cite{conway1990ACourseInFunctionalAnalysis_reference_for_completion_inner_product_space} now states that there exists an inner product $\la \cdot, \cdot \ra_{\cH}$ on $\cH$ such that
	\begin{align*}
	\la \iota(f),\iota(g) \ra_{\cH} = \la f , g \ra,
	\end{align*}
	for any $f,g\in H_{\text{pre}}$. We especially have that $\cH$ is an $\R$-Hilbert space  and we define $\phi:\cX\to \cH$ by $\phi(x) = \iota ( k(\cdot, x) )$,
	then we have that
	\begin{align*}
	\la \phi(x), \phi(y) \ra_\cH = \la k(\cdot,x) , k(\cdot,y) \ra = k(x,y). \\
	\end{align*}
	As regards the separability of $\cH$, we note that since $\cX$ is separable it has a countable dense subset $D$. We claim that the countable collection of maps
	\begin{align*}
	A =\lb  \sum_{i=1}^m a_ik(\cdot, x_i) : m\in \N, a\in \bQ^m, x\in D^m \rb,
	\end{align*}
	is dense in $H^{pre}$. To see this, let $f\in H^{pre}$ and note that $f$ has representation $f=\sum_{i=1}^m q_i k(\cdot ,x_i)$ for some $m\in \N,q\in \R^m,x\in\cX^m$. Since $D$ is dense in $\cX$ there exists a sequence $(x^n)\subset D^m$ such that $x^n_i \to_n x_i$ and since $\bQ$ is dense in $\R$ there exists a sequence $(q^n)\subset \bQ^m$ such that $q^n_i \to_n q_i$. Now define the sequence $(f_n)\subset A$   given by
	\begin{align*}
	f_n = \sum_{i=1}^m q_i^n k(\cdot , x_i^n) \in A,
	\end{align*}
	for all $n\in \N$. Note that if $\|f_n-f\|_{H^{pre}}^2\to_n 0$ we have that 	$A$ is a countable dense set in $H^{pre}$, rendering $H^{pre}$ separable.
	By writing the difference in norm as a linear combination of inner products one may realize that
	\begin{align*}
	\|f_n-f\|_{H^{pre}}^2 &= \sum_{i=1}^m \sum_{j=1}^m \Big[ q_i^n q_j^n k(x_i^n,x_j^n) +q_i q_j k(x_i,x_j) -2 \lp q_i^n q_j k(x_i^n,x_j)\rp \Big] \to_n 0,
	\end{align*}
	where we used the continuity of $k:\cX\times \cX\to \R$.  Thus $H^{pre}$ is a separable inner product space and since $\cH$ is the completion of $H^{pre}$, we  get that $\cH$ is a separable $\R$-Hilbert space (cf. theorem 3.25. \cite{caratheodoryfunctions_ref}).
\end{p}
The above lemma now allows us to establish a condition on metric spaces to be of negative type in terms of the isometric embeddability into Hilbert spaces.  \\ \\
First let us define what we mean by an isometric embedding between spaces, and show some important properties of such mappings. For any two metric spaces $(\cX,d_\cX)$ and $(\cY,d_\cY)$ we say that a map $f:(\cX,d_\cX) \to( \cY,d_\cY)$ is an isometry if it satisfies $d_\cX(x_1,x_2) = d_\cY(f(x_1),f(x_2)),
$ and if such a mapping exists, we say that $(\cX,d_\cX)$ is isometrically embeddable into $(\cY,d_\cY)$.
Whenever we do not specify the metric for  the isometric embeddings, it is equipped with the natural metric. For example, if we have a metric space $(\cX,d_\cX)$ and a Hilbert space $\cH$, then we say that $\phi:(\cX,d_{\cX}^{1/2})\to \cH$ is an isometric embedding if
\begin{align} \label{eq.temp.1233123333}
d_{\cX}(x_1,x_2) = \| \phi(x_1)-\phi(x_2) \|_{\cH}^2.
\end{align}
for all $x_1,x_2\in \cX$. We will apply the convention that, whenever we mention a isometric embedding from a metric space to a Hilbert space, it should be understood as in \cref{eq.temp.1233123333}.  Note that the isometric mapping is automatically injective and continuous.  \\ \\ 
An important property of such isometric embeddings into Hilbert spaces is that they are Pettis integrable with respect to every measure in $M^1(\cX)$. A fact that is used in the  main theorem below about the equivalence between metric spaces of negative type and their embeddability into Hilbert spaces, hence we briefly explain what we mean by Pettis integrable. 

A thorough introduction to the Pettis integral of Hilbert space valued mappings can be found in  \cref{Integration of Hilbert space valued mappings}. For the reader who is somewhat familiar with the Pettis integral the following one sentence recap can be given: If $\phi:\cX\to\cH$ is a $\bK$-Hilbert space valued mapping, which is Pettis integrable with respect to a finite positive measure $\mu$ on $(\cX,\cB(\cX))$, then the Pettis integral of $\phi$ over $\cX$ with respect to $\mu$, denoted $\int \phi \, d\mu$, is defined to be the unique element in $\cH$ satisfying
\begin{align*}
h^* \lp \int \phi \, d\mu\rp = \int h^* \circ \phi(x) \,  d\mu(x) \quad \quad \forall h^*\in \cH^*,
\end{align*}
where $\cH^*$ is the continuous dual space of $\cH$, that is \begin{align*}
\cH^*=\{h^*:\cH\to \bK: h^* \textit{ is continuous and linear}\}.
\end{align*}
Note that above and for the remainder of this thesis we let $\bK$ denote either $\R$ or $\bC$, and in any given scenario we only use $\bK$ whenever the result or statement holds for both $\R$ and $\bC$ respectively. E.g. a $\bK$-Hilbert space is either an $\R$-Hilbert space (real Hilbert space) or a $\bC$-Hilbert space (complex Hilbert space).

\begin{lemma} \label{lemma_isometric_embeddings_are_pettis_integrable_wrt_first_moment_measures}
	Any isometric embedding $\phi:(\cX,d_\cX^{1/2})\to \cH$ into a $\bK$-Hilbert space is Pettis integrable with respect to both Jordan-Hahn decompositions $\mu^+$ and $\mu^-$, for any $\mu\in M^1(\cX)$. 
\end{lemma}
\begin{p}
	Let $\mu\in M^1(\cX)$ be any finite signed Borel measure on $\cX$ with finite first moment. By \cref{Theorem_definition_Pettis_Integral} is suffices to show that $\phi$ is $\mu^\pm$-scalarly integrable. That is, it suffices to show that $h^*\circ \phi \in L^1(\mu^\pm)$ for all $h^*\in \cH$. \textit{Measurability}: $\phi$ is an  $(d_\cX^{1/2},d_\cH)$-isometry, that is
	\begin{align*}
	\sqrt{d_\cX(x_1,x_2)}=\|\phi(x_1)-\phi(x_2)\|_\cH \quad \quad  \quad \forall x_1,x_2\in \cX,
	\end{align*}
	from which continuity and hence $\cB(\cX)/\cB(\cH)$-measurability follows. By the same reasoning we get that every function $h^*\in \cH^*$ is continuous and hence $\cB(\cH)/\cB(\R)$-measurable, so composition $h^*\circ \phi$ is indeed $\cB(\cX)/\cB(\R)$-measurable. \textit{Integrability:} Note that for any $h^*\in \cH^*$ and $x\in\cX$ we have that $\|h^*\|_{\text{op}}<\i$ and $|h^*\circ \phi(x)| \leq \|h^*\|_{\text{op}}\|\phi(x)\|_{\cH}$, where $\|\cdot\|_{\mathrm{op}}$ denotes the operator norm on $\cH^*$. Thus we get that
	\begin{align*}
	\frac{1}{\|h^*\|_{\text{op}}}\int |h^* \circ \phi(x)| d \mu^\pm (x) &\leq \int  \|\phi(x) \|_\cH \, d\mu^\pm (x) \\
	&\leq \int \| \phi(x)\|_\cH \, d|\mu |(x) \\
	& \leq \int \| \phi(x)- \phi(o) \|_\cH \, d|\mu|(x) + \int \|\phi(o)\|_\cH \, d|\mu|(x) \\
	& =\int d_\cX(x,o)^{1/2} \, d|\mu|(x) +  |\mu|(\cX)\|\phi(o)\|_\cH \\
	& \leq \int d_\cX(x,o) \, d|\mu|(x) +  |\mu|(\cX)(\|\phi(o)\|_\cH  +1) <\i,
	\end{align*}
	for any $o\in \cX$, since $\mu$ has finite first moment and its total variation is finite (cf. Corollary 3.1.2 \cite{MeasureTheory2007bogachev} ). One might note that only half moments of $\mu$ were needed in order to ensure integrability.
\end{p}

For our purpose it does not suffice with the Pettis integral with respect to positive measures. Hence we will define the Pettis integral with respect to signed measures in the same way as one defines the Lebesgue integral with respect to signed measures. That is, for a finite signed measure $\mu$ on $(\cX,\cB(\cX))$ with Jordan-Hahn decompositions $\mu^+$ and $\mu^-$ and a Hilbert space valued mapping $f:\cX\to\cH$, we define the Pettis integral of $f$ over $\cX$ with respect to $\my$ by
\begin{align*}
\int f \, d\mu := \int f\, d\mu^+ - \int f \, d\mu^-,
\end{align*}
whenever both Pettis integrals exist. We may also realize that the unique defining property of the Pettis integral agrees with this definition. The linearity of $h^*\in \cH^*$ allows us to say
\begin{align*}
h^* \lp  \int f\, d\mu \rp &= h^* \lp  \int f\, d\mu^+ \rp - h^* \lp  \int f\, d\mu^- \rp\\
&=\int h^* \circ f \, d\mu^+ -\int h^*\circ f\, d\mu^+ = \int h^* \circ f \, d\mu,  \quad \quad \forall h^* \in \cH^*,
\end{align*}
and as a consequence the Pettis integral of $f$ over $\cX$ with respect to a finite signed measure $\mu$ is the unique element in $\cH$ satisfying
\begin{align*}
h^* \lp  \int f\, d\mu \rp= \int h^* \circ f \, d\mu,  \quad \quad \forall h^* \in \cH^*.
\end{align*}
Uniqueness is realized in the following way: assume for contradiction that \begin{align*}
h^*(a)=\int h^*\circ f \, d\mu = h^*(b), \quad \quad \forall h^*\in \cH^*,
\end{align*}
for two $a,b\in \cH$ with $a\not = b$.  Now note that  $\cH^*$ separates points in $\cH$ for any Hilbert space $\cH$ (cf. theorem 3.4 \cite{rudin1991functional}), meaning that $h^*(a)=h^*(b)$ for all $h^*\in \cH^*$ implies that $a=b$. This is a contradiction, so we have that $a=b$, proving uniqueness.\\ \\
Having defined all necessary Pettis integral terminology needed, we may now state one of the most important theorems of this section, namely the previously mentioned equivalence between a metric space of negative type and its isometric embeddability into a Hilbert space.
\begin{theorem} \label{theorem_equivalence_of_negative_type}
	Let $(\cX,d_\cX)$ be a metric space. Then the following statements are equivalent
	\begin{itemize}
		\item[1)] $(\cX,d_\cX)$ is a metric space of negative type.
		\item[2)] There exists an isometric embedding $\phi:(\cX,d_\cX^{1/2})\to \cH$ into a separable $\R$-Hilbert space.
		\item[3)] There exists an isometric embedding $\phi':(\cX,d_\cX^{1/2})\to \cH'$ into a separable $\bC$-Hilbert space.
	\end{itemize} 
\end{theorem}
\begin{p}
	\textit{2)}$\iff$\textit{3)}: To fully understand the equivalence between \textit{2)} and \textit{3)} we encourage the reader to understand how the realification of a complex Hilbert space and the complexification of a real Hilbert space are constructed (see \cref{Appendix_Complexification_and_Realification}). 
	
	Nevertheless if  $\phi:(\cX,d_\cX^{1/2})\to \cH$ is an isometric embedding into a separable $\R$-Hilbert space, then we may note that $\phi'=\mathrm{cpx}\circ \phi:(\cX,d_\cX^{1/2})\to \cH^\bC$ is an isometric embedding into a separable $\bC$-Hilbert space, where the complexification map $\mathrm{cpx}:\cH\to \cH^\bC$ into the complexification $\cH^\bC$ of $\cH$ is  an isometry. Conversely if $\phi':(\cX,d_\cX^{1/2})\to \cH'$ is an isometric embedding into a separable $\bC$-Hilbert space, then we may note that $\phi=\mathrm{re}\circ \phi:(\cX,d_\cX^{1/2})\to \cH'^\R$ is an isometric embedding into a separable $\R$-Hilbert space, where the realification map $\mathrm{re}:\cH' \to \cH'^\R$ into the realification $\cH'^{\R}$ of $\cH'$ is an isometry. \\ \\
	\textit{1)}$\implies$\textit{2)}: Assume that $(\cX,d_\cX)$ is of negative type. For some $o\in \cX$ we define the function $d_o:\cX \times \cX \to \R$ by
	\begin{align*}
	d_o(x_1,x_2)= \frac{d_\cX(x_1,o)+d_\cX(x_2,o)-d_\cX(x_1,x_2)}{2}.
	\end{align*}
	By \cref{theorem_negative_type_iff_positive_kernel} $d_o$  is a positive definite kernel (since $2d_o$ is). Furthermore $\cX$ is separable and $d_0:\cX\times \cX\to \R$ is continuous, by the continuity of $d_\cX:\cX\times \cX\to\R$. Hence by \cref{lemma_symmetric_positive_definite_kernel_implies_existence_of_hilbert_space_and_inner_prod_relation} we know that there exists a separable $\R$-Hilbert space $\cH$ and a mapping $\phi:\cX\to \cH$ such that
	\begin{align*}
	d_o(x_1,x_2) = \la \phi(x_1),\phi(x_2) \ra_\cH.
	\end{align*}
	We realize that
	\begin{align*}
	\| \phi(x_1)-\phi(x_2) \|_\cH^2 &= \la \phi(x_1),\phi(x_1)\ra_\cH + \la\phi(x_2),\phi(x_2)\ra -2\la\phi(x_1),\phi(x_2)\ra_\cH \\
	&= d_o(x_1,x_1)+d_o(x_2,x_2)-2d_o(x_1,x_2) \\
	&= d_\cX(o,x_1)+d_\cX(o,x_2)-d_\cX(x_1,o)-d_\cX(x_2,o)+d_\cX(x_1,x_2) \\
	&= d_\cX(x_1,x_2),
	\end{align*}
	proving that there exists an isometric embedding $\phi:(\cX,d_\cX^{1/2}) \to \cH$ into a separable $\R$-Hilbert space. \\ \\	
	\textit{2)}$\implies$\textit{1)}: Assume that there is an an isometric embedding $\phi:(\cX,d_\cX^{1/2})\to \cH$ into a separable $\R$-Hilbert space. 	 Let $n\in\N$, $x_1,...,x_n\in \cX$ and $\alpha_1,...,\alpha_n\in\R$ with $\sum_{i=1}^n \alpha_i =0$. Define a signed measure $\mu\in M(\cX)$ by 
	 \begin{align*}
	 \mu = \sum_{i=1}^n \alpha_i \delta_{x_i} = \underbrace{\sum_{i=1}^n (\alpha_i \lor 0) \delta_{x_i}}_{:=\mu_1} - \underbrace{\sum_{i=1}^n ((-\alpha_i) \lor 0) \delta_{x_i}}_{:=\mu_2}.
	 \end{align*}
	 Note that $\mu=\mu_1-\mu_2$ is not necessarily the Jordan-Hahn decomposition, but we know that $\mu^+ \leq \mu_1$ and $\mu^- \leq \mu_2$, thus $|\mu|\leq  \mu_1+\mu_2$, by the proof of \cref{lemma_M(X)_and_M^1(X)_are_vector_spaces}. Hence $
	 \int d_\cX(x,o) d|\mu|(x) \leq \sum_{i=1}^n d_\cX(x_i,o)|\alpha_i| < \i $
	 for any $o\in \cX$, so $\mu\in M^1(\cX)$. Thus $D(\mu)$ is well-defined and equals 
	 \begin{align*}
	 D(\mu) = \int \int d_\cX(x,y) \, d\mu(x) \, d\mu(y) = \sum_{i=1}^n \sum_{j=1}^n \alpha_i \alpha_j d_\cX(x_i,x_j),
	 \end{align*}
	 where we have used that $\int f \, d\mu := \int f \, d\mu^+ - \int f \, d\mu^- = \int f \, d\mu_1 - \int f\, d\mu_2$, see proof of \cref{lemma_beta_phi_is_linear} for further explanation.
	 Therefore it suffices to show that $D(\mu)\leq 0$. Since $\phi:(\cX,d_\cX^{1/2})\to \cH$ is an isometric embedding we have that
	 \begin{align*}
	 D(\mu) &= \int \int d_\cX(x,y) \, d\mu(x) \, d\mu(y)  \\
	 &= \int \int \| \phi(x)-\phi(y) \|_\cH^2 \,  d\mu(x) \, d\mu(y)  \\
	 &= \int\int \|\phi(x)\|_\cH^2 + \|\phi(y) \|_\cH^2 - 2\la \phi(x),\phi(y)\ra_\cH \,  d\mu(x) \, d\mu(y)  \\
	 &=  2\mu(\cX) \int \|\phi(x) \|_\cH^2 \, d\mu(x) -  2\int \int\la \phi(x),\phi(y)\ra_\cH \,  d\mu(x) \, d\mu(y) , 
	 \end{align*}
	 by integrability of each term. This is seen by noting that $x\mapsto \|\phi(x)\|_\cH^2\in \cL^1(\mu)$ (expand with triangle inequality around $o\in \cX$  to get integrable upper bound $d_\cX(x,o) + c_1 +c_2d_\cX(x,o)^{1/2}$) and that $(x,y)\mapsto\la \phi(x),\phi(y)\ra_{\cH}\in \cL^1(\mu\times \mu)$ (by Cauchy-Schwarz inequality $|\la \phi(x),\phi(y)\ra_{\cH}|\leq \|\phi(x)\|_\cH \|\phi(y)\|_\cH $, where the upper bound is integrable by proof of \cref{lemma_isometric_embeddings_are_pettis_integrable_wrt_first_moment_measures}).
	 It holds that  $\mu(\cX)=0$ (i.e. $\mu\in M^1_0(\cX)$) since $\sum_{i=1}^n \alpha_i = 0$, and therefore we get
	 \begin{align*}
	 D(\mu) = -2 \int \int \la \phi(x),\phi(y)\ra_\cH\, d\mu(x) \, d\mu(y).
	 \end{align*}
 Note that $\phi$ is Pettis integrable with respect to any $\mu\in M^1(\cX)$ (see \cref{lemma_isometric_embeddings_are_pettis_integrable_wrt_first_moment_measures}) so the Pettis integral $\int \phi \, d\mu\in\cH$ exists. By the unique defining property of the Pettis integral we can derive that
	 \begin{align*}
	 \int \int \la \phi(x),\phi(y)\ra_\cH\, d\mu(x) \, d\mu(y)&= \int \Big\la \int  \phi \, d\mu,\phi(y) \Big\ra_\cH  \, d\mu(y) \\
	 &=  \Big\la \int  \phi \, d\mu, \int\phi  \, d\mu \Big\ra_\cH   \\
	 &= \Big\| \int  \phi \, d\mu \Big\|_\cH^2,
	 \end{align*}
	 where we used that $x\mapsto \la x,\phi(y)\ra_\cH$ for every fixed $y\in\cX$ and that $y\mapsto \la \int \phi \, d\mu, y\ra_{\cH}$ are continuous linear mappings ($\cH$ is a $\R$-Hilbert space).  Hence
	 \begin{align} \label{eq_temp_D(mu)_equality_2_normbeta(mu)}
	 D(\mu) = -2\Big\| \int  \phi \, d\mu \Big\|_\cH^2 \leq 0.
	  \end{align}
	 This concludes the proof. As a closing remark, we may add that \cref{eq_temp_D(mu)_equality_2_normbeta(mu)} holds for general finite signed measures $\mu$ with finite first moment and $\mu(\cX)=0$, i.e. for any $\mu\in M^1_0(\cX)$.
\end{p}
The above theorem is crucial for our further development. Because if we assume that $(\cX,d_\cX)$ and $(\cY,d_\cY)$ are metric spaces of negative type, then  this theorem is mainly responsible for the derivation of the alternative Hilbert space representation of the distance covariance measure $dcov(\theta)$. \\ \\
The next order of business is to derive this alternative representation and this is done through the analysis of some rather abstract maps from $M^{1,1}(\cX\times \cY)$ to a tensor product of Hilbert spaces.  By \cref{lemma_isometric_embeddings_are_pettis_integrable_wrt_first_moment_measures} we may define the following mean embedding map, which is extensively used throughout the remainder of the thesis. 

\begin{definition}
	For any isometric embedding $\phi:(\cX,d^{1/2})\to \cH$ into a $\bK$-Hilbert space, we  define the mean embedding of $\phi$ as the map $\beta_\phi:M^1(\cX)\to \cH$  given by 
	\begin{align*}
	\beta_\phi(\mu) = \int \phi \, d\mu,
	\end{align*}
	for any $\mu \in M^1(\cX)$.
\end{definition}
We may also note that when $\mu\in M^1_1(\cX)\subset M^1(\cX)$, that is a probability measure with finite first moment then $\beta_\phi(\mu) = E(\phi(X))$ if $X\sim \mu$ (see \cref{defi_Hilbert_space_expectation_and_variance}). That is why we call $\beta_\phi$ the mean embedding map of $\phi$.  We will later use that this mapping is linear, so we start by showing this property.
\begin{lemma} \label{lemma_beta_phi_is_linear}
	For any isometric embedding $\phi:(\cX,d^{1/2})\to \cH$ into a $\bK$-Hilbert space,   it holds that 
		\begin{align*}
	\beta_\phi( a\mu+b\nu) = a\beta_\phi(\mu)+b\beta_\phi(\nu),
	\end{align*}
	for any $\mu,\nu\in M^1(\cX)$ and  $a,b\in \R$. Thus if $\cH$ is a $\R$-Hilbert space then $\beta_\phi$ is linear, and if $\cH$ is a $\bC$-Hilbert space then $\beta_\phi$ is linear if we view it as a map into the realification $\cH^\R$.
\end{lemma}
\begin{p}
	We refer the reader to appendix \cref{Appendix_Complexification_and_Realification} for the construction of the realification of a complex Hilbert space. \\ \\
	For any two signed measures $\mu,\nu\in M^1(\cX)$ and $a,b\in \R$, we have to show that $
	\beta_\phi(a\mu+b\nu) = a\beta_\phi(\mu)+b\beta_\phi(\nu).$
	That is, we need to show that
	\begin{align*}
	\int \phi \, d(a\mu+b\nu) 	= a \int \phi \, d\mu
	+ b \int \phi \, d\nu.
	\end{align*}
	By the unique defining property of the Pettis integral with respect to finite signed measures, we know that is suffices to show that
	\begin{align*}
	\int h^*\circ \phi \, d(a\mu+b\nu)= h^*\lp a\lp \int \phi \, d\mu\rp
	+ b\lp \int \phi \, d\nu\rp \rp ,
	\end{align*}
	for any $h^*\in\cH^*$. For any $h^*\in \cH^*$ we use the linearity of $h^*$ and the unique defining property for each of the individual Pettis integrals to rewrite the right-hand side
	\begin{align*}
	\int h^*\circ\phi \, d(a\mu+b\nu)= a \int h^*\circ\phi \, d\mu 
	+ b \int h^*\circ\phi \, d\nu.
	\end{align*}
	This equality indeed holds, and follows from standard integration theory  with respect to signed measures, proving that $\beta_\phi:M^1(\cX)\to \cH$ is a linear map. \\ \\
	To keep the thesis self-contained we sketch the proof of this equality. It holds by the standard approach of showing that the equality holds for characteristic, hence simple functions, followed by an approximation of measurable functions by simple functions. Lastly one goes to the limit of these approximations by for example Lebesgue's dominated convergence theorem for signed measures.
\end{p}
To summarize our previous findings: if both $(\cX,d_\cX)$ and $(\cY,d_\cY)$ are metric spaces of negative type, then we know that  there exist two isometric embeddings $\phi:(\cX,d_\cX^{1/2})\to \cH_1$ and $\psi:(\cY,d_\cY^{1/2})\to\cH_2$ into two separable $\bK$-Hilbert spaces $\cH_1$ and $\cH_2$. 

Hence the idea is now to construct the tensor product map of these two embeddings $\phi\otimes \psi:\cX\times \cY\to\cH_1\otimes \cH_2$ and show that it is Pettis integrable with respect to any measure in $M^{1,1}(\cX\times \cY)$. An introduction to the construction of the tensor product of Hilbert spaces can be found in \cref{Appendix_Tensor_product_of_Hilbert_spaces}. This allows for the construction of a mean embedding map $\beta_{\phi\otimes \psi}:M^{1,1}(\cX\times \cY) \to \cH_1\otimes \cH_2$ in the same fashion as we constructed $\beta_\phi$.  The purpose of constructing $\beta_{\phi\otimes \psi}:M^{1,1}(\cX\times \cY) \to \cH_1\otimes \cH_2$, is as mentioned in the introduction to \cref{section_negative_and_strong_negative_type}, that we can represent the distance covariance measure $dcov(\theta)$ in terms of this Hilbert space mapping . \\ \\
First we prove a lemma, which guarantees sufficient integrability needed to construct the mean embedding map $\beta_{\phi\otimes \psi}:M^{1,1}(\cX\times \cY) \to \cH_1\otimes \cH_2$.
\begin{lemma} \label{lemma_tensor_product_map_of_isometric_embeddings_is_pettis_integrable}
	Let $\phi:(\cX,d_{\cX}^{1/2})\to \cH_1$ and $\psi:(\cY,d_{\cY}^{1/2})\to \cH_2$ be any two isometric embeddings into two separable Hilbert spaces with the same scalar field. Then $\phi\otimes \psi:(\cX\times \cY) \to \cH_1 \otimes \cH_2$ defined by
	\begin{align*}
	\phi\otimes \psi (x,y) = \phi(x) \otimes \psi(y),
	\end{align*}
	is Pettis integrable with respect to any $\theta \in M^{1,1}(\cX\times \cY)$. 
\end{lemma}
\begin{p}
	Note $(\cH_1\otimes \cH_2,\la \cdot , \cdot \ra_{\cH_1\otimes \cH_2})$ is a Hilbert space, and its construction can be found in \cref{Appendix_Tensor_product_of_Hilbert_spaces}. This construction requires that both Hilbert spaces have the same scalar field, which is why we insist that $\cH_1$ and $\cH_2$ have the same scalar field. \\ \\
	The mapping $\phi\otimes \psi:(\cX\times \cY) \to \cH_1 \otimes \cH_2$ (with slight abuse of notation) given by $\phi\otimes \psi (x,y) = \phi(x) \otimes \psi(y)$ is realized to be the composition $ f\circ (\phi\times \psi)$, where $f:\cH_1\times \cH_2 \to \cH_1 \otimes \cH_2$ maps products into the embedding of simple tensors in $\cH_1 \otimes \cH_2$  by $f(h_1,h_2)=h_1\otimes h_2$ (formally $\iota(h_1\otimes h_2)$, see \cref{Appendix_Tensor_product_of_Hilbert_spaces}) and $\phi\times \psi:\cX\times \cY\to \cH_1\times \cH_2$ given by $\phi\times \psi(x,y)=(\phi(x),\psi(y))$. \\
	By \cref{appendix_lemma_tensor_map_is_measurable} we have that $f$ is $\cB(\cH_1)\otimes \cB(\cH_2)/\cB(\cH_1\otimes \cH_2)$-measurable and $(x,y)\mapsto \psi(x)$, $(x,y)\mapsto \psi(y)$ are both continuous and hence measurable ($\cB(\cX)\otimes \cB(\cY)=\cB(\cX\times \cY)$), and as a consequence the bundle map $\phi\times \psi:(x,y) \mapsto (\phi(x),\psi(y))$ is $\cB(\cX)\otimes \cB(\cY)/\cB(\cH_1)\otimes \cB(\cH_2)$-measurable (cf. theorem 13.10 \cite{schilling}). We conclude that the composition $\phi\otimes \psi:(\cX\times \cY) \to \cH_1 \otimes \cH_2$ is indeed $\cB(\cX)\otimes \cB(\cY)/\cB(\cH_1\otimes \cH_2)$-measurable. \\ \\
	By  \cref{lemma_f_is_pettis_integrabel_if_measurable_and_norm_of_f_is_integrable} it now suffices to show that
	\begin{align*}
	\int \| \phi\otimes \psi(x,y)\|_{\cH_1\otimes \cH_2} \, d|\theta|(x,y) <\i,
	\end{align*}
	for any $\theta \in M^{1,1}(\cX\times \cY)$, for which we know that $\pi_1(|\theta|)\in M^1(\cX)$ and $\pi_2(|\theta|)\in M^1(\cY)$ for any $\theta \in M^{1,1}(\cX\times \cY)$. 	Assume that $\theta \in M^{1,1}(\cX\times \cY)$ and note that  
	$$
	\|\phi(x)\|_{\cH_1} \leq \|\phi(x)-\phi(o)\|_{\cH_1}+\|\phi(o)\|_{\cH_1} = d_{\cX}(x,o)^{1/2}+\|\phi(o)\|_{\cH_1},
	$$
	where the upper bound is square integrable with respect to $\pi_1(|\theta|)$ (see proof of \cref{lemma_isometric_embeddings_are_pettis_integrable_wrt_first_moment_measures}). Similarly we also get that $\|\psi(y)\|_{\cH_2}$ is square integrable with respect to $\pi_2(|\theta|)$.  Thus we may conclude that
	$(x,y)\mapsto \|\phi(x)\|_{\cH_1}\in \mathcal{L}^2(\theta)$ and $(x,y)\mapsto \|\psi(y)\|_{\cH_2}\in \mathcal{L}^2(\theta)$. Recall from \cref{Appendix_Tensor_product_of_Hilbert_spaces} that the inner product $ \la \cdot , \cdot \ra_{\cH_1 \otimes\cH_2}$  on $\cH_1\otimes \cH_2$ satisfies the following equality for simple tensor products
	\begin{align*}
	\|\phi(x)\otimes \phi(y)\|_{\cH_1 \otimes\cH_2} &= \la\phi(x)\otimes \phi(y),\phi(x)\otimes \phi(y) \ra_{\cH_1 \otimes\cH_2}^{1/2} \\
	&= \la\phi(x),\phi(x) \ra_{\cH_1}^{1/2} \la\psi(y),\psi(y) \ra_{\cH_2}^{1/2} \\
	&= \|\phi(x)\|_{\cH_1} \|\psi(y)\|_{\cH_2} \in \mathcal{L}^1 (\theta),
	\end{align*}
	by the Cauchy-Schwarz inequality, proving that  $\phi\otimes \psi:(\cX\times \cY) \to \cH_1 \otimes \cH_2$ is Pettis integrable with respect to  any $\theta \in M^{1,1}(\cX\times \cY)$.
\end{p}
\begin{remark}
	In the above lemma, it is essential that the two Hilbert spaces $\cH_1$ and $\cH_2$ are separable if we are to prove measurability of $f:\cH_1\times \cH_2\to \cH_1 \otimes \cH_2$ by utilizing \cref{appendix_lemma_tensor_map_is_measurable}. Hence if we did not make the universal assumption of only considering separable metric spaces, then  \cref{theorem_equivalence_of_negative_type} would only  guarantee isometric embeddings into general Hilbert spaces. As a consequence it would not be sufficient to assume, that both marginal metric spaces are of negative type, in order for the Hilbert space representation of the distance covariance measure to hold, since this is constructed on the premise that $f$ is measurable. I do not postulate that measurability does not hold when $\cH_1$ or $\cH_2$ are non-separable, only that I failed to show this. Thus, this is yet another part of the thesis, where we directly use separability of the marginal metric spaces.
\end{remark}
A consequence of the above lemma is that the extension of the class of mean embedding mappings $\beta$ to tensor product of isometric embeddings is well-defined. \\ \\\\
\begin{definition} \label{defi_beta_tensor_mappings}
	For any two isometric embeddings $\phi:(\cX,d_{\cX}^{1/2})\to \cH_1$ and $\psi:(\cY,d_{\cY}^{1/2})\to \cH_2$ into two separable Hilbert spaces with the same scalar field, we may define the mapping  $\beta_{\phi\otimes\psi}:M^{1,1}(\cX\times \cY) \to \cH_1\otimes \cH_2$ 
	by the Pettis integral
	\begin{align*}
	\beta_{\phi\otimes\psi}(\theta) = \int \phi\otimes\psi \, d\theta,
	\end{align*}
	for any $\theta \in M^{1,1}(\cX\times \cY)$.
\end{definition} 
Furthermore,  this mean embedding map of the tensor product of isometries exhibits the same linear properties as the previously defined mean embedding map. Again this property is not needed for the derivation of the Hilbert space representation of the distance covariance measure, but it will be crucial when analysing the representation to find out for which marginal metric spaces we have the implication $dcov(\theta)=0\implies \theta=\mu\times \nu$.
\begin{corollary} \label{lemma_beta_tesor_is_linear}
	If $\phi:(\cX,d_{\cX}^{1/2})\to \cH_1$ and $\psi:(\cY,d_{\cY}^{1/2})\to \cH_2$ are isometric embeddings into two separable Hilbert spaces with the same scalar field $\bK$, then
	\begin{align*}
	\beta_{\phi\otimes \psi}(a\theta_1 +b \theta_2) = a\beta_{\phi\otimes \psi}(\theta_1)+b\beta_{\phi\otimes \psi}(\theta_2),
	\end{align*}
	 for any $\theta_1,\theta_2\in M^{1,1}(\cX\times \cY)$ and  $a,b\in \R$. Thus, if $\bK=\R$,  then $\beta_{\phi\otimes \psi}$ is linear, and if $\bK=\bC$, then $\beta_{\phi\otimes \psi}$ is linear when viewed as a map into the realification $(\cH_1\otimes \cH_2)^\R$.
\end{corollary}
\begin{p}
	The proof follows by arguments identical to those of  \cref{lemma_beta_phi_is_linear}.
\end{p}
Before proceeding with the proof of the alternative representation of the distance covariance measure in terms of the Hilbert space valued mapping $\beta_{\phi\otimes \psi}$, we state a lemma which will help facilitate the derivation of this representation This lemma  especially implies that, if both $(\cX,d_\cX)$ and $(\cY,d_\cY)$ are metric spaces of negative type as witnessed by isometric embeddings $\phi$ and $\psi$ respectively, then $d_\mu$ and $d_\nu$ have a Hilbert space representation in terms of $\phi,\beta_\phi$ and $\psi,\beta_\psi$ respectively.
\begin{lemma} \label{lemma_a_mu_and_D(mu)_and_d_mu_representations}
	Let $\phi:(\cX,d_\cX^{1/2})\to \cH$ be a isometric embedding into a separable $\bK$-Hilbert space and let $X\sim\mu\in M_1^1(\cX)$. It holds that 
	\begin{itemize}
		\item[1)] $a_\mu(x)= \|\phi(x)-\beta_\phi(\mu)\|^2_\cH + D(\mu)/2$,
		\item[2)] $D(\mu)= 2\text{Var}(\phi(X))$,
		\item[3)] $d_\mu(x,y) = - \la \phi(x)-\beta_\phi(\mu),\phi(y)-\beta_\phi(\mu) \ra_\cH- \la \phi(y)-\beta_\phi(\mu) ,\phi(x)-\beta_\phi(\mu)\ra_\cH,$
	\end{itemize}
for all $x,y\in \cX$.
\end{lemma}
\begin{p} 
	Let $\mu\in M^1_1(\cX)$ and let $X$ be a random Borel element in $\cX$ defined on some probability space $(\Omega,\F,P)$ with $X(P)=\mu$. Now note that $\beta_\phi(\mu)=E\phi(X)$ and hence  $$\int \|\phi(x)-\beta_\phi(\mu)\|^2_\cH  \, d\mu(x) =\int \|\phi(X)-E\phi(X)\|^2_\cH dP=E\|\phi(X)-E\phi(X)\|_\cH^2=:\text{Var}(\phi(X)),$$
	(see \cref{defi_Hilbert_space_expectation_and_variance}). Thus, using that $\phi:(\cX,d_\cX^{1/2})\to \cH$ is an isometric embedding, we get  
	\begin{align*}
	a_\mu(x) =& \int d(x,y)\, d\mu(y) = \int \| \phi(x)-\phi(y) \|^2_\cH \, d\mu(y) \\
	=& \int \|\phi(x)-\beta_\phi(\mu)+\beta_\phi(\mu)-\phi(y)\|^2_\cH\, d\mu(y)  \\
	=& \int \| \phi(x)-\beta_\phi(\mu) \|^2_\cH + \|\phi(y)-\beta_\phi(\mu)\|^2_\cH -  \la \phi(x)-\beta_\phi(\mu),\phi(y)-\beta_\phi(\mu) \ra_\cH \\
	&- \la \phi(y)-\beta_\phi(\mu),\phi(x)-\beta_\phi(\mu) \ra_\cH \, d\mu(y) \\
	=& \| \phi(x)-\beta_\phi(\mu) \|^2_\cH +\text{Var}(\phi(X)) - \overline{\int   \la \phi(y)-\beta_\phi(\mu),\phi(x)-\beta_\phi(\mu) \ra_\cH \, d\mu(y) } \\
	&- \int    \la \phi(y)-\beta_\phi(\mu),\phi(x)-\beta_\phi(\mu) \ra_\cH \, d\mu(y),
	\end{align*}
	where we used that each term is integrable.
	Since $h\mapsto \la h-\beta_\phi(\mu), \phi(x)-\beta_\phi(\mu) \ra_\cH$ is a continuous linear mapping for any $x\in \cX$, the unique defining property of the Pettis integral yields that \begin{align*}
	\int    \la \phi(y)-\beta_\phi(\mu),\phi(x)-\beta_\phi(\mu) \ra_\cH \, d\mu(y) &=    \Big\la  \int \phi \, d\mu -\beta_\phi(\mu),\phi(x)-\beta_\phi(\mu) \Big\ra_\cH  \\
	&=    \Big\la 0,\phi(x)-\beta_\phi(\mu) \Big\ra_\cH \\
	&=0.
	\end{align*}
	Hence $a_\mu(x)=\| \phi(x)-\beta_\phi(\mu) \|^2_\cH +\text{Var}(\phi(X))$ and
	\begin{align*}
	D(\mu)=\int a_\mu(x)\, d\mu(x) = 2\text{Var}(\phi(X)),
	\end{align*}
	proving $2)$, which implies $1)$ when inserting into $a_\mu(x)$. Lastly using the above proven equalities we get that 
	\begin{align*}
	d_\mu(x,y) &= d(x,y) - a_\mu(x)-a_\mu(y) + D(\mu) \\
	&= \|\phi(x)-\phi(y)\|^2_\cH - \| \phi(x)-\beta_\phi(\mu) \|^2_\cH - \| \phi(y)-\beta_\phi(\mu) \|^2_\cH.
	\end{align*}
	Furthermore for any $x,y,z\in \cH$, the linearity in the first argument and additivity in the second argument of $\la \cdot , \cdot \ra_\cH$ yield
	\begin{align*}
	\la x-z,y-z\ra_\cH &=\la x-z,x-z\ra_\cH +  \la x-z,y-x\ra_\cH \\
	&= \| x-z\|^2_\cH +  \la x-y,y-x\ra_\cH +\la y-z,y-x\ra_\cH \\
	&= \| x-z\|^2_\cH - \| x-y\|^2_\cH +\la y-z,y-z\ra_\cH+\la y-z,z-x\ra_\cH \\
	&=\| x-z\|^2_\cH - \| x-y\|^2_\cH + \| y-z \|^2_\cH - \la y-z , x-z\ra_\cH.
	\end{align*}
	As a consequence we have that \begin{align} \label{eq_temp_7648}
	 \| x-y\|^2_\cH -\| x-z\|^2_\cH -  \| y-z \|^2_\cH=-\la x-z,y-z\ra_\cH -\la y-z , x-z\ra_\cH,
	\end{align}
	allowing us to conclude that 
	\begin{align*}
	d_\mu(x,y) &= - \la \phi(x)-\beta_\phi(\mu),\phi(y)-\beta_\phi(\mu) \ra_\cH- \la \phi(y)-\beta_\phi(\mu) ,\phi(x)-\beta_\phi(\mu)\ra_\cH 
	\end{align*}
	and in the case that $\cH$ is a $\R$-Hilbert spaces we have that
	\begin{align*}
	d_\mu(x,y) = - 2\la \phi(x)-\beta_\phi(\mu),\phi(y)-\beta_\phi(\mu) \ra_\cH,
	\end{align*}
	which is what we wanted to prove.
\end{p}
Now we are ready to prove the alternative representation of $dcov(\theta)$ under the assumption that both $(\cX,d_\cX)$ and $(\cY,d_\cY)$ are metric spaces of negative type.

\begin{theorem} \label{lemma_beta_of_tensor_exists_and_representation_of_dcov}
	Let $(\cX,d_\cX)$ and $(\cY,d_\cY)$ have negative type as witnessed by the isometric embeddings $\phi:(\cX,d_\cX^{1/2})\to \cH_1$ and $\psi:(\cY,d_\cY^{1/2}) \to \cH_2$ into two separable Hilbert spaces with the same scalar field. If $\theta \in M_1^{1,1}(\cX\times \cY)$ have marginals $\mu\in M_1^1(\cX)$ and $\nu\in M_1^1(\cY)$, then it holds that
	\begin{align*}
	dcov(\theta)&= 4\|\beta_{\phi \otimes \psi}(\theta-\mu\times \nu) \|_{\cH_1\otimes \cH_2}^2 \\
	&=  4\Big\| \int \phi \otimes \psi \, d\theta - \int \phi \otimes \psi \, d\mu\times \nu  \Big\|_{\cH_1\otimes \cH_2}^2.
	\end{align*}
\end{theorem}
\begin{p}
	As in \cite{lyons2013distance} we show this for the case that $\cH_1$ and $\cH_2$ are separable $\R$-Hilbert spaces for simplicity. First note, by \cref{lemma_a_mu_and_D(mu)_and_d_mu_representations} identity $3)$ and the property of $\la \cdot , \cdot \ra_{\cH_1\otimes \cH_2}$ on simple tensors, we have that
	\begin{align*}
	\text{dcov}(\theta) &= \int d_\mu(x_1,x_2)d_\nu(y_1,y_2) \, d\theta \times \theta ((x_1,y_1),(x_2,y_2)) \\
	&= 4 \int \la \hat{\phi}(x_1),\hat{\phi}(x_2) \ra_{\cH_1} \la \hat{\psi}(y_1),\hat{\psi}(y_2) \ra_{\cH_2} \, d\theta \times \theta ((x_1,y_1),(x_2,y_2)) \\
	&= 4 \int \int \la \hat{\phi}(x_1) \otimes  \hat{\psi}(y_1) ,\hat{\phi}(x_2) \otimes \hat{\psi}(y_2) \ra_{\cH_1 \otimes \cH_2} \, d\theta (x_1,y_1) \, d\theta(x_2,y_2),
	\end{align*}
	where $\hat{\phi}(x)=\phi(x)-\beta_\phi(\mu)$ and $\hat{\psi}(y)=\psi(y)-\beta_\psi(\nu)$ for any $x\in\cX$ and $y\in\cY$. The mapping $\hat{\phi}\otimes \hat{\psi}:\cX\times \cY\to \cH_1\otimes \cH_2$ given by $\hat{\phi}\otimes \hat{\psi}(x,y)=\hat{\phi}(x)\otimes \hat{\psi}(y)$ is  Pettis integrable with respect to $\theta \in M^{1,1}(\cX\times \cY)$ so $\beta_{\hat{\phi}\otimes \hat{\psi}}:M^{1,1}(\cX\times \cY) \to \cH_1\otimes \cH_2$ is well-defined. This is seen by a simple replication of the arguments of \cref{lemma_tensor_product_map_of_isometric_embeddings_is_pettis_integrable} with the slight adjustment
	\begin{align*}
	\|\hat{\phi}(x)\otimes \hat{\psi}(y) \|_{\cH_1\otimes \cH_2} =& \| \phi(x)-\beta_\phi(\mu) \|_{\cH_1} \|\psi(y)-\beta_\psi(\nu)\|_{\cH_2} \\
	\leq&  \| \phi(x) \|_{\cH_1} \|\psi(y)\|_{\cH_2}+ \| \phi(x)\|_{\cH_1} \|\beta_\psi(\nu)\|_{\cH_2}\\
	&+ \| \beta_\phi(\mu) \|_{\cH_1} \|\psi(y)\|_{\cH_2}+ \| \beta_\phi(\mu) \|_{\cH_1} \|\beta_\psi(\nu)\|_{\cH_2},
	\end{align*}
	which is still integrable with respect to $\theta \in M^{1,1}(\cX\times \cY)$, since each term is.
	Using that the inner product $\la \cdot , \cdot \ra_{\cH_1 \otimes\cH_2}$ is a linear and continuous mapping when fixing one of its arguments (hence an element of $(\cH_1\otimes\cH_2)^*$), we get by the defining property of the Pettis integral that
	\begin{align*}
	\text{dcov}(\theta) &= 4 \int \bigg\la \int \hat{\phi} \otimes  \hat{\psi} \, d\theta  ,\hat{\phi}(x_2) \otimes \hat{\psi}(y_2) \bigg\ra_{\cH_1 \otimes \cH_2} \, d\theta(x_2,y_2) \\
	&=  4  \bigg\la \int \hat{\phi} \otimes  \hat{\psi} \, d\theta  , \int\hat{\phi} \otimes \hat{\psi}   \, d\theta  \bigg\ra_{\cH_1 \otimes \cH_2} \\
	&= 4 \bigg\| \int \hat{\phi} \otimes  \hat{\psi} \, d\theta  \bigg\|_{\cH_1 \otimes\cH_2 }^2 \\
	&=4 \|\beta_{\hat{\phi}\otimes \hat{\psi}}(\theta)\|_{\cH_1 \otimes \cH_2 }^2 .
	\end{align*}
	Now we simply need to show that $\beta_{\hat{\phi}\otimes \hat{\psi}}(\theta)$ has the wanted representation. To that end, we may expand $\beta_{\hat{\phi}\otimes \hat{\psi}}(\theta)$ is the following way
	\begin{align*}
	\beta_{\hat{\phi}\otimes \hat{\psi}}(\theta) &= \int (\phi-\beta_\phi(\mu)) \otimes  (\psi-\beta_\psi(\nu)) \, d\theta  \\
	&= \int \phi\otimes \phi - \phi\otimes \beta_\psi(\nu) - \beta_\phi(\mu)\otimes \psi + \beta_\phi(\mu)\otimes \beta_\psi(\nu) \, d\theta \\
	&= \int \phi\otimes \phi \,d\theta- \int \phi\otimes \beta_\psi(\nu) \, d\theta - \int\beta_\phi(\mu)\otimes \psi  \, d\theta + \int \beta_\phi(\mu)\otimes \beta_\psi(\nu) \, d\theta,
	\end{align*}
	by the linearity of the Pettis integral (cf. \cref{appendix_corollary_linearity_of_Pettis_integral}), since each Pettis integral exists (integrability follows from the same bounds as above). Now note that since $\cH_1$ and $\cH_2$ are Hilbert spaces we know that there exist orthonormal bases $\{\phi_i:i\in I\}$ and $\{\psi_j:j\in J\}$ for $\cH_1$ and $\cH_2$ respectively, where $I,J$ is at most infinitely countable (see \cref{remark_eigenvalues_of_S_and_limit_dist}). By \cref{appendix_theorem_orthonormal_basis_tensor_product} we have that $\{\phi_i\otimes \psi_j:i\in I,j\in J\}$ is an orthonormal basis for $\cH_1\otimes \cH_2$. As a consequence (cf. theorem 6.26 (2) \cite{applied_analysis_hunter2001}) we can expand elements of $\cH_1\otimes \cH_2$ in terms of this basis in the following way	\begin{align*}
	\int \phi \otimes \beta_\psi(\nu) \,  d\theta = \sum_{(i,j)\in I\times J} \Big\la \phi_i\otimes \psi_j, 	\int \phi \otimes \beta_\psi(\nu) \, d\theta \Big\ra_{\cH_1\otimes \cH_2} \phi_i\otimes \psi_j,
	\end{align*}
	where the equality is meant with respect to $\|\cdot\|_{\cH_1\otimes \cH_2}$-norm convergence. But for any $(i,j)\in I \times J$ we have by the unique defining property of the Pettis integral
	\begin{align*}
	\Big\la \phi_i\otimes \psi_j, 	\int \phi \otimes \beta_\psi(\nu) \, d\theta \Big\ra_{\cH_1\otimes \cH_2}&= \int	\la \phi_i\otimes \psi_j, 	 \phi(x) \otimes \beta_\psi(\nu) \ra_{\cH_1\otimes \cH_2}\, d\theta(x,y) \\
	&= \int	\la \phi_i , \phi(x) \ra_{\cH_1}  \la \psi_j, \beta_\psi(\nu) \ra_{\cH_2} \,d\mu(x) \\
	&=  \int	\la \phi_i , \phi(x) \ra_{\cH_1}  \,  d\mu(x) \la \psi_j, \beta_\psi(\nu) \ra_{\cH_2} \\
	&=  	\Big\la \phi_i , \int \phi \, d\mu \Big\ra_{\cH_1}   \la \psi_j, \beta_\psi(\nu) \ra_{\cH_2} \\
	&= \la \phi_i\otimes \psi_j , \beta_\phi(\mu)\otimes \beta_\psi(\nu) \ra_{\cH_1\otimes \cH_2}.
	\end{align*}
	Thus the (possibly countably infinite)  series converges in $\|\cdot\|_{\cH_1\otimes \cH_2}$-norm to both $\int \phi \otimes \beta_\psi(\nu) d\theta$  and $\beta_\phi(\mu)\otimes \beta_\psi(\nu)$, so we conclude that they coincide. That is
	\begin{align*}
	\int \phi \otimes \beta_\psi(\nu) \, d\theta &=  \beta_\phi(\mu)\otimes \beta_\psi(\nu),
	\end{align*}
	and by analogous arguments we also get that $\int \beta_\phi(\mu) \otimes \psi \, d\theta=  \beta_\phi(\mu)\otimes \beta_\psi(\nu)$. Finally for any $h\in\cH^*$ we have that
	\begin{align*}
	h^*\lp 	\int \beta_\phi(\mu)\otimes \beta_\psi(\nu) \, d\theta \rp &= \int h^* (\beta_\phi(\mu)\otimes \beta_\psi(\nu)) \, d\theta(x,y) \\
	&=h^* (\beta_\phi(\mu)\otimes \beta_\psi(\nu)),
	\end{align*}
	hence by the unique defining property of the Pettis integral we conclude that $\int \beta_\phi(\mu)\otimes \beta_\psi(\nu) \, d\theta =\beta_\phi(\mu)\otimes \beta_\psi(\nu) $. Thus
	\begin{align*}
	\beta_{\hat{\phi}\otimes \hat{\psi}}(\theta) &=  \int \phi\otimes \phi \,d\theta - \beta_\phi(\mu)\otimes \beta_\psi(\nu) \\
	&= \beta_{\phi\otimes \psi}(\theta) - \beta_\phi(\mu)\otimes \beta_\psi(\nu) \\
	&= \beta_{\phi\otimes \psi}(\theta) - \beta_{\phi \otimes \psi}(\mu\times\nu) \\
	&= \beta_{\phi\otimes \psi}(\theta-\mu\times \nu),
	\end{align*}
	where we used that $\beta_{\phi\otimes \psi}:M^1(\cX\times \cY)\to \cH_1\otimes \cH_2$ is linear (cf. \cref{lemma_beta_tesor_is_linear}) and that $\beta_{\psi\otimes\psi}(\mu\times\nu)=\beta_\phi(\mu)\otimes\beta_\psi(\nu)$, which is seen by an expansion in terms of the orthonormal basis. That is, 
	\begin{align*}
	\Big\la \phi_i\otimes \psi_j, \beta_{\phi\otimes \psi}(\mu\times \nu) \Big\ra_{\cH_1\otimes\cH_2} &=  \Big\la \phi_i\otimes \psi_j, \int \phi \otimes \psi \, d\mu\times\nu \Big\ra_{\cH_1\otimes\cH_2} \\
	&=  \int\int \la \phi_i\otimes \psi_j, \phi(x) \otimes \psi(y)  \ra_{\cH_1\otimes\cH_2} \, d\mu(x) \, d\nu(y)  \\
	&=  \int \la \phi_i ,\phi(x) \ra_{\cH_1} \, d\mu(x) \int \la \psi_j, \psi(y)  \ra_{\cH_2}  \, d\nu(y) \\
	&=  \Big\la \phi_i , \int  \phi  \, d\mu \Big\ra_{\cH_1}  \Big\la \psi_j,\int  \psi \, d\nu \Big\ra_{\cH_2}    \\
	&=  \la \phi_i\otimes \psi_j, \beta_\phi(\mu)\otimes\beta_\psi(\nu) \ra_{\cH_1\otimes\cH_2},
	\end{align*}
		by the unique defining property of the Pettis integral for any $(i,j)\in I \times J$, so they coincide by the same arguments as above. We conclude that
	\begin{align*}
	dcov(\theta) = 4 \|\beta_{\phi\otimes \psi}(\theta-\mu\times\nu)\|_{\cH_1 \otimes \cH_2 }^2.
	\end{align*}
\end{p}
	This marks an end of this section, as we proved the alternative representation of the distance covariance measure, we set out to find.
\newpage
\subsection{Metric spaces of strong negative type} \label{section_strong_negative_type}
 In the previous section we showed that, if $(\cX,d_\cX)$ and $(\cY,d_\cY)$ where of negative type, then
\begin{align*}
dcov(\theta) = 4 \|\beta_{\phi\otimes \psi}(\theta-\mu\times\nu)\|_{\cH_1 \otimes \cH_2 }^2,
\end{align*}
for some isometric embeddings $\phi:\cX\to \cH_1$ and $\psi:\cY\to\cH_2$. In this section we will define a subclass of negative type metric spaces called metric spaces of strong negative type. We will show that, if both  marginal spaces $(\cX,d_\cX)$ and $(\cY,d_\cY)$ are of this so-called strong negative type, then there exist isometric embeddings $\phi$ and $\psi$ into Hilbert spaces, such that the corresponding mean embedding maps are injective on a certain class of measures. This injective property is then used to prove that the linear mean embedding of the tensor product $\phi\otimes \psi$ given by $\beta_{\phi\otimes \psi}$ is injective on the whole of $M^{1,1}(\cX\times \cY)$. Now since  $\theta-\mu\times \nu\in M^{1,1}(\cX\times \cY)$ for any $\theta\in M^{1,1}_1(\cX\times \cY)$, we will argue that this injectivity together with the alternative representation of $dcov(\theta)$ stated above, yield the converse implication needed to use the distance covariance metric as a direct indicator of independence.  That is, if $(\cX,d_\cX)$ and $(\cY,d_\cY)$ are metric spaces of strong negative type, then we show that $dcov(\theta)=0$ if and only if $\theta=\mu\times \nu$ for any $\theta\in M^{1,1}_1(\cX\times \cY)$ with marginals $\mu$ and $\nu$.

In this section we will  furthermore show that, if at least one of the marginal spaces  $(\cX,d_\cX)$ and $(\cY,d_\cY)$ are not of strong negative type (and non-singleton spaces), then there exists a probability measure in $M^{1,1}_1(\cX\times \cY)$ for which $dcov(\theta)=0$ but $\theta\not = \mu\times \nu$. This renders $dcov$ unusable as a direct indicator of independence in metric spaces that are not of strong negative type. 

Lastly since the definition of a metric space of strong negative type is rather abstract and not easily recognizable, we will prove a theorem that identifies every separable Hilbert space as a metric space of strong negative type. That is, we prove that the distance covariance measure can directly determine whether random elements with values in two separable Hilbert spaces are independent or not.  \\ \\
The definition of metric spaces of strong negative type is given in terms of properties of $D:M^1(\cX)\to \R$, so before we continue we will present a lemma, connecting the previously analysed class of negative type metric spaces, with the mapping $D:M^1(\cX)\to\R$.
\begin{lemma} \label{lemma_strong_negative_type_D_leq_0}
	If $(\cX,d_\cX)$ is a metric space of negative type, then 
	\begin{align*}
	D(\mu_1-\mu_2)= \int d_\cX(x,y) \, d(\mu_1-\mu_2)\times (\mu_1-\mu_2)(x,y) \leq 0,
	\end{align*}
	for all $\mu_1,\mu_2\in M^1_1(\cX)$.
\end{lemma}
\begin{p}
	We initially note that $D(\mu_1-\mu_2)$ is well-defined for any $\mu_1,\mu_2\in M^1_1(\cX)\subset M^1(\cX)$, since $D$ is well-defined on $M^1(\cX)$, which is a vector space by \cref{lemma_M(X)_and_M^1(X)_are_vector_spaces}. \\ \\
	Let $\mu_1,\mu_2\in M_1^1(\cX)$ and let $(X_n)_{n\in\N}$ and $(Z_n)_{n\in\N}$ be two independent i.i.d. sequences of random elements in $\cX$ defined on a common probability space $(\Omega,\F,P)$ such that $X_1 \sim \mu_1$ and $Z_1\sim \mu_2$. This especially entails that $((X_n,Z_n))_{n\in \N}$ is an i.i.d. sequence of random elements in $\cX\times \cX$. 
	Fix $n\in \N$ and consider the  $n$-sample  empirical measures $\mu_1^{(n)},\mu_2^{(n)}:\Omega \to M_1^1{(\cX)}$ given by 
	$$
	\mu_1^{(n)}(\omega)= \frac{1}{n}\sum_{i=1}^n \delta_{X_i(\omega)} \quad \textit{ and } \quad \mu_2^{(n)}(\omega)=\frac{1}{n}\sum_{i=1}^n \delta_{Z_i(\omega)}.
	$$ 
	Let $K=(X_1,...,X_n,Z_1,...,Z_n)$ and note that for any $o\in \cX$
	\begin{align*}
	\int d_\cX(x,o) \, d(\mu_1^{(n)}(\omega)+\mu_2^{(n)}(\omega))(x)=	\frac{1}{n}\sum_{i=1}^{2n}d(K_i(\omega),o)  <\i,
	\end{align*}
	hence $\mu_1^{(n)}(\omega)-\mu_2^{(n)}(\omega)\in M^1(\cX)$ since $|\mu_1^{(n)}(\omega)-\mu_2^{(n)}(\omega)|\leq \mu_1^{(n)}(\omega)+\mu_2^{(n)}(\omega)$ for all $\omega \in \Omega$. Thus  \cref{lemma_metric_is_integral_product_of_M_with_finte_first_moment} yields that $d_\cX\in \mathcal{L}^1\lp (\mu_1^{(n)}-\mu_2^{(n)}) \times (\mu_1^{(n)}-\mu_2^{(n)})\rp
	$. Let $\beta=(1/n,...,1/n)\in \R^{n}$ and $\alpha=(\beta,-\beta)\in \R^{2n}$. Suppressing the $\omega$-notation, Fubini's theorem gives that
	\begin{align*}
	D(\mu_1^{(n)}-\mu_2^{(n)})&= \int d_\cX(x,y) \, d(\mu_1^{(n)}-\mu_2^{(n)})\times (\mu_1^{(n)}-\mu_2^{(n)})(x,y) \\
	&=\int \sum_{j=1}^{2n} \alpha_j d_\cX(y,K_j) \, d(\mu_1^{(n)}-\mu_2^{(n)})(y) \\
	&= \sum_{i=1}^{2n} \sum_{j=1}^{2n} \alpha_i \alpha_j d_\cX(K_i,K_j) \\
	&\leq 0,
	\end{align*}
	for all $\omega \in \Omega$, since $\cX$ is a metric space of negative type and $\sum_{i=1}^{2n} \alpha_i=0$. Now define $W_{1,n}=((X_1,Z_1),...,(X_n,Z_n))$ and note that
	\begin{align*}
	D(\mu_1^{(n)}-\mu_2^{(n)})&= \frac{1}{n^2}\sum_{i=1}^n \sum_{j=1}^n d_\cX(X_i,X_j)+d_\cX(Z_i,Z_j) - d_\cX(X_i,Z_j)-d_\cX(Z_i,X_j). \\
	&= \frac{1}{n^2}\sum_{i=1}^n \sum_{j=1}^n h((X_i,Z_i),(X_j,Z_j)) \\
	&= V_n^2(h,W_{1,n}),
	\end{align*}
	where $$
	V_n^2(h,W_{1,n})= \frac{1}{n^2} \sum_{i=1}^n \sum_{j=1}^n h(W_i,W_j),
	$$ is recognized as an $n$-sample V-statistic (see \cref{Appendix_V-statistics}) with symmetric kernel  $h:\cX^2 \times \cX^2 \to \R$ of degree $2$,  given by 
	$$
	h((x_1,z_1),(x_2,z_2)) = d_\cX(x_1,x_2)+d_\cX(z_1,z_2) - d_\cX(x_1,z_2)-d_\cX(z_1,x_2).
	$$ 
	The kernel $h$ is symmetric in the following sense 
	\begin{align*}
	h((x_1,z_1),(x_2,z_2))&=d_\cX(x_1,x_2)+d_\cX(z_1,z_2) - d_\cX(x_1,z_2)-d_\cX(z_1,x_2) \\
	&= d_\cX(x_2,x_1) +d_\cX(z_2,z_1) - d_\cX(x_2,z_1)- d_\cX(z_2,x_1) \\
	&= h((x_2,z_2),(x_1,z_1)).
	\end{align*}
	By the strong law of large numbers for $V$-statistics (\cref{theorem_SLLN_for_V-statistics}) we get that $$V_n^2(h,W_{1,n})\convas_n Eh((X_1,Z_1),(X_2,Z_2)),$$ if $E|h((X_1,Z_1),(X_2,Z_2))|<\i$ and $E|h((X_1,Z_1),(X_1,Z_1))|<\i$. To this end, note that 
	\begin{align*}
	E d_\cX (X_i,Z_j) &= \int d_\cX(x,z) \, d(X_i,Z_j)(P)(x,z) \\
	&= \int d_\cX(x,z) \, d\mu_1\times \mu_2 (x,z)\\
	&< \i ,
	\end{align*}
	for any $i,j\in\{1,2\}$  by \cref{lemma_metric_is_integral_product_of_M_with_finte_first_moment}. By the triangle inequality we get that
	\begin{align*}
	|h((x_1,z_1),(x_2,z_2))| &\leq d_\cX(x_1,x_2)+d_\cX(z_1,z_2) + d_\cX(x_1,z_2)+d_\cX(z_1,x_2), \\
	|h((x_1,z_1),(x_1,z_1))| &= 2d_\cX(x_1,z_1),
	\end{align*}
	for any $x_1,x_2,z_1,z_2 \in \cX$, hence for some $o\in\cX$ we have that
	\begin{align*}
	E|h((X_1,Z_1),(X_2,Z_2))| &\leq \int d_\cX \, d\mu_1^2 + \int d_\cX \, d\mu_2^2 +2\int d_\cX \, d\mu_1\times \mu_2 < \i,
	\\
	E|h((X_1,Z_1),(X_1,Z_1))|&= 2 \int d_\cX(x,y) \, d\mu_1\times \mu_2(x,y)
	< \i.	 
	\end{align*}
	We conclude that
	\begin{align*}
	D(\mu_1^{(n)}-\mu_2^{(n)})&=V_n^2(h,W_{1,n}) \\
	&\stackrel{a.s.}{\to}_n E[h((X_1,Z_1),(X_1,Z_2))]\\
	&= Ed_\cX(X_1,X_2)+Ed_\cX(Z_1,Z_2) - Ed_\cX(X_1,Z_2)-Ed_\cX(Z_1,X_2) \\
	&= \int d_\cX(x,y) \, d(\mu_1-\mu_2)\times (\mu_1-\mu_2)(x,y) \\
	&= D(\mu_1-\mu_2).
	\end{align*}
	Thus we have a non-positive sequence $(D(\mu_1^{(n)}(\omega)-\mu_2^{(n)}(\omega)))_{n\in\N}$, that for some $\omega\in\Omega$ converges to the constant $D(\mu_1-\mu_2)$, allowing us to conclude that 
	$$D(\mu_1-\mu_2)\leq 0,
	$$
	which proves the claim.
\end{p}
Now we state the definition of metric spaces of strong negative type.
\begin{definition}[Metric spaces of strong negative type] \label{defi_strong_negative_type}
	A metric space $(\cX,d_\cX)$ of negative type, is of strong negative type if
	\begin{align*}
	D(\mu_1-\mu_2) = 0 \iff \mu_1=\mu_2,
	\end{align*}
	for all $\mu_1,\mu_2\in M^1_1(\cX)$. 
	
\end{definition}

The next order of business is to prove an equivalence between the above definition and  injectivity of mean embedding maps $\beta_\phi:M^1(\cX)\to \cH$ on the subspace $M^1_1(\cX)\subset M^1(\cX)$.  Recall from \cref{theorem_equivalence_of_negative_type} that $(\cX,d_\cX)$ is a metric space of negative type if and only if there exist isometric embeddings into both a separable $\R$-Hilbert space and a separable $\bC$-Hilbert space.
\begin{lemma} \label{lemma_strong_negative_type_iff_beta_is_injective_on_prob_measures_of_finite_first_moment}
	The following statements are equivalent
	\begin{itemize}
		\item[1)] $(\cX,d_\cX)$ is a metric space of strong negative type.
		\item[2)] There exists an isometric embedding $\phi:(\cX,d_\cX^{1/2})\to \cH$ into a separable $\R$-Hilbert space, which induces a mean embedding map $\beta_\phi$ that is injective on $M^1_1(\cX)$.
		\item[3)] There exists an isometric embedding $\phi':(\cX,d_\cX^{1/2})\to \cH'$ into a separable $\bC$-Hilbert space, which induces a mean embedding map $\beta_{\phi'}$ that is injective on $M^1_1(\cX)$.
	\end{itemize} 
\end{lemma}
\begin{p}
	First note that, for any isometric embedding $\phi:(\cX,d_\cX^{1/2})\to \cH$ into a separable $\bK$-Hilbert space,  \cref{eq_temp_D(mu)_equality_2_normbeta(mu)} yields that
	\begin{align*}
	D(\mu) =  -2 \| \beta_\phi(\mu) \|_\cH^2,
	\end{align*}
	for any $\mu\in M^1_0(\cX)=\{\mu\in M^{1}(\cX):\mu(\cX)=0\}$. Now consider any two probability measures $\mu_1,\mu_2\in M^1_1(\cX)\subset M^1(\cX)$. By \cref{lemma_M(X)_and_M^1(X)_are_vector_spaces} we have that $M^1(\cX)$ is a vector space, so $\mu_1-\mu_2\in M^1(\cX)$. This signed measure  is obviously an element of $M^1_0(\cX)$, so the above applies, i.e.
	\begin{equation} \label{eq_temp_123321}
	D(\mu_1-\mu_2)=-2\|\beta_\phi(\mu_1-\mu_2)\|_\cH^2.
	\end{equation} By  \cref{lemma_beta_phi_is_linear} we have that $\beta_\phi$ is a linear  map on $M^1(\cX)$ (when viewing it as a map into the realification $\cH^\R$ if $\cH$ is a  $\bC$-Hilbert space), so $\beta_\phi(\mu_1-\mu_2) = \beta_\phi(\mu_1)-\beta_\phi(\mu_2)$. This allows us to  conclude that 
	\begin{align} \label{eq_temp_3223}
	D(\mu_1-\mu_2)=0 \iff \beta_\phi(\mu_1)=\beta_\phi(\mu_2).
	\end{align}
	\textit{1)}$\iff$\textit{2)}: Assume that $(\cX,d_\cX)$ is a metric space of strong negative type. Since $(\cX,d_\cX)$ is especially of negative type \cref{theorem_equivalence_of_negative_type} yields that there exists an isometric embedding $\phi:(\cX,d_\cX^{1/2})\to \cH$ into a $\R$-Hilbert space. By \cref{defi_strong_negative_type} and \cref{eq_temp_3223} we have  for any two probability measures $\mu_1,\mu_2\in M^1_1(\cX)$ that
	\begin{align*}
	\beta_\phi(\mu_1)=\beta_\phi(\mu_2) \iff D(\mu_1-\mu_2)=0 \iff \mu_1 = \mu_2,
	\end{align*}
	proving that $\beta_\phi$ is injective on $M^1_1 (\cX)$. Conversely assume that $\phi:(\cX,d_\cX^{1/2})\to \cH$ is an isometric embedding into a $\R$-Hilbert space such that $\beta_\phi$ is injective on $M^1_1 (\cX)$. Then for any two probability measures  $\mu_1,\mu_2\in M^1_1(\cX)$, \cref{eq_temp_3223} yields that
	\begin{align*}
	\mu_1=\mu_2 \iff \beta_\phi(\mu_1)=\beta_\phi(\mu_2) \iff D(\mu_1-\mu_2)=0,
	\end{align*}
	proving that $(\cX,d_\cX)$ is a metric space of strong negative type.\\ \\
	\textit{1)}$\iff$\textit{3)}: Start by invoking \cref{theorem_equivalence_of_negative_type} to get an isometric embedding $\phi':(\cX,d_\cX^{1/2})\to \cH'$ into a $\bC$-Hilbert space. Then the result follows by arguments identical to the previous equivalence.

\end{p}
This lemma is indeed very crucial for the following work. If we can identify a single isometric embedding into a separable $\bK$-Hilbert space, which induces a mean embedding map that is injective on $M^1_1(\cX)$, then $(\cX,d_\cX)$ is a metric space of strong negative type. This is the primary tool used, when we show that every separable Hilbert space is of strong negative type. In the case of finite-dimensional separable Hilbert spaces we identifying an isometric embedding into a separable $\bC$-Hilbert space $(L^2_\bC)$, that induces a mean embedding map that is injective on $M^1_1(\cX)$. Furthermore, infinite-dimensional separable Hilbert spaces are similarly shown to be of strong negative type; by identifying an isometric embedding into a separable $\R$-Hilbert space, which induces a mean embedding map that is injective on $M^1_1(\cX)$.

However, if we assume that $(\cX,d_\cX)$ is a metric space of strong negative type, then we know that there exists an isometric embedding $\phi:(\cX,d_\cX^{1/2})\to \cH$ into a separable $\R$-Hilbert space  that induces a mean embedding map $\beta_\phi$ that is injective on $M^1_1(\cX)$. This turns out to be very important in the following lemmas/theorems used to prove that, if both marginal metric spaces are of strong negative type, then $dcov(\theta)=0\iff\theta=\mu\times \nu$. To be perfectly clear on what is important about this: we only need to consider $\R$-Hilbert space valued isometric embeddings, a fact which greatly reduces to complexity of the following proofs (e.g. we do not have to consider Pettis integration of Hilbert space valued mappings with respect to complex measures.). \\ \\
Since we do not have to bother with embeddings into $\bC$-Hilbert spaces, the following statements are only concerned  with embeddings into  $\R$-Hilbert spaces, even though they may hold in both cases. The following lemma yields another domain on which our mean embeddings are injective.
\begin{lemma} \label{lemma_strong_negative_type_bet_phi_injective_on_M_0}
	If an isometric embedding $\phi:(\cX,d^{1/2})\to \cH$ into a separable $\R$-Hilbert space  induces a injective mean embedding map $\beta_\phi:M^1_1(\cX)\to \cH$, then $\beta_\phi:M^1_0(\cX)\to \cH$ is also injective
\end{lemma}
\begin{p}
	Note that $M^1_1(\cX)\cap M^1_0(\cX)=\emptyset$ and that the lemma does not state that $\beta_\phi$ is injective on  $M^1_1(\cX)\cup M^1_0(\cX)$, but that $\beta_\phi$ is injective when restricted the different domain $M^1_0(\cX)$. \\ \\
	Consider any isometric embedding $\phi:(\cX,d^{1/2})\to \cH$ into a separable $\R$-Hilbert space which induces a mean embedding map $\beta_\phi$ that is injective on $M^1_1(\cX)$. We obviously have that $M^1_0(\cX)$ is a $\R$-vector space, and since  $\beta_\phi:M^1_0(\cX)\to \cH$ is a linear map (\cref{lemma_beta_phi_is_linear}) it is injective if and only if the kernel $\mathrm{ker}(\beta_\phi)=\{\mu\in M^1_0(\cX):\beta_\phi(\mu)=0\}$ only contains the zero measure. Thus fix $\mu\in M^1_0(\cX)$ such that $\beta_\phi(\mu)=0$, and note that it suffices to show that $\mu=0$. 
	
	The Jordan-Hahn decomposition theorem yields that $\mu=\mu^+ - \mu^-$ where $\mu^+,\mu^-$ are non-negative finite singular measures on $(\cX,\cB(\cX))$, so it suffices to show that $\mu^+=\mu^-$.  Let $a:=\mu^+(\cX)=\mu^-(\cX)$. If $a=0$, then the non-negativity of $\mu^\pm$ implies that $\mu^+=\mu^-$. So it only remains to check the case where $a=\mu^\pm(\cX)>0$. We note that $	\mu^\pm_p:=a^{-1}\mu^\pm$,	are probability measures satisfying that $\mu^\pm= a\mu^\pm_p$. By the linearity of $\beta_\phi$ we have that
	\begin{align*}
0=	\beta_\phi(\mu)= \beta_\phi(a\mu^+_p - a\mu^-_p) = a\beta_\phi(\mu^+_p) -a\beta_\phi( \mu^-_p) \iff \beta_\phi(\mu^+_p)=\beta_\phi( \mu^-_p).
	\end{align*}
	Since $\beta_\phi$ is injective on probability measures, we conclude that $\mu^+_p=\mu^-_p$. As a consequence $\mu^+=\mu^-$, proving that $\beta_\phi$ is injective on $M^1_0(\cX)$.
\end{p}
This above lemma can now be used in conjunction with \cref{lemma_strong_negative_type_iff_beta_is_injective_on_prob_measures_of_finite_first_moment} to prove the following result.
\begin{theorem} \label{theorem_strong_negative_type_implies_existence_of_phi_such_that_beta_phi_is_injective_on_entire_domain}
	If $(\cX,d_\cX)$ is a metric space of strong negative type, then there exists an isometric embedding $\phi:(\cX,d_\cX^{1/2})\to \cH$ into a separable $\R$-Hilbert space which induces a mean embedding map $\beta_\phi$ that is injective on the whole domain $M^1(\cX)$. This isometric embedding might be different from the one guaranteed to exist by \cref{theorem_equivalence_of_negative_type}.
\end{theorem}
\begin{p} 	Assume that $(\cX,d_\cX)$ is a metric space of strong negative type. By \cref{lemma_strong_negative_type_iff_beta_is_injective_on_prob_measures_of_finite_first_moment} we know that there exists an isometric embedding $\phi:(\cX,d_\cX^{1/2})\to \cH$ into a separable $\R$-Hilbert space, which induces a mean embedding map $\beta_\phi$ that is injective on $M^1_1(\cX)$. By \cref{lemma_strong_negative_type_bet_phi_injective_on_M_0} we get that the mean embedding map $\beta_\phi$  is injective on  $M^1_0(\cX)$. \\ \\
	If  $\beta_\phi$ is injective on the whole of $M^1(\cX)$ we are done. On the other hand, if $\beta_\phi$ is not injective on the whole on $M^1(\cX)$, then we can construct another isometric embedding into a different $\R$-Hilbert space inducing a mean embedding map, which is. \\ \\
	To see how this is done we assume that $\beta_\phi$ is not injective on $M^1(\cX)$.  Since $\beta_\phi:M^1(\cX)\to \cH$ is linear and linear maps are injective if and only if the kernel only contains the zero element (in our case the zero measure), we may conclude that $\mathrm{ker}(\beta_\phi)=\{\mu\in M^1(\cX):\beta_\phi(\mu)=0\}\not =\{0\}$. That is, there exists at least one non-zero measure $\mu\in M^1(\cX)$ such that $\beta_\phi(\mu)=0$. 
	
	Next we will realize that every measure in $\mathrm{ker}(\beta_\phi)$ has a distinct measurement on the entire space $\cX$. That is, for any two distinct finite signed measures $\mu_1,\mu\in \mathrm{ker}(\beta_\phi)$, it holds that $\mu_1(\cX)\not = \mu_2(\cX)$. To see this, assume for contradiction that there are two distinct finite signed measures $\mu_1,\mu_2\in \mathrm{ker}(\beta_\phi)$ with $\mu_1(\cX)=\mu_2(\cX)$. But note that $\mu_1-\mu_2\in M^1_0(\cX)$ such that
	\begin{align*}
	\beta_\phi(\mu_1)-\beta_\phi(\mu_2)=0 \iff \beta_\phi(\mu_1-\mu_2)=0\iff \mu_1-\mu_2=0 \iff \mu_1=\mu_2,
	\end{align*}
	by linearity and the injectivity of $\beta_\phi$ on $M^1_0(\cX)$, a contradiction. \\ \\
	Now we turn our attention to the  construction of the $\R$-Hilbert space $\cH'$ and isometric embedding $\phi'$ into $\cH'$, which induces a mean embedding map, that is injective on the whole of $M^1(\cX)$. Let $\cH'=\cH\oplus \R$ be the direct sum of $\cH$ and $\R$. This space is an $\R$-Hilbert space given by the Cartesian product of $\cH$ and $\R$ with addition and scalar multiplication operations working  coordinate-wise. That is,
	\begin{align*}
	\cH \oplus \R = \{ x \oplus y : x\in \cH, y\in \R\}=\cH\times \R,
	\end{align*}
	with addition and scalar multiplication defined by $(x\oplus y)+ (x'\oplus y') = (x+x') \oplus (y+y')$ and $ a(x\oplus y) = (ax)\oplus(ay)$, $a\in\R$.	Furthermore we define an inner product on $\cH\oplus \R$,  $\la \cdot , \cdot \ra:\cH\oplus \R \times \cH \oplus \R \to \R$ given by $	\la x\oplus y , x' \oplus y' \ra = \la x,x' \ra_{\cH} + \la y,y' \ra_\R$. 	This is easily seen to be an inner product: symmetry, linearity in the first argument and positive-definiteness are all inherited by the same properties of the marginal inner product spaces. Lastly, $\cH\oplus \R$ is also complete with respect to the naturally induced metric (cf. section 1.6 \cite{conway1990ACourseInFunctionalAnalysis_reference_for_completion_inner_product_space}) rendering it a $\R$-Hilbert space.
	
	Separability of $\cH\oplus \R$ is realized by noting that, if $D_1$ and $D_2$ are countable dense sets in $\cH$ and $\R$ respectively, then $D_1\times  D_2$ is countable and dense in $\cH\oplus \R$. More specifically for any $x\otimes y\in \cH\otimes \R$ there exist sequences $(x_n)\subset D_1$ with $x_n\to_n x$ and $(y_n)\subset D_2$ with $y_n\to_n y$, hence the sequence $(x_n\oplus y_n)\subset D_1 \times D_2$ satisfies
	$\|x_n\oplus y_n - x\oplus y\|_{\cH\oplus \R} = \|(x_n-x)\oplus (y_n-y)\|_{\cH\oplus \R}= \sqrt{\|x_n-x\|_{\cH}^2+\|y_n-y\|_{\R}^2}\to_n 0$, proving that $D_1 \times D_2$ is dense in $\cH\oplus \R$. \\ \\
	Now we construct a candidate for $\phi':\cX \to \cH'$ and subsequently show that it is indeed an isometric embedding of $(\cX,d_\cX^{1/2})$ into $\cH'$, which induces an injective mean embedding map $\beta_{\phi'}:M^1(\cX)\to \cH'$. 
	
	Let $\phi':\cX\to \cH'$ be given by $	\phi'(x)=\phi(x)\oplus 1$, and note that for any $x,y\in \cX$
	\begin{align*}
	d_\cX(x,y)&= \|\phi(x)-\phi(y)\|_\cH^2 \\
	&= \|(\phi(x)-\phi(y))\oplus 0\|_{\cH'}^2  \\
	&= \| \phi(x)\oplus 1-\phi(y)\oplus 1 \|_{\cH'} \\
	&= \| \phi'(x)-\phi'(y) \|_{\cH'}^2,
	\end{align*}
	proving that $\phi':(\cX,d^{1/2})\to\cH'$ is an isometric embedding into an $\R$-Hilbert space.   Now note that the linear map $\beta_{\phi'}:M^1(\cX)\to  \cH'$ given by the Pettis integral of $\phi'$ with respect to the argument, satisfies that $\beta_{\phi'}(\mu)$ is the  unique element in $\cH'$ such that
	\begin{align*}
    h^*(\beta_{\phi'}(\mu)) = \int h^* \circ \phi' \, d\mu, \quad \quad \forall h^*\in (\cH')^*.
	\end{align*}
	Fix $h^*\in \cH'^*:=(\cH')^*$ and note that Riesz's representation theorem yields that there exists a unique $(x'\oplus y')\in \cH'$ such that $h^*( x\oplus y)=\la x\oplus y, x'\oplus y' \ra_{\cH'}$, hence
	\begin{align*}
	\int h^* \circ \phi' \, d\mu &= \int \la \phi'(x),x'\oplus y' \ra_{\cH'} \, d\mu(x) \\
	&= \int \la \phi(x)\oplus 1, x'\oplus y' \ra_{\cH'} \, d\mu(x)\\
	&= \int \la \phi(x),x' \ra_\cH + \la 1,y' \ra_\R \, d\mu(x) \\
	&=  \Big\la\int \phi \, d\mu,x' \Big\ra_\cH  + \mu(\cX) \la 1,y' \ra_\R \\
	&=  \la\beta_\phi(\mu),x' \ra_\cH  +  \la \mu(\cX),y' \ra_\R \\
	&=   \la \beta_\phi(\mu) \oplus \mu(\cX) , x' \oplus y' \ra_{\cH'} \\
	&= h^*(\beta_\phi(\mu) \oplus \mu(\cX) ).
	\end{align*}
	Thus 
	\begin{align*}
	\beta_{\phi'}(\mu) = \beta_\phi(\mu) \oplus \mu(\cX). 
	\end{align*}
	By the linearity of $\beta_{\phi'}:M^1(\cX)\to \cH'$ we have that it is injective if and only if the kernel $\mathrm{ker}(\beta_{\phi'})$ only contains the zero element of $M^1(\cX)$, i.e. the zero measure. But we note that the zero element of $\cH'=\cH\oplus \R$ is $0\oplus 0$, hence
	\begin{align*}
	\mathrm{ker}(\beta_{\phi'}) &= \{ \mu\in M^1(\cX):\beta_{\phi'}(\mu)=0\oplus 0 \} \\
	&= \{ \mu\in M^1(\cX):\beta_{\phi}(\mu)=0 \text{ and }\mu(\cX)=0 \} \\
	&=\{ \mu\in M^1(\cX):\mu\in\mathrm{ker}(\beta_\phi) \text{ and } \mu(\cX)=0 \}.
		\end{align*}
 We previously showed that every element of $\mathrm{ker}(\beta_\phi)$ has a distinct measure on the entire space. Therefore only one measure in  $\mathrm{ker}(\beta_\phi)$ measures the entire space $\cX$ to zero. We easily realize that the only measure in $\mathrm{ker}(\beta_\phi)$ that assigns the entire space to zero is the zero measure. That is, $\mathrm{ker}(\beta_{\phi'})=\{0\}$, proving that $\beta_{\phi'}:M^1(\cX)\to \cH'$ is an injective mapping. We conclude that there exists an $\R$-Hilbert space isometric embedding $\phi':(\cX,d_\cX^{1/2})\to \cH'$, which induces a mean embedding map $\beta_{\phi'}$ that is injective on the whole of $M^1(\cX)$.
\end{p}
Hence we have that, if both $(\cX,d_\cX)$ and $(\cY,d_\cY)$ are metric spaces of strong negative type, then we know that there exist two isometric embeddings $\phi:\cX\to \cH_1$ and $\psi:\cY\to \cH_2$ into two separable $\R$-Hilbert spaces that induce mean embeddings $\beta_\phi:M^1(\cX)\to\cH_1$ and $\beta_\psi:M^1(\cX)\to\cH_2$ that are injective. \Cref{lemma_strong_negative_type_beta_tensor_injective} below will furthermore show that, if this is indeed the case, then the tensor mean embedding $\beta_{\phi\otimes \psi}:M^{1,1}(\cX\times \cY)\to \cH_1\otimes \cH_2$ is injective. From here it is easy to prove that $dcov(\theta)=0\iff \theta=\mu\times \nu$, using the alternative representation of the distance covariance. \\ \\
The proof of \cref{lemma_strong_negative_type_beta_tensor_injective} will utilize a specific continuous linear map, hence we start by proving that such a map is indeed unique and well-defined.

\begin{lemma} \label{lemma_existence_of_continuous_linear_map_T_h}
	Let $\cH_1\otimes \cH_2$ be the tensor product of two separable $\R$-Hilbert spaces. For any $h\in \cH_1$, there exists a unique continuous and linear map $T_h:\cH_1\otimes \cH_2\to \cH_2$, such that
	\begin{align*}
	T_h (h_1\otimes h_2) = \la h_1,h\ra_{\cH_1} h_2,
	\end{align*}
	for any $h_1\in \cH_1$ and $h_2\in \cH_2$.
\end{lemma}
\begin{p}
	Since $\cH_1$ and $\cH_2$ are separable Hilbert spaces we know that there exist two orthonormal bases $\{e_{1,i}:i\in I\}$ and $\{e_{2,j}:j\in J\}$ for $\cH_1$ and $\cH_2$ respectively, where $I$ and $J$ are two at most infinitely countable index sets (see \cref{remark_eigenvalues_of_S_and_limit_dist}). By \cref{appendix_theorem_orthonormal_basis_tensor_product}, we have that an orthonormal basis for $\cH_1\otimes \cH_2$ is given by $B=\{e_{1,i}\otimes e_{2,j}:(i,j)\in I \times J\}$. Now for any $h\in \cH_1$, we define the mapping $T_{h}:\mathrm{span}(B)\to \cH_2$ by letting 
	\begin{align*}
	T_{h}(e_{1,i}\otimes e_{2,j}) = \la e_{1,i} ,h \ra_{\cH_1} e_{2,j},
	\end{align*}
	and extending by linearity to $\mathrm{span}(B)$ (finite linear combinations of elements of $B$). That is,  for any $v\in\mathrm{span}(B)$, there exist $n,m\in\N$ and $I_1\subset I$, $J_1\subset J$ with $|I_1|=n$, $|J_1|=m$ and a family of $\R$-scalars $\{a_{i,j}\}$ such that $v = \sum_{(i,j)\in I_1\times J_1} a_{i,j} e_{1,i} \otimes e_{2,j}$, and then let
	\begin{align*}
	T_{h}(v) = \sum_{(i,j)\in I_1\times J_1} a_{i,j} T_{h}(e_{1,i} \otimes e_{2,j}) = \sum_{(i,j)\in I_1\times J_1} a_{i,j} \la e_{1,i}, h \ra_{\cH_1} e_{2,j}.
	\end{align*}
	We equip $\mathrm{span}(B)$  with the inner product $\la \cdot , \cdot \ra$ (and induced norm $\| \cdot \|$) given by the restriction of $\la\cdot , \cdot \ra_{\cH_1\otimes \cH_2}$ to $\mathrm{span}(B)$. Now recall that for a finite sum of mutually orthogonal elements it holds that $	\| \sum_{i=1}^n u_i  \|^2 = \sum_{i=1}^n \|u_i\|^2$, hence
	\begin{align*}
	\|T_{h}(v)\|_{\cH_2}^2&=	\Big\| \sum_{(i,j)\in I_1\times J_1} a_{i,j} \la e_{1,i}, h \ra_{\cH_1} e_{2,j} \Big\|_{\cH_2}^2 = \Big\| \sum_{j\in J_1} e_{2,j}\Big(\sum_{i\in I_1} a_{i,j} \la e_{1,i}, h \ra_{\cH_1} \Big)
	\Big\|_{\cH_2}^2 \\
	&= \sum_{j\in J_1} \| e_{2,j} \|_{\cH_2}^2 \Big|\sum_{i\in I_1} a_{i,j} \la e_{1,i}, h \ra_{\cH_1} 
	\Big|^2 \leq  \sum_{j\in J_1} \lp \sum_{i\in I_1} a_{i,j}^2 \rp \lp \sum_{i\in I_1}  \la e_{1,i}, h \ra_{\cH_1}^2 \rp ,
	\end{align*}
	by Cauchy-Schwarz's inequality. 	Each of the latter factors can be bounded from above by $\|h\|_{\cH_1}$; by noting that $\sum_{l\in I_1} \la e_{1,i} , h \ra_{\cH_1}^2\leq  \sum_{l\in I} \la e_{1,l} , h \ra_{\cH_1} ^2=\|h\|_{\cH_1}^2$, by Parseval's identity since $\{e_{1,i}:i\in I\}$ is an orthonormal basis for $\cH_1$. Thus
	\begin{align*}
	\|T_{h}(v)\|_{\cH_2}^2 &\leq  \|h\|_{\cH_1}^2\sum_{(i,j)\in I_1\times J_1} a_{i,j}^2  = \|h\|_{\cH_1}^2\sum_{(i,j)\in I_1\times J_1} \|a_{i,j}e_{1,i}\otimes e_{2,j}\|^2 \\
	&= \|h\|_{\cH_1}^2 \Big\|\sum_{(i,j)\in I_1\times J_1} a_{i,j}e_{1,i}\otimes e_{2,j}\Big\|^2 = \|h\|_{\cH_1}^2 \|v\|^2,
	\end{align*}
	where we used that $e_{1,i}\otimes e_{2,j}\perp e_{1,l}\otimes e_{2,p}$ for any $(i,j)\not = (l,p)$ and that $e_{1,i}\otimes e_{2,j}$ are normalized. This proves that $T_h$ is a bounded linear map. By the bounded linear transformation (BLT) theorem (see theorem 5.19 \cite{applied_analysis_hunter2001}) there exists a unique bounded linear extension of $T_h$, to the closure of the original domain $\overline{\mathrm{span}(B)}$ (seen as a subset of $\cH_1\otimes \cH_2$). Since  $B$ is a basis for $\cH_1\otimes \cH_2$, we know that $\overline{\mathrm{span}(B)}=\cH_1\otimes \cH_2$. Hence the theorem gives a unique bounded linear map $\overline{T_h}:\cH_1\otimes \cH_2\to \cH_2$, with the property that
	\begin{align*}
	\overline{T_h}(v)=T_h(v),
	\end{align*}
	for all $v\in\mathrm{span}(B)$. Assume without loss of generality, that $\cH_1$ and $\cH_2$ are infinite-dimensional, such that the respective bases are countably infinite and enumerated by the natural numbers. Any $h_1\in\cH_1$ and $h_2\in\cH_2$, may be expanded in terms of the bases $h_1=\sum_{i=1}^\i \la e_{1,i} , h_1 \ra_{\cH_1} e_{1,i}$ and $h_2=\sum_{j=1}^\i \la e_{2,j} , h_2 \ra_{\cH_2} e_{2,j}$. Thus we realize (see \cref{Appendix_Tensor_product_of_Hilbert_spaces}) that \begin{align*}
	h_1\otimes h_2= \sum_{i=1}^\i \sum_{j=1}^\i \la e_{1,i} , h_1 \ra_{\cH_1} \la e_{2,j} , h_2 \ra_{\cH_2} e_{1,i}\otimes e_{2,j}.
	\end{align*}
	The linearity and continuity of $\overline{T_h}$ yield that
	\begin{align*}
	\overline{T_h}\lp h_1\otimes h_2\rp &= \sum_{i=1}^\i \sum_{j=1}^\i \overline{T_h}\lp   (\la e_{1,i} , h_1 \ra_{\cH_1}e_{1,i})\otimes ( \la e_{2,j} , h_2 \ra_{\cH_2}e_{2,j})\rp \\
	&= \sum_{i=1}^\i \sum_{j=1}^\i T_h\lp   (\la e_{1,i} , h_1 \ra_{\cH_1}e_{1,i})\otimes ( \la e_{2,j} , h_2 \ra_{\cH_2}e_{2,j})\rp \\
	&= \sum_{i=1}^\i \sum_{j=1}^\i \la \la e_{1,i} , h_1 \ra_{\cH_1}e_{1,i} ,h \ra_{\cH_1} \la e_{2,j} , h_2 \ra_{\cH_2}e_{2,j} \\
	&=   \Big\la \sum_{i=1}^\i\la e_{1,i} , h_1 \ra_{\cH_1}e_{1,i} ,h \Big\ra_{\cH_1}  \sum_{j=1}^\i \la  e_{2,j} , h_2 \ra_{\cH_2}e_{2,j} \\
	&= \la h_1 , h \ra_{\cH_1} h_2,
	\end{align*}
	which is what we wanted to show.
\end{p}
\newpage 
\begin{lemma} \label{lemma_strong_negative_type_beta_tensor_injective}
	If $(\cX,d_\cX)$ and $(\cY,d_\cY)$ are metric spaces of strong negative type, then there exist two isometric embeddings $\phi:(\cX,d_\cX^{1/2})\to \cH_1$ and $\psi:(\cY,d_\cY^{1/2})\to \cH_2$ into two $\R$-Hilbert spaces, such that the mean embedding of the tensor map $\beta_{\phi \otimes \psi}:M^{1,1}(\cX\times \cY)\to\cH_1\otimes \cH_2$ is injective.
\end{lemma}
\begin{p}
	Since $(\cY,d_\cY)$ is a metric space of strong negative type we invoke \cref{theorem_strong_negative_type_implies_existence_of_phi_such_that_beta_phi_is_injective_on_entire_domain} to get an isometric embedding $\psi:(\cY,d_\cY^{1/2})\to \cH_2$ into a separable  $\R$-Hilbert space, which induces a mean embedding $\beta_{\psi}:M^1(\cY)\to \cH_2$ that is injective. \\ \\ 
	Furthermore since $(\cX,d_\cX)$ is a metric space of strong negative type we invoke  \cref{theorem_equivalence_of_negative_type} to say that there exists an isometric embedding $\phi'':(\cX,d_\cX^{1/2}) \to \cH_1'$  into  a separable  $\R$-Hilbert space. Now translate this embedding, such that the translation has the origin of $\cH_1'$ in its image. That is define a new embedding $\phi':x\mapsto \phi''(x)-\phi''(o)\in \cH_1'$ for some $o\in \cX$. This is obviously still an isometric embedding in the same fashion as $\phi''$ is, but it indeed contains the origin of $\cH_1'$ in its images.
	
	Now since the $(\cX,d_\cX)$ is of strong negative type  \cref{lemma_strong_negative_type_iff_beta_is_injective_on_prob_measures_of_finite_first_moment} yields that there exists an isometric embedding $\phi''':\cX\to \cH_1'''$  with mean embedding map that is injective on $M^1_1(\cX)$. Now note that, for any $\mu_1,\mu_2\in M^1_1(\cX)$ \cref{eq_temp_123321} gives that
	\begin{align*}
-2\|\beta_{\phi'}(\mu_1)-\beta_{\phi'}(\mu_2)\|_{\cH'_1}=	D(\mu_1-\mu_2)=-2\|\beta_{\phi'''}(\mu_1)-\beta_{\phi'''}(\mu_2)\|_{\cH'''_1},
	\end{align*}
	yielding that $\beta_{\phi'}:M^1_1(\cX)\to \cH_1'$ is injective, since $\beta_{\phi'''}:M^1_1(\cX)\to \cH_1'''$ is (note that this actually proves that if a metric space is of strong negative type, then all mean embeddings of isometric embeddings are injective in this manner). \Cref{lemma_strong_negative_type_bet_phi_injective_on_M_0} now yield that  $\beta_{\phi'}:M^1_0(\cX)\to \cH_1'$ is injective. By the proof of \cref{theorem_strong_negative_type_implies_existence_of_phi_such_that_beta_phi_is_injective_on_entire_domain} we get, that $\phi:x \mapsto \phi'(x)\oplus 1$ is an isometric embedding into $\cH_1:=\cH_1'\oplus \R$, which induces a mean embedding $\beta_\phi:M^1(\cX)\to \cH_1$ that is injective. \\ \\
	This explicit construction is necessary, since we later need that $\phi$ is constructed in the above manner in terms of $\phi'$, which contains the origin of $\cH_1'$ in its image. More specifically we will use that
	\begin{align*}
	\la \phi(x) , \phi(z) \ra_{\cH_1}  &= \la \phi'(x)\oplus 1 , \phi'(z) \oplus 1\ra_{\cH_1} \\ &= \la \phi'(x), \phi'(z) \ra_{\cH_1'} + \la 1,1\ra_{\R} \\
	&= \la \phi'(x), \phi'(z) \ra_{\cH_1'} + 1,
	\end{align*}
	and that $\phi'(\cX)$ contains the origin of $\cH_1'$. \\ \\ Since $\beta_{\phi\otimes \psi}:M^{1,1}(\cX\times \cY)\to \cH_1\otimes \cH_2$ is linear it suffices to show that $\mathrm{ker}(\beta_{\phi\otimes \psi})=\{0\}$. Hence consider any $\theta\in \mathrm{ker}(\beta_{\phi\otimes \psi})$. That is, a $\theta\in M^{1,1}(\cX\times\cY)$ with $\beta_{\phi\otimes \psi}(\theta)=0$ and note that we have to show that $\theta=0$. \newpage
	\textbf{Step 1): Reformulating the problem.} \\It suffices to show that
	\begin{align*}
	\theta(A\times B)= 0, \quad \quad \forall A\in\cB(\cX),B\in\cB(\cY),
	\end{align*}
	since such sets form an intersection stable generator for $\cB(\cX)\otimes \cB(\cY)$. Fix an arbitrary set $B\in\cB(\cY)$ and define the finite signed measure (cf. \cref{lemma_help_injectivity_proof}) $\mu_B:\cB(\cX)\to\R$ by
	\begin{align*}
	\mu_B(A):=\theta(A\times B) = \int 1_{A}(x)1_{B}(y) \, d\theta(x,y),
	\end{align*}
	for any $A\in\cB(\cX)$. The problem has now been reduced to showing that $\mu_B=0$.
	
	 For any $A\in \cB(\cX)$, it holds that $\mu_B(A)=\theta^+(A\times B)-\theta^-(A\times B)$, hence
	\begin{align*}
	|\mu_B|(A) &\leq  \theta^+(A\times B) + \theta^-(A\times B) = |\theta|(A\times B) \leq |\theta|(A\times \cY) = \pi_1(|\theta|)(A).
	\end{align*}
	Now note that $\theta\in M^{1,1}(\cX\times \cY)$, which by definition means that $\pi_1(|\theta|)\in M^{1}(\cX)$, so the above inequality also yields that $\mu_B\in M^1(\cX)$. We especially know $\phi$ is Pettis integrable with respect to such measures (see \cref{lemma_isometric_embeddings_are_pettis_integrable_wrt_first_moment_measures}), meaning that $h^*\circ \phi \in\cL^1(\mu_B)$ for all $h^*_1\in\cH_1^*$. We also note that $(x,y) \mapsto \phi(x)1_B(y)$ is Pettis integrable with respect to $\theta$, since it is jointly measurable and 
	$$
		\int |h_1^* (\phi(x)1_B(y))| \, d|\theta|(x,y) \leq \int |h_1^* \circ \phi(x)| \, d|\theta|(x,y) = \int |h_1^*\circ \phi(x)|\, d\pi_1(|\theta|)(x) <\i,
	$$
	 for any $h_1^*\in \cH_1^*$.
	Hence the Pettis integrals $\int \phi(x)1_B(y) \,d\theta(x,y)$ and $ \int \phi \, d\mu_B$ 
	 exist, but we also have that $h_1^*\circ \phi(x) 1_B(y) \in \cL^1(\theta)$, so  
	$$
		\int h_1^*\circ \phi(x) \, d\mu_B(x) = \int h_1^*\circ \phi(x) 1_B(y) \, d\theta(x,y)= \int  h_1^* (\phi(x) 1_B(y) ) \, d\theta(x,y),
		$$
	for any $h_1^*\in\cH^*_1$; by \cref{lemma_help_injectivity_proof}. By the unique defining property of the Pettis integral we get that 
	$$
	\beta_\phi(\mu_B)= \int \phi \, d\mu_B = \int \phi(x)1_B(y) \,d\theta(x,y).
	$$
	The linearity and injectivity of $\beta_\phi:M^1(\cX)\to \cH_1$, yield that $\mu_B=0$ if and only if $\beta_\phi(\mu_B)=0$, which by the unique defining property of the Pettis integral and the above equality happens if and only if \begin{align*}
	 h_1^*(\beta_\phi(\mu_B))=\int h_1^*(\phi(x)1_B(y)) \,d\theta(x,y)=0,
	 	\end{align*}
	 for all $h_1^*\in\cH^*_1$. In the remainder of this proof we use the Riesz representation theorem to uniquely connect every $h_1\in \cH_1$ with $h_1^*\in \cH^*_1$, by the identity $h^*_1(x)=\la x, h_1\ra_{\cH_1}$ for all $x\in \cX$. \\ \\
	 Now for any $h_1\in \cH_1$, we define the finite signed measure (\cref{lemma_help_injectivity_proof})  $\nu_{h_1}:\cB(\cY)\to\R$ by
	\begin{align*}
	\nu_{h_1}(B) :=\int \la \phi(x),h_1 \ra_{\cH_1} 1_B(y) \, d\theta(x,y) =\int h_1^*(\phi(x)1_B(y)) \,d\theta(x,y),
	\end{align*}
	for any $B\in \cB(\cY)$. \textit{Note that $\nu_{h_1}$ does not necessarily have finite first moment for all $h_1\in \cH_1$, rendering the original proof in \cite{lyons2013distance} incorrect}. We realize that, it suffices to show that $\nu_{h}=0$ for all $h\in \cH_1$. \\ \\
	\textbf{Step 2.1): Proving that $\mathbf{\nu_{h}=0}$ for all $\mathbf{h\in \overline{\mathrm{span}(\mathrm{Im}(\phi))}}$.}\\
	Recall the continuous linear map $T_{h}:\cH_1\otimes \cH_2\to \cH_2$ from \cref{lemma_existence_of_continuous_linear_map_T_h} and that $\phi\otimes \psi:(x,y)\mapsto \phi(x)\otimes \psi(y)\in \cH_1\otimes \cH_2$ is Pettis integrable with respect to $\theta$. \Cref{appendix_lemma_linear_map_inside_pettis_integral_also_pettis_integrable} yields that the composition  
	$T_{h_1}\circ \phi\otimes \phi:\cX\times \cY\to \cH_2$ is Pettis integrable with respect to $\theta$, and
	\begin{align*}
	\int \la \phi(x),h_1 \ra_{\cH_1}\psi(y) \, d\theta(x,y)=\int T_{h_1}( \phi\otimes \phi ) \, d\theta = T_{h_1} \lp  \int \phi \otimes \psi \, d\theta \rp = T_{h_1}(\beta_{\phi\otimes \psi }(\theta)) =0,
	\end{align*}
	for any $h_1\in\cH_1$, by the linearity of $T_h$ and the assumption that $\beta_{\phi\otimes \psi}(\theta)=0$. \\ \\
	Fix $h_1\in \mathrm{Im}(\phi)\subset \cH_1$, where $\mathrm{Im}(\phi)=\phi(\cX)$ is the image of $\phi$. and note that there exists a $z\in \cX$ such that $\phi(z)=h_1$.  Recall from \cref{eq_temp_7648}, that
	\begin{align} \label{temp_eq_21321323}
	 \| x-y\|^2_{\cH_1'} -\| x-w\|^2_{\cH_1'} -  \| y-w \|^2_{\cH_1'}=-2\la x-w,y-w\ra_{\cH_1'},
	\end{align}
	for any $x,y,w\in \cH_1'$. Let $w\in \cH_1$, $x-w=\phi(x)$ and $y-w=\phi(z)$. \Cref{temp_eq_21321323} together with the initial construction of $\phi$, yield that
	\begin{align*}
2|\la \phi(x),h_1 \ra_{\cH_1}| &=	2|\la \phi(x),\phi(z) \ra_{\cH_1}| \\ &= 2|\la \phi'(x)\oplus 1,\phi'(z)\oplus 1 \ra_{\cH_1'\oplus \R}|\\ 
	&\leq 2|\la \phi'(x),\phi'(z) \ra_{\cH_1'}|+2 \\
	&= \left| \|\phi'(x)\|^2_{\cH_1'} +\|\phi'(z)\|^2_{\cH_1'}-\|\phi'(x)-\phi'(z)\|_{\cH_1'}^2  \right|+2 \\
	&= \left| \|\phi'(x)-\phi'(o)\|^2_{\cH_1'} +\|\phi'(z)-\phi'(o)\|^2_{\cH_1'}-\|\phi'(x)-\phi'(z)\|_{\cH_1'}^2  \right|+2 \\
	&= |d_\cX(x,o)+d_\cX(z,o)-d_\cX(x,z)| +2\\
	&\leq d_\cX(z,o) + |d_\cX(x,o)-d_\cX(x,z)|+2 \\
	& \leq 2d_\cX(z,o)+2,
	\end{align*}
	where $\phi'(o)=0$ by construction of $\phi'$; by the triangle inequality and the reverse triangle inequality. Now note that for any $y'\in \cY$, 
	\begin{align*}
	\int d_\cY(y,y') \, d|\nu_{h_1}|(y)  &\leq \int d_\cY(y,y)|\la \phi(x),h_1\ra_{\cH_1}| \, d|\theta|(x,y)  \\
	&\leq \lp d_\cX(z,o)+1\rp
	 \int d_\cY(y,y')\, d|\theta|(x,y) <\i ,
	\end{align*}
	so $\nu_{h_1}\in M^1(\cY)$. 
	 A consequence of the above mentioned Pettis integrability of $T_{h_1}\circ \phi\otimes \psi$ is that,  $$h^*_2(\la \phi(x),h_1 \ra_{\cH_1}\psi(y))=\la \phi(x),h_1 \ra_{\cH_1}h^*_2(\psi(y))\in \cL^{1}(\theta),$$ for all $h_2^*\in\cH_2^*$. By \cref{lemma_help_injectivity_proof} we have that $h^*_2(\psi(y))\in \cL^1(\nu_{h_1})$  and
	\begin{align*}
	\int h^*_2(\psi(y)) \, d\nu_{h_1}(y) = \int \la \phi(x),h_1 \ra_{\cH_1}h^*_2(\psi(y)) \, d\theta(x,y) = \int  h^*_2(\la \phi(x),h_1 \ra_{\cH_1}\psi(y)) \, d\theta(x,y),
	\end{align*}
	for all $h_2^*\in \cH_2^*$. By the unique defining property of the Pettis integral we now have that
	\begin{align*}
	\beta_\psi(\nu_{h_1}) = \int \psi \, d\nu_{h_1} = \int \la \phi(x),h_1 \ra_{\cH_1}\psi(y) \, d\theta(x,y)=0,
	\end{align*}
	so $\nu_{h_1}=0$ by the linearity and injectivity of $\beta_\phi$ on $M^1(\cX)$. \\ \\
	Now fix $h_1\in \mathrm{span}(\mathrm{Im}(\phi))\subset \cH_1$. For any fixed $B\in \cB(\cY)$ we obviously have that $\cH_1\ni h_1\mapsto \nu_{h_1}(B)$ is linear, but it is indeed also continuous. To see this, note that
	\begin{align*}
	|\nu_{h_1}(B)| &\leq \int |\la \phi(x),h_1 \ra_{\cH_1}|  \, d|\theta|(x,y) \\
	&\leq \|h_1\|_{\cH_1} \int \|\phi(x)\|_{\cH_1} \, d|\theta|(x,y),
	\end{align*}
	by the Cauchy-Schwarz inequality, proving that $h_1\mapsto \nu_{h_1}(B)$ is a bounded linear map since $\|\phi(x)\|_{\cH_1}\in \cL^1(\pi_1(|\theta|))$. 	 Since $h_1\in \mathrm{span}(\mathrm{Im}(\phi))\subset \cH_1$ can be written as $h_1=a_1h_{1,1}+\cdots +a_nh_{1,n}$ for $a_1,...,a_n\in \R$ and $h_{1,1},...,h_{1,n}\in \mathrm{Im}(\phi)$ for some $n\in \N$, we get that
	\begin{align*}
	\nu_{h_1}(B)= a_1 \nu_{h_{1,1}}(B) + \cdots +a_n \nu_{h_{1,n}}(B)=0,
	\end{align*}
	since $\nu_{h_{1,1}},...,\nu_{h_{1,n}}=0$. Hence $\nu_{h_1}=0$, since $\nu_{h_1}(B)=0$ for all $B\in \cB(\cY)$. \\ \\
	Now fix $h_1\in \overline{\mathrm{span}(\mathrm{Im}(\phi))}\subset \cH_1$. Note that any point in the closure can be written as the limit of elements in $\mathrm{span}(\mathrm{Im}(\phi))$. That is, $h_1=\lim_{n\to\i }h_{1,n}$ for some sequence $(h_{1,n})_{n\in \N}\subset \mathrm{span}(\mathrm{Im}(\phi))$. By the continuity of $h\mapsto \nu_{h}(B)$ we get that
	\begin{align*}
	\nu_{h_1}(B) = \lim_{n\to \i} \nu_{h_{1,n}}(B) = 0 ,
	\end{align*}
	for any $B\in \cB(\cY)$, proving that $\nu_{h_1}=0$ for all $h_1\in \overline{\mathrm{span}(\mathrm{Im}(\phi))}\subset \cH_1$.\\ \\
	\textbf{Step 2.1): Proving that $\mathbf{\nu_{h}=0}$ for all $\mathbf{h\in \cH_1}$.}\\
	Simply note that, since $F:=\overline{\mathrm{span}(\mathrm{Im}(\phi))}$ is a closed linear subspace of $\cH_1$, we have that $\cH_1= F+F^{\perp}$ (cf. corollary 6.15 \cite{applied_analysis_hunter2001}). That is, for any $h\in \cH_1$ there exist unique elements $u \in F $ and $w\in F^{\perp}$, such that $h=u+w$. Hence
	\begin{align*}
	\nu_{h}= \nu_{u}+\nu_{w}= \nu_{w}.
	\end{align*}
	Lastly, note that since $w\in \overline{\mathrm{span}(\mathrm{Im}(\phi))}^{\perp}$, then $w\perp \phi(x)$ for any $x\in \cX$. Thus we have that
	\begin{align*}
	\nu_{w}(B) = \int \la \phi(x),w\ra_{\cH_1}1_B(y)\, d\theta(x,y)  =0,
	\end{align*}
	for any $B\in \cB(\cY)$. We conclude that $\nu_{w}=0$, hence $\nu_{h}=0$, for any $h\in \cH_1$. Which is what we wanted to show.
\end{p}
\begin{remark}
	\Cref{lemma_strong_negative_type_beta_tensor_injective} is lemma 3.8 in \cite{lyons2013distance}. However the proof is vastly different from the original, since the original proof erroneously used that $\nu_h$ is a measure with finite first moment for every $h\in \cH_1$, without proving so. In personal communication with Russell Lyons, he acknowledges that the original proof is indeed not correct. However, he took a look at the proof again and provided me with the idea that resulted in  the above proof, where we circumvent the problem with $\nu_h$ not necessarily having first moment for all $h\in \cH_1$.
\end{remark}
We have now proved every essential lemma, which we will use to prove the main result of this section.
\begin{theorem} \label{theorem_dcov_zero_iff_independence}
	Let $(\cX,d_\cX)$ and $(\cY,d_\cY)$ be metric spaces of strong negative type. For any $\theta \in M^{1,1}_1(\cX\times \cY)$ with marginals $\mu\in M^1_1(\cX)$ and $\nu\in M^1_1(\cY)$, it holds that
	\begin{align*}
	dcov(\theta)=0\iff \theta=\mu\times \nu.
	\end{align*}
\end{theorem}
\begin{p}
	Since $(\cX,d_\cX)$ and $(\cY,d_\cY)$ are metric spaces of strong negative type, \cref{lemma_strong_negative_type_beta_tensor_injective} yields that there exist isometric embeddings $\phi:(\cX,d_\cX^{1/2}) \to \cH_1$ and $\psi:(\cY,d_\cY^{1/2})\to \cH_2$ into separable $\R$-Hilbert spaces such that the linear mapping $\beta_{\phi\otimes \psi}:M^{1,1}(\cX\times \cY) \to \cH_1\otimes \cH_2$ is injective. Let $\theta\in M^{1,1}_1(\cX\times \cY)$ and note that by \cref{lemma_beta_of_tensor_exists_and_representation_of_dcov} we have that
	\begin{align*}
	dcov(\theta) = 4\| \beta_{\phi\otimes \psi}(\theta-\mu\times \nu) \|_{\cH_1\otimes \cH_2}^2.
	\end{align*}
	If $\theta=\mu\times \nu$ then we obviously have that $dcov(\theta)=0$. For the converse statement, note that if  $dcov(\theta)=0$ then $\beta_{\phi\otimes \psi}(\theta-\mu\times \nu)=0$. Furthermore, as $\theta,\mu\times \nu\in M^{1,1}_1(\cX\times \cY)\subset M^{1,1}(\cX\times \cY)$ we get that $\theta-\mu\times \nu \in M^{1,1}(\cX\times \cY)$, since $M^{1,1}(\cX\times \cY)$ is a vector space. By the linearity and injectivity of $\beta_{\phi\otimes \psi}:M^{1,1}(\cX\times \cY) \to \cH_1\otimes \cH_2$ we see that the kernel $\mathrm{ker}(\beta_{\phi\otimes \psi})=\{m\in M^{1,1}(\cX\times \cY):\beta_{\phi\otimes\psi}(m)=0\}$ only contains the zero measure, hence
		\begin{align*}
		\theta-\mu\times \nu=0 \iff \theta=\mu\times \nu,
		\end{align*}
	which concludes the proof.	
\end{p}
This concludes the most important part of this section - namely that the distance covariance measure can be used as a direct measure of independence, in the case where both marginal metric spaces are of strong negative type. \\ \\
\textit{An important question to ask is whether metric spaces of strong negative type is the smallest class of metric spaces where distance covariance works as a direct measure of independence.}

The answer to this question is a partial yes. If we only consider metric spaces consisting of two or more points, then the answer is yes  by \cref{theorem_need_strong_negative_type_in_spaces_with_more_and_two_points} below. On the other hand, if one of the marginal spaces  consist only of a singleton then distance covariance can be used as a direct measure of independence regardless of the properties of the other marginal metric space. However such scenarios are not important for the independence problem, because every measure in $M^{1,1}_1(\cX\times \cY)$ is a product measure as explained by the below remark.
\begin{remark}
	Note that in independence testing, we are only interested in metric spaces $(\cX,d_\cX)$ and $(\cY,d_\cY)$ which consist of two or more points because, if one of the spaces consists only of a singleton, then the null-hypothesis is always satisfied. 
	That is, if $\cX$ is a space consisting only of a singleton $\cX=\{x\}$, then any probability measure $\theta\in M^{1,1}_{1}(\cX\times \cY)$ satisfies the null-hypothesis $\theta=\mu\times \nu$. This is seen by noting that any metric on $\cX$ generates the trivial Borel sigma algebra $\cB(\cX)=\{\emptyset,\cX\}$. As a consequence we have that $\theta(\cX\times B) = \pi_2(\theta)(B) = \nu(B)$ and $\theta(\emptyset \times B) = \theta(\emptyset)=0$, for any $B\in\cB(\cY)$. This proves that $\theta=\mu\times \nu$ since $\theta(A\times B)=\mu(A) \nu(B)$ for all $A\in\cB(\cX)$ and $B\in\cB(\cY)$ (intersection stable generator for $\cB(\cX)\otimes\cB(\cY)$).
\end{remark}
As mentioned above the distance covariance measure cannot  directly be used to verify independence if the (non-singleton) marginal metric spaces are not of strong negative type. To this end, we simply need to show that the distance covariance measure is flawed as a direct indicator of independence in metric spaces with at least two points that is not of strong negative type.
\begin{theorem} \label{theorem_need_strong_negative_type_in_spaces_with_more_and_two_points}
	Consider two metric spaces $(\cX,d_\cX)$ and $(\cY,d_\cY)$ both consisting of two or more points. If at least one  of these metric spaces is not of strong negative type, then there exists a measure $\theta\in M^{1,1}_1(\cX\times \cY)$ with marginals $\mu\in M^1_1(\cX)$ and $\nu\in M^1_1(\cY)$, such that $dcov(\theta)=0$ but 
	\begin{align*}
	\theta\not =\mu\times \nu.
	\end{align*}

\end{theorem}
\begin{p}
	Assume without loss of generality, that $(\cX,d_\cX)$ is not of strong negative type, meaning that there exist two distinct probability measures $\mu_1,\mu_2\in M^1_1(\cX)$ such that $D(\mu_1-\mu_2)=0$. Fix any such $\mu_1,\mu_2\in M^1_1(\cX)$ and two  distinct elements $y_1,y_2\in\cY$ (possible since $\cY$ is a non-singleton space). Now define the measure $\theta\in M^{1,1}_1(\cX\times \cY)$ by
	\begin{align*}
	\theta= \frac{\mu_1 \times \delta_{y_1} +\mu_2\times \delta_{y_2} }{2}.
	\end{align*}
	The marginals of this measure are easily seen to be $\mu=(\mu_1+\mu_2)/2\in M^1_1(\cX)$ and $\nu=(\delta_{y_1}+\delta_{y_2})/2\in M^1_1(\cY)$. Since $\mu_1\not=\mu_2$ we know that there exists a set $A\in \cB(\cX)$ such that $\mu_1(A)\not = \mu_2(A)$. Let $B\in\cB(\cY)$ such that $y_1\in B\not \ni y_2$ and note that
	\begin{align*}
	\theta(A\times B) = \frac{\mu_1(A)\delta_{y_1}(B) + \mu_2(A)\delta_{y_2}(B)}{2}= \frac{\mu_1(A)}{2},
	\end{align*}
	but
	\begin{align*}
	\mu\times \nu(A\times B) = \frac{\mu_1(A)+\mu_2(A)}{2}\frac{\delta_{y_1}(B)+\delta_{y_2}(B)}{2} \not = \frac{\mu_1(A)}{2},
	\end{align*}
	proving that $\theta\not = \mu\times \nu$. Now it suffices to show that this measure satisfies that $dcov(\theta)=0$ and to this end note that $D(\mu_1-\mu_2)=D(\mu_1)+D(\mu_2)-2D\mu_1\mu_2$, where we denoted $D\mu_1\mu_2:=\int d_\cX(x,x')\,d\mu_1\times\mu_2(x,x')$, which allows us to say that
	\begin{align*}
	D\lp \frac{\mu_1+\mu_2}{2}\rp &= \int \int d_\cX(x,x') \, d\lp \frac{\mu_1+\mu_2}{2}\rp(x) \, d\lp \frac{\mu_1+\mu_2}{2}\rp(x')\\	
	&=\frac{1}{4}\lp D(\mu_1)+D(\mu_2)+2D\mu_1\mu_2\rp=\frac{1}{4} D(\mu_1-\mu_2)+D\mu_1\mu_2
	, \\
	D\lp \frac{\delta_{y_1}+\delta_{y_2}}{2}\rp &=  \frac{1}{4}\lp 2\int\int d_\cY(y,y')\, d\delta_{y_1}(y)\,d\delta_{y_1}(y')\rp = \frac{1}{2}d_\cY(y_1,y_2).
	\end{align*}
	We also have that $a_{\mu}(x)=(a_{\mu_1}(x)+a_{\mu_2}(x))/2$ and $a_{\nu}(y)=(a_{\delta_{y_1}}(y)+a_{\delta_{y_2}}(y))/2$, hence
	\begin{align*}
	d_\mu(x,x') &= d_\cX(x,x') - \frac{a_{\mu_1}(x)+a_{\mu_2}(x)}{2} - \frac{a_{\mu_1}(x')+a_{\mu_2}(x')}{2}+\frac{1}{4} D(\mu_1-\mu_2)+D\mu_1\mu_2, \\
	d_\nu(y,y') &=d_\cY(x,x') - \frac{a_{\delta_{y_1}}(y)+a_{\delta_{y_2}}(y)}{2} - \frac{a_{\delta_{y_1}}(y')+a_{\delta_{y_2}}(y')}{2}+\frac{1}{2}d_\cY(y_1,y_2).
	\end{align*}
	Now note that by using the symmetry of $d_\mu$ and $d_\nu$, Fubini's theorem yields that
	\begin{align*}
	dcov(\theta)=  \frac{1}{4}\Big[&\int d_{\mu}(x,x') \, d\mu_1^2( x,x') \int d_\nu(y,y') \, d\delta_{y_1}^2(y,y')\\
	&+\int d_{\mu}(x,x') \,d\mu_2^2(x,x') \int d_\nu(y,y') \, d\delta_{y_2}^2(y,y')  \\
	&+ 2\int d_{\mu}(x,x') \, d\mu_1\times \mu_2(x,x') \int d_\nu(y,y') \,  d\delta_{y_1}\times \delta_{y_2}(y,y')\Big].
	\end{align*}
	Now by direct calculation we may derive all of the following equalities
	\begin{align*}
	\begin{array}{ll}
	\int d_\nu(y,y') \, d \delta_{y_1}^2(y,y') = -\frac{1}{2}d_\cY(y_1,y_2), &\int d_\nu(y,y') \, d \delta_{y_2}^2(y,y') = -\frac{1}{2}d_\cY(y_1,y_2),\\
	\int d_\nu(y,y') \, d \delta_{y_1}\times\delta_{y_2}(y,y') = \frac{1}{2}d_\cY(y_1,y_2),	&\int d_\mu(x,x') \, d \mu_1\times \mu_2(x,x') = -\frac{1}{4}D(\mu_1-\mu_2),\\
	\int d_\mu(x,x') \, d \mu_1^2(x,x') = \frac{1}{4}D(\mu_1-\mu_2),	&\int d_\mu(x,x') \, d \mu_2^2(x,x') = \frac{1}{4}D(\mu_1-\mu_2),
	\end{array}
	\end{align*}
	and as a consequence we get that
	\begin{align*}
	dcov(\theta)= \frac{1}{4} \lp -\frac{1}{8}-\frac{1}{8}-\frac{2}{8}\rp D(\mu_1-\mu_2)d_\cY(y_1,y_2) =-\frac{1}{8}	D(\mu_1-\mu_2)d_\cY(y_1,y_2)=0,
	\end{align*}
	proving that the distance covariance measure is zero, yet $\theta\not = \mu\times \nu$.
\end{p}

Having established that it is sufficient and necessary to assume strong negativity of the (non-singleton) marginal sample spaces, in order for the distance covariance to be used as a direct indicator of independence, we will now explore what kind of metric spaces are indeed of strong negative type. \\ \\
It is not easy to directly verify that a metric space is of strong negative type by the definition, so the last agenda of this section is to find a subclass of strong negative type metric spaces, with easily verifiable conditions or at the very least are more well-known. To this end, we have the following theorem, which tells us that every separable Hilbert space is a metric space of strong negative type.

\begin{theorem} \label{theorem_separable_hilbert_spaces_are_of_strong_negative_type}
	Every separable Hilbert space is  a metric space of strong negative type.
\end{theorem}
\begin{p}
	We prove this theorem in the case that the Hilbert space $\cX$ is an $\R$-Hilbert space, but we note that the arguments are similar for $\bC$-Hilbert spaces, by using the intermediate space $l^2_\bC(A)$ instead of $l^2_\R(A)$ used below. Alternatively, instead of considering $l^2_\bC(A)$, we could  start by embedding $\cX$ into its realification $\cX^{\R}$ with the additive isometric isomorphism $\mathrm{re}:\cX\to \cX^\R$ (see \cref{Appendix_Complexification_and_Realification}) and repeating the steps below. That is, we could use the embedding scheme $\cX\stackrel{\mathrm{re}}{\to}\cX^\R \stackrel{T}{\to}l^2_\R(A)$ and continue as below. \\ \\ 
    First we need some initial considerations about the $l^2(A)$ space.  Recall that $l^2(A)$ for $A\in\{\{1\},\{1,2\},...,\N\}$ is the separable Hilbert space of square summable  $A$-length real sequences. That is, 
    \begin{align*}
    l^2(A):=l^2_\R(A) = \lb (x_i)_{i\in A} \in \R^{|A|}: \sum_{i\in A} x_i^2 <\i  \rb\subset \R^{|A|},
    \end{align*} with coordinate-wise addition and $\R$-scalar multiplication when equipped with the inner product $\la (x_n),(y_n) \ra_{2}= \sum_{i=1}^{|A|} x_ny_n$ is a separable $\R$-Hilbert space (cf. example 21.2 \cite{schilling}). Also note that if $x^n \to_n x$ in $l^2(A)$, then $\pi_i(x^n)=x^n_i\to x_i=\pi_i(x)$ in $\R$ for all $i\in A$, so the coordinate projections are continuous.  If $\cX$ is an $|A|$-dimensional separable $\R$-Hilbert space then there exists a linear isometric homeomorphism $T:\cX\to l^2(A)$. Any orthonormal basis for $\cX$ has cardinality $|A|$, hence we may enumerate such a basis by $\{e_i:i\in A\}$. Now note that
    \begin{align*}
    T(x)= (\la x, e_i \ra_{\cX})_{i\in A},
    \end{align*}
    is a well-defined mapping from $\cX$ to $l^2(A)$ by Bessel's inequality $\sum_{i\in A} \la x, e_i \ra_{\cX}^2\leq \|x\|_\cX$. It is furthermore linear by the linearity of $\la \cdot, x\ra_\cX$ and the coordinate-wise addition and scalar multiplication on $l^2(A)$. We also note that, if $0=\|T(x)\|_2= \sqrt{\sum_{i\in A} \la x, e_i \ra_{\cX}^2}$, then $\la x, e_i \ra_{\cX}=0$ for all $i\in A$. This implies that $x=0$, since $\{e_i:i\in A\}$ is am orthonormal basis (cf. theorem 6.26 (a) \cite{applied_analysis_hunter2001}). This proves that $T$ is injective. Surjectivity of $T$ follows by noting that for any $(x_i)_{i\in A}\in l^2 (A)$ then $x=\sum_{i\in A} x_ie_i\in \cX$. This holds since $\sum_{i\in A} x_ie_i$ converges if and only if $\sum_{i\in A} \|x_ie_i\|_\cX^2<\i$ by the orthogonality of the terms and Lemma 6.23 \cite{applied_analysis_hunter2001}. Hence $x\in \cX$, since $ \sum_{i\in A} \|x_ie_i\|_\cX^2 =\sum_{i\in A} x_i^2 <\i$, where we used the normality of $(e_i)_{i\in A}$ and that $(x_i)_{i\in A}\in l^2(A)$. Now fix any $(x_i)_{i\in A}\in l^2(A)$ and note that
    \begin{align*}
    T\lp \sum_{i\in A} x_ie_i \rp =   \sum_{i\in A} T\lp x_ie_i \rp =  \sum_{i\in A} (\la x_ie_i , e_j\ra  )_{j\in A} =  (\la x_ie_i , e_i\ra  )_{i\in A} = (x_i  )_{i\in A},
      \end{align*}
      proving that $T$ is surjective, hence a bijective linear transformation (isomorphism). By Parseval's identity it furthermore holds that $     \|T(x)\|_2^2 = \sum_{i\in A} \la x, e_i \ra_{\cX}^2 = \|x\|_\cX^2$,   for any $x\in \cX$, proving that $T:\cX\to l^2(A)$ is a linear isometry - Thus $T$ is an linear isometric isomorphism, and as a consequence also a linear isometric homeomorphism. \\ \\
      We now consider finite dimensional and infinite dimensional separable Hilbert spaces separately, since arguments from the case of finite dimensional Hilbert spaces, will be used later to prove that the distance covariance measure in metric spaces (which we have derived) coincides with the distance covariance measure in Euclidean spaces introduced in $\cite{szekely2007measuring}$. \\ \\
     \textit{Finite dimensional separable Hilbert spaces:}\\
	If $(\cX,\la \cdot , \cdot\ra_{\cX})$ is a  finite dimensional $\R$-Hilbert space, then $\cX$ is isometric homeomorphic to $\R^n$ for some $n=\mathrm{dim}(\cX)\in\N$. This is seen by noting  $F:l^2(\{1,...,n\})\to \R^n$ given by $F((x_i)_{i\in \{1,...,n\}} ) =(x_1,...,x_n)\in \R^n$, is an obvious linear isometric isomorphism. The composition $\iota= F\circ T:\cX \to \R^n$ is therefore a linear isometric homeomorphism.  \\ \\	
	Let $f:\R^n\to \R$ be given by $f(s)=\|s\|^{-(n+1)}$ (zero in zero) and define  $\phi:\R^n\to L_\bC^2(\R^n , f\cdot \lambda^n)$ by
	\begin{align*}
	\phi(x)(s) = c(1-e^{is^\t x}),
	\end{align*}
	where $\lambda^n$ is the Lebesgue measure on $(\R^n,\cB(\R^n))$ and $c$ is an appropriate constant (chosen below). This map is called the Fourier  embedding in \cite{lyons2013distance}. That $\phi(x)\in L_\bC^2(\R^n , f\cdot \lambda^n)$ follows from the fact that $1-e^{is^\t x}+1-\overline{e^{is^\t x}}=2(1-\cos(s^\t x)) $,
 hence	\begin{align*}
	\|\phi(x)\|_{L_\bC^2(\R^n , f\cdot \lambda^n)}^2 &= c^2\int |1-e^{is^\t x}|^2 \,d (f\cdot \lambda^n)(s) \\
	&= c^2\int \frac{(1-e^{is^\t x})(1-\overline{e^{is^\t x})}}{\|s\|^{n+1}} \,d  \lambda^n(s) \\
	&= c^2\int \frac{1-e^{is^\t x}+1-\overline{e^{is^\t x}}}{\|s\|^{n+1}} \,d  \lambda^n(s) \\
		&= 2c^2 \int \frac{1-\cos(s^\t x)}{\|s\|^{n+1}} \,d  \lambda^n(s) \\
	&= 2c^2  c_n \|x\| \\
	&<\i,
	\end{align*}
	where the fifth equality and constant $c_n\in\R$ are found in lemma 1 \cite{szekely2007measuring}. We may also note that
	\begin{align*}
	|c(1-e^{is^\t x})-c(1-e^{is^\t x'})|^2 &= c^2|e^{is^\t x'}-e^{is^\t x}|^2 \\
	&= c^2(e^{is^\t x'}-e^{is^\t x})(\overline{e^{is^\t x'}}-\overline{e^{is^\t x}}) \\
	&=c^2\lp (e^{is^\t x'}\overline{e^{is^\t x'}}+e^{is^\t x}\overline{e^{is^\t x}}- e^{is^\t x'}\overline{e^{is^\t x}}-e^{is^\t x}\overline{e^{is^\t x'}})\rp \\
	&=c^2 \lp 2-e^{is^\t (x'-x)}-e^{is^\t (x-x')}\rp \\
	&=2c^2 \lp 1- \cos (s^\t (x-x'))\rp.
	\end{align*}
	As a consequence of the same lemma used above, we also get that
	\begin{align*}
	\|\phi(x)-\phi(x')\|_{L_\bC^2(\R^n , f\cdot \lambda^n)}^2 &= \int \frac{|c(1-e^{is^\t x})-c(1-e^{is^\t x'})|^2}{\|s\|^{n+1}} \, d\lambda^n(s) \\
	&= 2c^2\int \frac{1- \cos (s^\t (x-x'))}{\|s\|^{n+1}} \, d\lambda^n(s) \\
	&= 2c^2c_n\|x-x'\|.
	\end{align*}
	for any $x,x'\in \R^n$. Hence if we set $c^2=1/(2c_n)$ we get that
	\begin{align*}
	d_\cX(x,x') = \|\iota(x)-\iota(x')\| = \|\phi(\iota(x))-\phi(\iota(x'))\|_{L_\bC^2(\R^n , f\cdot \lambda^n)}^2,
	\end{align*}
	proving that $\phi'\equiv\phi\circ\iota=\phi \circ F\circ T:(\cX,d_\cX^{1/2})\to L_\bC^2(\R^n , f\cdot \lambda^n)$ is an isometric embedding into a $\bC$-Hilbert space, which proves that $(\cX,d_\cX)$ is a metric space of negative type by \cref{theorem_equivalence_of_negative_type}. We need to show  that $(\cX,d_\cX)$ is of strong negative type, and by \cref{lemma_strong_negative_type_iff_beta_is_injective_on_prob_measures_of_finite_first_moment} it suffices to show that the mean embedding map $\beta_{\phi'}:M^1(\cX)\to L_\bC^2(\R^n , f\cdot \lambda^n)$ is injective on $M^1_1(\cX)$. To this end, let $\mu\in M^1_1(\cX)$ and note that $\beta_{\phi'}(\mu)$ is the unique element in $L_\bC^2(\R^n , f\cdot \lambda^n)$ satisfying
	\begin{align*}
	g^*(\beta_{\phi'}(\mu)) &= \int g^*(\phi'(x)) \, d\mu(x) \\
	&= \int \la \phi'(x),g\ra_{L_\bC^2(\R^n , f\cdot \lambda^n)} \, d\mu(x) \\
	&= \int \int \phi'(x)(s)\overline{g(s)} \, d(f\cdot \lambda^n)(s) \mu(x)\\
	&= \int \int c(1-e^{is^\t \iota(x)})\overline{g(s)} \, d(f\cdot \lambda^n)(s) \, d\mu(x) \\
	&= \int  c\Big(1- \int e^{is^\t x} \, d\iota(\mu)(x)\Big)   \, \overline{g(s)} \, d(f\cdot \lambda^n)(s) \\
	&= \int c(1-\widehat{\iota(\mu)}(s))\overline{g(s)} \, d(f\cdot \lambda^n)(s) \\
	&= \la c(1-\widehat{\iota(\mu)}), g \ra_{L_\bC^2(\R^n , f\cdot \lambda^n)} \\
	&= g^*(c(1-\widehat{\iota(\mu)})),
	\end{align*}
	by Fubini's theorem, where $\widehat{\iota(\mu)}:\R^n \to \bC$ is the characteristic function corresponding to the push-forward measure $\iota(\mu)$ on $(\R^n,\cB(\R^n))$ (see \cref{Appendix_Characteristic_fucntions} for more details on characteristic functions). The application of Fubini's theorem is justified by noting that $(f\cdot \lambda^n$) is a $\sigma$-finite measure and that $|\phi'(x)|,|g|\in L^2(\R^n,f\cdot \lambda^n)$, implying that
	\begin{align*}
	\int \int |\phi'(x)(s) \overline{g(s)}| \, d(f\cdot \lambda^n)(s) \, d\mu(x) &= \int \la |\phi'(x)|,|g|\ra_{L^2(\R^n,f\cdot \lambda^n)} \, d\mu(x) \\
	&\leq \int |\la |\phi'(x)|,|g|\ra_{L^2(\R^n,f\cdot \lambda^n)}| \, d\mu(x) \\
	&\leq \int \|\phi'(x)\|_{L^2(\R^n,f\cdot \lambda^n)} \|g\|_{L^2(\R^n,f\cdot \lambda^n)} \, d\mu(x) \\
	&= \|g\|_{L^2(\R^n,f\cdot \lambda^n)} \int d_\cX(x,0)^{1/2} \, d\mu(x) \\
	&<\i,
	\end{align*}
	where we used Cauchy-Schwarz's inequality and that $\mu$ has finite first moment. We furthermore used that $\phi'(0)$ is the zero element of $L^2(\R^n,f\cdot \lambda^n)$  such that $$\|\phi'(x)\|_{L^2(\R^n,f\cdot \lambda^n)}=\|\phi'(x)-\phi'(0)\|_{L^2(\R^n,f\cdot \lambda^n)} = d_\cX(x,0)^{1/2}.$$ To see the $\sigma$-finitesness, note that $(B_j)_{j\geq 1}$ given by $B_j=\{0\}\cup (B(0,j)\setminus B(0,1/j))$ for $j\geq 1$, is a exhausting sequence of Borel sets with $(f\cdot \lambda^n)(B_j)<\i$. This is seen by noting that $f(s)\leq (1/j)^{n+1}$ for all $s\in B_j$ (recall that we set $f(0)=0$), such that $(f\cdot \lambda^n) ( B_j ) = \int_{B_j} f(s) \, d\lambda^n(s) \leq (1/j)^{n+1} \lambda^n ( B(0,j)) <\i$, proving that $(f\cdot \lambda^n)$ is a $\sigma$-finite measure.
	
	This proves that $\beta_{\phi'}(\mu)=c(1-\widehat{\iota(\mu)})$ in $L_\bC^2(\R^n , f\cdot \lambda^n)$. \\ \\ Assume for contradiction that $\beta_{\phi'}$ is not injective on $M^1_1(\cX)$. Then there exists two distinct Borel probability measures $\mu_1,\mu_2\in M^1_1(\cX)$ such that $\beta_{\phi'}(\mu_1)=\beta_{\phi'}(\mu_2)$ in  $L_\bC^2(\R^n , f\cdot \lambda^n)$, which by the above derivation happens if and only if $c(1-\widehat{\iota(\mu_1)}(s))=c(1-\widehat{\iota(\mu_2)}(s))$ for $(f\cdot \lambda^n)$-almost all $s\in\R^n$. Since $f(s)> 0$ for all $s\in\R^n \setminus \{0\}$ we may conclude that $\widehat{\iota(\mu_1)}(s)= \widehat{\iota(\mu_2)}(s)$ for $\lambda^n$-almost all $s\in\R^n$. By \cref{appendix_theorem_properties_char} item \cref{char_prooerties_unique} and \cref{char_prooerties_relax} we have that $\iota(\mu_1)=\iota(\mu_2)$, since their characteristic functions coincide. 
	
	Now recall that $\iota$ is a isometric homeomorphism, i.e. it has a well-defined inverse $\iota^{inv}:\R^n \to \cX$ that is continuous and measurable. As a consequence we have that  $\iota^{-1}(\iota(A))=A$ and  the pre-image of the inverse coincides with the image $\iota(A)=(\iota^{inv})^{-1}(A)$, such that $\iota(A)\in\cB(\R^n)$ for any $A\in\cB(\cX)$. Hence for any $A\in \cB(\cX)$
	\begin{align*}
	\mu_1(A) = \mu_1(\iota^{-1}(\iota(A))) = \iota(\mu_1)(\iota(A))=\iota(\mu_2)(\iota(A))=\mu_2(\iota^{-1}(\iota(A)))=\mu_2(A),
	\end{align*}
	a contradiction. We conclude that $\beta_{\phi'}$ is injective on $M^1_1(\cX)$, proving that $\cX$ is a metric space of negative type. \\ \\
	\textit{Infinite dimensional separable Hilbert spaces:}\\
	Let $(\cX,\la \cdot , \cdot\ra_{\cX})$ be a infinite dimensional separable $\R$-Hilbert space. We will show that $\cX$ equipped with the naturally induced metric $d_\cX$ is of strong negative type by showing that there exists an isometric embedding $\phi:(\cX,d_\cX^{1/2})\to L^2(\R^\i \times \R, \rho \times \lambda)$, into a separable $\R$-Hilbert space, which induces an injective mean embedding map $\beta_\phi:M^{1}_1(\cX)\to L^2(\R^\i \times \R, \rho \times \lambda)$, where $\rho$ is a probability measure on $(\R^\i,\cB(\R^\i))$ and $\lambda$ is the Lebesgue measure on $(\R,\cB(\R))$. This will indeed imply that $(\cX,d_\cX)$ is a metric space of strong negative type by \cref{lemma_strong_negative_type_iff_beta_is_injective_on_prob_measures_of_finite_first_moment}. \\ \\
    The isometric embedding we are going to construct requires the definition of another map, which we start by defining.     Let $(Z_n)_{n\in \N}$ be an independent and identically standard normal distributed sequence of random variables defined on some probability space $(\Omega,\F,P)$. Let $Z:\Omega \to \R^\i$ denote its natural embedding into $\R^\i$ equipped with the product $\sigma$-algebra $\cB(\R^\i)$ (see \cref{Appendix_Product_spaces_section}). $Z$ is obviously measurable and we may denote its law $\rho=Z(P)$ on $(\R^\i,\cB(\R^\i))$.   For any $u=(x_n)_{n\in\N}\in l^2(\N)$ one may show that
    \begin{align*}
    \sum_{n=1}^k u_nZ_n \convd_k Z(u)\sim \cN\lp 0, \sum_{n=1}^\i u_n^2\rp = \cN\lp 0, \|u\|_2^2\rp,
    \end{align*}
    by the use of characteristic functions and Levy's continuity theorem. But we also have that $\sum_{n=1}^k u_nZ_n \stackrel{P-a.s.}{\longrightarrow}_k \sum_{n=1}^\i u_nZ_n$, so  $\sum_{n=1}^\i u_nZ_n\sim \cN\lp 0, \|u\|_2^2\rp$ by uniqueness of limits. Hence, $\Big|\sum_{n=1}^\i u_nZ_n\Big|<\i $ $P$-almost surely. We may also note that $|Z(u)|$ is half-normal distributed, so \begin{align} \label{ep_temp_99786}
    E|Z(u)|=c\|u\|_2,
    \end{align}
    where $c=\sqrt{2}/\sqrt{\pi}$. Let $g:\R^\i \times \R^\i \to \bar{\R}$ be given by
    \begin{align*}
    g(u,v) = \limsup_{N\to\i } \sum_{n=1}^N u_nv_n,
    \end{align*}
    and note by the above convergence in distribution we have that for fixed $u\in l^2(\N)$, $g(u,v)$ converges in $\R$ for $\rho$-almost all $v\in\R^\i$. Furthermore let $g|_{l^2(\N)}$ denote the restriction of $g$ to $l^2(\N)\times \R^\i$. Also note that $g$ and $g|_{l^2(\N)}$ are respectively $\cB(\R^\i)\otimes \cB(\R^\i)$- and $\cB(l^2(\N))\otimes \cB(\R^\i )$-measurable  mappings, since they are given by the limit supremum of measurable mappings (coordinate projections are continuous, hence measurable). In general, when we write $f|_A$ we mean the restriction of $f$ to $A$. \\ \\
    Let $T:\cX \to  l^2(\N)$  denote the linear isometric homeomorphism defined in the beginning of the proof. We define our embedding $\phi:\cX \to L^2(\R^\i \times \R, \rho\times \lambda)$ into the separable  $\R$-Hilbert space  $L^2(\R^\i \times \R, \rho\times \lambda)$ by 
    \begin{align*}
    \phi(x)(v,s)= 1_{[g(T(x),v)/c,\i)}(s) - 1_{[0,\i)}(s).
    \end{align*}
    To see that $\phi(x)\in L^2(\R^\i \times \R, \rho\times \lambda)$ for all $x\in\cX$, we simply note that
    \begin{align*}
    \int| \phi(x)(v,s)|^2 \, d\rho \times \lambda (v,s) =&\int_{B\times \R} 1_{[g(T(x),v)/c,0)}(s) \, d\rho \times \lambda(v,s) \\
    &+ \int_{B^c\times \R} 1_{[0,g(T(x),v)/c)}(s) \, d\rho \times \lambda(s) \\
    =& \int_B -g(T(x),v)/c \, d\rho(v) + \int_{B^c} g(T(x),v)/c \, d\rho(v)  \\
    =& \frac{1}{c} \int |g(T(x),v)| \, d\rho(v) = \frac{1}{c}E|(g(T(x),Z)| \\
    =& \frac{1}{c}E\left|  Z(T(x))\right| = \|T(x)\|_2 = \|x\|_{\cX},
    \end{align*}
   which is finite, where $B=\{v\in \R^\i:g(T(x),v) < 0\}$. Furthermore by similar reasoning
    \begin{align*}
    \|\phi(x)-\phi(x')\|_{L^2(\R^\i \times \R, \rho\times \lambda)}^2 =& \int |1_{[g(T(x),v)/c,\i)}(s)-1_{[g(T(x'),v)/c,\i)}(s)|^2 \, d\rho \times \lambda(v,s)  \\
    =&\int_{D\times \R} 1_{[g(T(x),v)/c,g(T(x'),v)/c))}(s) \, d\rho \times \lambda (v,s) \\
    &+ \int_{D^c\times \R} 1_{[g(T(x'),v/c),g(T(x),v)/c)}(s) \, d\rho \times \lambda (v,s) \\
    =& \frac{1}{c}\int |g(T(x),v)-g(T(x'),v)| \, d\rho(v) \\
    =&\frac{1}{c}\int |g(T(x-x'),v)| \, d\rho(v) = \frac{1}{c} E|Z(T(x-x')| \\
    =& \|T(x-x')\|_2 = \|x-x'\|_\cX,
    \end{align*}
    where $D=\{v\in\R^\i  :g(T(x),v)<g(T(x'),v) \}$, proving that $\phi:(\cX,d_\cX^{1/2})) \to L^2(\R^\i \times \R, \rho\times \lambda)$ is an isometric embedding. It remains to be shown that the mean embedding map $\beta_{\phi}:M^1(\cX) \to L^2(\R^\i \times \R, \rho\times \lambda)$ given by
    \begin{align*}
    \beta_{\phi}(\mu)= \int \phi \, d\mu,
    \end{align*}
    is injective on $M^1_1(\cX)\subset M^1(\cX)$. That is, $\beta_{\phi}(\mu_1)=\beta_{\phi}(\mu_2)$ implies $\mu_1=\mu_2$ for any two measures $\mu_1,\mu_2\in M^1_1(\cX)$. Thus fix $\mu_1,\mu_2\in M^1_1(\cX)$ with $\beta_{\phi}(\mu_1)=\beta_{\phi}(\mu_2)$ and denote $\mu=\mu_1-\mu_2$. By the same arguments as in the above proof for finite-dimensional Hilbert spaces  it suffices to show that $T(\mu_1)=T(\mu_2)$  as measures on $(l^2(\N),\cB(l^2(\N)))$. \\ \\
    Now we will show that these measures on $(l^2(\N),\cB(l^2(\N)))$ are uniquely determined by their finite-dimensional distributions pushed forward to $\R$ by dual mappings, what this exactly entails will become clear below.  
    
    First we show that $ \cB(l^2(\N))=\cB(\R^\i) \cap l^2(\N)$. Note that $l^2(\N)\in \cB(\R^\i)$, which follows by noting that all coordinate projections  $\pi_i:\R^\i\to \R$ are continuous mappings, by the definition of the product topology that induced the Borel $\sigma$-algebra $\cB(\R^\i)$. As a consequence they are also Borel measurable, rendering the map $f:\R^\i \to [0,\i]$ defined by $f(x)=\sum_{i=1}^\i |\pi_i(x)|^2$ Borel measurable. Now note that $f(x)=\|x\|_2^2$ for any $x\in l^2(\N)$ and $f(x)<\infty \iff x \in  l^2(\N)$, so $l^2(\N)=f^{-1}([0,\i))\in \cB(\R^\i)$. We now show that $\cB(l^2(\N))=l^2(\N)\cap \cB(\R^\i)$ by showing each is included in the other. A generator for the trace $\sigma$-algebra $\cB(\R^\i) \cap l^2(\N)$ is given by $\{l^2(\N)\cap \pi_i^{-1}((a,b)):a,b\in\R,i\in\N\}=\{ (\pi_i|_{l^2(\N)})^{-1}((a,b)):a,b\in\R,i\in\N\}$ (see \cref{Appendix_Product_spaces_section}) and all these sets are open in $l^2(\N)$ since $\pi_i|_{l^2(\N)}$ is continuous for all $i\in\N$. This proves that $\cB(\R^\i) \cap l^2(\N) \subset \cB(l^2(\N))$, since we have inclusion of their generators. For the converse we show that every open ball in $l^2(\N)$ lies in $\cB(\R^\i)$, proving that $\cB(\R^\i) \cap l^2(\N) \supset \cB(l^2(\N))$ since the open balls in $l^2(\N)$ generate $\cB(l^2(\N))$. To this end, note that every open ball in $l^2(\N)$ has the form $\{x\in l^2(\N): \|x-y\|_2<\delta\}=l^2(\N)\cap \{x\in\R^\i:f(x-y)<\delta^2\}$ for some $y\in l^2(\N)$ and $\delta >0$. We realize that $x\mapsto g(x):=f(x-y)$ is Borel measurable, hence the above set becomes $l^2(\N)\cap g^{-1}([0,\delta^2))\in \cB(\R^\i).$ Thus $ \cB(l^2(\N))=\cB(\R^\i) \cap l^2(\N)$. 
    
    Let $\tilde{T}(\mu_1)$ and $\tilde{T}(\mu_2)$  denote the extensions of $T(\mu_1)$ and $T(\mu_2)$  to $(\R^\i,\cB(\R^\i))$ given by
	\begin{align*}
	\tilde{T}(\mu_1)(A) = T(\mu_1)(A\cap l^2(\N)) \quad \text{ and } \quad  \tilde{T}(\mu_2)(A)=T(\mu_2)(A\cap l^2(\N)),
	\end{align*}
    for any $A\in\cB(\R^\i)$. Since $ \cB(l^2(\N))=\cB(\R^\i) \cap l^2(\N)$, it suffices to show that $\tilde{T}(\mu_1)=\tilde{T}(\mu_2)$. Furthermore it is well know that probability measures of $(\R^\i,\cB(\R^\i))$ are uniquely determined by their finite dimensional distributions meaning that it suffices to show that $\pi_{1\cdots n}\circ \tilde{T}(\mu_1)=\pi_{1\cdots n}\circ \tilde{T}(\mu_2)$ for all $n\in \N$. This equality can be expressed through equality of their respective characteristic functions, see \cref{appendix_theorem_properties_char} item \cref{char_prooerties_unique} and \cref{char_prooerties_relax}. We note that
    \begin{align*}
    \int_{\R^n} e^{i\la s, x\ra_{\R^n} } \, d\pi_{1\cdots n}\circ \tilde{T}(\mu_i)(x) =  \int_{\R} e^{ix} \, ds^* \circ \pi_{1\cdots n}\circ \tilde{T}(\mu_i)(x),
    \end{align*}
    showing that it suffices to show that $s^* \circ \pi_{1\cdots n}\circ \tilde{T}(\mu_1)=s^* \circ \pi_{1\cdots n}\circ \tilde{T}(\mu_2)$ on $(\R,\cB(\R))$ for all $n\in\N$ and  $s\in\R^n$, where we use Riesz's  representation theorem to uniquely connect dual mappings $s^*\in (\R^n)^*$ with a corresponding element $s\in\R^n$ such that $s^*(x)=\la x , s \ra_{\R^n}$.

    Now we show that it suffices to prove that \begin{align*}
    s^*\circ \lp \pi_{1\cdots n}|_{l^2(\N)}\rp
     \circ T(\mu_1)=s^*\circ \lp \pi_{1\cdots n}|_{l^2(\N)}\rp
      \circ T(\mu_2),
    \end{align*}
     for all $s\in \R^n$ and $n\in\N$ in order for $T(\mu_1)=T(\mu_2)$, which is exactly what we meant by saying that probability measures on $(l^2(\N),\cB(l^2(\N))$ are uniquely determined by their finite dimensional distributions pushed forward to $\R$ by dual mappings. Simply note for any $s^*\in(\R^n)^*$, $n\in \N$ and $y\in \R$ it holds that
    \begin{align*}
    s^*\circ \lp \pi_{1\cdots n}|_{l^2(\N)}\rp
     \circ T(\mu_i) ((-\i, y]) &= \lp \pi_{1\cdots n}|_{l^2(\N)} \rp
     \circ T(\mu_i) ((s^*)^{-1}(\i, y]) \\
    &=   T(\mu_i) ((\pi_{1\cdots n}|_{l^2(\N)})^{-1} \circ (s^*)^{-1}((-\i, y])) \\ 
    &=   T(\mu_i) ((\pi_{1\cdots n})^{-1} \circ (s^*)^{-1}((-\i, y]) \cap l^2(\N)) \\
    &=   \tilde{T}(\mu_i) ((\pi_{1\cdots n})^{-1} \circ (s^*)^{-1}((-\i, y]) ) \\
    &= s^* \circ \pi_{1\cdots n}\circ \tilde{T}(\mu_i)((-\i,y]),
    \end{align*}
    proving the claim, since $\{(-\i,y]:y\in\R\}$ is an intersection stable generator for $\cB(\R)$. \\ \\ 
    Hence with a little manipulation of the above representation, we see that it suffices to show that the following cdfs coincide. That is,
    \begin{align*}
    F_{1,n,s}(y):=T(\mu_1)(u\in l^2(\N): \la u_{\leq n},s \ra \leq y) =T(\mu_2)(u\in l^2(\N): \la u_{\leq n},s \ra \leq y)=:F_{2,n,s}(y),
    \end{align*}
    for all  $s\in \R^n$, $n\in \N$ and $y\in \R$.
    \\ \\
    To this end we need a property implied by the assumption that $\beta_\phi(\mu)=0$. Note that $\beta_\phi(\mu)$ is the unique element in $L^2(\R^\i \times \R, \rho\times \lambda)$ satisfying 
    \begin{align*}
    h^*(\beta_{\phi}(\mu))&= \int h^* \circ \phi(x) \, d\mu(x) \\
        &= \int \int h(v,s)\phi(x)(v,s) \, d\rho \times \lambda (v,s) \, d\mu(x) \\
    &= \int h(v,s) \int \phi(x)(v,s) \, d\mu(x) \, d\rho\times \lambda(v,s) \\
    &= \int h(v,s) \lp \int 1_{(-\i,cs]}(g(T(x),v))\, d\mu(x) -1_{[0,\i)}(s)\mu(\cX) \rp \, d\rho\times \lambda(v,s) \\
    &= \int h(v,s)  \mu(x\in \cX: g(T(x),v)\leq cs) )  \, d\rho\times \lambda(v,s) \\
    &= h^*[(v,s) \mapsto \mu(x\in\cX:g(T(x),v)\leq cs)],
    \end{align*}
    for all $h^*\in L^2(\R^\i \times \R, \rho\times \lambda)^*$, where we used Riesz's representation theorem and Fubini's theorem, which is justified since \begin{align*}
    \int \int |h(v,s)\phi(x)(v,s)| \, d \rho\times \lambda (v,s) \, d|\mu|(x)& \leq \|h\|_{L^2(\R^\i\times \R,\rho\times \lambda)} \int d_\cX(x,0)^{1/2} \, d|\mu|(x)  <\i,
    \end{align*}
    by Cauchy-Schwarz's inequality and that $\mu=\mu_1-\mu_2$ has finite first moment.  We conclude that the mean embedding  $\beta_\phi(\mu)$ is the map in $L^2(\R^\i \times \R, \rho\times \lambda)$ given by
    \begin{align*}
    \beta_\phi(\mu):(v,s) &\mapsto \mu(x\in\cX:g(T(x),v)\leq cs)\\
    &=\mu(g(T(\cdot),v)^{-1}((-\i,cs])) \\
    &=\mu(T^{-1}\circ g|_{l^2(\N)}(\cdot,v)^{-1}((-\i,cs]))\\
    &= T(\mu)( g|_{l^2(\N)}(\cdot,v)^{-1}((-\i,cs])) \\
   &= T(\mu)( u\in l^2(\N) : g(u,v) \leq cs).
       \end{align*}
    Let $L:\R^\i \to \R^\i$ denote the left shift operator, i.e. $  L((x_n)_{n\in \N}) = (x_{n+1})_{n\in \N}$ for any $(x_n)_{n\in\N}\in \R^\i$. For any $k\in\N$ and $u\in\R^\i$, we write $u_{\leq k}=(u_1,...,u_k)$ and $u_{>k}=L^k(u)$ such that $u=(u_{\leq k}, u_{>k})$. Note that $\rho = \pi_{1\cdots k}(\rho)\times L^k(\rho)=\pi_{1\cdots k}(\rho)\times \rho$, since $L$ is $\rho$-measure preserving ($(Z_n)$ is stationary). Furthermore we note that $\pi_{1\cdots k}(\rho)=\cN(0,I_k)$, such that $\pi_{1\cdots k}(\rho) \ll \lambda^k$.  Since we assumed that $\beta_\phi(\mu)=0$,  Tonelli's theorem yields that for any $k\in \N$
    \begin{align*}
    0&= \|\beta_{\phi}(\mu)\|_{L^2(\R^\i \times \R, \rho\times \lambda)}^2 \\
    &= \int | T(\mu)(u\in l^2(\N) :g(u,v)\leq cs)|^2 \, d \rho \times \lambda(v,s) \\
    &= \int\int\int | T(\mu)(u\in l^2(\N) :g(u,(v_{\leq k},v_{>k}))\leq cs)|^2 \, d\lambda(s) \, d\pi_{1\cdots k}(\rho)(v_{\leq k}) \, d\rho(v_{\geq k}).
    \end{align*}
Hence for every $k\in \N$ there exists a $\rho$-almost sure set $H_k$, such that for every $v\in H_{k}$ there exists a $\lambda^k$-almost everywhere set $H_{k,v}$, such that for every $s\in H_{k,v}$ there exists a $\lambda$-almost everywhere set $H_{k,v,s}$, with the following properties. For any $k\in \N,v\in H_k,s\in H_{k,v},y\in H_{k,v,s}$ it holds that
\begin{align*}
0=&T(\mu)(u\in l^2(\N): \la u_{\leq k}, s\ra_{\R^k} +g(u_{>k},v) \leq y ) \\
=& T(\mu_1)(u\in l^2(\N): \la u_{\leq k}, s\ra_{\R^k} +g(u_{>k},v) \leq y )\\
&-T(\mu_2)(u\in l^2(\N): \la u_{\leq k}, s\ra_{\R^k} +g(u_{>k},v) \leq y ) \\
\iff& \\
&T(\mu_1)(u\in l^2(\N): \la u_{\leq k}, s\ra_{\R^k} +g(u_{>k},v) \leq y ) \\
=&T(\mu_2)(u\in l^2(\N): \la u_{\leq k}, s\ra_{\R^k} +g(u_{>k},v) \leq y ). \notag
\end{align*}

First we show that for any $k\in\N,v\in H_k$ this statement can be strengthened to all $s\in \R^k,y\in \R$. Fix any $k\in \N$ and $v\in H_k$ and let $s\in \R^k$. Note that we can find a sequence $(s_n)\subset H_{k,v}$ such that $s_n \to_n s$. 
Now let $X_1 \sim T(\mu_1)$ and $X_2\sim T(\mu_2)$ be defined on some probability space \pspace, and note that \begin{align*}
\la \pi_{1\cdots k}(X_i) , s_n \ra_{\R^k}  +g(L^k(X_i),v)
\convas_n &\la \pi_{1\cdots k}(X_i) , s \ra_{\R^k}  +g(L^k(X_i),v)
\end{align*}

for $i=1,2$ by continuity of the inner product, so we also have convergence in distribution. Hence also point-wise convergence of the cdfs for every continuity point of the limit distribution cdfs. Let $D_{k,v,s}^1$ and $D_{k,v,s}^2$ denote the corresponding discontinuity points of the limit cdfs, which are  at most countably infinite ($\lambda$-nullsets). Thus for every $y$ in the $\lambda$-almost everywhere  set
\begin{align*}
(\R\setminus D_{k,v,s}^1 )\cap (\R\setminus D_{k,v,s}^2) \cap \lp \bigcap_{n=1}^\i H_{k,v,s_{n}}\rp,
\end{align*}
we have that
\begin{align*}
&T(\mu_1)(u\in l^2(\N): \la u_{\leq k}, s\ra_{\R^k} +g(u_{>k},v) \leq y ) \\
=& P(\la \pi_{1\cdots k}(X_1) , s \ra_{\R^k}  +g(L^k(X_1),v)\leq y) \\
=& \lim_{n\to \i } P(\la \pi_{1\cdots k}(X_1) , s^n \ra_{\R^k} +g(L^k(X_1),v)\leq y)\\
=& \lim_{n\to \i }T(\mu_1)(u\in l^2(\N): \la u_{\leq k}, s^n\ra_{\R^k} +g(u_{>k},v) \leq y ) \\
=& \lim_{n\to \i }T(\mu_2)(u\in l^2(\N): \la u_{\leq k}, s^n\ra_{\R^k} +g(u_{>k},v) \leq y ) \\
=& \lim_{n\to \i } P(\la \pi_{1\cdots k}(X_2) , s^n \ra_{\R^k} +g(L^k(X_2),v)\leq y)\\
=& P(\la \pi_{1\cdots k}(X_2) , s \ra_{\R^k}  +g(L^k(X_2),v)\leq y) \\
=&T(\mu_2)(u\in l^2(\N): \la u_{\leq k}, s\ra_{\R^k} +g(u_{>k},v) \leq y ).
\end{align*}
So we have two mappings that are cadlag in $s$ (they are cdfs) which coincide $\lambda$-almost all $s\in \R$. By trivial $\ep/\delta$-arguments (similar to those below) they must coincide for all $s\in \R$.

We conclude that for every $k\in \N$ there exists a $\rho$-almost sure set $H_k$ such that for all $v\in H_k$ it holds that
\begin{align} \label{eq_g1}
&T(\mu_1)(u\in l^2(\N): \la u_{\leq k}, s\ra_{\R^k} +g(u_{>k},v) \leq y ) \\
=&T(\mu_1)(u\in l^2(\N): \la u_{\leq k}, s\ra_{\R^k} +g(u_{>k},v) \leq y ) \notag
\end{align}
for all $s\in \R^k$ and $y\in \R$  \\ \\ 
    With this in mind, we fix $\ep>0$. Note for any $u\in l^2(\N)$ we have that $\|u\|_2^2 = \sum_{n=1}^\i u_n^2 <\i$ and this entails that $\|u_{>k}\|_2 \leq \|u\|_2$ for all $k\in \N$  and
    \begin{align*}
    \lim_{k\to \i} \|u_{>k}\|_2^2 =\lim_{k\to \i}\sum_{i=1}^\i u_{i+k}^2 = \lim_{k\to \i}\sum_{i=k+1}^\i u_{i}^2=0.
    \end{align*}
    Furthermore since $T:\cX\to l^2(\N)$ is an linear isometry, an application of the abstract change of variable theorem gives us that
    \begin{align*}
    \int \|u\|_2 \, dT(\mu_1)(u) &\leq \int \|T(x)\|_2 \, d\mu_1(x) = \int \|x\|_\cX \, d\mu_1(x) = \int d_\cX(x,0) \, d\mu_1(x)<\i,    \end{align*}
    for all $k\in \N$, since $\mu_1\in M^1_1(\cX)$ (finite first moment). As a consequence the Lebesgue dominated convergence theorem yields that $\lim_{k\to\i} c\int \|u_{>k}\|_2 \, dT(\mu_1)(u) = 0$, implying that there exists a $K\in \N$ (dependent on $\ep$) such that
    \begin{align} \label{eq_temp_2442}
    c\int \|u_{>k}\|_2 \, dT(\mu_1)(u) < \ep^2, \quad \quad \forall k \geq K.
    \end{align}
    Now let $A(\ep)$ be the set, given by $ A(\ep)= \{ (u,v) \in l^2(\N) \times \R^\i : |g(u_{>K},v)|\geq \ep \},$ and note that $A(\ep)\in \cB(l^2(\N))\otimes \cB(\R)$, since $g|_{l^2(\N)}$ is $\cB(l^2(\N))\otimes \cB(\R)$-measurable. Furthermore
    \begin{align}
    T(\mu_1) \times \rho (A(\ep)) &\leq \ep^{-1}\|g(u_{>K},v)\|_{L^1(T(\mu_1)\times \rho)} 
    = \ep^{-1} \int E|g(u_{>K},Z)| \, dT(\mu_1)(u)  \notag\\
    &= \ep^{-1} c \int \|u_{>K}\|_2 \, dT(\mu_1)(u) < \ep, \label{eq_temp_908}
    \end{align}
    where we used Markov's inequality, Tonelli's theorem, \cref{ep_temp_99786} and \cref{eq_temp_2442}. For any $v\in\R^\i$ we denote the section set  $A(\ep,v) = \{ u\in l^2(\N) : |g(u_{>K},v)|\geq \ep\}$. By theorem 3.4.1 \cite{MeasureTheory2007bogachev} we have that $v\mapsto T(\mu_i1)(A(\ep,v))$ is $\cB(\R^\i)/\cB(\R)$-measurable and
    \begin{align*}
    T(\mu_1) \times \rho (A(\ep)) = \int T(\mu_1)(A(\ep,v)) \, d\rho(v).
    \end{align*}
    If $\rho  (v\in \R^\i : T(\mu_1)(A(\ep,v))\geq \ep)=1$, then we have that $T(\mu_1) \times \rho (A(\ep)) \geq \int \ep \, d\rho(v) =\ep$. This is in contradiction with \cref{eq_temp_908}, hence $\rho  (v\in \R^\i : T(\mu_1)(A(\ep,v))< \ep)>0$. Thus there exists a $\rho$-positive probability set $E\in \cB(\R^\i)$ such that $T(\mu_1)(A(\ep,v))<\ep$ for all  $v\in E$. Now note that
    \begin{align*} 
    \rho(E\cap H_K) = \rho(E) > 0,
    \end{align*}
    proving that $E\cap H_K$ is non-empty. Hence there exists a $v\in E\cap H_K\subset \R^\i$ such that \begin{align}
    \label{eq_g} T(\mu_1)(\{ u\in l^2(\N) : |g(u_{>K},v)|\geq \ep\})<\ep
    \end{align}
    and 
    \begin{align} \label{eq_temp_23337}
    &T(\mu_1)(u\in l^2(\N) :\la u_{\leq K}, s \ra_{\R^k} + g(u_{>K},v)\leq y) \\
    =&T(\mu_2)(u\in l^2(\N) :\la u_{\leq K}, s \ra_{\R^k} + g(u_{>K},v)\leq y), \notag
    \end{align}
    for all $y\in \R$ and $s\in \R^K$. Moreover these two properties also implies that 
    $$T(\mu_2)(\{ u\in l^2(\N) : |g(u_{>K},v)|\geq \ep\})<\ep.$$
    To see this, note that we may take any sequence $(\ep_n)_{n\geq 1}\subset \R$ such that $\ep_n\uparrow \ep$. It holds that $( u\in l^2(\N) : g(u_{>K},v)\leq \ep_n) \subset ( u\in l^2(\N) : g(u_{>K},v)\leq \ep_{n+1})$ for any $n\in \N$, so continuity from below and \cref{eq_temp_23337} (with $s=0$ and $y=\ep_n$) yield that
    \begin{align}
    T(\mu_1)(u\in l^2(\N) : g(u_{>K},v)< \ep) &= T(\mu_1)\lp \bigcup_{n=1}^\i ( u\in l^2(\N) : g(u_{>K},v)\leq \ep_n)  \rp \label{eq_ggg} \\
    &= \lim_{n\to \i } T(\mu_1)( u\in l^2(\N) : g(u_{>K},v)\leq \ep_n)\notag \\
    &= \lim_{n\to \i } T(\mu_2)( u\in l^2(\N) : g(u_{>K},v)\leq \ep_n)  \notag\\
    &= T(\mu_2)(u\in l^2(\N) : g(u_{>K},v)< \ep). \notag
    \end{align}
    Thus by \cref{eq_g,eq_temp_23337,eq_ggg}  we get that
    \begin{align*}
    T(\mu_2)(u\in l^2(\N): |g(u_{>K},v)|\geq \ep ) =& 1- T(\mu_2)(u\in l^2(\N): g(u_{>K},v)< \ep ) \\
    &+T(\mu_2)(u\in l^2(\N): g(u_{>K},v)\leq - \ep ) \\
    =& 1- T(\mu_1)(u\in l^2(\N): g(u_{>K},v)< \ep ) \\
    &+T(\mu_1)(u\in l^2(\N): g(u_{>K},v)\leq - \ep ) \\
    =& T(\mu_1)(u\in l^2(\N): |g(u_{>K},v)|\geq \ep ) \\
    <& \ep,
    \end{align*}
      Now note that, when suppressing $u\in l^2(\N)$,  we get
    \begin{align*}
    F_{i,K,s}(y-\ep)-\ep=
    &T(\mu_i)(u: \la u_{\leq K}, s\ra \leq y-\ep)-\ep \\
    =  &T(\mu_i)([u: \la u_{\leq K}, s\ra \leq y-\ep]\cap  [u : |g(u_{>K},v)| <\ep]) \\
    &+  T(\mu_i)([u: \la u_{\leq K}, s\ra \leq y-\ep]\cap[u : |g(u_{>K},v)| \geq \ep])-\ep \\
    \leq & T(\mu_i)([u: \la u_{\leq K}, s\ra+ g(u_{>K},v) \leq y])  \\
    &+  T(\mu_i)([u : |g(u_{>K},v)| <\ep])-\ep \\
    <& T(\mu_i)([u: \la u_{\leq K}, s\ra+ g(u_{>K},v) \leq y]) \\
    = & T(\mu_i)([u: \la u_{\leq K}, s\ra+ g(u_{>K},v) \leq y]\cap  [u : |g(u_{>K},v)| <\ep]) \\
    &+ T(\mu_i)([u: \la u_{\leq K}, s\ra+ g(u_{>K},v) \leq y]\cap  [u : |g(u_{>K},v)| \geq \ep]) \\
    = & T(\mu_i)([u: \la u_{\leq K}, s\ra- \ep \leq y]) \\
    &+ T(\mu_i)( [u : |g(u_{>K},v)| \geq \ep]) \\
    <& F_{i,K,s}(y+\ep) + \ep,
    \end{align*}
    for any $y\in \R$ and $s\in \R^K$. We note that the expression after the first strict inequality is interchangeable in $i$ for all $s\in \R^K$; by \cref{eq_temp_23337}. As a consequence we have that
    \begin{align*}
     F_{1,K,s}(y-\ep)-\ep <  F_{2,K,s}(y+\ep)+\ep \quad \quad \text{and} \quad \quad  F_{2,K,s}(y-\ep)-\ep <  F_{1,K,s}(y+\ep)+\ep,
    \end{align*}
    for all $y\in \R$ and $s\in\R^K$. Note that $K$ depends on $\ep$, but for any $n\leq K$ it holds that
    \begin{align*}
    T(\mu_i)(u: \la u_{\leq K}, (s,(0,...,0))\ra \leq y-\ep) = T(\mu_i)(u: \la u_{\leq n}, s\ra \leq y-\ep),
    \end{align*}
    for any $s\in\R^n$ and $y\in \R$, so the inequalities also hold for any $n\leq K$. Now fix $n\in \N$  and note that for any $\ep>0$ we may choose $K\geq n$, implying that
   \begin{align*}
   F_{1,n,s}(y-\ep)-\ep <  F_{2,n,s}(y+\ep)+\ep \quad \quad \text{and} \quad \quad  F_{2,n,s}(y-\ep)-\ep <  F_{1,n,s}(y+\ep)+\ep,
   \end{align*}
   for all $y\in \R$, $\ep>0$ and  $s\in \R^n$.  Since $F_{1,n,s}$ and $F_{2,n,s}$ are càdlàg functions, we may let $\ep\downarrow 0$ and get that
    \begin{align} \label{eq_temp_77663}
    F_{1,n,s}(y-) \leq  F_{2,n,s}(y) \quad \quad \text{and} \quad \quad F_{2,n,s}(y-) \leq F_{1,n,s}(y),
    \end{align}
    for all $y\in \R$ and  $s\in \R^n$. Now fix $s\in \R^n$, and let $D_{1,n,s}$ and $D_{2,n,s}$ be the sets of discontinuities of $F_{1,n,s}$ and $F_{2,n,s}$ on $\R$   respectively. For any $y\in \R\setminus (D_{1,n,s}\cup D_{2,n,s})$, \cref{eq_temp_77663} yields that
    \begin{align*}
     F_{1,n,s}(y) & = F_{2,n,s}(y).
    \end{align*}
    For any $y\in D_{1,n,s}\cup D_{2,n,s}$ there exists a sequence $(y_k) \downarrow y$ such that $y_k\not \in D_{1,n,s}\cup D_{2,n,s}$, hence $F_{1,n,s}(y_k)  = F_{2,n,s}(y_k)$ for all $k \in \N$, since there is at most a countable number of discontinuities of càdlàg functions over a finite interval (e.g. $[y,y+1]$). As a consequence of right-continuity of the cumulative distribution functions we get that
    \begin{align*}
    F_{1,n,s}(y) = \lim_{k\to \i}  F_{1,n,s}(y_k) =\lim_{k\to \i}  F_{2,n,s}(y_k) =F_{2,n,s}(y),
    \end{align*}
    proving that $F_{1,n,s}(y)=F_{2,n,s}(y)$ for all $y\in\R$. Note $n\in\N$ and $s$ were arbitrarily chosen, so we conclude that $F_{1,n,s}(y)=F_{2,n,s}(y)$ for all $y\in \R, n\in \N$ and  $s\in \R^n$. We conclude that $\beta_{\phi}$ is injective on $M^1_1(\cX)$, such that $(\cX,d_\cX)$ is of strong negative type.
\end{p}
This concludes the section on distance covariance in metric spaces of strong negative type.
\section{Properties of distance covariance in metric spaces} \label{section_properties_of_dcov}
In this section, we will derive some rudimentary properties of the distance covariance measure in metric spaces. These properties include absolute bounds on the distance covariance measure, and when these bounds are attained.  We will also show that, when the marginal spaces are finite-dimensional Euclidean spaces, then our distance covariance measures coincide with the squared distance covariance measure from \cite{szekely2007measuring}. \\ \\
We \textbf{stress} that distance covariance measure cannot be used to measure any kind of dependence degree. It only serves as a direct indicator of independence or the alternative in metric spaces of strong negative type. \\ \\
As previously, let $(\cX,d_\cX)$ and $(\cY,d_\cY)$ be separable metric spaces and  $X,Y$ be random Borel elements with values in $\cX$ and $\cY$, with simultaneous distribution $(X,Y)\sim \theta\in M^{1,1}_1(\cX\times \cY)$ and marginal distributions $X\sim \mu\in M^1_1(\cX)$ and $Y\sim \nu\in M^1_1(\cY)$. Now recall the alternative representation of $dcov(X,Y):=dcov(\theta)$ from \cref{section_distance_covariance} given by
\begin{align*}
dcov(X,Y) =& Ed_\mu(X,X')d_\nu(Y,Y')\\
=& E\Big[ \Big( d_\cX(X,X') -a_{\mu}(X) - a_{\mu}(X')+ D(\mu) \Big) \\
&\times \Big(  d_\cY( Y \, ,  Y' ) -a_{\nu}(\,Y\, ) - a_{\nu}(Y')+ D(\nu) \Big)\Big],
\end{align*}
where $(X',Y')$ is an independent copy of $(X,Y)$. Similarly we may define $dcov(X,X):=dcov(\mu\times \mu)$. By reasoning similar to that of the derivation of the above representation we get that
\begin{align*}
dcov(X,X) &=Ed_\mu(X,X')^2\\
&= E\Big[ \Big( d_\cX(X,X') -a_{\mu}(X) - a_{\mu}(X')+ D(\mu) \Big)^2 \Big],
\end{align*}
where $X'$ is an independent copy of $X$.
\begin{theorem} \label{theorem_bounds_on_dcov}
	For any random element $(X,Y) \sim \theta \in M_1^{1,1}(\cX\times \cY)$ with marginals $\mu\in M^1_1(\cX)$ and $\nu\in M^1_1(\cY)$, it holds that \begin{align*}
	|dcov(X,Y)| \leq \sqrt{dcov(X,X)dcov(Y,Y)} \leq D(\mu)D(\nu).
	\end{align*}
	Furthermore, if both metric spaces $(\cX,d_\cX)$ and $(\cY,d_\cY)$ are of negative type, then
	\begin{align*}
	dcov(X,Y) \geq 0.
	\end{align*}
\end{theorem}
\begin{p}
	Recall that by \cref{lemma_d_mu_is_squre_integrabel_with_respect_to_product_measure_in_M11} we have that $d_\mu\in \mathcal{L}^2(\mu\times \mu)$ and vice versa for $d_\nu$. As we did below \cref{defi_Distance_covariance} we may view $d_\mu$ and $d_\nu$ as mappings from $(\cX\times \cY)^2$ and say that $d_\mu\in \mathcal{L}^2((\cX\times \cY)^2, \theta \times \theta)$. This allows us to use Cauchy-Schwarz' inequality to get
	\begin{align*}
	|dcov(X,Y)| &\leq \int |d_\mu d_\nu| \, d \theta \times \theta \leq \|d_\mu\|_{2} \|d_\nu \|_{2} = \left(\int d_\mu^2 \, d\theta \times \theta \right)^{\frac{1}{2}} \left(\int d_\nu^2 \, d\theta \times \theta \right)^{\frac{1}{2}}.
	\end{align*}
	Utilizing Tonelli's theorem on both integrals and the fact that the $y$-coordinates in the first integral and $x$-coordinates in the second integral are superfluous, the last expression equals
	\begin{align*}	
	 \left(\int d_\mu^2 \, d\mu \times \mu \right)^{\frac{1}{2}} \left(\int d_\nu^2 \, d\mu \times \mu \right)^{\frac{1}{2}} &=\sqrt{E(d_\mu(X,X')d_\mu(X,X'))E(d_\nu(Y,Y')d_\nu(Y,Y')))} \\&= \sqrt{dcov(X,X)dcov(Y,Y)},
	\end{align*}
	where $(X',X')$ is an independent copy of $(X,X)$ and $(Y',Y')$ is an independent copy of $(Y,Y)$. This proves the first inequality of the theorem. \\ \\
	Recall the two inequalities from \cref{eq_showing_that_|d(xy)-a_mu(x)|<=a_mu(y)} in \cref{lemma_d_mu_is_squre_integrabel_with_respect_to_product_measure_in_M11}: $		d_\cX(x,y) \leq a_\mu(x) + a_\mu(y)$ and $a_\mu(x) \leq d_\cX(x,y) + a_\mu(y)$. These inequalities show that $|d_\cX(x,y)-a_\mu(x)|\leq a_\mu(y)$, and therefore we have that
	\begin{align*}
	E|(d_\cX(X,X')-a_\mu(X))a_\mu(X)|\leq E|a_\mu(X')a_\mu(X)|=E a_\mu(X')E a_\mu(X) =D(\mu)^2<\i.
	\end{align*}
	Thus Fubini's theorem yields that
	\begin{align*}
	E[(d_\cX(X,X')-a_\mu(X))a_\mu(X)] &= \int \int (d_\cX(x,y)-a_\mu(y))a_\mu(y)\, d\mu(x) \, d\mu(y) =0,
	\end{align*}
	and analogously  $E[(d_\cX(X,X')-a_\mu(X'))a_\mu(X')]=0$.
	Thus using the representation from above and expanding, we get that 
	\begin{align*}
	dcov(X,X) =& E \Big[ d_\cX(X,X')^2+a_\mu(X)^2+a_\mu(X')^2+D(\mu)^2 \\
	&- 2d_\cX(X,X')a_\mu(X) - 2d_\cX(X,X')a_\mu(X') )-2D(\mu)a_\mu(X)-2D(\mu)a_\mu(X') \\
	&+ 2d_\cX(X,X')D(\mu) +2a_\mu(X)a_\mu(X') \Big]\\
	=& E\Big[ d_\cX(X,X')^2-d_\cX(X,X')a_\mu(X) -d_\cX(X,X')a_\mu(X')+D(\mu)^2  \\
	&- (d_\cX(X,X')-a_\mu(X))a_\mu(X)) - (d_\cX(X,X')-a_\mu(X'))a_\mu(X'))\\
	&-2D(\mu)a_\mu(X)-2D(\mu)a_\mu(X')+2d_\cX(X,X')D(\mu) +2 a_\mu(X)a_\mu(X')\Big].
	\end{align*} 
	Since $d_\cX(X,X')\leq a_\mu(X)+a_\mu(X')$ we have that $d_\cX(X,X')^2\leq d_\cX(X,X')(a_\mu(X)+a_\mu(X'))$, proving that the first three terms have an upper bound of zero. The fifth and sixth term was shown above to have expectation zero. Thus using linearity of the expectation (all individual terms are integral)
	\begin{align*}
	dcov(X,X) \leq&  D(\mu^2) -2D(\mu)(Ea_\mu(X)+Ea_\mu(X')-Ed_\cX(X,X'))+2Ea_\mu(X)a_\mu(X')\\
	=&D(\mu)^2,
	\end{align*}
	where we have used that $X\independent X'$, so $D(\mu)=Ed_\cX(X,X')=Ea_\mu(X)=\sqrt{Ea_\mu(X)a_\mu(X')}$. By similar arguments we also get that $dcov(Y,Y)\leq D(\nu)^2$, proving the second inequality of the theorem.\\ \\
	As regards the last inequality, assume that $(\cX,d_\cX)$ and $(\cY,d_\cY)$ are of negative type. By \cref{theorem_equivalence_of_negative_type} there exist isometric embeddings $\phi:(\cX,d_\cX^{1/2})\to \cH_1$ and $\psi:(\cY,d_\cY^{1/2})\to \cH_2$ into Hilbert spaces. Then by the representation of $dcov$ found in \cref{lemma_beta_of_tensor_exists_and_representation_of_dcov}, we have that
	\begin{align*}
	dcov(X,Y) = 4\|\beta_{\phi\otimes\psi}(\theta-\mu\times \nu)\|_{\cH_1\otimes\cH_2}^2\geq 0.
	\end{align*}
\end{p}

As seen in the above theorem, we have that 
\begin{align*}
dcov(X,Y)\leq\sqrt{dcov(X,X)dcov(Y,Y)}.
\end{align*}
Hence, if either $dcov(X,X)=0$ or $dcov(Y,Y)=0$, then $dcov(X,Y)=0$. Thus for metric spaces of strong negative type, we might have that information only about  $X$ or $Y$, would be sufficient to conclude that $X\independent Y$. There is only one scenario where this is possible, and that is when the marginal distributions are degenerate (concentrated on a singleton), which automatically implies independence. 

However, as we shall see in the below theorem, this is also the case for arbitrary separable metric spaces. This next theorem also entails that $|dcov(X,Y)|$ can only attain the upper bound $D(\mu)D(\nu)$ if both $X$ and $Y$ are concentrated at two points in $\cX$ and $\cY$ respectively. \\ \\
Before stating the above mentioned theorem, we need to prove a lemma regarding the support of Borel probability measures on separable metric spaces.
\begin{lemma} \label{lemma_borel_prob_measure_on_separable_metric_support}
	Every Borel probability measure $\mu$ on a separable metric space has support of full measure. That is, $\mu(\mathrm{supp}(\mu))=1$ and as a consequence $\mu^2(\mathrm{supp}(\mu)^2)=1$ since $\mathrm{supp}(\mu)^2=\mathrm{supp}(\mu^2)$.
\end{lemma}
\begin{p}
	The support $\textrm{supp}(\mu)$ of a Borel probability measure $\mu$ on a metric space $(\cX,d)$ is defined by
	\begin{align*}
	\textrm{supp}(\mu)=\{x\in \cX \, | \, \forall N\in \mathcal{N}_x : \mu(N)>0 \},
	\end{align*}
	where $\mathcal{N}_x$ is the set of all open neighbourhoods of $x$. Hence we also have that the complement of the support of $\mu$ is given by
	\begin{align*}
	\textrm{supp}(\mu)^c = \{x\in \cX \, | \,  \exists N\in \cN_x : \mu(N)=0  \} = \bigcup_{O\in \mathcal{O}(\cX) : \mu(O)=0} O,
	\end{align*}
	where $\mathcal{O}(\cX)$ are the open sets of $\cX$ (this latter representation coincides with the definition in \cite{MeasureTheory2007bogachev2} p. 77).  We obviously have that $\{O\in \mathcal{O}(\cX) : \mu(O)=0\}$ is an open cover of $\textrm{supp}(\mu)^c$ and by separability (see \cite{Billingsley_convergence_1999} section M3) of $\cX$ we know that it has a countable sub-cover $\{O_n: n \in \mathbb{N}\}\subset\{O\in \mathcal{O}(\cX) : \mu(O)=0\} $. Hence by the countable sub-additivity of $\mu$ we get that
	\begin{align*}
	\mu(\textrm{supp}(\mu)^c ) = \mu \lp \bigcup_{n=1}^\i O_n \rp \leq \sum_{n=1}^\i \mu(O_n) = 0,	\end{align*}
	proving that $\textrm{supp}(\mu)$ has full measure. \\ \\
	 We recall from \cref{Appendix_Product_spaces_section} that $\cX\times \cX$ equipped with the product topology is a metrizable topological space - the maximum metric $\rho^{\max}:\cX^2\times \cX^2 \to [0,\i)$ given by $\rho^{\max}(x,y)=d_\cX(x_1,y_1)\lor d_\cX(x_2,y_2)$ induces the product topology. By \cref{theorem_product_sig_alg_is_tensor_when_separable} it also holds that $\cX\times \cX$ is separable and that $\cB(\cX)\otimes \cB(\cX)=\cB(\cX\times \cX)$, implying that $\mu^2=\mu\times \mu$ is a Borel probability measure on the separable metric space $(\cX\times \cX,d_{\rho^{\max}})$, so we also have that $\mu^2(\text{supp}(\mu^2))=1$.  \\ \\
	 Lastly we show that $\textrm{supp}(\mu^2)= \textrm{supp}(\mu)^2$.  The inclusion $\subset$ easily follows from contraposition. Let $x=(x_1,x_2) \in  \cX^2\setminus \text{supp}(\mu)^2$ and assume without loss of generality that $x_1\not \in \text{supp}(\mu)$. Then there exists an $N\in\mathcal{N}_{x_1}$ such that $\mu(N)=0$. We furthermore have that $x\in N\times \cX\in \cB(\cX\times \cX)$  with $\mu^2(N\times \cX)=\mu(N)=0$, proving that $x \in \cX\setminus \text{supp}(\mu^2)$. 
	 
	 The converse inclusion $\supset$ follows by noting that for any $x=(x_1,x_2)\in \text{supp}(\mu)^2$ we have $\forall N_1\in \mathcal{N}_{x_1},N_2\in\mathcal{N}_{x_2}$ that $\mu(N_1),\mu(N_2)>0$. Hence fix $x=(x_1,x_2)\in \text{supp}(\mu)^2$ and note that for any open neighbourhood $N\in \mathcal{N}_x$, the definition of open sets in metric spaces, yields there exists a $\delta>0$ such that the open ball $B_{\rho^{\max}}(x,\delta)\subset N$. It is obvious by the definition of the maximum metric that $ B_{\rho^{\max}}(x,\delta) = B_{d_\cX}(x_1,\delta) \times B_{d_\cX}(x_2,\delta)$,	 hence we have that
	 \begin{align*}
	 \mu^2(N)\geq \mu^2(B_{\rho^{\max}}(x,\delta)) 
	 = \mu (B_{d_\cX}(x_1,\delta) ) \mu ( B_{d_\cX}(x_2,\delta) ) 
	 >0,
	 \end{align*}
	 since $B_{d_\cX}(x_1,\delta)\in \cN_{x_1}$ and $ B_{d_\cX}(x_2,\delta)\in \cN_{x_2}$, proving that $x\in \text{supp}(\mu^2)$.
\end{p}
\begin{remark}
	The assumption that $(\cX,d_\cX)$ is a separable metric space is essential to the above proof. In fact  there exist Borel probability measures on a topological space which have no support; see example 7.1.3 \cite{MeasureTheory2007bogachev2}. Whether or not the topological space considered in example 7.1.3 \cite{MeasureTheory2007bogachev2} is metrizable is not investigated further, but it serves as an indicator that we might run into further trouble if we did not restrict ourselves to separable metric spaces.
\end{remark}
\begin{theorem} \label{theorem_dcov(XX)=0_iff_X_degenerate}
	Let $(\cX,d_\cX)$ be a metric space. For any random element $X\sim \mu\in M^1_1(\cX)$ we have that
	\begin{align*}
	dcov(X,X)=0 \iff X \textit{ is degenerate},
	\end{align*}
	and 
	\begin{align*}
	dcov(X,X)=D(\mu)^2 &\iff X \textit{ is concentrated on at most two points}.
	\end{align*}
\end{theorem}
\begin{p}
		Let $X'$ be an independent copy of $X$ and recall that $dcov(X,X)=dcov(\mu\times \mu)=Ed_\mu(X,X')^2$.\\ \\
		\textit{First equivalence:}	Since $(\cX,d_\cX)$ is separable we have by \cref{lemma_borel_prob_measure_on_separable_metric_support} that $\mu^2(\text{supp}(\mu)^2)=1$, and as a consequence we have that \begin{align*}
	dcov(X,X)=0 \iff 	d_\mu(x_1,x_2)=0 \quad \forall (x_1,x_2)\in \text{supp}(\mu)^2.
	\end{align*}
	To see this note that
	\begin{align*}
	dcov(X,X)=0\iff d_\mu(X,X')\stackrel{a.s.}{=}0 \iff d_\mu(x_1,x_2)=0 \, \, \text{ for $\mu^2$-almost all }(x_1,x_2)\in \cX^2.
	\end{align*}
	The latter set for which the equality must hold may obviously be intersected with another almost sure set at no cost. That is, for $\mu^2$-almost all $(x_1,x_2)\in \text{supp}(\mu)^2$. But for contradiction assume that there exists an $(x_1,x_2)\in \text{supp}(\mu^2)$ where $d_\mu(x_1,x_2)\not = 0$. We recall by definition of the support of $\mu^2$ (see previous lemma) that $\mu^2(N)>0$ for every  open neighbourhood $N\in \mathcal{N}_{(x_1,x_2)}$ of $(x_1,x_2)$. The function $d_\mu$ is continuous since $(x_1,x_2) \mapsto d_\cX(x_1,x_2)$ and $(x_1,x_2)\mapsto a_\mu(x_1),a_\mu(x_2)$ are continuous. Hence there exists a $\delta>0$ such that $d_\mu(x_1,x_2)\not =0$ for all $(x_1,x_2)\in B_{\rho^{\max}}((x_1,x_2),\delta)$. But since $B_{\rho^{\max}}((x_1,x_2),\delta)$ is an open neighbourhood of $(x_1,x_2)$ we have that $\mu^2(B_{\rho^{\max}}((x_1,x_2),\delta))>0$, proving that $d_\mu(x_1,x_2) \not = 0$ with positive probability - a contradiction. \\ \\
	Assume that $X$ is degenerate. That is,  $X=c$ almost surely or equivalently $\text{supp}(\mu)=\{c\}$, for some $c\in \cX$. We obviously have that $a_\mu(x)=d_\cX(c,x)$ for all $x\in \cX$ and $D(\mu)=\int\int  d_\cX(x,y) \, d\mu(x) \, d\mu(y)=d_\cX(c,c)=0$. Thus
	\begin{align*}
	d_\mu(x_1,x_2)=-d_\cX(c,x_1) -d_\cX(c,x_2) =0,
	\end{align*}
	for all $(x_1,x_2)\in \text{supp}(\mu)^2= \{c\}^2$, so $dcov(X,X)=0$.\\ \\
	 Conversely if $dcov(X,X)=0$ or equivalently $d_\mu(x_1,x_2)=0$ for all $(x_1,x_2)\in \text{supp}(\mu)^2$, then we have that $	0=d_\mu(x,x)=-2a_\mu(x)+D(\mu)$,
	proving that $a_\mu(x)=D(\mu)/2$ for all $x\in \text{supp}(\mu)$. As a consequence, we have that $d_\mu(x_1,x_2)=d_\cX(x_1,x_2)$, hence $d_\cX(x_1,x_2)=0$ for all $(x_1,x_2)\in \text{supp}(\mu)^2$. In other words, the distance between any two points in the non-empty support $\text{supp}(\mu)$ is zero, proving that 	$\text{supp}(\mu)$ is a singleton or equivalently that $X$ is degenerate.	\\ \\
	\textit{Second equivalence:} First note that by the above proof we have that, if $X$ is concentrated on a single point, then $dcov(X,X)=0=D(\mu)^2$. Additionally note that, if $X$ is concentrated on two points  $x_1,x_2\in \cX $ or equivalently $\mathrm{supp}(\mu)=\{x_1,x_2\}$, then by direct calculation we see that
	\begin{align*}
	dcov(X,X) =&Ed_\mu(X,X')^2\\
	=& 2p(1-p)d_\mu(x_1,x_2)^2  + p^2d_\mu(x_1,x_1)^2+(1-p)^2d_\mu(x_2,x_2)^2 \\
		=&D(\mu)^2,
	\end{align*}
	where we used the symmetry of $d_\mu$ and assumed that $P(X=x_1)=p$ and $P(X=x_2)=1-p$ for some $p\in(0,1)$. \\ \\
	For the converse, assume for contradiction that $dcov(X,X)=D(\mu)^2$ and $\mathrm{card}(\mathrm{supp}(\mu))\geq 3$. Recall from the proof of \cref{theorem_bounds_on_dcov} that $dcov(X,X)\leq D(\mu)^2$ where the only upper bound we used was 
	\begin{align*}
	d_\cX(X,X')^2\leq d_\cX(X,X')( a_\mu(X) + a_\mu(X')).
	\end{align*}
	Hence we have  $dcov(X,X)=D(\mu)^2$ if and only if $d_\cX(X,X')^2 =d_\cX(X,X')( a_\mu(X) + a_\mu(X'))$ almost surely. In other words, we have equality if and only if $d_\cX(x,y)=a_\mu(x)+a_\mu(y)$ for all $(x,y)\in \mathrm{supp}(\mu^2)=\mathrm{supp}(\mu)^2$ with $x\not = y$. Now fix any $(x,y)\in \mathrm{supp}(\mu)^2$ with $d_\cX(x,y)\not =0$ and note that
	\begin{align*}
	 d_\cX(x,y)= \int   \lp d_\cX(x,z) +d_\cX(y,z) \rp
	    \, d\mu(z) \quad \quad \text{and} \quad \quad d_\cX(x,y)\leq d_\cX(x,z)+d_\cX(y,z) \, \, \forall z\in \cX,
	\end{align*}
	so we must have that $d_\cX(x,y)=d_\cX(x,z)+d_\cY(y,z)$ for $\mu$-almost all $z\in \cX$. As a consequence, it must especially hold that $d_\cX(x,y)=d_\cX(x,z)+d_\cX(y,z)$ for all $z\in \mathrm{supp}(\mu)$, by continuity of $z\mapsto d_\cX(x,z)+d_\cX(y,z)$. Thus
	\begin{align*}
	d_\cX(x,y)=d_\cX(x,z)+d_\cX(y,z) \quad \quad \forall x,y,z\in \mathrm{supp}(\mu) \text{ with } x\not=y\not =z,
	\end{align*}
	and since we assumed that $\mathrm{card}(\mathrm{supp}(\mu))\geq 3$, there exist three such points. That is, there exist three distinct points $x,y,z\in \mathrm{supp}(\mu)$ such that
	\begin{align*}
	d_\cX(x,y)=d_\cX(x,z)+d_\cX(y,z)  \quad \quad \text{and} \quad \quad d_\cX(x,z)=d_\cX(x,y)+d_\cX(z,y). 
	\end{align*}
	By inserting the second equation in the first, we get that
	\begin{align*}
	d_\cX(x,y) =  d_\cX(x,y) + d_\cX(z,y)+ d_\cX(y,z) \iff 2d_\cX(y,z) =0,
	\end{align*}
	a contradiction, proving that $\mathrm{card}(\mathrm{supp}(\my))\leq 2$.
\end{p}
Now to the last item on the agenda of this section, namely proving that the distance covariance measure in metric spaces coincides with the distance covariance from \cite{szekely2007measuring}, when the marginal spaces are finite-dimensional Euclidean spaces.
\begin{theorem} \label{theorem_equivalence_distance_covariance_in_metric_spaces_and_Euclidean}
	Let $(X,Y)\sim \theta\in M^{1,1}_1(\R^n \times \R^m)$ have marginals $X\sim \mu\in M^1_1(\R^n)$ and $Y\sim \nu \in M^1_1(\R^m)$ for some $n,m\in \N$. It then holds that the square root of the  distance covariance measure in metric spaces coincides with the distance covariance in Euclidean spaces from \cite{szekely2007measuring}. That is,
	\begin{align*}
	\sqrt{dcov(X,Y)}=\dcov(X,Y):= \sqrt{\frac{1}{c_nc_m} \int_{\R^{n}\times \R^m} \frac{|\hat{\theta}((t,s))-\hat{\mu}(t)\hat{\nu}(s)|^2}{\|t\|_{\R^n}^{n+1} \|s\|_{\R^m}^{m+1}} d\lambda^n\times \lambda^m (t,s)}.
	\end{align*}
	Here  $\hat{\theta}:\R^{n+m}\to \bC$, $\hat{\mu}:\R^n\to \bC$ and $\hat{\nu}:\R^m \to \bC$ are the characteristic functions corresponding to the probability measures $\theta$, $\mu$  and $\nu$ respectively. $\lambda^k$ is the Lebesgue measure on $\R^k$ and  $c_k= \pi^{(1+k)/2}/\Gamma((1+k)/2)$ for any $k\in \N$.  
\end{theorem}
\begin{p}
First note that $\R^n$ and $\R^m$ are indeed separable Hilbert spaces for any $n,m\in \N$ and therefore \cref{theorem_separable_hilbert_spaces_are_of_strong_negative_type} yields that they are of strong negative type. As a consequence  of \cref{theorem_bounds_on_dcov} we have that $dcov(X,Y)\geq 0$, so $\sqrt{dcov(X,Y)}$ is indeed well-defined. \\ \\
In the proof of theorem 3.27 we saw that with $f_n:\R^n\to \R$ be given by $f_n(s)=\|s\|_{\R^n}^{-(n+1)}$ (set $f_n(0)=0$), then $\phi:(\R^n, d_{\R^n}^{1/2})\to L_\bC^2(\R^n , f_n\cdot \lambda^n)$ given by
\begin{align*}
\phi(x)(t) = \frac{1}{\sqrt{2c_n}}\lp 1-e^{it^\t x}\rp ,
\end{align*}
is a well-defined isometric embedding into the $\bC$-Hilbert space $L_\bC^2(\R^n , f_n\cdot \lambda^n)$. We also saw that the mean embedding of $\phi$ was given by $\beta_\phi(\mu) = \lp  1- \hat{\mu}\rp/\sqrt{2c_n}$ for any $\mu\in M^1_1(\R^n)$. We define the isometric embedding on the other marginal space in an identical fashion, $\psi:(\R^m,d_{\R^m}^{1/2})\to L^2_\bC(\R^m, f_m\cdot \lambda^m)$ given by $\psi(y)(s)=\lp 1-e^{is^\t y}\rp/\sqrt{2c_m}$ which has mean embedding given by $\beta_\psi(\nu)=(1-\hat{\nu})/\sqrt{2c_m}$ for any $\nu\in M^1_1(\R^m)$. Now note that $\cH_1:=L_\bC^2(\R^n , f_n\cdot \lambda^n)$ and $\cH_2:=L^2_\bC(\R^m, f_m\cdot \lambda^m)$ are separable since $\R^n$ and $\R^m$ are separable (cf. theorem 4.13 [Bre10]). Hence we  know that the map $U$ taking simple tensors from  $\cH_1\otimes \cH_2$ to $\cH_3:=L^2_\bC(\R^n\times \R^m , (f_n\cdot \lambda^n )\times (f_m\cdot \lambda^m))$ by
\begin{align*}
U(x\otimes y)(t,s)= x(t)y(s),
\end{align*}
for any $x\in \cH_1$ and $y\in \cH_2$, extends uniquely to a  unitary isomorphism of $L_\bC^2(\R^n , f_n\cdot \lambda^n)\otimes L^2_\bC(\R^m, f_m\cdot \lambda^m)$ onto $L^2_\bC(\R^n\times \R^m , (f_n\cdot \lambda^n )\times (f_m\cdot \lambda^m))$ (see p.  51 [RS72] and theorem 7.16 [Fol95]). That is, we have that 
\begin{align*}
U:&L_\bC^2(\R^n , f_n\cdot \lambda^n)\otimes L^2_\bC(\R^m, f_m\cdot \lambda^m)=\cH_1\otimes \cH_2\\
\to& L^2_\bC(\R^n\times \R^m , (f_n\cdot \lambda^n )\times (f_m\cdot \lambda^m)) = \cH_3,
\end{align*}
is an isometry satisfying
\begin{align*}
\la U(z),U(w) \ra_{\cH_3} = \la z,w \ra_{\cH_1\otimes \cH_2},
\end{align*}
for any $z,w\in \cH_1\otimes \cH_2$ and it has a well-defined inverse $U^{-1}:\cH_3\to\cH_1\otimes \cH_2$. Now note that 
\begin{align*}
\phi(x)(t)\psi(y)(s)= \frac{\lp 1-e^{it^\t x}\rp \lp 1-e^{is^\t y}\rp}{2\sqrt{c_nc_m}}= \frac{ 1-e^{it^\t x}-e^{is^\t y} + e^{i(t,s)^\t (x,y)}}{2\sqrt{c_nc_m}},
\end{align*}
for any $x,t\in \R^n$ and $y,s\in \R^m$. Hence we get that
\begin{align*}
\int \phi(x)(t)\psi(y)(s) \, d\theta(x,y) &= \frac{1-\hat{\mu}(t)-\hat{\nu}(s)+\hat{\theta}(t,s)}{2\sqrt{c_nc_m}},
\end{align*}
for any $t\in \R^n$ and $s\in \R^m$. For notational simplicity in the following arguments, we define the maps $\xi,\tilde{\xi}:\R^n \times \R^m \to \bC$ by
\begin{align*}
\xi(t,s) =  \frac{1-\hat{\mu}(t)-\hat{\nu}(s)+\hat{\theta}(t,s)}{2\sqrt{c_nc_m}} \quad \quad \text{and} \quad \quad \tilde{\xi}(t,s) =  \frac{1-\hat{\mu}(t)-\hat{\nu}(s)+\widehat{\mu\times \nu}(t,s)}{2\sqrt{c_nc_m}}.
\end{align*}
 Now fix any $g\in (\cH_1\otimes \cH_2)^*$ and note that by the Riesz representation theorem, there exists a unique $g\in \cH_1\otimes \cH_2$ such that $g^*(x)=\la x, g \ra_{\cH_1\otimes \cH_2}$. As a consequence we have that the mean embedding $\beta_{\phi\otimes \psi}:M^{1,1}_1(\R^n\times \R^m)\to \cH_1\otimes \cH_2$ fulfils
\begin{align*}
g^*(\beta_{\phi\otimes \psi}(\theta)) &= \int g^*(\phi(x)\otimes \psi(y)) \, d\theta(x,y) \\
&= \int \la \phi(x)\otimes \psi(y) , g \ra_{\cH_1\otimes \cH_2} \, d\theta(x,y) \\
&= \int \la U(\phi(x)\otimes \psi(y)) ,U(g) \ra_{\cH_3} \, d\theta(x,y) \\
&= \int \int  U(\phi(x)\otimes \psi(y))(t,s)  \overline{U(g)(t,s)} \, d(f_n\cdot \lambda^n)\times (f_m\cdot \lambda^m)(t,s) \, d\theta(x,y) \\
&= \int \int  \phi(x)(t)\psi(y)(s) \, d\theta(x,y)  \, \overline{U(g)(t,s)} \, d(f_n\cdot \lambda^n)\times (f_m\cdot \lambda^m)(t,s)  \\
&= \int \frac{1-\hat{\mu}(t)-\hat{\nu}(s)+\hat{\theta}(t,s)}{2\sqrt{c_nc_m}} \overline{U(g)(t,s)} \, d(f_n\cdot \lambda^n)\times (f_m\cdot \lambda^m)(t,s)  \\
& =\left\la  \xi, U(g) \right\ra_{\cH_3} \\
&= \la U(U^{-1}(\xi)),U(g) \ra_{\cH_3} \\
&= \la U^{-1}(\xi),g \ra_{\cH_1\otimes \cH_2} \\
&= g^*(U^{-1}(\xi)).
\end{align*}
Since $g^*\in (\cH_1\otimes \cH_2)^*$ was arbitrarily chosen, the unique defining property of the Pettis integral yields that $\beta_{\phi\otimes \psi}(\theta) =U^{-1}(\xi)\in \cH_1\otimes \cH_2$. In an identical fashion we deduce that $\beta_{\phi\otimes \psi}(\mu\times \nu)= U^{-1}(\tilde{\xi})\in \cH_1\otimes \cH_2$. It obviously holds that $\widehat{\mu\times \nu}(t,s)=\hat{\mu}(t)\hat{\nu}(s)$ for any $(t,s)\in \R^n \times \R^m$, and since $U$ is an isometry, \cref{lemma_beta_of_tensor_exists_and_representation_of_dcov} yields
\begin{align*}
\sqrt{dcov(X,Y)}&= 2\|\beta_{\phi \otimes \psi}(\theta-\mu\times \nu) \|_{\cH_1\otimes \cH_2} \\
&= 2\|\beta_{\phi \otimes \psi}(\theta)-\beta_{\phi \otimes \psi}(\mu\times \nu) \|_{\cH_1\otimes \cH_2} \\
&=2\|U^{-1}(\xi)- U^{-1}(\tilde{\xi}) \|_{\cH_1\otimes \cH_2} \\
&= 2\|U(U^{-1}(\xi))- U(U^{-1}(\tilde{\xi})) \|_{\cH_3} \\
&= 2\|\xi-\tilde{\xi} \|_{\cH_3} \\
&=\sqrt{\frac{1}{c_nc_m} \int \left|\hat{\theta}(t,s)-\hat{\mu}(t)\hat{\nu}(s) \right|^2 \, d(f_n\cdot \lambda^n)\times (f_m\cdot \lambda^m)(t,s)} \\
&=\sqrt{\frac{1}{c_nc_m} \int_{\R^{n}\times \R^m} \frac{|\hat{\theta}((t,s))-\hat{\mu}(t)\hat{\nu}(s)|^2}{\|t\|_{\R^n}^{n+1} \|s\|_{\R^m}^{m+1}} d\lambda^n\times \lambda^m (t,s)},
\end{align*}
which is what we wanted to show.
\end{p} 
\newpage
\begin{remark}
	The above theorem becomes rather trivial, if we assume that $(X,Y)\sim \theta\in M^{2,2}_1(\cX\times \cY)$. In this case we have that
	\begin{align*}
	\int d_\cX(x,o)d_\cY(y,o) \, d\theta(x,y) \leq \int d_\cX(x,o)^2 d\mu(x) \int d_\cY(y,o)^2 \, d\nu(y)<\i,
	\end{align*}
	by the Cauchy-Schwarz inequality. Now note  that (see next section)
	\begin{align*}
	dcov(X,Y)=& E \Big[ \Big( d_\cX(X_1,X_2)-d_\cX(X_1,X_3)-d_\cX(X_2,X_4)+d_\cX(X_3,X_4)\Big) \\
	& \times \Big(d_\cY(Y_1\, ,Y_2\, )-d_\cY(Y_1\, ,Y_5\, )-d_\cY(Y_2\, ,Y_6\, )+d_\cY(Y_5\, ,Y_6\, )  \Big)\Big],	
	\end{align*}
	where $(X_1,Y_1),...,(X_6,Y_6)$ are mutually independent random elements with distribution $\theta$.
	By expanding, we see that it is the expectation of a sum of individually integrable terms (use triangle inequality to get terms on the above form). As a consequence, we may split up the expectation into a sum of expectations. By tirelessly reducing the expression of 16 expectations, one gets that
	\begin{align} \label{eq_temp_229975}
	dcov(X,Y)= Ed_\cX(X,X')d_\cY(Y,Y')+Ed_\cX(X,X')Ed_\cY(Y,Y')-2Ed_\cY(X,X')d_\cY(Y,Y''),
	\end{align}
	where $(X,Y),(X',Y')$ and $(X'',Y'')$ are mutually independent random elements with distribution $\theta$. Hence in the case that $\cX$ and $\cY$ are finite-dimensional  Euclidean spaces we get that $dcov(X,Y)$ coincides with, the square of the distance covariance $\dcov(X,Y)^2$ from \cite{szekely2007measuring}, and the square of the Brownian distance covariance $\mathcal{W}(X,Y)^2$ from \cite{szekely2009brownian} (compare with the expressions in theorem 7 and 8 in \cite{szekely2009brownian}).
\end{remark}
With this remark, we end the section on basic properties of the distance covariance measure.

\section{Asymptotic consistent tests of independence} \label{section_asyp_test_main}
In the previous sections we defined the distance covariance measure $$dcov:M^{1,1}_1(\cX\times \cY)\to \R,$$ and showed that it can be used as a direct indicator of independence, whenever the marginal spaces $\cX$ and $\cY$ are metric spaces of strong negative type. That is, 
\begin{align*}
dcov(\theta)=0 \iff \theta=\mu\times \nu,
\end{align*}
for all $\theta\in M^{1,1}_1(\cX\times \cY)$. However as $dcov(\theta)$ is only know when $\theta$ is, we can not directly use it in the non-parametric independence problem stated in the introduction of this thesis. \\ \\ 
\textit{Let us recall the probabilistic set-up and the non-parametric independence problem}:\\  $(\cX,d_\cX)$ and $(\cY,d_\cY)$ denotes two generic metric spaces, and  $(Z_i)_{i\in \N}=((X_i,Y_i))_{i\in \N}$ is an independent and identically distributed sequence of random Borel elements, defined on a probability space $(\Omega,\F,P)$. It is assumed that, each pair of random elements $Z_i=(X_i,Y_i)$ takes values in the product space $\cX\times \cY$, such that
\begin{align*}
(X_i,Y_i)\sim \theta\in M_1(\cX\times \cY), \quad \quad X_i\sim \mu\in M_1(\cX), \quad \quad Y_i\sim \nu \in M_1(\cY),
\end{align*}
for each $i\in \N$. Now suppose that we are given a finite collection of paired sample points $z_{1,n}=[(x_i,y_i)]_{1\leq i \leq n}$, where each pair $(x_i,y_i)$ is a realization of $(X_i,Y_i)$. Given this collection of samples how can we, without restricting $\theta$ to a specific parametric class of distributions, draw inference on whether to reject the null-hypothesis of independence
\begin{align*}
H_0: \theta = \mu \times \nu,
\end{align*}
in favor of the alternative hypothesis of dependence 
\begin{align*}
H_1:\theta \not = \mu\times \nu.  \\ 
\end{align*}
In this section, we will finally provide an answer this question by constructing estimators of the distance covariance measure $dcov(\theta)$ and utilizing their asymptotic properties to create asymptotically consistent tests of independence.

 The asymptotic properties of the estimators holds for general separable metric spaces, but when constructing the asymptotically consistent tests of independence we will restrict the both marginal spaces to be metric spaces of strong negative type (e.g. separable Hilbert spaces), since we need to utilize that $dcov(\theta)=0\iff \theta=\mu\times \nu$. \\ \\
The construction of these asymptotically consistent tests of independence, is split up into the three following subsections:  \\ \\
\begin{tabularx}{1\linewidth}{>{\raggedleft}p{2cm}X}
	\textit{\Cref{section_estimators}} & We introduce two different estimators for $dcov(\theta)$, which will yield two different statistical tests of independence. It is seen that, $dcov$ is a so-called regular functional, and one may recall that such functionals are the building blocks of the so-called $U$- and $V$-statistic estimators. Our choice of estimators for $dcov(\theta)$ are therefore given by a $U$- and a $V$-statistic estimator.   \vspace{0.01cm}\\
	\textit{\Cref{section_asymp_prop_of_estimators}} & We show that the estimators from \cref{section_estimators} are both strongly consistent and if scaled correctly also possess rather complicated asymptotic distributions, under certain moment conditions of the underlying distribution $\theta$. \vspace{0.5cm}\\
	\textit{\Cref{section_tests}} & We formally describe the statistical models for which the asymptotic properties from \cref{section_asymp_prop_of_estimators} yield asymptotically consistent tests of independence. These tests turns out to have non-traceable rejection thresholds, so we end this last section by describing how one may reasonably bootstrap the rejection thresholds.   \vspace{0.2cm}\\
\end{tabularx}

\newpage 
\subsection{Estimators for the distance covariance measure} \label{section_estimators}
This section is dedicated to defining two different estimators of the distance covariance measure. First we derive an alternative representation of $dcov(\theta)$, and to that extent  define $f_\cX:\cX^4 \to \R$, $f_\cY:\cY^4 \to \R$ and $h:(\cX\times \cY)^6\to\R$ by
\begin{align*}
f_\cX(x_1,x_2,x_3,x_4)&= d_\cX(x_1,x_2) -d_\cX(x_1,x_3) - d_\cX(x_2,x_4) +d_\cX(x_3,x_4), \\
f_\cY(y_1\, ,y_2\, ,y_3,y_4)&= d_\cY(y_1,y_2) -d_\cY(y_1,y_3) - d_\cY(y_2,y_4) +d_\cY(y_3,y_4),
\end{align*}
and
\begin{align*}
&h((x_1,y_1),...,(x_6,y_6)) = f_\cX(x_1,x_2,x_3,x_4)f_\cY(y_1,y_2,y_5,y_6),
\end{align*}
for any $x_1,...,x_6\in \cX$ and $y_1,...,y_6\in \cY$.
\begin{lemma} \label{lemma_h_is_theta^6_integrable}
	For any $\theta \in M_1^{1,1}(\cX\times \cY)$ it holds that $h\in \mathcal{L}^1((\cX\times \cY)^6,\theta^6)$ and that
	\begin{align*}
	d_\mu(x_1,x_2)&= \int f_\cX(x_1,x_2,x_3,x_4) \, d\theta^2((x_3,y_3),(x_4,y_4)), \\
	d_\nu(y_1,y_2)&= \int f_\cY(y_1,y_2,y_5,y_6) \, d\theta^2((x_5,y_5),(x_6,y_6)),	\end{align*}	
	for any $x_1,x_2\in \cX$ and $y_1,y_2\in \cY$.
\end{lemma}
\begin{p}
First let $i$ denote either $\cX$ or $\cY$ and let $w$ be arguments in the corresponding space. Note that the two first inequalities of \cref{lemma_inequality_f} yield
	\begin{align*}
	\frac{|f_i(w_1,w_2,w_3,w_4)|}{2}\leq \min \{d_i(w_1,w_4),d_i(w_2,w_3) \}.
	\end{align*}
	The fact that $h$ is $\bigotimes_{i=1}^6 (\cB(\cX)\otimes \cB(\cY))/\cB(\R)$-measurable  is seen by noting that $h$ is the product of sums where each term is measurable. Using the above upper bound of $|f_i|$ we get  that
	\begin{align*}
	\int |h| \, d\theta^6 &= \int |f_\cX(x_1,x_2,x_3,x_4)f_\cY(y_1,y_2,y_5,y_6)| \, d\theta^6 ((x_1,y_1),...,(x_6,y_6)) \\
	&\leq 4\int d_\cX(x_1,x_4)d_\cY(y_2,y_5) d\theta^6((x_1,y_1),...,(x_6,y_6)) \\
	&= 4 \int d_\cX(x_1,x_4) \, d\theta^2((x_1,y_1),(x_4,y_4))  \int d_\cY(x_2,x_5) \, d\theta^2((x_2,y_2),(x_5,y_5)) \\
	&=4 \|d_\cX\|_{\mathcal{L}^1(\mu\times \mu)}\|d_\cY\|_{\mathcal{L}^1(\nu\times \nu)} \\
	&<\i,
	\end{align*}
	where we used Tonelli's theorem in the second equality, abstract change of variable in the third and \cref{lemma_metric_is_integral_product_of_M_with_finte_first_moment} to bound the last expression. As regards to the two equalities one can easily show, using similar arguments as above, that the two integrals exists. Thus
	\begin{align*}
	\int  f_\cX(x_1,x_2,x_3,x_4) \, d\theta^2((x_3,y_3),(x_4,y_4)) &= \int  f_\cX(x_1,x_2,x_3,x_4) \, d\mu\times \mu(x_3,x_4) \\
	&= d_\cX(x_1,x_2)- a_\mu(x_1)-a_\mu(x_2)+D(\mu) \\
	&= d_\mu(x_1,x_2),
	\end{align*}
	using linearity, Fubini's theorem and that $\mu$ is a probability measure. Analogous arguments yield the equality for $d_\nu(y_1,y_2)$.
\end{p}
An immediate consequence of the above lemma is that 
\begin{align*}
	dcov(\theta)&= \int d_\mu(x_1,x_2)d_\nu(y_1,y_2) \, d\theta\times \theta((x_1,y_1),(x_2,y_2)) \\
	&= \int \Big( \int f_\cX(x_1,x_2,x_3,x_4) \, d\theta^2((x_3,y_3),(x_4,y_4)) \\
	&\quad \times \int f_\cY(y_1,y_2,y_5,y_6) \, d\theta^2((x_5,y_5),(x_6,y_6)) \Big) \, d\theta^2((x_1,y_1),(x_2,y_2)) \\
	&= \int h((x_1,y_1),...,(x_6,y_6))\,  d\theta^6((x_1,y_1),...,(x_6,y_6)),
\end{align*}
that is, $dcov:M_1^{1,1}(\cX\times \cY)\to \R$ is a regular functional with  kernel $h$ of degree $6$. Regular functionals are the building blocks of $U$- and $V$-statistics so it seems quite intriguing to create our estimators using such statistics.
\begin{remark}
	The kernel $h$ is in general not symmetric. To see this, let $\cX=\cY=\R$ be equipped with the Euclidean metric. By insertion we see that
	\begin{align*}
	h((1,1),(3,0),(2,0),(4,0),(0,0),(0,0))=4, 
	\end{align*}
	  but if we permutate the third and fourth argument pairs, we get that\begin{align*}
	    h((1,1),(3,0),(4,0),(2,0),(0,0),(0,0))=0,
	  \end{align*}
	   proving that is $h$ is not symmetric.
\end{remark}
Now we may define the estimators of $dcov(\theta)$. Before doings so, recall that the random empirical measure $\theta_n:\Omega \to M_1^{1,1}(\cX\times \cY)$ of $\theta$, based on the first $n$ samples $Z_{1,n}=((X_1,Y_1),...,(X_n,Y_n))$ of our sample sequence, is defined in the usual way as
\begin{align*}
\theta_n(\omega)(A) = \frac{1}{n}\sum_{k=1}^n \delta_{(X_k(\omega),Y_k(\omega))}(A),
\end{align*}
for any $\omega\in \Omega$ and $A\in \cB(\cX)\otimes \cB(\cY)$
\begin{definition}[Estimators for the distance covariance measure]
	We define the empirical distance covariance as the (in general biased) plug-in estimator $dcov(\theta_n)$. It is easily seen that this estimator is a V-statistic  with non-symmetric kernel $h$ of degree 6 given by
	\begin{align*}
	V_n^6(h,Z_{1,n}) = \frac{1}{n^6} \sum_{i_1=1}^n \cdots \sum_{i_6=1}^n h((X_{i_1},Y_{i_1}),...,(X_{i_6},Y_{i_6})) .
	\end{align*}
	In addition to the V-statistic estimator, we may also consider the corresponding U-statistic with kernel $h$. That is, for a sample size of $n>6$,  the unbiased estimator given by
\begin{align*}
\tilde{U}_n^6(h,Z_{1,n})= \frac{1}{n_{(6)}} \sum_{1\leq i_1 \not = \cdots \not = i_6 \leq n}h((X_{i_1},Y_{i_1}),...,(X_{i_6},Y_{i_6})),
\end{align*}
where $n_{(6)}=n(n-1)\cdots(n-5)$.
\end{definition}
\begin{remark} \label{remark_symmetrized_kernels}
	In \cite{lyons2013distance}, Russell Lyons only considers the V-statistic plug-in estimator $V_n^6(h,Z_{1,n})$. We introduce a second estimator given by the above U-statistic. This new estimator was deviced after discovering that the original moment assumptions in \cite{lyons2013distance} were insufficient to guarantee strong consistency of $V_n^6(h,Z_{1,n})$ (explained in detail in the next section). As we shall see later, the U-statistic estimator $\tilde{U}_n^6(h,Z_{1,n})$ is guaranteed to be strongly consistent, under weaker moment conditions than those for $V_n^6(h,Z_{1,n})$.We refer the reader to \Cref{Appendix_U-statistics,Appendix_V-statistics} for quick introductions to the theory of U- and V-statistics. We will however note that the theory for U- and V-statistics typically works from the outset of symmetric kernels. \\ \\In the case of the above U-statistics we can write it in the regular form with a symmetric kernel (see \cref{Appendix_U-statistics} for explanation). That is,
	\begin{align*}
	\tilde{U}_n^6(h,Z_{1,n})=U_n^6(\bar{h},Z_{1,n}) = {n \choose 6}^{-1} \sum_{1\leq i_1 < \cdots <i_6 \leq n} \bar{h}((X_{i_1},Y_{i_1}),...,(X_{i_6},Y_{i_6})),
	\end{align*}
	where $\bar{h}$ is the symmetrized version of $h$ given by
	\begin{align*}
	\bar{h}(z_1,...,z_6)=\frac{1}{6!} \sum_{\sigma\in \Pi_6} h(z_{\sigma(1)},...,z_{\sigma(6)}),
	\end{align*}
	with $\Pi_6$ being the set of all permutations of $\{1,...,6\}$. So instead of working with the U-statistic $\tilde{U}_n^6(h,Z_{1,n})$ with an (in general) non-symmetric kernel, we can work with $U_n^6(\bar{h},Z_{1,n})$ with symmetric kernel for which most theorems regarding U-statistics are formulated. \\ \\
	 Likewise in the case of V-statistics we may note that, for any $\sigma \in \Pi_6$
	\begin{align*}
	\int h(z) \, d\theta^6 (z)  =\int h(z) \, d\theta^6 (z_{\sigma}) =\int h(z_{\sigma}) \, d\theta^6 (z),
	\end{align*}
	by the abstract change of variable theorem, where $z_\sigma=(z_{\sigma(1)},...,z_{\sigma(6)})$. This shows that the integral of $h$ with respect to $\theta^6$ is invariant under permutations of the integrand's arguments. By Minkowski's inequality this especially implies integrability of $\bar{h}$ with respect to $\theta^6$, but it also shows that 
	\begin{align*}
	dcov(\theta)&=  \int h((x_{1},y_{1}),....,(x_{6},y_{6})) \, d\theta^6((x_1,x_2),...,(x_6,y_6)) \\
	&=  \int \bar{h}((x_{1},y_{1}),....,(x_{6},y_{6})) \, d\theta^6((x_1,x_2),...,(x_6,y_6)). 
	\end{align*}
	for any $\theta \in M_1^{1,1}(\cX\times \cY)$. Since $\theta_n(\omega)\in M_1^{1,1}(\cX\times \cY)$ for any $\omega\in \Omega$ (it is a finitely supported probability measure) we get that 
	\begin{align*}
	V_n^6(h,Z_{1,n})= dcov(\theta_n)=V_n^6(\bar{h},Z_{1,n}),
	\end{align*}
	which proves that instead of working with the V-statistic $V_n^6(h,Z_{1,n})$ with an (in general) non-symmetric kernel, we can work with $V_n^6(\bar{h},Z_{1,n})$ with symmetric kernel for which most theorems regarding V-statistics are formulated. 
	
	Furthermore, the $V$-statistic can easily be rewritten in a form similar to the estimator from \cite{szekely2007measuring}. To this end, note that it obviously also hold that $\theta_n(\omega)\in M_1^{2,2}(\cX\times \cY)$ for any $\omega \in \Omega$, so by utilizing the expression from \cref{eq_temp_229975} we have that
	\begin{align*}
		V_n^6(h,Z_{1,n})=& dcov(\theta_n) \\
		=&  \int d_\cX(x_1,x_2)d_\cY(y_1,y_2) \, d\theta_n\times \theta_n ((x_1,y_1),(x_2,y_2))\\
		+& \int d_\cX(x_1,x_2) \,d\theta_n\times \theta_n ((x_1,y_1),(x_2,y_2)) \int d_\cY(y_1,y_2) \, d\theta_n\times \theta_n ((x_1,y_1),(x_2,y_2))\\
		-&2 \int d_\cX(x_1,x_2)d_\cY(y_1,y_3) \,d\theta_n\times \theta_n \times \theta_n ((x_1,y_1),(x_2,y_2),(x_3,y_3)) \\
		=&  \frac{1}{n^2} \sum_{i=1}^n\sum_{j=1}^n d_\cX(X_i,X_j)d_\cY(Y_i,Y_j) \\
		+&  \lp \frac{1}{n^2}\sum_{i=1}^n \sum_{j=1}^n d_\cX(X_i,X_j) \rp
		 \lp \frac{1}{n^2}\sum_{i=1}^n \sum_{j=1}^n  d_\cY(Y_i,Y_j)\rp
		   \\
		-& \frac{2}{n^3} \sum_{i=1}^n \sum_{j=1}^n  \sum_{k=1}^n d_\cX(X_i,X_j)d_\cY(Y_i,Y_k),
	\end{align*}
	for any $\omega\in \Omega$.
	
\end{remark} \newpage
\subsection{Asymptotic properties of the estimators} \label{section_asymp_prop_of_estimators}
Now we start by justifying the choice of the above estimators for distance covariance. It turns out both estimators are strongly consistent and we are able to derive asymptotic distributions under the null-hypothesis, which allows for the construction of asymptotically consistent statistical tests for independence. 
\begin{theorem}[Strong consistency of estimators] \label{theorem_strong_consistency_of_emperical_dcov}
	If $\theta\in M_1^{1,1}(\cX\times \cY)$, then
	\begin{align*}
	\tilde{U}_n^6(h,Z_{1,n}) \longrightarrow_n dcov(\theta) \quad P\text{-almost surely.}
	\end{align*}
	If it furthermore holds that $E([d_\cX(X_{1},x)d_\cY(Y_{1},y)]^{5/6})<\i$ for any $x\in \cX$ and $y\in \cY$, then
		\begin{align*}
	V_n^6(h,Z_{1,n}) \longrightarrow_n dcov(\theta) \quad P\text{-almost surely.}
	\end{align*}
\end{theorem}
\begin{p}
	We start by showing the almost sure convergence of $V_n^6(h,Z_{1,n})$. By the above \cref{remark_symmetrized_kernels} it suffices to show that 
	 \begin{align*}
	 V_n^6(\bar{h},Z_{1,n})\convas_n E\bar{h}(Z_1,...,Z_6)=dcov(\theta),
	 \end{align*}
	 where $\bar{h}$ is the symmetrized version of $h$ defined in
	  \cref{remark_symmetrized_kernels}. By the strong law of large numbers for V-statistics -  \cref{theorem_SLLN_for_V-statistics} - it suffices to show that
	\begin{align*}
	E|\bar{h}((X_{i_1},Y_{i_1}),...,(X_{i_6},Y_{i_6})|^{\frac{\#\{i_1,...,i_6\}}{6}}<\i,
	\end{align*}
	for all $1\leq i_1 \leq \cdots \leq i_6 \leq 6$. Since $\bar{h}$ is the sum of all permutations of the given indices the above integrability conditions are especially satisfied if $$
	E|h((X_{i_1},Y_{i_1}),...,(X_{i_6},Y_{i_6})|^{\frac{\#\{i_1,...,i_6\}}{6}}<\i,$$
	for all $(i_1,...,i_6)\in \{1,...,6\}^6$; by the sub-additivity of $x\mapsto |x|^p$ for $0<p\leq 1$. Hence let $(i_1,...,i_6)\in \{1,...,6\}^6$ and set $p=\#\{i_1,...,i_6\}/6\in (1/6,1]$. Note that
		\begin{align*}
		E|h(Z_{i_1},...,Z_{i_6})|^p= &E|f_\cX(X_{i_1},X_{i_2},X_{i_3},X_{i_4})f_\cY(Y_{i_1},Y_{i_2},Y_{i_5},Y_{i_6})|^p \\
		\leq& E[d_\cX(X_{i_1},X_{i_4})d_\cY(Y_{i_2},Y_{i_5})]^p \\
		=& E([(d_\cX(X_{i_1},x)+d_\cX(X_{i_4},x))(d_\cY(Y_{i_2},y)+d_\cY(Y_{i_5},y))]^p) \\
		\leq&  E[d_\cX(X_{i_1},x)d_\cY(Y_{i_2},y)]^p+E[d_\cX(X_{i_1},x)d_\cY(Y_{i_5},y)]^p \\
		&+ E[d_\cX(X_{i_4},x)d_\cY(Y_{i_2},y)]^p+E[d_\cX(X_{i_4},x)d_\cY(Y_{i_5},y)]^p,
			\end{align*}
			for some $x\in \cX$ and $y\in \cY$, where we used the first two inequalities of \cref{lemma_inequality_f}. In all terms above  we either have that the indices of the random elements are distinct or coincide. Lets consider any of the above terms separately. If the indices are distinct the assumption of $\theta\in M_1^{1,1}(\cX\times \cY)$ is sufficient for finiteness since $p\leq 1$; e.g. $Ed_\cX(X_{i_1},x)d_\cY(Y_{i_2},y)=a_\mu(x)a_\nu(y)<\i$. However, in the case that they coincide the additional assumption that $$E([d_\cX(X_{1},x)d_\cY(Y_{1},y)]^{5/6})<\i,$$ implies that the $p$'th moments are finite, since $p=\#\{i_1,...,i_6\}/6\leq 5/6$ whenever one or more indices coincide. We conclude that under the assumptions of the theorem, the claimed almost sure convergence of the V-statistic hold. \\ \\
			As regards to the almost sure convergence of the U-statistic, we also note that by the above \cref{remark_symmetrized_kernels},  we have to show almost sure convergence
			\begin{align*}
			\tilde{U}_n^6(h,Z_{1,n})={n\choose 6}^{-1} \sum_{1\leq i_1 < \cdots <i_6\leq n} \bar{h}((X_{i_1},Y_{i_1}),...,(X_{i_6},Y_{i_6})\convas_n E\bar{h}(Z_1,...,Z_6),
			\end{align*}
			where $\bar{h}$ is the symmetrized version of $h$ mentioned before. The strong law of large numbers for U-statistics -  \cref{theorem_SLLN_for_U-statistics} - yields the wanted almost sure convergence if
			\begin{align*}
			E|\bar{h}((X_1,Y_1),...,(X_6,Y_6))|<\i.
			\end{align*}
			We note that in the case of U-statistics opposed to V-statistics, we only need integrability of the kernel when all arguments are independent and this is the reason for weaker moment assumptions. By the triangle inequality, it suffices to show that each term of $\bar{h}((X_1,Y_1),...,(X_6,Y_6))$ is integrable. That is, 
			\begin{align*}
			E|h((X_{i_1},Y_{i_1}),...,(X_{i_6},Y_{i_6})|<\i,
			\end{align*}
			for all $1\leq i_1\not = \cdots \not= i_6 \leq 6$. We note that for any such $i_1,...,i_6$ the arguments in the above expectation are mutually independent copies of $(X_1,Y_1)$ implying that their simultaneous distribution is given by the six-fold product measure $\theta^6$. Hence for any $1\leq i_1\not = \cdots \not= i_6 \leq 6$ we have that
			\begin{align*}
			E|h((X_{i_1},Y_{i_1}),...,(X_{i_6},Y_{i_6})| &= \int |h(z_1,...,z_6)| \, d\theta^6(z_1,...,z_6)=\|h\|_{\mathcal{L}^1((\cX\times \cY)^6,\theta^6)} < \i,
			\end{align*}
			where we applied \cref{lemma_h_is_theta^6_integrable} to ensure finiteness. We conclude that the claimed almost sure convergence of the U-statistic hold.
			
\end{p}
\begin{remark}
	The almost sure convergence of the empirical distance covariance $V_n^6(h,Z_{1,n})$ in the previous theorem is proposition 2.6 in \cite{lyons2013distance} by Russell Lyons. In that paper the almost sure convergence is claimed to hold whenever  $\theta\in M_1^{1,1}(\cX\times \cY)$, i.e. the same conditions  for which we showed the almost sure convergence of the U-statistic estimator. In \cite{lyons2013distance} it is proved that $E|h((X_1,Y_1),...,(X_6,Y_6))|<\i$, after which it is stated that the almost sure convergence follows. The weakest conditions (that I am aware of) under which the SLLN for V-statistics applies are those of \cite{SLLN_for_V-statistics} (see \cref{theorem_SLLN_for_V-statistics}). 
	
	However I am unable to verify the conditions of \cref{theorem_SLLN_for_V-statistics} under the sole assumption that $\theta\in M_1^{1,1}(\cX\times \cY)$. The problem lies within showing sufficient integrability of the kernel $h(Z_{i_1},...,Z_{i_6})$ whenever two or more indices coincide. In \cite{lyons2013distance} $|h|$ is bounded from above by using the triangle inequality on the  factors $|f_\cX|$ and $|f_\cY|$. The two inequalities for $|f_i|$ used in \cite{lyons2013distance} are a subset of all such inequalities, which are as follows
			\begin{align*}
			\frac{|f_i(z_1,z_2,z_3,z_4)|}{2}\leq \left\{
			\begin{array}{lll}
			d_i(z_1,z_4), & d_i(z_2,z_3), & d_i(z_1,z_2) \lor d_i(z_1,z_3), \\
			d_i(z_1,z_2) \lor d_i(z_1,z_4), & d_i(z_1,z_2) \lor d_i(z_2,z_3), & d_i(z_1,z_2) \lor d_i(z_2,z_4), \\
			d_i(z_1,z_4)  \lor d_i(z_1,z_3), & d_i(z_1,z_4)  \lor d_i(z_2,z_3), & d_i(z_1,z_4)  \lor d_i(z_2,z_4), \\
			d_i(z_2,z_3) \lor d_i(z_1,z_3), & d_i(z_2,z_3) \lor d_i(z_2,z_4),  & d_i(z_3,z_4) \lor d_i(z_1,z_3), \\ 
			d_i(z_3,z_4) \lor d_i(z_1,z_4), & d_i(z_3,z_4) \lor d_i(z_2,z_3), &d_i(z_3,z_4) \lor d_i(z_2,z_4).  	\end{array}
			\right.,
			\end{align*}
	 where $i=\cX$ or $i=\cY$ and $z_1,...,z_4$ are elements in the corresponding metric space (see \cref{lemma_inequality_f}). It is easy to see that, whenever the arguments in $h(Z_{i_1},...,Z_{i_6})$ have distinct indices, then the 1st and 2nd (as we used above) or the 12th and 15th (as Lyons used) inequality used on $|f_\cX|$ and $|f_\cY|$ respectively yield independent factors. Hence allowing us to conclude finiteness of $E|h(Z_{i_1},...,Z_{i_6})|$ whenever the marginal distributions of $X_1$ and $Y_1$ have finite first moments. However, in the case that the indices are not distinct, e.g. when $i_1=i_2\not=i_3\not= \cdots \not=i_6$, the above inequalities do not yield independent factors. To see this, note that any combination of the above inequalities used on $|f_\cX|$ and $|f_\cY|$ gives upper bounds dependent on $X_{i_1}$ and $Y_{i_1}$ respectively. That is, the upper bounds for $|f_\cX|$ and $|f_\cY|$ are dependent on $X_{i_1}$ or $X_{i_2}$ and $Y_{i_1}$ or $Y_{i_2}$ respectively. Since $i_1=i_2$ we get that all possible combinations of the above inequalities result in upper bound factors that are mappings of $X_{i_1}$ and $Y_{i_1}$ respectively. Since $X_{1}$ and $Y_1$ are  in general not independent  we get upper bounds which possibly consists of dependent factors. I have been unsuccessful in resolving this matter, hence I assumed the ad-hoc condition  $E([d_\cX(X_{1},x)d_\cY(Y_{1},y)]^{5/6})<\i$.
	 
	 The ad-hoc condition is obviously satisfied if $Ed_\cX(X_{1},x)d_\cY(Y_{1},y)<\i$, which in the case when $\cX=\cY=\R$ are equipped with the Euclidean metric is equivalent to the existence of the covariance $\text{Cov}(X_1,Y_1)$. By the virtue of Cauchy-Schwarz inequality this stronger integrability condition is  also satisfied if $\theta\in M_1^{2,2}(\cX\times \cY)$ \\ \\
	 In personal communication with Russell Lyons he acknowledges that his original conditions are insufficient and recommended that they should be replaced by second moments. However, as seen above the slightly weaker ad-hoc condition which we used suffices. 
\end{remark}
The next order of business is to show results concerning the asymptotic distribution of our estimators under the null-hypothesis, when the sample size $n$ tends to infinity. Before we proceed with this, we introduce some lemmas which will facilitate the following theorem regarding the asymptotic distributions of the estimators.
\begin{lemma}
	For any $\theta\in M_1^{1,1}(\cX\times \cY)$ satisfying  the null-hypothesis $\theta=\mu\times \nu$, we have that the kernel $h$ and its symmetrized version $\bar{h}$ are square integrable with respect to $\theta^6$. That is, we have that $h,\bar{h}\in\mathcal{L}^2((\cX\times \cY)^6,\theta^6)$.
\end{lemma}
\begin{p}
	First note that under the null hypothesis we can factorize the following expectation
	\begin{align*}
	\|h((X_1,Y_1),...,(X_6,Y_6))\|_2 &= [Ef_\cX(X_1,X_2,X_3,X_4)^2 f_\cY(Y_1,Y_2,Y_5,Y_6)^2]^{1/2}  \\ 
	&= \|f_\cX(X_1,X_2,X_3,X_4)\|_2\|f_\cY(Y_1,Y_2,Y_5,Y_6)\|_2.
	\end{align*}
	These two factors are finite, and one can realize this by either using equality one and two from \cref{lemma_inequality_f} on both $|f_\cX|$ and $|f_\cY|$. Alternatively, we note that
	\begin{align*}
	f_\cX(x_1,x_2,x_3,x_4) &= d_\cX(x_1,x_2)-d_\cX(x_1,x_3)-d_\cX(x_2,x_4)+d_\cX(x_3,x_4) \\
	&= d_\mu(x_1,x_2) - d_\mu(x_1,x_3)-d_\mu(x_2,x_4)+d_\mu(x_3,x_4),
	\end{align*}
	for any $x_1,...,x_4\in \cX$, such that
	\begin{align*}
	\|f_\cX(X_1,X_2,X_3,X_4)\|_2 \leq& \|d_\mu(X_1,X_2)\|_2 + \|d_\mu(X_1,X_3)\|_2 \\
	& +\|d_\mu(X_2,X_4)\|_2+\|d_\mu(X_3,X_4)\|_2 <\i,
	\end{align*}
	by Minkowski's inequality, where we used \cref{lemma_d_mu_is_squre_integrabel_with_respect_to_product_measure_in_M11} for finiteness of each of the four terms. This combined with analogous arguments for the finiteness of $\|f_\cY(Y_1,Y_2,Y_5,Y_6)\|_2$, proves that $h\in\mathcal{L}^2(\theta^6)$. \\ \\
	As regards  the symmetrized version $\bar{h}$, we get by Minkowski's inequality that
	\begin{align*}
\|\bar{h}((X_1,Y_1),...,(X_6,Y_6))\|_2 &\leq \frac{1}{6!}\sum_{\sigma\in \Pi_6}\left\| h\lp (X_{\sigma(1)},Y_{\sigma(1)}),...,(X_{\sigma(6)},Y_{\sigma(6)})\rp \right\|_2 \\
&=\|h((X_1,Y_1),...,(X_6,Y_6))\|_2 <\i,
\end{align*}
where we used that all of the terms in the sum over all permutations $\sigma$ of $\{1,..,6\}$ are identically equal to the $\mathcal{L}^2((\cX\times \cY)^6,\theta^6)$-norm of $h$. Hence both $h$ and $\bar{h}$ are square integrable kernels under the null-hypothesis.
\end{p}
The limit distribution of $U$- and $V$-statistics in the case of non-degenerate kernels is given by a normal distribution. However, as we shall see in the following lemma, the kernel for both the $U$- and $V$-statistic estimators for distance covariance measure, is degenerate of order 1, under the null-hypothesis. An implication of this is that, rather than having a nice normal distribution as a limit, we instead get rather complex limit distributions called  Gaussian chaos distributions.
\begin{lemma} \label{theorem_bar_h_is_degenerate_and_identity}
For any $ \theta \in M_{1,nd}^{1,1}(\cX\times \cY) $ satisfying the null-hypothesis $\theta=\mu\times \nu$, it holds that
	\begin{align*}
	\bar{h}^{(2)}((x_1,y_1),(x_2,y_2)) = \frac{1}{15} d_\mu(x_1,x_2)d_\nu(y_1,y_2),
	\end{align*}
	 and as a consequence we have that the symmetric kernel $\bar{h}$   is $\theta$-degenerate of order $1$. 	 
\end{lemma}
\begin{p}
		
	Assume that $\theta = \mu\times \nu\in M_1^{1,1}(\cX\times \cY)$ such that the marginals $\mu$ and $\nu$ are non-degenerate. As a consequence  $E\bar{h}(Z_1,...,Z_n)=dcov(\theta)=0$ and the (so-called $\theta$-canonical) mappings $\bar{h}^{(k)}:(\cX\times \cY)^k \to \R$ for $1\leq k\leq 6$ from \cref{Appendix_various_results_and_defitions_for_U-_and_V-statistics} become $	 \bar{h}^{(1)}(z_1)=\bar{h}_1(z_1)
	$
	and then recursively
	\begin{align*}
	\bar{h}^{(k)}(z_1,..,z_k)= \bar{h}_k(z_1,..,z_k)- \sum_{j=1}^{k-1}\sum_{1\leq i_1 < \cdots <i_{j} \leq k } \bar{h}^{(j)}(z_{i_1},...,z_{i_{j}}),
	\end{align*}
	 for $k=2,...,6$. Here the subscript functions $\bar{h}_k:(\cX\times \cY)^k\to \R$  are conditional expectations
	\begin{align*}
	\bar{h}_k(z_1,...,z_k)&=E\bar{h}(z_1,...,z_k,Z_{k+1},...,Z_6)
	\\
	&\stackrel{\theta^k\text{-}a.e.}{=}E(\bar{h}(Z_1,...,Z_6)|(Z_1,...,Z_k)=(z_1,...,z_k)),
	\end{align*}
	for $k=1,...,6$. \\ \\
	We start by showing that the first $\theta$-canonical mapping is identically zero, i.e.  $\bar{h}^{(1)}(z)=0$ for all $z\in \cX \times \cY$. 	Note that $\bar{h}^{(1)}(z)=E\bar{h}(z,Z_2,...,Z_6) $ and that $\bar{h}$ is the symmetrized version of $h$. Hence all terms of $\bar{h}(z,Z_2,...,Z_6)$ are given by $h$ with $z$ in argument number $1,2,3,4,5$ or $6$ while $Z_2,...,Z_6$ will be placed in one of the $5!$ possible permutation of the remaining arguments. The placement of the random elements $Z_2,...,Z_6$ in the remaining arguments does not matter since $Z_2,...,Z_6$ are independent and identically distributed. That is, we realize that \begin{align*}
	\bar{h}^{(1)}(z) 	&=\frac{1}{6!} \sum_{k=1}^6 \sum_{\sigma\in \Pi_{2,6}} Eh(Z_{\sigma(2)}...,Z_{\sigma(k)},z,Z_{\sigma(k+1)},...,Z_{\sigma(6)})  \\
	&=  \frac{1}{6!} \sum_{k=1}^6  5!E h(Z_{1},...,Z_{k-1},z,Z_{k+1},...,Z_{6}),
	\end{align*}
	where $ \Pi_{2,6}$ is the set of all permutations of $\{2,...,6\}$. Let $\Lambda:(\cX\times \cY)\times \{1,...,6\}$ be  given by	$\Lambda(z,k)= E h(Z_{1},...,Z_{k-1},z,Z_{k+1},...,Z_{6})$ 
	for all $z\in \cX\times \cY$ and $k=1,...,6$. We see that
	\begin{align*}
	\Lambda(z,1)=& E[d_\cX(x,X_2)-d_\cX(x,X_3)-d_\cX(X_2,X_4)+d_\cX(X_3,X_4)] \\
	&\times E[d_\cY(y,Y_2)-d_\cY(y,Y_5)-d_\cY(Y_2,Y_6)+d_\cY(Y_5,Y_6)] \\
	=&[a_\mu(x)-a_\mu(x)-D(\mu)+D(\mu)][a_\nu(y)-a_\nu(y)-D(\nu)+D(\nu)] =0,
	\end{align*}
	and similarly $\lambda(z,2)=0$ for all $z=(x,y)\in (\cX\times \cY)$. For $k=3$ we have that
	\begin{align*}
	\Lambda(z,3)=& E[d_\cX(X_1,X_2)-d_\cX(X_1,x)-d_\cX(X_2,X_4)+d_\cX(x,X_4)] \\
	&\times E[d_\cY(Y_1,Y_2)-d_\cY(Y_1,Y_5)-d_\cY(Y_2,Y_6)+d_\cY(Y_5,Y_6)]\\
	=& [D(\mu) - a_\mu(x)-D(\mu)+a_\mu(x)][D(\nu)-D(\nu)-D(\nu)+D(\nu)]=0,
	\end{align*}
	and similarly $\lambda(z,4)=0$ for all $z=(x,y)\in(\cX\times \cY)$. For $\lambda(\cdot,5)$ and $\lambda(\cdot,6)$ we get mirrored expressions of the cases $k=3$ and $k=4$, all resulting in zero. We conclude that $\Lambda(z,k)=0$ for all $z\in (\cX\times \cY)$ and $k=1,...,6$. Hence $
	\bar{h}^{(1)}=0$, which means that the kernel $\bar{h}$ is at least  $\theta$-degenerate of first order.	As a consequence we have that
	\begin{align*}
	\bar{h}^{(2)}(z_1,z_2) = \bar{h}_2(z_1,z_2) - \bar{h}^{(1)}(z_1) - \bar{h}^{(1)}(z_1)  =\bar{h}_2(z_1,z_2) ,
	\end{align*}
	for all $z_1,z_2\in \cX\times \cY$. With $\delta_{z}$ denoting the Dirac measure at $z\in\cX\times \cY$, we may realize that for any $z_1,z_2\in \cX\times \cY$ 
	\begin{align*}
	\bar{h}^{(2)}(z_1,z_2)  =& E\bar{h}(z_1,z_2,Z_3,Z_{4},Z_5,Z_6) \\
	=&\frac{1}{6!}\sum_{\sigma\in \Pi_6} \int h(u_{\sigma(1)},...,u_{\sigma(6)})  \, d \delta_{z_1} \times \delta_{z_2}\times \theta^4 (u_1,...,u_6)\\
	=&\frac{1}{6!}\sum_{\stackrel{\sigma\in \Pi_6:}{\sigma(1),\sigma(2)\in\{1,2\}}} \int h(u_{\sigma(1)},...,u_{\sigma(6)})  \, d \delta_{z_1} \times \delta_{z_2}\times \theta^4 (u_1,...,u_6) \\
	&+ \frac{1}{6!}\sum_{\stackrel{\sigma\in \Pi_6:}{\lnot (\sigma(1),\sigma(2)\in\{1,2\}})} \int h(u_{\sigma(1)},...,u_{\sigma(6)})  \, d \delta_{z_1} \times \delta_{z_2}\times \theta^4 (u_1,...,u_6).
	\end{align*}
	The latter sum vanishes since every term is zero. 
	
	To see this, note that under the null-hypothesis $\theta=\mu\times \nu$, the fact that $\delta_{z_1}=\delta_{x_1}\times \delta_{y_1}$,  allows Fubini's theorem to factorize the integral of $h$ into two integrals of $f_\cX$ and $f_\cY$ with respect to measures depending on the specific choice of permutation $\sigma$. What we will realize is that in any of the permutations in the latter sum, one of the factor integrals will always be zero. \textit{The arguments are trivial so, if the reader can take the fact that the latter sum vanishes at face value, then the following wall-of-text can be skipped.} 
	
	To that extent, we may note that the permutations in question ($\sigma\in \Pi_6 $ for which it does not hold that $ \sigma(1),\sigma(2)\in\{1,2\}$) will at most allow one of the Dirac measures to act on argument 1 or 2 of the integrand $h$.
	
	First we consider permutations where $\sigma(1)\in\{1,2\}$ or $\sigma(2)\in\{1,2\}$, that is the cases where one of the Dirac measures acts on argument 1 or 2.	Assume that $\sigma(1)\in\{1,2\}$ such that $\delta_{z_1}$ or $\delta_{z_2}$ will act on the first argument of $h$. To further clarify, note that since $\sigma\in \Pi_6$ such that $\lnot(\sigma(1),\sigma(2)\in\{1,2\})$, we know that the Dirac measure not acting on argument 1 must act on arguments $3,4,5$ or $6$. If this latter Dirac measure acts on arguments 3 or 4 we have that the factor integral with integrand $f_\cY$ becomes $
	\int f_\cY(u_1,u_2,u_5,u_6) \, d\delta_{y_{\sigma(1)}} \times \nu^3 (u_1,u_2,u_5,u_6) =a_\nu(y_{\sigma(1)})-a_\nu(y_{\sigma(1)})-D(\nu)+D(\nu )=0$.
	On the other hand, if it instead acts on argument 5 or 6, then the factor integral with integrand $f_\cX$ becomes $
\int f_\cX(u_1,u_2,u_3,u_4) \, d\delta_{x_{\sigma(1)}} \times \mu^3 (u_1,u_2,u_3,u_4) =a_\mu(x_{\sigma(1)})-a_\mu(x_{\sigma(1)})-D(\mu)+D(\mu )=0$. 
	By similar arguments one can realize that one of the two factor integrals is also always zero, if we instead assume that $\sigma(2)\in\{1,2\}$, e.g. in the case that $\sigma(2)\in\{1,2\}$ and the other Dirac measure acts on 3 or 4 then the factor integral with integrand $f_\cY$ becomes $ a_\nu(y_{\sigma(2)})-D(\nu)-a_\nu(y_{\sigma(2)})-D(\nu )=0$.

	It remains to be shown that one of the factor integrals is always zero in the case that both Dirac measures act on argument number 3,4,5 or 6. It suffices to consider two different scenarios: Either both Dirac measures acts on the same argument pair $\{3,4\}$ or $\{5,6\}$ or both Dirac measures acts on different arguments - one from each pair $\{3,4\}$ and $\{5,6\}$.	If both act on the same argument pair $\{3,4\}$ or $\{5,6\}$ we get that the factor integral with integrand $f_\cY$ or $f_\cX$  becomes zero respectively, since they would equal  $Ef_\cY(Y_1,Y_2,Y_5,Y_6)=0$ or $Ef_\cX(X_1,X_2,X_3,X_4)=0$.	Lastly, if both Dirac measures act on different argument pairs, one from each pair $\{3,4\}$ and $\{5,6\}$ we still get zero. Assume that $\delta_{z_1}$ acts on argument 3,  while $\delta_{z_2}$ acts on argument 5 or 6. Then the factor integral of $f_\cX$ becomes $
	\int f_\cX(u_1,u_2,u_3,u_4) \, d\mu^2 \times \delta_{x_1} \times \mu (u_1,u_2,u_3,u_4) = D(\mu)-a_\mu(x_1)-D(\mu)+a_\mu(x_1)=0$ 
	and if $\delta_{z_1}$ instead acted on argument 4 the factor integral would become $D(\mu)-D(\mu)-a_\mu(x_1)+a_\mu(x_1)=0$. Interchanging $\delta_{z_1}$ with $\delta_{z_2}$ in the above considerations, one obtains that one of the factor integrals are zero in the remaining cases (simply interchange $x_1$ with $x_2$ in the above expressions). \\ \\
	Thus we have that
\begin{align*}
\bar{h}_2(z_1,z_2)  =& \frac{1}{6!}\sum_{\stackrel{\sigma\in \Pi_6:}{\sigma(1),\sigma(2)\in\{1,2\}}} \int h(u_{\sigma(1)},...,u_{\sigma(6)})  \, d \delta_{z_1} \times \delta_{z_2}\times \theta^4 (u_1,...,u_6) \\
=& \frac{1}{6!}\bigg( \sum_{\stackrel{\sigma\in \Pi_6:}{\sigma(1)=1,\sigma(2)=2}} Eh(z_1,z_2,Z_3,Z_4,Z_5,Z_6) +\sum_{\stackrel{\sigma\in  \Pi_6:}{\sigma(1)=2,\sigma(2)=1}} Eh(z_2,z_1,Z_3,Z_4,Z_5,Z_6) \bigg),
\end{align*}
Now note that in each of the above sums there are $4!$ identical terms and that
\begin{align*}
Eh(z_1,z_2,Z_3,Z_4,Z_5,Z_6)=&Eh(z_2,z_1,Z_3,Z_4,Z_5,Z_6) \\
=&\lp d_\cX(x_1,x_2)-Ed_\cX(x_1,X_3)-Ed_\cX(x_2,X_4)+Ed_\cX(X_3,X_4) \rp \\
&\times \lp d_\cY(y_1,y_2)-Ed_\cY(y_1,Y_5)-Ed_\cY(y_2,Y_6)+Ed_\cY(Y_5,Y_6)\rp \\
=& d_\mu(x_1,x_2)d_\nu(y_1,y_2),
\end{align*}
implying that 
\begin{align*}
\bar{h}^{(2)}(z_1,z_2) &= \frac{1}{15}d_\mu(x_1,x_2)d_\nu(y_1,y_2),
\end{align*}
proving the wanted equality.\\ \\
Hence it only remains to be shown that the kernel $\bar{h}$ is degenerate or order 1. 	By \cref{Appendix_definition_degenerate_kernel} we need to show that  $\bar{h}^{(1)}(Z_1)=0$ almost surely and  $\bar{h}^{(2)}(Z_1,Z_2)\not = 0$ with positive probability.	We have already shown that $\bar{h}^{(1)}=0$, meaning that $\bar{h}$ is degenerate of at least order 1, so it only remains to be shown that $\bar{h}^{(2)}(Z_1,Z_2)\not = 0$ with positive probability, as to guarantee that $\bar{h}$ is exactly degenerate of order 1. Note that by the equalities shown above
	\begin{align*}
	\bar{h}^{(2)}(Z_1,Z_2)&= \bar{h}_2(\lp X_1,Y_1\rp,\lp X_2,Y_2\rp) =\frac{1}{15}d_\mu(X_1,X_2)d_\nu(Y_1,Y_2),
	\end{align*}
	and for contradiction assume that $\bar{h}^{(2)}(Z_1,Z_2)= 0$ almost surely. By the independence $d_\mu(X_1,X_2)\independent d_\nu(Y_1,Y_2)$ under the null-hypothesis we have that
	\begin{align*}
0&= P(d_\mu(X_1,X_2)d_\nu(Y_1,Y_2)\not=0)\\
&=  P([ d_\mu(X_1,X_2)\not = 0 ] \cap [ d_\nu(Y_1,Y_2)\not=0]) \\
&= P\lp d_\mu(X_1,X_2)\not = 0 \rp P\lp d_\nu(Y_1,Y_2)\not=0\rp.
	\end{align*}
	This shows that at least one of the two factors must be zero almost surely, i.e. $d_\mu(X_1,X_2)=0$ almost surely or $d_\nu(Y_1,Y_2)=0$ almost surely. Assume without loss of generality that $d_\mu(X_1,X_2)=0$ almost surely. By the proof of \cref{theorem_dcov(XX)=0_iff_X_degenerate} this implies that $X_1$ is degenerate, which is a contradiction. We conclude that $\bar{h}^{(2)}(Z_1,Z_2)\not = 0$ with positive probability, proving that $\bar{h}$ is degenerate of order 1.
\end{p}
Before proceeding with the theorem regarding the asymptotic distribution of the estimators, we will continue with a remark containing thorough explanations and analysis of the limiting distribution.
\begin{remark} \label{remark_eigenvalues_of_S_and_limit_dist}
	Assume that the null-hypothesis is satisfied. Now define the linear operator $S:L^2(\cX\times \cY,\cB(\cX)\otimes \cB(\cY),\theta)\to L^2(\cX\times \cY,\cB(\cX)\otimes \cB(\cY),\theta)$  by	\begin{align*}
	S(f)(x,y) =\int d_\mu(x,x')d_\nu(y,y')f(x',y') \, d\theta(x',y'),
	\end{align*}
	 For notational simplicity denote $L^2(\theta)=L^2(\cX\times \cY,\cB(\cX)\otimes \cB(\cY),\theta)$ with norm $\|\cdot \|_2$ induced by the inner product $\la \cdot , \cdot \ra_2:L^2(\theta)\times  L^2(\theta)\to\R$ given by
		\begin{align*}
		\la f,g \ra_2 =  \int f(z)g(z)\, d\theta(z).
		\end{align*}
		First we show that the obviously linear map $S$ is in fact an operator between $L^2(\theta)$ and  $L^2(\theta)$. Under the null-hypothesis $\theta=\mu\times \nu$ we have that
	\begin{align*}
	\|d_\mu d_\nu\|^2_{L^2(\theta^2)}:=& \int |d_\mu(x,x')d_\nu(y,y')|^2 \, d\theta^2((x,y),(x',y')) \\
	=& \int |d_\mu(x,x')|^2 \, d \mu^2(x,x') \int |d_\nu(y,y')|^2 \, d\nu^2(y,y')\\
	<&\i,
	\end{align*}
	by Tonelli's theorem and \cref{lemma_d_mu_is_squre_integrabel_with_respect_to_product_measure_in_M11}. The proof of \cref{lemma_d_mu_is_squre_integrabel_with_respect_to_product_measure_in_M11} can easily be adjusted to show that also $d_\mu(x,\cdot):x'\mapsto d_\mu(x,x')\in \mathcal{L}^2(\mu)$ and $d_\nu(y,\cdot):y'\mapsto d_\nu(y,y')\in \mathcal{L}^2(\mu)$ for all $x\in\cX$ and $y\in\cY$. This especially implies that $d_\mu(x,\cdot)d_\nu(y,\cdot):(x',y')\mapsto d_\mu(x,x')d_\nu(y,y')\in L^2(\theta)$ for all $(x,y)\in\cX\times\cY$ since $\theta=\mu\times \nu$ (formally when looking at $d_\mu(x,\cdot)d_\nu(y,\cdot)$ as an element in $L^2(\theta)$ we consider its equivalent class, but this will not create any confusion). Hence the square integrability of $S(f)$ for any $f\in L^2(\theta)$ follows by noting that
	\begin{align*}
	\|S(f)\|_2^2 &= \int \lv \int d_\mu(x,x')d_\nu(y,y')f(x',y') \, d\theta(x',y')\rv^2 \, d\theta(x,y) \\
	&\leq \int \lp \int |d_\mu(x,x')d_\nu(y,y')f(x',y')| \,d\theta(x',y')\rp^2 \, d\theta(x,y) \\
	&\leq \int  \|d_\mu(x,\cdot)d_\nu(y,\cdot)\|_2^2 \,\|f\|_2^2  \, d\theta(x,y) \\
	&=\|f\|_2^2\int |d_\mu(x,x')d_\nu(y,y')|^2 \, d\theta^2((x,y),(x',y)) \,\\
	&<\i,
	\end{align*}
	by the Cauchy-Schwarz inequality, proving that $S$ is indeed a linear operator on $L^2(\theta)$. In fact we have that $S$ is a bounded linear operator, since $\|S(f)\|_2 \leq \|d_\mu d_\nu\|_{L^2(\theta^2)} \|f\|_2$, by the above inequality. Since $L^2(\theta)$ is a Hilbert space we know that it has an orthonormal basis $\{x_\alpha :\alpha \in A\}$ for some index-set $A$ (cf. theorem 6.29 \cite{applied_analysis_hunter2001}). Furthermore it holds that $L^2(\theta)$ is separable,  since $(\cX\times \cY,\rho^{{\max}})$ is a separable metric space and $\theta$ is a Borel measure on it (cf. theorem 4.13 \cite{brezis2010functional}). Since a separable Hilbert space admits a countable orthonormal basis, we conclude that $A$ is either finite or countably infinite (cf. proposition 2.3.8 \cite{sunder1998functionalanalyssis_spectral_theory}). Note that the bounded linear operator $S$ is a Hilbert-Schmidt operator, if the Hilbert-Schmidt norm $\|S\|_{HS}^2:= \sum_{\alpha\in A} \|S(x_\alpha)\|_2^2<\i$ (definition 1 section 10.6 \cite{dunford1963linear_part2}). Since $A$ is at most countably infinite we can use Tonelli's theorem to interchange the summation and integration in the following way
	\begin{align*}
	\|S\|_{HS}^2 &= \sum_{\alpha\in A} \int |S(x_\alpha)(z)|^2\, d\theta(z)= \int \sum_{\alpha\in A} |S(x_\alpha)(z)|^2\, d\theta(z)\\
	&=  \int\sum_{\alpha\in A} \lv \int d_\mu(x,x')d_\nu(y,y')x_\alpha(x',y') \, d\theta(x',y') \rv^2 \, d\theta(x,y)\\
	&= \int\sum_{\alpha\in A} |\la d_\mu(x,\cdot)d_\nu(y,\cdot),x_\alpha \ra_2|^2 \, d\theta(x,y) = \int \|d_\mu(x,\cdot)d_\nu(y,\cdot)\|_2^2\, d\theta(x,y)\\
	&= \int |d_\mu(x,x')d_\nu(y,y')|^2\, d\theta^2((x,y),(x',y') =\|d_\mu d_\nu\|_{L^2(\theta^2)}^2 <\i,
	\end{align*}
	where in the fourth equality we used Parseval's identity. We conclude that the integral operator $S$ is a Hilbert-Schmidt operator and hence also compact (cf. theorem 6 section 10.6 \cite{dunford1963linear_part2}). Moreover the integral operator is self-adjoint
	\begin{align*}
	\la S(f),g\ra &= \int \int d_\mu(x,x')d_\nu(y,y') f(x'.y') \, d\theta(x',y') g(x,y) \, d\theta(x,y) \\
	&= \int \int d_\mu(x,x')d_\nu(y,y') g(x,y) \, d\theta(x,y) f(x'.y') \, d\theta(x',y') \\
	&= \la f,S(g)\ra.
	\end{align*}
	Thus $S$ is a self-adjoint compact linear operator on a separable Hilbert space, and by the Hilbert-Schmidt theorem (see theorem 6.2.3. \cite{eidelman2004functional} or theorem 8.94 \cite{An_introduction_to_partial_differential_equations_renardy2006}) the set of non-zero eigenvalues counted according to multiplicity $(\lambda_i)$ of $S$ is either finite or countably infinite. Furthermore the eigenvalues may be indexed in absolute descending order  $|\lambda_i|\geq |\lambda_{i+1}|$ and they possess the property that $\lim_{i\to \i}\lambda_i =0$. The set of corresponding eigenfunctions $(e_i)$, i.e. $S(e_i)=\lambda_i e_i$, may be assumed orthonormal. Lastly, the theorem also states that the orthonormal set $(e_i)$ is actually a orthonormal basis for $\text{Range}(S)$. Assume without loss of generality that there are infinitely many non-zero eigenvalues and note that since  $S$ is a self-adjoint Hilbert-Schmidt integral operator on $L^2(\theta)$ with kernel $d_\mu d_\nu\in L^2(\theta\times \theta)$, it satisfies the conditions of exercise 56 \cite{dunford1963linear_part2}, which then states that $
	d_\mu d_\nu((x,y),(x',y'))=d_\mu(x,x')d_\nu(y,y') = \sum_{i\geq 1}\lambda_i e_i(x,y)e_i(x',y'),
$
	where the convergence of the series happens in $L^2(\theta\times \theta)$. That is,
	\begin{align*}
	K_n =E\lp d_\mu(X_1,X_2)d_\nu(Y_1,Y_2)-\sum_{i=1}^n \lambda_i e_i(Z_1)e_i(Z_2)\rp^2 \to_n 0.
	\end{align*}
	The eigenfunctions $e_i$ obviously satisfy the following properties
	\begin{align*}
	Ee_i(Z_1)^2=1, \, \quad Ee_i(Z_1)e_j(Z_1)=0 \text{ for } i\not=j,
	\end{align*}
	and
	\begin{align*}
	E(d_\mu(X_1,X_2)d_\nu(Y_1,Y_2)e_i(Z_1)|Z_2=(x,y)) &= \int d_\mu(x,x')d_\nu(y,y')e_i(x',y') \, d\theta(x',y') \\
	&= S(e_i)(x,y) = \lambda_ie_i(x,y),
	\end{align*}
	$\theta$-almost surely, hence $E(d_\mu(X_1,X_2)d_\nu(Y_1,Y_2)e_i(Z_1)|Z_2)=\lambda_ie_i(Z_2)$ almost surely.	
	Thus when expanding $K_n$ we get	
		\begin{align*}
K_n =& \|d_\mu d_\nu\|_{L^2(\theta^2)}^2 
+ \sum_{i=1}^n\Big(\lambda_i^2Ee_i(Z_1)^2Ee_i(Z_2)^2 
-2  \lambda_iEd_\mu(X_1,X_2)d_\nu(Y_1,Y_2) e_i(Z_1)e_i(Z_2) \Big) \\
=&\|d_\mu d_\nu\|_{L^2(\theta^2)}^2+\sum_{i=1}^n\lambda_i^2 -2 \lambda_i E[E(d_\mu(X_1,X_2)d_\nu(Y_1,Y_2) e_i(Z_1)|Z_2)e_i(Z_2)] \\
=& \|d_\mu d_\nu\|_{L^2(\theta^2)}^2+\sum_{i=1}^n\lambda_i^2 -2\lambda_i^2 Ee_i(Z_2)^2 
=\|d_\mu d_\nu\|_{L^2(\theta^2)}^2-\sum_{i=1}^n\lambda_i^2 ,
	\end{align*}
	where we used that $\lambda_i^2Ee_i(Z_1)e_i(Z_2)e_j(Z_1)e_j(Z_2)=\lambda_i^2Ee_i(Z_1)e_j(Z_1)Ee_i(Z_2)e_j(Z_2)=0$ for $i\not=j$. Since $K_n\to_n 0$ we get that
	\begin{align*}
	\sum_{i=1}^\i \lambda_i^2 = \|d_\mu d_\nu\|_{L^2(\theta^2)}^2<\i.
	\end{align*}
	That is, the sequence of non-zero eigenvalues $(\lambda_i )$ of $S$ repeated according to multiplicity is square summable.\\ \\
	Now, for an independent and identically distributed sequence $(W_i)_{1\geq 1}$ of standard normal distributed random variables, define $L_n = \sum_{i=1}^{n} \lambda_i (W_i^2-1) \in \mathcal{L}^2(\Omega,\F,P)$ 
	for all $n\geq 1$, and let 
	\begin{align*}
	(\Omega,\F,P)\ni \omega \mapsto \sum_{i=1}^{\infty} \lambda_i (W_i(\omega)^2-1) \in \overline{\R},
	\end{align*}
	 denote the pointwise limit as $n$ tends to infinity. This pointwise limit is welldefined and almost surely finite by Khinchin-Kolmogorov's convergence theorem. That is, $L_n$ converges almost surely and in $\mathcal{L}^2$, since 
	 $ \sum_{i=1}^\i \mathrm{Var} (\lambda_i(W_i^2-1) ) = 2\sum_{i=1}^\i \lambda_i^2 <\i $ and $\lambda_i (W_i^2-1)$ has mean zero and finite variance. This also entails that the pointwise limit 
	\begin{align*}
	\sum_{i=1}^{\infty} \lambda_i (W_i^2-1)\in \mathcal{L}_{\bar{\R}}^2(\Omega,\F,P),
	\end{align*} is almost surely finite.  Thus the following limit distribution of $nU_n^6(h,Z_{1,n})$ and $V_n^6(h,Z_{1,n})$ are well-defined distributions on $\R$. By a standard convolution argument, we also see that the distribution is absolutely continuous with respect to the Lebesgue measure, hence the corresponding cumulative distribution function is continuous.	
\end{remark}
\begin{theorem}[Limiting distribution of estimators under the null hypothesis] \label{theorem_asymp_dist_of_estimators}
	If $ \theta \in M_{1,nd}^{1,1}(\cX\times \cY)$ satisfies the null-hypothesis $\theta=\mu\times \nu$, then
	\begin{align*}
	nV_n^6(h,Z_{1,n}) \convd \sum_{i=1}^\i \lambda_i (W_i^2-1) + D(\mu)D(\nu),
	\end{align*}
	and 
	\begin{align*}
	n\tilde{U}_n^6(h,Z_{1,n}) \convd \sum_{i=1}^\i \lambda_i (W_i^2-1),
	\end{align*}
	as $n\to \i$.	Where $(W_n)_{n\in\N}$ is a sequence of independent and identically standard normal distributed random variables, and $(\lambda_i)$ are the eigenvalues counted with multiplicity of the linear operator $S:L^2(\cX\times \cY,\cB(\cX)\otimes \cB(\cY),\theta)\to L^2(\cX\times \cY,\cB(\cX)\otimes \cB(\cY),\theta)$ given by
	\begin{align*}
	S(f)(x,y) =\int d_\mu(x,x')d_\nu(y,y')f(x',y') \, d\theta(x',y').
	\end{align*}
\end{theorem}
\begin{p}
The convergence of the U-statistics follows quite effortlessly from the well-documented limit theorem of U-statistics with 1st order degenerate kernel  - see \cref{appendix_theorem_asymp_dist_of_2_deg_degenerate_U_stat}. The convergence in distribution of the V-statistics is a little more complicated, since this is not a theorem explicitly found in the literature we have referenced. Such limit theorems can be found in e.g. \cite{Ustatistics_in_banach_spaces_Borovskikh}, where the limit distribution is stated in terms of multiple stochastic integrals. In order to avoid the theory of multiple stochastic integrals, we can with a little more work derive the limit distribution of $nV_n^6(h,Z_{1,n})$, using various decomposition theorems and asymptotic properties of U-statistics. \\ \\
First we show the wanted convergence in distribution of the scaled U-statistic $n\tilde{U}_6^6(h,Z_{1,n})=nU_n^6(\bar{h},Z_{1,n})$. Under the null-hypothesis this is a centered U-statistic and since $\bar{h}\in \mathcal{L}^2(\theta^6)$ (by \cref{lemma_h_is_theta^6_integrable}) is a symmetric kernel with  $\theta$-degeneracy of first order, we get that
\begin{align*}
n U_n^6(\bar{h},Z_{1,n})=n (U_n^6(\bar{h},Z_{1,n})-E\bar{h}(Z_1,...,Z_6)) \convd \frac{6(6-1)}{2} \sum_{i=1}^\i \lambda^*_i (W_i^2-1),
\end{align*}
as $n$ tends to infinity, by \cref{appendix_theorem_asymp_dist_of_2_deg_degenerate_U_stat}, where $(\lambda_i^*)$ are the eigenvalues of 
\begin{align*}
S^*:L^2(\theta)&\to L^2(\theta), \\
S^*(f)(x,y) &= \int \bar{h}^{(2)}((x,y),(x',y'))f(x',y') d\theta(x',y) \\
&= \frac{1}{15}\int d_\mu(x,x')d_\mu(y,y')f(x',y') \, d\theta(x',y') \\
&= \frac{1}{15}S(f)(x,y),
\end{align*}
by \cref{theorem_bar_h_is_degenerate_and_identity}. Let $(\lambda_i)$ be all the non-zero eigenvalues of $S$ counted according to its multiplicity and let $(e_i)$ be the corresponding eigenfunctions descriped in the above \cref{remark_eigenvalues_of_S_and_limit_dist}. We note that $S^*(e_i)=\frac{1}{15}S(e_i)=\frac{\lambda_i}{15}e_i$, so if we enumerate $\lambda_i^*=\lambda_i/15$ for all $i\geq 1$, every non-zero eigenvalue  for $S^*$ repeated according to multiplicity will be given by $(\lambda_i^*)$. This is easily seen by observing that if a non-zero eigenvalue of $S^*$ is missing from the list $(\lambda_i^*)$, then there will also be missing a non-zero eigenvalue of $S$ in the list $(\lambda_i)$ - a contradiction. \Cref{remark_eigenvalues_of_S_and_limit_dist} also showed that $\sum_{i\geq 1}\lambda_i (W_i^2-1)$ is almost surely convergent, hence
\begin{align*}
\frac{6(6-1)}{2} \sum_{i=1}^\i \lambda^*_i (W_i^2-1) = 15 \sum_{i=1}^\i \frac{\lambda_i}{15}(W_i^2-1) = \sum_{i=1}^\i \lambda_i(W_i^2-1),
\end{align*}
almost surely, proving the wanted convergence in distribution of the scaled U-statistic $nU_n^6(\bar{h},Z_{1,n})$. \\ \\
Now we will show the claimed convergence in distribution of the scaled V-statistics $nV_n^6(h,Z_{1,n})$. We note that $V_n^6(h,Z_{1,n})=V_n^6(\bar{h},Z_{1,n})$ under the null-hypothesis, is a centered V-statistic with symmetric kernel $\bar{h}$ of degree $6$. Hence by \cref{appendix_lemma_V-statistic_decomposed_into_V-statistics} we decompose it into a linear combination of six V-statistics. The last four of these we furthermore decompose into a linear combination of U-statistics using \cref{appendix_lemma_V-statistic_decomposed_into_U-statistics}. That is,
	  \begin{align}
	  V_n^6(\bar{h},Z_{1,n}) &= \sum_{c=1}^6 {6 \choose c} V_n^c(\bar{h}^{(c)},Z_{1,n}). \notag\\
	  &=6V_n^1(\bar{h}^{(1)},Z_{1,n})+15V_n^2(\bar{h}^{(2)},Z_{1,n})+ \sum_{c=3}^6 {6 \choose c} \sum_{d=1}^c {n \choose d}n^{-c} U_n^d(\bar{h}^{(c)}_{cd},Z_{1,n}), \label{theorem_asymptotic_distribution_eq_1}
	  \end{align}
	   where $\bar{h}^{(c)}:(\cX\times \cY)^c\to \R$ is defined in \cref{Appendix_various_results_and_defitions_for_U-_and_V-statistics}  (or the previous theorem)  and $\bar{h}^{(c)}_{cd}:(\cX\times \cY)^d\to \R$ (which is also found in \cref{Appendix_various_results_and_defitions_for_U-_and_V-statistics}), is a symmetric kernel of degree $d$ defined by \begin{align*}
	  \bar{h}^{(c)}_{cd}(z_1,...,z_d)= \sum_{\stackrel{v_1+\cdots+ v_d=c}{v_j\geq 1}} \frac{c!}{v_1!\cdots v_d!} \bar{h}^{(c)}(\underbrace{z_1,...,z_1}_{v_1\textit{ times}},...,\underbrace{z_d,...,z_d}_{v_d\textit{ times}}),
	  \end{align*}
	  for all $c=3,..,6$ and $d=1,...,c$, e.g. $\bar{h}^{(c)}_{cc}(z_1,...,z_c)=c!\bar{h}^{(c)}(z_1,...,z_c)$ and $\bar{h}^{(c)}_{c1}(z_1)=\bar{h}^{(c)}(z_1,...,z_1)$ We used the convention that the superscript is read first, i.e.  $\bar{h}^{(c)}_{cd}=(\bar{h}^{(c)})_{cd}$. \\ \\
	  We need to find the limiting distribution of $nV_n^6(\bar{h},Z_{1,n})$ and we do this by multiplying $n$ on both sides of \cref{theorem_asymptotic_distribution_eq_1} and showing that the right hand side converges in distribution to the claimed limiting distribution. In the proof of \cref{theorem_bar_h_is_degenerate_and_identity} we showed that $\bar{h}^{(1)}=0$, so more specifically we only need to prove that
	  \begin{itemize}
	  	\item[\textbf{\textit{(1)}}] $15nV_n^2(\bar{h}^{(2)},Z_{1,n}) \convd \sum_{i} \lambda_i (W_i^2-1) + D(\mu)D(\nu)$, as $n\to\i$.
	  	\item[\textbf{\textit{(2)}}] ${n \choose d}n^{1-c} U_n^d(\bar{h}^{(c)}_{cd},Z_{1,n})\convp 0$, as $n\to\i$,  for all $c=3,...,6$ and $d=1,...,c$.
	  \end{itemize}
  which will yield the wanted convergence of $nV_n^6(\bar{h},Z_{1,n})$; by Slutsky's theorem. \\ \\
	  \textbf{\textit{(1)}:}   By the identity of $\bar{h}^{(2)}$ in \cref{theorem_bar_h_is_degenerate_and_identity} we have that
	  \begin{align}
	  15nV_n^2(\bar{h}^{(2)},Z_{1,n}) &=  \frac{1}{n} \sum_{i_1=1}^n \sum_{i_2=1}^n d_\mu(X_{i_1},X_{i_2})d_\nu(Y_{i_1},Y_{i_2}) \notag \\
	   &= \frac{2}{n} \sum_{1\leq i_1 < i_2\leq n}  d_\mu(X_{i_1},X_{i_2})d_\nu(Y_{i_1},Y_{i_2})+ \frac{1}{n} \sum_{i=1}^n  d_\mu(X_{i},X_{i})d_\nu(Y_{i},Y_{i}), \label{temp.eq.decomp.of.15V_n}
	  \end{align}
	  by the symmetry $d_\mu(x,x')d_\nu(y,y')= d_\mu(x',x)d_\nu(y',y)$.  Now note that under the null-hypothesis $ E|d_\mu(X_1,X_1)d_\nu(Y_1,Y_1)| = E|d_\mu(X_1,X_1) |E|d_\nu(Y_1,Y_1) | \leq 6D(\mu)D(\nu) <\i$  by the triangle inequality, and \begin{align*}
	  Ed_\mu(X_1,X_1)d_\nu(Y_1,Y_1)=\lp -2D(\mu)+D(\mu)\rp \lp -2D(\nu)+D(\nu)\rp=D(\mu)D(\nu).
	  \end{align*}
	  Hence   the last term in \cref{temp.eq.decomp.of.15V_n} converges almost surely
	  \begin{align*}
	  \frac{1}{n} \sum_{i=1}^n  d_\mu(X_{i},X_{i})d_\nu(Y_{i},Y_{i}) \convas_n D(\mu)D(\nu),
	  \end{align*}
	  by the regular strong law of large numbers. By Slutsky's theorem it now suffices to show the following convergence in distribution 
	  \begin{align*}
	  \frac{2}{n} \sum_{1\leq i_1 < i_2\leq n}  d_\mu(X_{i_1},X_{i_2})d_\nu(Y_{i_1},Y_{i_2}) \convd_n  \sum_{i} \lambda_i (W_i^2-1),
	  \end{align*}
	  and to that end Slutsky's theorem also yields that it suffices to show the convergence in distribution of the expression in question multiplied by a factor that tends to one in probability,
	  \begin{align*}
	  \lp \frac{n}{n-1}\rp
	  \frac{2}{n} \sum_{1\leq i_1 < i_2\leq n}  d_\mu(X_{i_1},X_{i_2})d_\nu(Y_{i_1},Y_{i_2}) &= n {n\choose 2}^{-1} \sum_{1\leq i_1 < i_2\leq n}  d_\mu(X_{i_1},X_{i_2})d_\nu(Y_{i_1},Y_{i_2}) \\
	  &= n U_n^2(d_\mu d_\nu,Z_{1,n}).
	  \end{align*}
	  Now note that $d_\mu(x,x') d_\nu(y,y')=15\bar{h}^{(2)}((x,y),(x',y'))$ is a completely $\theta$-degenerate (cf. \cref{appendix_corollary_h^(k)_is_completely_degenerate}) symmetric kernel of degree 2 with $Ed_\mu(X_1,X_2) d_\nu(Y_1,Y_2)=dcov(\theta)=0$, hence
	  \begin{align*}
	   n U_n^2(d_\mu d_\nu,Z_{1,n}) \convd_n \sum_{i} \lambda_i (W_i^2-1),
	  \end{align*}
	  by \cref{appendix_theorem_asymp_dist_of_2_deg_degenerate_U_stat}, where $(\lambda_i)$ are the eigenvalues of $S$ counted according to multiplicity. We conclude that the convergence in distribution of \textit{(1)} holds. \\ \\
	  \textbf{\textit{(2)}:}
 First we note that the factor multiplied with the U-statistics has different asymptotic properties depending on the $c$'s and $d$'s. We have that
 	\begin{align*}
 	\lim_{n\to \i }{n\choose d}n^{1-c}  = \lim_{n\to\i  }  \frac{n^{1+d-c}}{d!} =\left\{ \begin{array}{ll}
 	\infty & \textit{if } d=c \\
 	1/d! & \textit{if } d=c-1 \\
 	0 & \textit{if } d<c-1
 	\end{array}\right. .	\end{align*}
 	Hence for any $c=3,...,6$ it suffices to show that

 	\begin{itemize}[leftmargin=+.6in]
 	\item[\textbf{\textit{(2.1)}}] If $d=c$ then $	{n \choose c} n^{1-c} U_n^c(\bar{h}^{(c)}_{cc},Z_{1,n})  \convp_n 0$.
	\item[\textbf{\textit{(2.2)}}] If $d=c-1$ then $
		U_n^{c-1}(\bar{h}^{(c)}_{c(c-1)},Z_{1,n})\convp_n E\bar{h}^{(c)}_{c(c-1)}(Z_1,...,Z_{c-1}) =0$.
	\item[\textbf{\textit{(2.3)}}] If $d<c-1$ then $
 			U_n^d(\bar{h}^{(c)}_{cd},Z_{1,n}) \convp_n E\bar{h}_{cd}^{(c)}(Z_1,...,Z_d)\in \R$. \\
 	\end{itemize}
 		\textbf{\textit{(2.1):}} We note that for any $c\in\{3,...,6\}$	\begin{align} \label{eq_asymptotic_dist_3i}
 		{n \choose c} n^{1-c} U_n^c(\bar{h}^{(c)}_{cc},Z_{1,n}) = c! n^{-(c-1)} \sum_{1\leq i_1< \cdots <i_c\leq n} \bar{h}^{(c)}(Z_{i_1},...,Z_{i_c}),
 		\end{align}
 		and realize that the factor $n^{-(c-1)}$ tends to zero much slower than ${n \choose c}^{-1}$ which is the normalization factor on regular U-statistic type sums. Hence the regular SLLN for U-statistics is insufficient for our purpose. Luckily it turns out that $\bar{h}^{(c)}$ is a degenerate kernel which comes to our aid as centered U-statistics with degenerate kernel are, under certain conditions, guaranteed to converge to zero much faster than regular centered U-statistics. This SLLN for centered U-statistics with degenerate kernels can be found in \cite{SLLN_for_V-statistics} and is also stated in the appendix under \cref{theorem_SLLN_for_degenerate_U-statistics}.  \\ \\
 		Fix $c\in\{3,...,6\}$ and note that by \cref{appendix_corollary_h^(k)_is_completely_degenerate} the kernel $\bar{h}^{(c)}$ of degree $c$ is completely degenerate, i.e. degenerate of order $c-1$ or has rank $r=c$. The reader is encouraged to read the conditions and statement of \cref{theorem_SLLN_for_degenerate_U-statistics} - the SLLN for centered U-statistics with degenerate kernels. Firstly we note that the  order of normalization $c-1$ in \cref{eq_asymptotic_dist_3i} lies within allowed interval $(c-r/2,c)=(c/2,c)$. Hence by \cref{theorem_SLLN_for_degenerate_U-statistics} we have that
 		\begin{align} \label{eq_asymp_dist_degenerate_SLLN_convergence}
 		n^{-(c-1)} \sum_{1 \leq i_1 <\cdots <i_{c}\leq n}[\bar{h}^{(c)}(Z_{i_1},...,Z_{i_{c}}) -E\bar{h}^{(c)}(Z_1,...,Z_c)] \convp_n 0,
 		\end{align}
 		if $E|\bar{h}^{(c)}|^{\frac{r}{(c-1)+r-c}}= E|\bar{h}^{(c)}|^{\frac{c}{c-1}}<\i$. Since $c/(c-1)<2$ it suffices to show that $E|\bar{h}^{(c)}|^{2}<\i$ to ensure the convergence in \cref{eq_asymp_dist_degenerate_SLLN_convergence} . By examining the recursive nature of $\bar{h}^{(c)}$ we see that it entirely consists of linear combinations of $\bar{h}_j$ for $1\leq j\leq c$ with argument spanning over all subsets of $(Z_1,...,Z_c)$ of cardinality $j$. Hence by Minkowski's inequality it sufficies to show that all of the aforementioned terms of the linear combination are square integrable. We note that any such $j$-cardinality subset has distribution $\theta^j$, so we only need to show that  $\bar{h}_j \in\mathcal{L}^2((\cX\times \cY)^j,\theta^j)$ for all $j\in\{1,...,c\}$. To this extend we simply note that 
 			\begin{align*}
 			E\bar{h}_{j}(Z_1,...,Z_j)^{2} &=  EE(\bar{h}(Z_1,...,Z_6)|Z_1,...,Z_j)^2 \\
 			&\leq EE(\bar{h}(Z_1,...,Z_6)^2| Z_1,...,Z_j) \\
 			&= E\bar{h}(Z_1,...,Z_6)^2,
 			\end{align*}
 			by Jensen's conditional  inequality ($x\mapsto x^2$ is convex). By \cref{lemma_h_is_theta^6_integrable} we have that the right-hand side is finite, so convergence in probability stated in \cref{eq_asymp_dist_degenerate_SLLN_convergence} holds. \\ \\ 			
 			To conclude the wanted convergence it suffices to show that $E\bar{h}^{(c)}(Z_1,...,Z_c)=0$. By similar considerations as we initially did  with the above square integrability of $\bar{h}^{(c)}(Z_1,...,Z_c)$ regarding the recursive nature of $\bar{h}^{(c)}$, we note that it suffices to show that  $E\bar{h}_j(Z_1,...,Z_j)=0$ for all $1\leq j\leq c$. Thus note that for any $1\leq j\leq c$, that $\bar{h}_j(z_1,...,z_j)$ is a conditional expectation of $h(Z_1,...,Z_6)$ given $(Z_1,...,Z_j)=(z_1,...,z_j)$, such that
 		\begin{align*}
 		E\bar{h}_j(Z_1,...,Z_j) = EE(\bar{h}(Z_1,...,Z_6)|Z_1,...,Z_j) =E\bar{h}(Z_1,...,Z_6)=0,
 		\end{align*}
 		proving that the claimed convergence in statement (2.1) holds. \\ \\
 		\textbf{\textit{(2.2):}} We will show this by using \cref{theorem_SLLN_for_U-statistics} - the regular strong law of large numbers for U-statistics, to establish that $
 		U_n^{c-1}(\bar{h}^{(c)}_{c(c-1)},Z_{1,n}) \convas_n E\bar{h}^{(c)}_{c(c-1)}(Z_1,...,Z_{c-1})$ and hereafter showing that $E\bar{h}^{(c)}_{c(c-1)}(Z_1,...,Z_{c-1})=0$ for all $c\in\{3,...,6\}$, proving the wanted convergence. \\ \\
 		Fix $c\in \{3,...,6\}$ and note that in order to apply the SLLN for U-statistics  it suffices to show that $E|\bar{h}^{(c)}_{c(c-1)}(Z_1,...,Z_{c-1})|<\i$, since the kernel $\bar{h}^{(c)}_{c(c-1)}$ is symmetric. Recall 
 		\begin{align*}
 		\bar{h}^{(c)}_{c(c-1)}(z_1,...,z_{c-1})= \sum_{\stackrel{v_1+\cdots+ v_{c-1}=c}{v_j\geq 1}} \frac{c!}{v_1!\cdots v_{c-1}!} \bar{h}^{(c)}(\underbrace{z_1,...,z_1}_{v_1\textit{ times}},...,\underbrace{z_{c-1},...,z_{c-1}}_{v_{c-1}\textit{ times}}),
 		\end{align*}
 		and note that any solution to $v_1+\cdots +v_{c-1}=c$ with $v_j\geq 1$ will consist of $v_i=2$ and $v_1,...,v_{i-1},v_{i+1},...,v_{c-1}=1$ for some $i\in \{1,...,c-1\}$. Now for any sequence $z=(z_k)_{k\in \N}$ we define the projection onto the $j$-first without the $i$'th coordinate as
 		\begin{align*}
 		z_{j\setminus i} =\pi_{1,...,i-1,i+1,...,j}(z)= (z_k:k\in \{1,...,j\}\setminus \{i\}),
 		\end{align*}
 		  e.g. $z_{n\setminus 1}=(z_2,...,z_n)$. Since $\bar{h}^{(c)}$ is a symmetric mapping we can always move the two identical arguments up to the first two argument positions, that is
 		\begin{align*}
 		\bar{h}^{(c)}_{c(c-1)}(z_1,...,z_{c-1})&= \frac{c!}{2!} \sum_{i=1}^{c-2} \bar{h}^{(c)}(z_1,...,z_{i-1},z_i,z_i,z_{i+1},...,z_{c-1})\\
 		&=\frac{c!}{2!} \sum_{i=1}^{c-2} \bar{h}^{(c)}(z_i,z_i,z_{(c-1)\setminus i}).
 		\end{align*}
 		 Using the fact that $Z_1,...,Z_{c-1}$ are independent and identically distributed we get that
 		\begin{align*}
 		E|\bar{h}^{(c)}_{c(c-1)}(Z_1,...,Z_{c-1})|&\leq\frac{c!}{2!} \sum_{i=1}^{c-2}E| \bar{h}^{(c)}(Z_i,Z_i,Z_{(c-1)\setminus i})| \\
 		&= \frac{c!(c-2)}{2!}E|\bar{h}^{(c)}(Z_1,Z_1,Z_2,...,Z_{c-1})|,
 		\end{align*}
 		by the triangle inequality and linearity of the expectation. As argued in (2.1) we have that $\bar{h}^{(c)}(Z_1,Z_1,Z_2,...,Z_{c-1})$ is a linear combination of $\bar{h}_j$ for $1\leq j \leq c$ with arguments spanning over all sublists of $(Z_1,Z_1,Z_2,...,Z_{c-1})$ of cardinality $j$. Whenever those sublists of cardinality $j$ have only one occurrence of $Z_1$ the square integrability of $\bar{h}_j(Z_1,...,Z_j)$ shown in (2.1) implies integrability in particular. Hence the only terms of the aforementioned linear combination needing attention are those, where the sublist of cardinality $j$ have both occurrences of $Z_1$. Again by the i.i.d. property of $Z_1,...,Z_{c-1}$ the particular composition of these ordered sublists is not important, implying that we only need to show that $E|\bar{h}_1(Z_1)|,E|\bar{h}_2(Z_1,Z_1)|, E|\bar{h}_3(Z_1,Z_1,Z_2)|,... ,E|\bar{h}_c(Z_1,Z_1,Z_2,...,Z_{c-1})|<\i $.	The first of these expectation is finite since  $\bar{h}_1\in\mathcal{L}^2(\theta)$. For $2\leq j\leq c$ we have that
 		\begin{align*}
 		&E|\bar{h}_j(Z_1,Z_1,Z_2,...,Z_{j-1})| \\
 		=&\int \Big|\int \bar{h}(z_1,z_1,z_2,...,z_{j-1},z_{j+1},...,z_6)\, d\theta^{6-j}(z_{j+1},...,z_6) \Big| \, d\theta^{j-1}(z_1,...,z_{j-1})\\
 		\leq& \int \int |\bar{h}(z_1,z_1,z_2,...,z_{j-1},z_{j+1},...,z_6)|\, d\theta^{6-j}(z_{j+1},...,z_6) \, d\theta^{j-1}(z_1,...,z_{j-1}) \\
 		=& \int |\bar{h}(z_1,z_1,z_2,...,z_{6-1})| \, d\theta^{6-1}(z_1,z_2,...,z_{6-1}) \\
 		=& E|\bar{h}(Z_1,Z_1,Z_2,...,Z_5)|,
 		\end{align*}
 		where we used the triangle inequality for integrals and Tonelli's theorem.
 		Now this upper bound is easily seen finite, by using the triangle inequality on all $6!$ terms of $\bar{h}$. That is, we get finiteness if $		E|h(Z_{i_1},...,Z_{i_6})| <\i$,
 		for all $(i_1,...,i_6)\in\{1,...,5\}^6$ where all but two indices are distinct. We can actually show even stronger integrability, which becomes useful in the proof of statement (2.3). To this extend take any indices $(i_1,...,i_6)\in\{1,...,6\}^6$ and note that under the null-hypothesis the expectation factorizes 
 		\begin{align*}
 		E|h(Z_{i_1},...,Z_{i_6})| &= E|f_\cX(X_{i_1},X_{i_2},X_{i_3},X_{i_4})|E|f_\cY(Y_{i_1},Y_{i_2},Y_{i_5},Y_{i_6})| \\
 		&\leq [ 8 Ed_\cX(X_1,x) ] [ 8Ed_\cY(Y_1,y)] \\
 		&= 16 a_\mu(x)a_\mu(y)<\i,
 		\end{align*}
 		for any $x\in\cX$ and $y\in \cY$, where we used the triangle inequality to say that $|f_\cX(x_1,x_2,x_3,x_4)|\leq 2\sum_{i=1}^4d_\cX(x_i,x)$ with a similar inequality for $|f_\cY|$. Thus we have argued that the SLLN for U-statistics applies and it remains to be shown that the limit, given by the  expectation of the kernel, is zero. That is, \begin{align*}
 		E\bar{h}^{(c)}_{c(c-1)}(Z_1,...,Z_{c-1})=0.
 		\end{align*}
 	 By the above discussion about $\bar{h}^{(c)}_{c(c-1)}$ we have that
 	 \begin{align*}
 	 E\bar{h}^{(c)}_{c(c-1)}(Z_1,...,Z_{c-1})&=\frac{c!}{2!} \sum_{i=1}^{c-2} E\bar{h}^{(c)}(Z_i,Z_i,Z_{(c-1)\setminus i}) = \frac{c!(c-2)}{2!}E\bar{h}^{(c)}(Z_1,Z_1,Z_2,...,Z_{c-1}),
 	 \end{align*}
 	 so it suffices to show that the last expectation is zero. Hence we note that 
 	 \begin{align*}
 	 E\bar{h}^{(c)}(Z_1,Z_1,Z_2,...,Z_{c-2})&= EE(\bar{h}^{(c)}(Z_1,Z_1,Z_2,...,Z_{c-2})|Z_1) = E (\bar{h}^{(c)})_2(Z_1,Z_1) =0,
 	 \end{align*}
 	 where in the last equality we used \cref{appendix_theorem_properties_of_h^(k)} since $2<c$. This concludes the proof of statement (2.2). 
 	  \\ \\
 		\textbf{\textit{(2.3):}}
Fix $c\in\{3,...,6\}$ and let $1\leq d < c-1$. We realize that the wanted convergence  $
U_n^d(\bar{h}^{(c)}_{cd},Z_{1,n}) \convp_n E\bar{h}_{cd}^{(c)}(Z_1,...,Z_d)\in \R$, holds if we can show that the conditions of the SLLN for U-statistics are satisfied. Since  $\bar{h}^{(c)}_{cd}$ is a symmetric kernel of degree $d$, we only need to show that the kernel is integrable, i.e. $E|\bar{h}^{(c)}_{cd}(Z_1,...,Z_d)|<\i$. Hence note that
\begin{align*}
E|\bar{h}^{(c)}_{cd}(Z_1,...,Z_d)|& \leq \sum_{\stackrel{v_1+\cdots+ v_{d}=c}{v_j\geq 1}} \frac{c!}{v_1!\cdots v_{d}!} E|\bar{h}^{(c)}(Z_1^{(v_1)},...,Z_d^{(v_d)})|,
\end{align*}
where $Z_j^{(v_j)}=(Z_j,...,Z_j)\in \lp \cX\times \cY\rp^{v_j}$ for all $1\leq j \leq d$.
By similar considerations as we have done previously, we may note that any of the above terms $\bar{h}^{(c)}(Z_1^{(v_1)},...,Z_d^{(v_d)})$ can be written as a linear combination of $\bar{h}_j$ with arguments spanning over all ordered sublists of $(Z_1^{(v_1)},...,Z_d^{(v_d)})$ with cardinality $j$ (i.e. $j$-size sublists) for all $1\leq j\leq c$. \\  \\
Consider any of the terms in the above sum: $E|\bar{h}^{(c)}(Z_1^{(v_1)},...,Z_d^{(v_d)})|$ for some $v_1+\cdots +v_d=c$ with $v_1,...,v_d\geq 1$, and realize that the term is finite by the triangle inequality, if all individual terms in its linear combination have finite expectation. Thus for any $1\leq j\leq c$ we fix an arbitrary  ordered sublist of $(Z_1^{(v_1)},...,Z_d^{(v_d)})$ of cardinality $j$. We note that this ordered sublist can be written as $(Z_{i_1},...,Z_{i_j})$ for some $1\leq i_1 \leq \cdots \leq i_j\leq d$. It furthermore holds that $\#\{i_1,...,i_j\}=k$ for some $k\in\{1,...,d\}$. Thus we may establish that
\begin{align*}
E|\bar{h}_j(Z_{i_1},...,Z_{i_j})| &= \int |\bar{h}_j(z_{i_1},...,z_{i_j})|\, d\theta^{k}(z_1,...,z_k) \\
&\leq\int \int  |\bar{h}(z_{i_1},...,z_{i_j},w_1,...,w_{6-j})| \, d\theta^{6-j}(w_1,...,w_{6-j}) \, d\theta^{k}(z_1,...,z_k) \\
&= E|\bar{h}(Z_{i_1},...,Z_{i_j},Z_{d+1},...,Z_{6+d-j})|.
\end{align*}	
The kernel  $\bar{h}$ is the symmetrized version of $h$, that is it is a linear combination of $h$ with arguments spanning over every possible permutation of the list $(Z_{i_1},...,Z_{i_j},Z_{d+1},...,Z_{6+d-j})$. We realize that it suffices to show that $E|h(Z_{i_1},...,Z_{i_6})|<\i$ for any $(i_1,...,i_6)\in\{1,...,6\}^6$, which was done in the proof of statement (2.2) above. We conclude that the wanted convergence in statement (2.3) holds. Hence we have argued that statements (2.1), (2.2) and (2.3) hold, implying  the wanted convergence in distribution of our $V$-statistics.
\end{p}
\begin{remark}
	As regards the limit distribution of $nV_n^6(h,Z_{1,n})$, we note that it indeed differs from the claimed limit distribution from Theorem 2.7 \cite{lyons2013distance}. In the proof of that theorem, it is stated (without proof) that $\sum \lambda_i = \int d_\mu(x,x)d_\nu(y,y)\, d\theta(x,y)$, and we have been unable to prove this equality. In case that the equality holds we obviously have that
	\begin{align*}
	\sum_{i=1}^\i \lambda_i &= \int d_\mu(x,x)d_\nu(y,y)\, d\theta(x,y) = E\lp d_\mu(X,X)d_\nu(Y,Y)\rp \\&=	 E(-2a_\mu(X_1)+D(\my))E(-2a_\nu(Y_1)+D(\nu)) \\
	&= \lp -D(\mu)\rp \lp-D(\nu) \rp= D(\mu)D(\nu) <\i,
	\end{align*}
	such that
	\begin{align*}
	\sum_{i=1}^\i \lambda_i (W_i^2-1)+E\lp d_\mu(X,X)d_\nu(Y,Y)\rp = \sum_{i=1}^n \lambda_iW_i^2,
	\end{align*}
	almost surely, showing why $nV_n^6(h,Z_{1,n})\convd_n  \sum_{i=1}^n \lambda_iW_i^2$ in \cite{lyons2013distance}. \\ \\
	Let us try to examine a possible way to arrive at the above equality. If $S$ is a trace class operator, i.e. $\sum_{i=1}^\i |\lambda_i|<\i$, then the trace of $S$ is given by $\text{Tr}(S)=\sum_{i=1}^\i \lambda_i$. Under certain conditions a trace class operator has  trace given by integral of the kernel over the diagonal; see for example exercise 49 in \cite{dunford1963linear_part2} or \cite{essay_compact_operators_bill_casselman}. In the affirmative of the previous conditions we have that $\sum \lambda_i = \int d_\mu(x,x)d_\nu(y,y)\, d\theta(x,y)$. However, we have not even been  successful in affirming that $S$ is of trace class. 
	
	Furthermore exercise 49 in $\cite{dunford1963linear_part2}$ gives conditions for which $S$ is of trace class and has trace given be the integral of kernel over the diagonal. This exercise specifically requires that $S$ is a composition of two Hilbert-Schmidt integral operators, i.e.  our kernel needs to satisfy
	\begin{align*}
	d_\mu(x,x')d_\nu(y,y') = \int A_1((x,y),(s,t))A_2((s,t),(x',y')) \, d\theta(s,t)
	\end{align*}
	for two Hilbert-Schmidt integral operator kernels $A_1$ and $A_2$. We have not been able to prove such a factorization, so we are not able to justify the conditions of this exercise. In the proof of theorem 2.7 \cite{lyons2013distance} there is a reference to \cite{serfling2009approximation} and within this book there is a remark on p. 227 stating a similar trace formula. This remark  refers to exercise 49 in $\cite{dunford1963linear_part2}$, hence this might be what motivated the equality \cite{lyons2013distance} (only speculation). 
	
	In personal communication with Russell Lyons he acknowledges that it is not evident that $S$ is of trace class, so it remains an open problem.
\end{remark}
With this remark, we end this subsection about the asymptotic properties of our estimators.
\newpage
\subsection{Asymptotically consistent tests of independence} \label{section_tests}
In this section we will discuss how to construct an asymptotically consistent statistical test of independence, using the theory derived in the previous sections. The tests we construct have rejection thresholds given by quantiles of unknown distributions, so we will finally show how these thresholds can be consistently bootstrapped.
\subsubsection{Statistical models and specification of tests}

A statistical test of significance level $\alpha\in(0,1)$ is said to be asymptotically consistent at level $\alpha$ if (1) the probability of rejecting a true hypothesis (Type I error) tends to $\alpha$ and (2) the probability of failing to reject a wrong hypothesis (Type II error) tends to zero as the sample size tends to infinity. \\ \\ First we present the general setup of the statistical models in which we can  test the null-hypothesis against its general alternative, using the theory of distance covariance in metric spaces examined in the previous sections. 
\begin{definition}
	Let $(\cX,d_\cX)$ and $(\cY,d_\cY)$ be separable metric spaces of strong negative type and consider the following three  statistical models
	\begin{itemize}
		\item[1)]  The first statistical model is given by the sample space $(\cX\times \cY,\cB(\cX)\otimes \cB(\cY))$ and the non-parametric family of probability measures $\mathcal{P}_1=M_{1,nd}^{1,1}(\cX\times \cY)$.
		\item[2)] The second statistical model is given by the sample space $(\cX\times \cY,\cB(\cX)\otimes \cB(\cY))$ and the non-parametric family of probability measures  $\mathcal{P}_2$  given by the subset of $M_{1,nd}^{1,1}(\cX\times \cY)$ such that every $\theta\in \mathcal{P}_2$ satisfies $\int [d_\cX(x,x')d_\cY(y,y')]^{5/6} \, d\theta(x,y)<\i$ for some $x'\in \cX$ and $y'\in \cY$. 
		\item[3)] The third statistical model is given by the sample space $(\cX\times \cY,\cB(\cX)\otimes \cB(\cY))$ and the non-parametric family of probability measures $\mathcal{P}_3=M_{1,nd}^{2,2}(\cX\times \cY)$.
	\end{itemize}
\end{definition}
Having established the statistical models we now focus on devising an asymptotically consistent statistical test, which can test the null-hypothesis
\begin{align*}
H_0 : \theta =\mu\times \nu \quad \quad \text{ against  the general alternative } \quad \quad 
H_1: \theta \not = \mu\times \nu.
\end{align*}
For the first statistical model we will construct a statistical test with test statistic given by the $U$-statistic estimator of $dcov$ and for the second statistical model we will construct a statistical test with test statistic given by the $V$-statistic estimator of $dcov$. However, note that the three models are nested, $\mathcal{P}_3\subset \mathcal{P}_2 \subset \mathcal{P}_1$. Hence, every test that is asymptotically consistent in the first model is also asymptotically consistent in the second model, and every test that is asymptotically consistent in the second model is also asymptotically consistent in the third model. \\ \\

For all statistical models we assume that $Z=((X_k,Y_k))_{k\in \N}$ are independent pairs of random elements with values in $\cX\times \cY$, all defined on a common probability space $(\Omega,\F,P)$, such that each pair has simultaneous probability distribution $(X_1,Y_1)(P)=\theta\in \mathcal{P}_i$. 

As previously, we denote the  first $n$ sample pairs by $Z_{1,n}=((X_1,Y_1),...,(X_n,Y_n))$ but now we also let $\gamma$ and $\eta$ denote the probability distributions on $\R$ of the limiting variables of the scaled estimators 	$n\tilde{U}_n^6(h,Z_{1,n})$ and $nV_n^6(h,Z_{1,n})$ respectively.   That is,
\begin{align*}
 \sum_{i} \lambda_i (W_i^2-1) \sim \gamma \quad \quad \text{and} \quad \quad \sum_{i} \lambda_i (W_i^2-1) + D(\mu)D(\nu) \sim \zeta,
\end{align*}
and let $F_\gamma$ and $F_\eta$ denote the respective cumulative distribution functions.
Here $(W_n)_{n\in\N}$ is a sequence of independent and identically standard normal distributed random variables, and $(\lambda_i)$ are the eigenvalues counted with multiplicity of the linear operator $S:L^2(\cX\times \cY,\cB(\cX)\otimes \cB(\cY),\theta)\to L^2(\cX\times \cY,\cB(\cX)\otimes \cB(\cY),\theta)$ given by
\begin{align*}
S(f)(x,y) =\int d_\mu(x,x')d_\nu(y,y')f(x',y') \, d\theta(x',y').
\end{align*}
\begin{theorem} 
	Consider the following two statements
	\begin{itemize}
		\item[1)] For any fixed significance level $\alpha\in (0,1)$, we have that the statistical test that rejects the null-hypothesis if
		$$
		n\tilde{U}_n^6(h,Z_{1,n})> q_{1-\alpha}^{(\gamma)}=\inf\{x\in \R: 1-\alpha \leq F_\gamma(x)\},
		$$
		is an asymptotically consistent test of independence at level $\alpha$.
		\item[2)] For any fixed significance level $\alpha\in (0,1)$, we have that the statistical test that rejects the null-hypothesis if
		$$
		nV_n^6(h,Z_{1,n})> q_{1-\alpha}^{(\eta)}=\inf\{x\in \R: 1-\alpha \leq F_\eta(x)\},
		$$		
		is an asymptotically consistent test of independence at level $\alpha$.
	\end{itemize}
Statement \textit{1)} is true in all three of the considered statistical models, but statement \textit{2)} is only guaranteed to be true in the second and third statistical model.
\end{theorem}
\begin{p}
	Let us consider test \textit{1)} in the first statistical model. Let the   test statistic based on the first $n$ sample pairs $Z_{1,n}=((X_1,Y_1),...,(X_n,Y_n))$ be given by the scaled $U$-statistic estimator of the distance covariance measure. That is, the $n$'th test statistic is given by
	\begin{align*}
	n\tilde{U}_n^6(h,Z_{1,n}),
	\end{align*}
	for every $n\in \N$. Under the null-hypothesis $\theta=\mu\times \nu$,   \cref{theorem_asymp_dist_of_estimators} yields that
	\begin{align*}
	n\tilde{U}_n^6(h,Z_{1,n}) \convd_n \sum_{i} \lambda_i (W_i^2-1) \sim \gamma,
	\end{align*}
	where the limiting distribution $\gamma$ is a well-defined probability distribution on $\R$ with continuous distribution function (see \cref{remark_eigenvalues_of_S_and_limit_dist}). On the other hand, if $\theta\not =\mu\times \nu$, then 
	\begin{align*}
	\tilde{U}_n^6(h,Z_{1,n}) \convas_n dcov(\theta)> 0,
	\end{align*}
	by \cref{theorem_strong_consistency_of_emperical_dcov}, \cref{theorem_bounds_on_dcov} and \cref{theorem_dcov_zero_iff_independence}, implying that $n\tilde{U}_n^6(h,Z_{1,n}) \convas_n \i$. Thus we realize that large values of our test statistic $n\tilde{U}_n^6(h,Z_{1,n})$ are in disagreement with the null-hypothesis. It is therefore reasonable to devise a test that rejects the null-hypothesis if the test statistic is observed to be larger than a certain threshold. If we let this threshold be the $(1-\alpha)$-quantile  $q_{1-\alpha}^{(\gamma)}$ of the limit distribution $\gamma$, then we see that
	\begin{align*}
	P\lp n\tilde{U}_n^6(h,Z_{1,n})> q_{1-\alpha}^{(\gamma)}\rp = 1-P\lp n\tilde{U}_n^6(h,Z_{1,n})\leq  q_{1-\alpha}^{(\gamma)}\rp \to_n  1-F_\gamma \lp q_{1-\alpha}^{(\gamma)}\rp =\alpha,
	\end{align*}
	under the assumption that the null-hypothesis $\theta=\mu\times \nu$ is true. This is seen by noting that the above convergence in distribution implies convergence of the cumulative distribution functions in every point (since the limit distribution has a continuous cdf.). Thus the probability of rejecting the null-hypothesis, even though it is true, is asymptotically $\alpha$. Furthermore we see that
	\begin{align*}
	P(n\tilde{U}_n^6(h,Z_{1,n}) \leq  q_{1-\alpha}^{(\gamma)}) = 1 - P(n\tilde{U}_n^6(h,Z_{1,n}) >  q_{1-\alpha}^{(\gamma)}) \to_n 0,
	\end{align*}
	under the assumption that the null-hypothesis is false. This follows from the fact that the above almost sure convergence implies convergence in probability towards infinity. Thus the probability of accepting the null-hypothesis, even though it is false,  is asymptotically zero. We conclude that for any level $\alpha\in (0,1)$,  the test that rejects the null-hypothesis if
	\begin{align*}
	n\tilde{U}_n^6(h,Z_{1,n})> q_{1-\alpha}^{(\gamma)},
	\end{align*}
	where $q_{1-\alpha}^{(\gamma)}$ is the $(1-\alpha)$-quantile of $\gamma$, is an asymptotically consistent test at level $\alpha$ in the first statistical model.

	The asymptotically consistency of the test proposed in \textit{2)}, follows by identical arguments. However, the $V$-statistic estimator is only guaranteed to be strongly consistent in the second and third models, because of the additional moment condition from  \cref{theorem_strong_consistency_of_emperical_dcov}. Thus the test in \textit{2)} is only guaranteed to be asymptotically consistent in the second and third statistical models.
\end{p}
At a first glance one might think we devised statistical tests that are directly usable in practice, but unfortunately one may realize that the proposed thresholds for rejection depends on the specific underlying distribution $\theta$. That is, the eigenvalues $(\lambda_i)$ of the integral operator $S$ are dependent on the specific choice of $\theta$. Thus without knowing $\theta$ we cannot find the eigenvalues $(\lambda_i)$ analytically and as a consequence we have no idea how the $\gamma$ and $\eta$ distributions behave and we especially do not know where the $(1-\alpha)$-quantiles are located.
\subsubsection{Bootstrapping of test thresholds}
Fortunately for us Miguel A. Arcones and Evarist Giné proved in 1992 \cite{arcones1992bootstrap} that the limiting distribution of both degenerate $U$- and $V$-statistics can be consistently bootstrapped. However one needs to be careful when doing this, since the naive bootstrap approach of simply sampling with replacement from the empirical distribution and inserting into $n\tilde{U}_n^6(h,\cdot)$ and $nV_n^6(h,\cdot)$ fails to be consistent in general. In our case the $U$- and $V$-statistic estimators are both $\theta$-degenerate of order 1 (see \cref{theorem_bar_h_is_degenerate_and_identity}) and \cite{arcones1992bootstrap} proves consistency of a bootstrapping approach which utilizes that the asymptotic distribution of such degenerate statistics is solely determined by the leading terms in the Hoeffding decomposition. In the proof of \cref{theorem_asymp_dist_of_estimators} we only saw this this for the $V$-statistic estimator since we referred to the literature for the $U$-statistics estimator. Nevertheless, we saw that the specific asymptotic distribution $\eta$ of $nV_n^6(h,Z_{1,n})$ was derived solely from the decomposition term $V_n^2(\bar{h}^{(2)},Z_{1,n})$, as every other decomposition term converged to zero. This is the reason for the bootstrapping approach proposed in \cite{arcones1992bootstrap}, instead samples with replacement from the empirical distribution and inserts these samples into $U$- and $V$-statistics with empirically modified kernels based on $\bar{h}^{(2)}$.

We go into detail on how to bootstrap the limit distribution $\gamma$, of our scaled $U$-statistics $n\tilde{U}_n^6(h,Z_{1,n})$ under the null-hypothesis, but refer the reader to \cite{arcones1992bootstrap} for a similar approach for limit the distribution $\eta$ of our scaled $V$-statistics $nV_n^6(h,Z_{1,n})$. \\ \\
To this end, let $(Z_n)=((X_n,Y_n))$ be an i.i.d. sequence defined on a common probability space $(\Omega,\F,P)$ such that each pair is distributed according to a  $\theta\in M^{1,1}_{1,nd}(\cX\times \cY)$ that satisfies the null-hypothesis $\theta=\mu\times \nu$. Furthermore let $\theta_n(\omega)=n^{-1}\sum_{i=1}^n \delta_{(X_i(\omega),Y_i(\omega))}$ denote the $n$'th empirical measure given a realization $\omega\in \Omega$. We assume that $Z_{n1}^{*\omega},...,Z_{nn}^{*\omega}$ denotes i.i.d. random elements in $\cX\times \cY$ with distribution function $\theta_n(\omega)$ for any $\omega\in \Omega$ and $n\in \N$. Recall $\bar{h}^{(2)}:(\cX\times \cY)^2\to \R$ defined in \cref{appendix_definition_h^(k)}, and note that it can be written as
\begin{align*}
\bar{h}^{(2)}(z_1,z_2) =& \int \bar{h}(z_1,...,z_2,w_{3},...,w_6) \, d\theta^{4}(w_{3},...,w_6) - \int \bar{h}(z_1,w_2,...,w_6) \, d\theta^{5}(w_{2},...,w_6) \\
&- \int \bar{h}(z_2,w_2,...,w_6) \, d\theta^{5}(w_{2},...,w_6w) + \int \bar{h}(w_1,...,w_6) \, d\theta^{6}(w_{1},...,w_6) .
\end{align*}
The previously mentioned empirically modified version of $\bar{h}^{(2)}$ is given by the above expression, but where we interchange the true distribution $\theta$ with the realized empirical distribution $\theta_n(\omega)$. That is, the empirically modified version of $\bar{h}^{(2)}$ is given by 
\begin{align*}
&\bar{h}^{(2)}_{\theta_n(\omega)}(z_1,z_2)\\ =& \int \bar{h}(z_1,z_2,w_{3},...,w_6) \, d\theta_n(\omega)^{4}(w_{3},...,w_6) - \int \bar{h}(z_1,w_2,...,w_6) \, d\theta_n(\omega)^{5}(w_{2},...,w_6) \\
&- \int \bar{h}(z_2,w_2,...,w_6) \, d\theta_n(\omega)^{5}(w_{2},...,w_6) + \int \bar{h}(w_1,...,w_6) \, d\theta_n(\omega)^{6}(w_{1},...,w_6) \\
=& \frac{1}{n^4}\sum_{i_3=1}^n \cdots \sum_{i_6=1}^n \bar{h}(z_1,z_2,Z_{i_3}(\omega),...,Z_{i_6}(\omega))  -\frac{1}{n^5} \sum_{i_2=1}^n \cdots \sum_{i_6=1}^n\bar{h}(z_1,Z_{i_2}(\omega),...,Z_{i_6}(\omega))  \\
&- \frac{1}{n^5}\sum_{i_2=1}^n \cdots \sum_{i_6=1}^n\bar{h}(z_2,Z_{i_2}(\omega),...,Z_{i_6}(\omega))  + \frac{1}{n^6}\sum_{i_1=1}^n \cdots \sum_{i_6=1}^n\bar{h}(Z_{i_1}(\omega),...,Z_{i_6}(\omega)) .
\end{align*}
The bootstrap consistency theorem of \cite{arcones1992bootstrap} (theorem 2.4) states that, if the symmetric kernel $\bar{h}:(\cX\times \cY)^6 \to \R$ satisfies the integrability condition
\begin{align*}
E|\bar{h}(Z_{i_1},...,Z_{i_6})|^{\frac{2}{6}\#\{i_1,..,i_6\}}<\i \quad \quad \text{ for all } (i_1,...,i_6)\in\{1,..,6\}^6,
\end{align*}
then
\begin{align*}
\frac{2!{ 6 \choose 2}}{n} \sum_{1\leq i_1  <i_2\leq n} \bar{h}^{(2)}_{\theta_n(\omega)}(Z_{ni_{1}}^{*\omega},Z_{ni_2}^{*\omega}) \convd_n \gamma \quad \quad \textit{for } P\textit{-almost all }\omega\in\Omega.
\end{align*}
 As we have argued before, $\bar{h}$ is the symmetrized version of $h$ so Minkowski's inequality (or see below if exponent is less than one) yields that the above integrability holds if $E|h(Z_{i_1},...,Z_{i_6})|^{\frac{2}{6}\#\{i_1,...,i_6\}}<\i $ for any  $(i_1,...,i_6)\in\{1,..,6\}^6$. In the case that all indices are distinct, we note that the requirement is square integrability of $h$ with respect to $\theta^6$, which is guaranteed by \cref{lemma_h_is_theta^6_integrable}. Hence denote $p= \frac{2}{6}\#\{i_1,..,i_6\}$  and fix any indices $(i_1,...,i_6)\in\{1,..,6\}^6$ such that $p\leq \frac{2}{6}5=\frac{5}{3}$. We see that 
	\begin{align*}
&E|h(Z_{i_1},...,Z_{i_6})|^p\\
= &E|f_\cX(X_{i_1},X_{i_2},X_{i_3},X_{i_4})f_\cY(Y_{i_1},Y_{i_2},Y_{i_5},Y_{i_6})|^p \\
\leq& E[d_\cX(X_{i_1},X_{i_4})d_\cY(Y_{i_2},Y_{i_5})]^p \\
=& E[(d_\cX(X_{i_1},x)+d_\cX(X_{i_4},x))(d_\cY(Y_{i_2},y)+d_\cY(Y_{i_5},y))]^p \\
\leq& \big[   E(d_\cX(X_{i_1},x)^p)^{1/p}E(d_\cY(Y_{i_2},y)^p)^{1/p}+E(d_\cX(X_{i_1},x)^p )^{1/p} E(d_\cY(Y_{i_5},y)^p)^{1/p} \\
&+ E(d_\cX(X_{i_4},x)^p )^{1/p} E(d_\cY(Y_{i_2},y)^p)^{1/p}+E(d_\cX(X_{i_4},x)^p )^{1/p} E(d_\cY(Y_{i_5},y)^p)^{1/p} \big]^p,
\end{align*}
for some $x\in \cX$ and $y\in \cY$, where we used Minkowski's inequality, the first two inequalities of \cref{lemma_inequality_f} and that $\theta=\mu\times \nu$ (if $p<1$ create similar upper bounds by the inequality $\|f+g\|_p^p \leq \|f\|_p^p+\|g\|_p^p$). From this we see,  it is sufficient that $\theta\in M^{5/3,5/3}_{1,nd}(\cX\times \cY)$ in order to guarantee that the bootstrap consistency theorem holds. 

These arguments entail that we are only guaranteed to have convergence in distribution of the bootstrap statistics towards  the $\gamma$ distribution  in the third statistical model where $\mathcal{P}_3=M^{2,2}_{1,nd}(\cX\times \cY)$. Note that we do not state, that the integrability condition is not satisfied in the first and second statistical models, but that the above upper bounds are only sufficiently tight in third statistical model. \\ \\
Now let us describe the heuristics behind bootstrap approach to approximate the $(1-\alpha)$-quantile of the $\gamma$ distribution. Let $F_n^{*\omega}$ denote the cumulative distribution function of the random variable $\frac{30}{n} \sum_{1\leq i_1  <i_2\leq n} \bar{h}^{(2)}_{\theta_n(\omega)}(Z_{ni_{1}}^{*\omega},Z_{ni_2}^{*\omega})$ for any $\omega\in \Omega$ and $n\in \N$. Since $\gamma$ has a continuous cumulative distribution function $F_\gamma$, the bootstrap consistency theorem yields that
\begin{align*}
F_n^{*\omega}(x) \to_n F_\gamma(x) \quad \quad \textit{for all } x\in \R,
\end{align*}
for $P$-almost all $\omega\in \Omega$. This is of course equivalent to the convergence of the quantile functions (lemma 21.2 \cite{van2000asymptotic}), i.e.
\begin{align*}
q_{p}^{(F_n^{*\omega})}:=&\inf\{x\in \R : p\leq F_n^{*\omega}(x) \} \\\to_n& \inf \{x\in \R :p\leq F_\gamma(x) \}= q^{(\gamma)}_{p} \quad \quad \textit{for all } p\in (0,1),
\end{align*}
for $P$-almost all $\omega\in \Omega$. Now the bootstrap approach for approximating $q_{1-\alpha}^{(\gamma)}$ makes the approximation $q_{1-\alpha}^{F_n^{*\omega}} \approx q_{1-\alpha}^{(\gamma)}$ for any realization $\omega\in \Omega,\alpha\in(0,1)$ and $n\in \N$, which is deemed reasonable if $n$ is large by the above quantile convergence. 

Hence fix $\omega\in \Omega$ and $n\in \N$ (denoting the given sample-size) and note that we have reduced the problem of finding the rejection threshold to finding the quantile $q_{1-\alpha}^{F_n^{*\omega}}$. Let  $[(Z_{1n}^{*\omega}(i),...,Z_{nn}^{*\omega}(i))]_{i=1}^\i$ be independent copies of $n$ i.i.d. random variables $Z_{1n}^{*\omega}(1),...,Z_{nn}^{*\omega}(1)$ each distributed according to the $n$'th empirical measure  $\theta_n(\omega)=n^{-1}\sum_{i=1}^n \delta_{(X_i(\omega),Y_i(\omega))}$ given the realization $\omega$. By the Glivenko-Cantelli theorem we have that
\begin{align*}
&\sup_{x\in \R}|F_{mn}^{*\omega}(x)-F_n^{*\omega}(x)|\\
=&\sup_{x\in \R}\left|\frac{1}{m} \sum_{j=1}^m 1_{ (-\infty  ,x]}\lp \frac{30}{n} \sum_{1\leq i_1  <i_2\leq n} \bar{h}^{(2)}_{\theta_n(\omega)}(Z_{ni_{1}}^{*\omega}(j),Z_{ni_2}^{*\omega}(j))\rp
 -F_n^{*\omega}(x)\right| \to_m 0,
\end{align*}
almost surely. Hence for large $m$ we may reasonably approximate the unknown distribution function $F_n^{*\omega}$ by  $F_{mn}^{*\omega}$. Since we know the empirical measure $\theta_n(\omega)$ we may generate realizations of
$$
\lp Z_{1n}^{*\omega}(1),...,Z_{nn}^{*\omega}(1)\rp ,..., \lp Z_{1n}^{*\omega}(m),...,Z_{nn}^{*\omega}(m) \rp,
$$
 for some arbitrarily large $m\in \N$. Based on these samples we may calculate $$\frac{30}{n} \sum_{1\leq i_1  <i_2\leq n} \bar{h}^{(2)}_{\theta_n(\omega)}(Z_{ni_{1}}^{*\omega}(j),Z_{ni_2}^{*\omega}(j)) \quad \quad \text{for }j=1,...,m ,
 $$
and find corresponding the $(1-\alpha)$-quantile $q^*_{1-\alpha}(m)$ of the resulting empirical distribution.  We say that this quantile $q^*_{1-\alpha}(m)$ approximates the true $(1-\alpha)$-quantile of the $\gamma$ distribution, through the above reasoning of the approximations $q^*_{1-\alpha}(m) \approx q_{1-\alpha}^{F_n^{*\omega}} \approx q^{(\gamma)}_{1-\alpha}$, for some large $m\in \N$.  \\ \\
We can summarize this approach in the following bootstrap and test algorithm, where we are given empirical samples $[z_i]_{i=1}^n = [(x_i,y_i)]_{i=1}^n$ assumed to be a realization of $[(X_i,Y_i)]_{i=1}^n$.
\begin{itemize}
	\item[1)] Choose a large $m\in \N$ and sample with replacement $m\times n$ times from the observed empirical distribution $\theta_n$ placing $1/n$ point-mass at $(x_i,y_i)$ for all $1\leq i \leq n$, yielding an $m\times n$ array of samples
	\begin{align*}
	\begin{array}{ccc}
	z_{n1}^*(1) & \cdots & z_{nn}^*(1)  \\
	\vdots & \vdots & \vdots \\ 
	z_{n1}^*(m) & \cdots & z_{nn}^*(m)
	\end{array}.
	\end{align*}
	\item[2)] For each $1\leq j \leq m$ calculate 
	\begin{align*}
	\tau^*_n(j)=\frac{30}{n} \sum_{1\leq i_1  <i_2\leq n} \bar{h}^{(2)}_{\theta_n}(z_{ni_{1}}(j),z_{ni_2}(j)),
	\end{align*}
	and denote the resulting empirical cdf by $
	F_m^{*}(x)= \frac{1}{m}\sum_{j=1}^m 1_{(-\infty,x]}(\tau^*_n(j))
$
	\item[3)] Calculate the corresponding empirical $(1-\alpha)$-quantile $q_{1-\alpha}^*(m)=\inf\{x\in \R: (1-\alpha)\leq F_m^*(x) \}$ and reject the null-hypothesis if
	\begin{align*}
	n\tilde{U}_n^6(h,((x_1,y_1),...,(x_n,y_n))) > q_{1-\alpha}^*(m).
	\end{align*}\\ 
\end{itemize}

This concludes the last section of the thesis. We have provided a solution to the non-parametric independence problem, by constructing asymptotically consistent statistical tests for testing the null-hypothesis, and it was argued how one reasonably can bootstrap approximate the rejection thresholds of the aforementioned tests.

\section{Summary and future work}
{\large\textit{\textbf{Summary:}}} In this thesis we proposed a solution to the non-parametric independence problem, whenever the marginal metric spaces were separable and of strong negative type. We did this by introducing the distance covariance measure in metric spaces $dcov:M^{1,1}_1(\cX\times \cY)\to \R$. In order to ensure that $dcov$ is well-defined, we made the additional restriction only to consider marginal metric spaces that  are separable. 

 The distance covariance measure is however, not  a direct indicator of independence, for general marginal metric spaces. Thus we embarked on searching for further conditions on the marginal metric spaces, which would guarantee this property. To this end, we showed that whenever the marginal metric spaces are of negative type, we can represent the distance covariance measure in terms of mean embeddings of certain isometries into separable Hilbert spaces.  This representation resulted in the definition of the subset of negative type metric spaces, called metric spaces of strong negative type. With some effort we were able to show, that the distance covariance measure is a direct indicator of independence, whenever the marginal spaces are metric spaces of strong negative type. Additionally, it was shown that every separable Hilbert space is a metric space of strong negative type.
 
 Then we constructed two estimators for the distance covariance measure, a $V$-statistic estimator as in \cite{lyons2013distance}, but also a new one given by the corresponding $U$-statistic. We proved, that both estimators are strongly consistent and possesses well-defined asymptotic distributions. We also argued that the moment conditions in \cite{lyons2013distance} for strong consistency of the $V$-statistic estimators are not sufficient. However, our $U$-statistic estimator only needed the weaker moment assumptions in order to guarantee strong consistency. As regards to the asymptotic distribution of the $V$-statistics, it is still  unresolved whether the asymptotic distribution indeed can be written as in \cite{lyons2013distance}. 
 
 Nevertheless the aforementioned asymptotic properties of the estimators was combined with the developed theory of the distance covariance measure, to construct  statistical tests of independence. These tests were guaranteed to be asymptotically consistent under certain conditions, one of which was that the marginal spaces must be of strong negative type. Lastly as the tests were constructed with rejection thresholds given by non-traceable quantiles, we argued that they could be reasonably bootstrapped. \\ \\
{\large\textit{\textbf{Future work:}}}  There are a few  things, which would be very interesting to explore and examine. First of all, it would be interesting to do a simulation experiment, to see how the two different statistical tests compare to each other. For example, it would be interesting to examine the statistical power "$1-P(\text{Type II error})$"  of the two tests for varying sample sizes, to see how many sample points each test would need to yield a reasonably low  frequency of Type II errors. It could also be of interest to actually apply the statistical tests, to real sample-data with values in a non-Euclidean space. E.g. functional data where each realization is seen as a sample-path with valued in a $L^2$-space.
 
 In  \cite{sejdinovic2013equivalence}, published in The Annals of Statistics, it is stated that the theory of distance covariance in metric spaces extends to semi-metric spaces. They furthermore state an equivalence between independence testing using distance covariance in semi-metric spaces and something called the Hilbert-Schmidt independence criterion. Since we did not have the time to pursue these claims, it could serve as a very interesting continuation of the thesis.
\section{Appendix}
	\subsection{Product spaces - metrics, topologies and $\sigma$-algebras} \label{Appendix_Product_spaces_section}
Let $I=\{1,...,n\}$ for some $n\in\N$, or $I=\N$ and consider any family of measurable spaces $((\cX_i,\cE_i))_{i\in I}$. The (Cartesian) product space $\cX = \prod_{i\in I} \cX_i$ has the following representations.
	\begin{itemize}	
		\item[$\cdot$] \textit{If $I=\{1,...,n\}$ then $\cX = \{(x_1,x_2,...,x_n):  x_k \in \cX_k, \, \,\forall k\in \{1,...,n\}\},$
		that is the set of all ordered $n$-tuples, with $x_k\in\cX_k$ for all $1\leq k \leq n$}.
		\item[$\cdot$] \textit{If $I=\N$ then $\cX = \{(x_1,x_2,x_3,...):  x_k \in \cX_k, \, \,\forall k\in I \},$ that is the set of all infinite sequences $(x_1,x_2,x_3,...)$, with $x_k\in\cX_k$  for all $k\in\N$}.
	\end{itemize}
	For any $i\in I$ we define the coordinate projection $\pi_i:\cX \to \cX_i$ by $\pi_i(x)=x_i$ for any $x\in \cX $. Furthermore for any $k\in I$ and $t\in I^k$ with $t_1<\cdots <t_k$, we define the simultaneous coordinate projection $\pi_{t_1 \cdots t_k}:\cX \to \prod_{i=1}^k \cX_{t_i}$ by 
	\begin{align} \label{defi_simultaneous_coordinate_projection}
	\pi_{t_1\cdots t_k}(x)=(\pi_{t_i}(x),...,\pi_{t_k}(x))=(x_{t_1},...,x_{t_k}),
	\end{align}
	and for simplicity we may also denote this simultaneous coordinate projection by $\pi_t$ whenever it is clear that $t\in I^k$ with $t_1< \cdots <t_k$.
	\begin{definition}[Product $\sigma$-algebra] \label{defi_product_sig_alg}
		 We define the product \sigalg $\bigotimes_{i\in I} \cE_i$ on  $ \prod_{i\in I} \cX_i$ as the smallest \sigalg  making every coordinate projection measurable . That is
		\begin{align} \label{defi_product_sig_alg_eq}
		\bigotimes_{i\in I} \cE_i = \sigma \lp \lb \pi^{-1}_i(E_i) : E_i \in \cE_i ,i\in I  \rb \rp.
		\end{align}
	\end{definition}
	It is fairly easy to show that the above \sigalg  is also generated by
		\begin{align*}
		\lb \prod_{i\in I}E_i : E_i \in \cE_i   \rb \quad \textit{ and } \quad \lb \pi^{-1}_i(E_i) : E_i \in \bE_i ,i\in I   \rb \quad \textit{ and } \quad \lb \prod_{i\in I} E_i: E_i\in \bE_i \rb,
		\end{align*}
		 when $\cE_i=\sigma(\bE_i)$ with $\cX_i\in \bE_i$ for all $i\in I$ (see proposition 1.3 and 1.4 \cite{folland1999real}). It is also evident that the simultaneous coordinate projections $\pi_{t_1 \cdots t_k}$ defined in \cref{defi_simultaneous_coordinate_projection} are $\bigotimes_{i\in I} \cE_i / \bigotimes_{i=1}^k \cE_{t_i}$-measurable. Moreover when $I=\N$ we have an intersection stable generator for $\bigotimes_{i\in \N} \cE_i$ given by
		 \begin{align*}
		 \left\{\bigcap_{i=1}^n \pi_i^{-1}(E_i): n\in\N, E_i\in \bE_i \right\},
		 \end{align*}
		 (see \cite{VS} lemma 2.3.2 for proof).\\ \\
	Let $((\cX_i,d_i))_{i\in I}$ be a  family of metric spaces. Whenever a norm or metric is introduced for a space, we always work with the corresponding metric topology $\cT_i$ and the Borel $\sigma$-algebra $\cB(\cX_i)=\sigma(\cT_i)$ induced by these, unless otherwise stated. \\ \\
	When considering the product space $\cX=\prod_{i\in I}\cX_i$ we always equip it with the product topology $\bigodot_{i\in I}\cT_i$, unless we introduce a metric in which case the above comment applies. The product topology is defined as the topology generated by (smallest topology containing) \begin{align}
	\mathcal{A}=\{\pi_i^{-1}(U): U \in \cT_i, i\in I\}, \label{defi_subbase_for_prod_topology}
	\end{align}
	that is the above family of sets is a subbase for the product topology on $\cX$. In other words the product topology is the smallest/coarsest topology making all coordinate projections continuous. We may also note that the Borel $\sigma$-algebra of the product space $\sigma(\bigodot_{i\in I}\cT_i)$ (which we also write as $\cB \big(\prod_{i\in I}\cX_i)\big)$ always contains the product Borel $\sigma$-algebra, that is $$\bigotimes_{i\in I} \cB(\cX_i)\subset \cB \bigg(\small\prod_{i\in I}\cX_i \bigg)
	.$$ To see this simply note that  $\lb \pi^{-1}_i(B_i) : B_i \in \cT_i ,i\in I  \rb=\cA \subset \cT$. By the remark below \cref{defi_product_sig_alg} yields that the former family of sets is a generator for the product Borel $\sigma$-algebra. Hence 
	$$\bigotimes_{i\in I} \cB(\cX_i)=\sigma(\lb \pi^{-1}_i(B_i) : B_i \in \cT_i ,i\in I  \rb)\subset \sigma\big( \bigodot_{i\in I}\cT_i\big)=\cB \big(\prod_{i\in I}\cX_i)\big).$$
		A natural question is whether the Borel  \sigalg $\cB(\cX)$ induced by product topology  coincides with the product \sigalg $\bigotimes_{i\in I} \cB(\cX_i)$. Nice properties follow if we indeed have equality, for instance every continuous function  becomes measurable with respect to the product $\sigma$-algebra. Unfortunately, equality does not always hold,  as shown in example 6.4.3 \cite{MeasureTheory2007bogachev2}. Though for some sufficiently nice spaces the equality does indeed hold. Below we prove that they coincide in the case of the case where all marginal spaces are separable.  \\ \\ 
	We note that if $I$ is finite then $\cX$ is metrizable, in the sense that the maximum/product metric $\rho^{\max}:\cX \times \cX \to [0,\i)$ given by
	$
	\rho^{\max}(x,y)= \max_{i\in I} d_i(x_i,y_i)
	$
	is a metric on $\cX$ which induces product topology (see remark below next theorem). Hence for finite $I$ it is evident that we have convergence in $\cX$ if and only if we have convergence in $\cX_i$ of each coordinate, and realize that this implies that every coordinate projection is continuous.

	\begin{theorem} \label{theorem_product_sig_alg_is_tensor_when_separable}
		If every metrizable topological space in the family $((\cX_i,\cT_i))_{i \in I}$ is separable, then $(\prod_{i\in I}\cX_i,\bigodot_{i\in I}\cT_i)$ is separable and the Borel \sigalg $\cB\lp\prod_{i\in I}\cX_i\rp$ induced by the product topology coincides with the product \sigalg $\bigotimes_{i\in I} \cB(\cX_i)$.
	\end{theorem}
	\begin{p}
		We prove it for $I=\{1,...,n\}$ for some $n\in\N$, but the proof when $I$ is a countably infinite index set follows by analogous steps (see for example Lemma 1.2 \cite{kallenberg1997foundationsofmodernprobability}). \\ \\
		First we show separability of $\cX$: Let $D_i=\{x_{i,1},x_{i,2},x_{i,3},...\}$ be a countable dense subset of $\cX_i$, for each $1\leq i \leq n$. Note that $ D= \prod_{i=1}^n D_i\subset \cX$ is the Cartesian product of $n$ countable sets, hence itself countable. Now fix an arbitrary $x=(x_1,...,x_n)\in \cX$ and note that since $D_i\subset \cX_i$ is dense, there exists a sequence $(x^k_i)_{k\in \N}$ in $D_i$ converging to $x_i$, for all $1\leq i \leq n$. 
		
		Now construct  the sequence $(x^k)_{k\in \N}$ in $D$ by setting $x^k=(x_1^k,...,x_n^k)$ for every $k\in \N$. Lastly, note that since convergence in $\cX$  is equivalent to convergence in each coordiante, we by construction of $(x^k)_{k\in \N}$ have that $x^k\to_k x$ since $x^k_i \to_k x_i$, for all $1\leq i \leq n$. Thus every point in $\cX$ is a limit point of the dense subset $D$, proving separability.\\ 
		
		As regards the claim about the $\sigma$-algebras it suffices to show that $\bG_1 \subset \sigma(\bG_2)$ and that $\sigma(\bG_1)\supset \bG_2$ for any two generators $\bG_1$ and $\bG_2$ such that $\sigma(\bG_1)=\cB\lp \cX \rp$ and $\sigma(\bG_2)=\bigotimes_{i=1}^n \cB(\cX_i)$ respectively. Since the open sets of $\cX_i$ generate $\cB(\cX_i)$, we get by the remark below \cref{defi_product_sig_alg}, that the generator $\bG_2$ can be choosen as
		\begin{align*}
		\bG_2=\lb \pi_i^{-1}(O_i) : O_i \text{ is open in } \cX_i,\forall 1\leq i \leq n  \rb,
		\end{align*}
		but by  \cref{defi_subbase_for_prod_topology} this is  exactly the family of sets generating the product topology. Hence $\bG_2$ must be a subset of the corresponding Borel \sigalg $\cB\lp \cX \rp$. 
		
		Conversely, let $\bG_1$ be the entire family of open sets in $\cX$ and note that since $\cX$ is separable, we know that it has a countable base (cf. Theorem M3 \cite{Billingsley_convergence_1999}) given by 
		\begin{align*}
		\bB =\lb \prod_{i=1}^n B_{d_i}(x_i,q) : x \in  D,q\in \bQ \rb = \lb \bigcap_{i=1}^n \pi_i^{-1}(B_{d_i}(x_i,q) ) : x \in  D,q\in \bQ \rb.
		\end{align*}
		Now realize that each set in the above base is a finite intersection of sets in $\bG_2$ since $B_{d_i}(x_i,q)$ is open in $\cX_i$ for any $x\in D$ and $q\in \bQ $, and therefore  $\bB \subset \sigma(\bG_2)=\bigotimes_{i=1}^n\cB(\cX_i) $. Lastly, we note that the countability of the base implies that any open set in $\cX$, i.e. any element of $\bG_1$ is a countable union of elements in $\bigotimes_{i=1}^n \cB(\cX_i)$. Since $\bigotimes_{i=1}^n \cB(\cX_i)$ contains countable unions of its members we conclude that, $\bG_1 \subset \bigotimes_{i=1}^n \cB(\cX_i)$.
	\end{p}

Moreover we have that a metrizable topological space is completely determined by its convergent sequences (see \cite{SpacesInWhichSequencesSuffice}). Furthermore if $((\cX_i,\cT_i))_{i \in I}$ are metrizable topological spaces then $(\Pi_{i\in I}\cX_i, \bigodot_{i\in I}\cT_i)$ is a metrizable topological space if $\# I\leq \aleph_0$ (see corollary 7.3 \cite{topology-dugundji}) and in general we have that the product spaces have coordinatewise convergence, i.e.  $x_n \to x$ in $(\Pi_{i\in I}\cX_i, \bigodot_{i\in I}\cT_i)$ if and only if $\pi_i(x_n) \to \pi_i(x)$ in $(\cX_i,\cT_i)$ for all $i\in I$ (see lemma 43.3 $\cite{munkres2000topology}$). A consequence of these facts is: If $((\cX_i,d_{\cX_i}))_{i\in I }$ is a family of metric spaces and $d$ is a metric on the product space $\Pi_{i\in I}\cX_i$ for which it holds that $d(x_n,x)\to 0 $ if and only  $d(\pi_i(x_n),\pi_i(x))\to 0$ for all $i\in I$ then $d$ induces the product topology on $\Pi_{i\in I}\cX_i$. 

\newpage
\subsection{U- and V-statistics}
In this section we introduce $U$- and $V$-statistics. We prove the Hoeffding decomposition theorem for $U$-statistics, and we will define various mappings used in the thesis. Lastly, we draw on the litterature to state some asymptotic properties of $U$- and $V$-statistics, which we use in the thesis.
\subsubsection{Hoeffding decomposition} \label{Appendix_Hoeffding_decomposition_section}
Let $(H_i)_{1\leq i \leq k}\subset H$ be a family of non-empty subspaces of a Hilbert space $H=L^2(\Omega,\F,P)$ each mutually orthogonal to each other, that is $H_i \perp H_j$ for $i\not = j$. If all orthogonal projections $P_{H_i}(T)$ for $1\leq i \leq k$ exist for some $T\in H$ then
\begin{align*}
P_{H_1+\cdots +H_k}(T)= P_{H_1}(T)+ \cdots + P_{H_k}(T),
\end{align*}
where $H_1+\cdots +H_k$ is the sum of subspaces defined by
\begin{align*}
H_1+\cdots +H_k = \{h_1+\cdots +h_k: h_1\in H_1,...,h_k\in H_k\}.
\end{align*}
This is easily established by an induction argument. For two spaces, let $P_{H_1}(T)\in H_1$ and $P_{H_2}(T)\in H_2$ with $\la T-P_{H_1}(T), h_1\ra=0$ and $\la T-P_{H_2}(T),h_2 \ra =0$ for all $h_1\in H_1, h_2\in H_2$. Note that $P_{H_1}(T)+P_{H_2}(T)\in H_1+H_2$ and for every $h_1+h_2\in H_1+H_2$ we have that
\begin{align*}
\la T-P_{H_1}(T)-P_{H_2}(T),h_1+h_2\ra =& \la T-P_{H_1}(T),h_1 \ra -\la P_{H_2}(T),h_1 \ra  \\
&+\la T-P_{H_2}(T),h_2 \ra -\la P_{H_1}(T),h_2 \ra \\
=& 0,
\end{align*}
since $P_{H_2}(T) \perp  h_1$ and $ P_{H_1}(T) \perp h_2$ because $H_1\perp H_2$, proving that the projection onto the sum space $H_1+H_2$ is given by the sum of the projections by standard equivalence of orthogonal projections. Now assume that $P_{H_1+\cdots +H_{n-1}}(T)= P_{H_1}(T)+ \cdots + P_{H_{n-1}}(T)$ for some $1<n<k$ and note that
$H_1+\cdots +H_{h-1} \perp H_n$ such that by the above arguments $P_{H_1+\cdots + H_{n-1}+H_n}(T)= P_{H_1+\cdots +H_{n-1}}(T)+P_{H_n}(T)$. Hence by the induction assumption we have that $P_{H_1+\cdots + H_n}(T)=P_{H_1}(T)+ \cdots + P_{H_n}(T)$, proving the general claim by induction. We will use this property so the next step is to define some mutually orthogonal spaces. \\ \\
Consider independent random elements $X_1,...,X_n$ in the measurable space $(\cY,\bK)$ defined on some probability space $(\Omega,\F,P)$, and let $A\subset \{1,...,n\}$. Furthermore let $|A|$ denote the number of elements in the set $A$ and let $H_A$ denote the subset of  $L^2(\Omega,\F,P)$ where each element can be written in the form
\begin{align*}
g_A(X_i:i\in A),
\end{align*}
where $g_A(X_i:i\in A) \in L^2(\Omega,P)$ with
\begin{align} \label{defi_conditioning_property_of_H_A}
E \lp g_A(X_i:i\in A)|X_j:j\in B\rp=0,
\end{align}
for all $B\subset \{1,...,n\}$ with $|B|<|A|$ (With the convention that conditioning on an empty index set, is the conditioning on the trivial $\sigma$-algebra $\{\Omega,\emptyset \}$, resulting in every random variable in $H_A$ for $|A|>0$ having zero expectation). Furthermore it is interpreted that $H_\emptyset$ is the set of (almost surely) constant functions, and it is easily verified that $H_A$  is indeed a subspace of $L^2(\Omega,\F,P)$ for any $A\subset \{1,...,n\}$. \\ \\
Now we show the family $(H_A)_{A\in\mathcal{P}(\{1,...,n\})}$ of subspaces are mutually orthogonal. Thus consider any two distinct sets $A,B \in \mathcal{P}:=\mathcal{P}(\{1,...,n\})$, and take any $g_A(X_i:i\in A)\in H_A$ and $g_B(X_i:i\in B) \in H_B$. If $A\cap B = \emptyset$ then by  the mutual independence of the $X_i$'s we have that
\begin{align*}
\la g_A(X_i:i\in A),g_B(X_i:i\in B) \ra= E \lp   g_A(X_i:i\in A) \rp E \lp  g_B(X_i:i\in B) \rp  =0,
\end{align*}
where we used that since $A\not = B$ so one of them has at least one element which by the above remark yields a zero expectation. If $A\cap B=C \not = \emptyset$ then since $A\not = B$ we have that  $|C|<|A|$ or $|C|<|B|$ or both. Without loss of generality we may assume that $|C|<|A|$ (otherwise interchange $A$ with $B$ below) and by removing redundant information on the conditioning (that is $E(X|Y,Z)=E(X|Z)$ if $X\independent Y$) we have that
\begin{align*}
\la g_A(X_i:i\in A),g_B(X_i:i\in B) \ra &= E \lp  E \lp g_A(X_i:i\in A) g_B(X_i:i\in B) \big| X_j:j\in B\rp  \rp \\
&= E \lp  E \lp g_A(X_i:i\in A)  \big| X_j:j\in B\rp g_B(X_i:i\in B)  \rp \\
&= E \lp  E \lp g_A(X_i:i\in A)  \big| X_j:j\in C\rp g_B(X_i:i\in B)  \rp \\
&=0,
\end{align*}
where we used  \cref{defi_conditioning_property_of_H_A}. This proves that  $(H_A)_{A\in \mathcal{P}}$ is a family of mutually orthogonal subspaces of $L^2(\Omega,\F,P)$.
\begin{lemma} \label{lemma_projection_onto_H_A_spaces}
	Let $T\in L^2(\Omega,\F,P)$ be an arbitrary square integrable random variable and let $A\subset \{1,...,n\}$. 
	\begin{itemize}
		\item[(1)] The orthogonal projection onto $H_A$ exists and is given by
		\begin{align} \label{projection_candidate_in_lemma_for_projection}
		P_A(T)=\sum_{B\subset A} (-1)^{|A|-|B|}E(T|X_i:i\in B).
		\end{align}
		\item[(2)] If $T \perp H_B$ for all $B\subset A$, then $E(T|X_i:i\in A)=0$.
		\item[(3)] For any measurable map $g:(\cX^{|A|},\otimes_{i=1}^{|A|}\F)\to(\R,\cB(\R))$ such that $g (X_i:i\in A)\in \mathcal{L}^2(\Omega,\F,P)$, it holds that $g (X_i:i\in A)\in  \sum_{B\subset A}H_B$.
	\end{itemize}
\end{lemma}
\begin{p}
	Let $A\subset \{1,...,n\}$ be any non-empty subset and  $B\subset A$  any subset hereof. Now note that by the independence of $X_1\independent \cdots \independent X_n$ we may remove redundant information as follows
	\begin{align*}
	E(E(T|X_i:i\in A)|X_j:j\in B) &= E(E(T|X_i:i\in A)|X_j:j\in A\cap B) \\
	&= E(T|X_j:j\in A\cap B),
	\end{align*}
	where we also used that only the smallest conditioning $\sigma$-algebra remains in an iterated conditional expectation. With $P_A(T)$ defined in $\cref{projection_candidate_in_lemma_for_projection}$ and $C\subsetneq A$ being a proper subset we have that
	\begin{align*}
	E(P_A(T)|X_i:i\in C) &= \sum_{B\subset A} (-1)^{|A|-|B|}E(T|X_i:i\in B\cap C).
	\end{align*}
	First realize when summing over all possible subsets $B\subset A$, we exactly hit all conditioning indexes of the form $B\cap C = D \subset C$ at least once in the sum. Hence we may instead sum over $D \subset  C$ and change the conditional expectation to $E(T|X_i:i\in D)$, but in order to do this we must count how many different subsets $B\subset A$ result in the same $D=B\cap C$. Or to put it differently when considering any subset  $D\subset C$ which and how many subsets $B\subset A$ result in $B\cap C=D$. 
	
	This only holds for sets of the form $B=D\cup K$, for any $K\subset A\setminus C$. To see this note that in order for $B\cap C = D$ it is necessary that $B\supset  D$, so it can only be true for $B$ of the the form $B=D\cup K$ for some set $K\subset A\setminus D$. Let $K\subset A\setminus D$ such that $K\cap C$ is non-empty, but then as $K\cap C \not \supset D$ we get
	$$
	B\cap C = (D\cup K)\cap C=(D\cap C) \cup(K\cap C)= D \cup (K\cap C) \not = D,
	$$
	so as stated above we need to restrict $K$ to be a subset of $A\setminus C$ in order for $B=D\cup K$ to fulfil $B\cap C = D$. 
	
In general for each $D\subset C$ we can chose $K$ to consist of $j\in \{0,...,|A|-|C|\}$ elements of $A\setminus C$, and there are exactly ${|A|-|C|}\choose{j}$ distinct ways of choosing $j$ distinct elements from $A\setminus C$. Hence by using the binomial formula we get that 
	\begin{align*}
	\sum_{B\subset A} (-1)^{|A|-|B|}E(T|X_i:i\in B\cap C) &=\sum_{D\subset C} \sum_{j=0}^{|A|-|C|} \, { |A|-|C|\choose j} (-1)^{|A|-(|D|+j)}E(T|X_i:i\in D) \\
	&=\sum_{D\subset C} (1-1)^{|A|-|C|} E(T|X_i:i\in D) \\
	&=0.
	\end{align*}
	Now realize that for any subset $B\subset \{1,...,n\}$ with $|B|<|A|$ we have that $B\cap A \subsetneq A$ and thus by removing redundant information from the conditioning we get
	\begin{align*}
	E(P_A(T)|X_i:i\in B) &= E(P_A(T)|X_i:i\in B\cap A) \\
	&= 0,
	\end{align*}
	by using the above. It is furthermore not hard to realize that $P_A(T)$ can be written as a measurable function composed with $(X_i:i\in A)$, which satisfies
	\begin{align*}
	\| P_A(T) \|_2 \leq \sum_{B\subset A} ||E(T|X_i:i\in B)||_2 =\sum_{B\subset A} \|T\|_2<\i,
	\end{align*}
	by Minkowski's inequality, proving that $P_A(T)\in H_A$. \\ \\
	Hence in order to check that $P_A(T)$ is indeed the projection of $T$ onto $H_A$ it remains to verify that $T-P_A(T) \perp h$ for all $h\in H_A$. Take any $h\in H_A$ and note that $h=g_A(X_i:i\in A)$ for some $g_A\in \mathcal{L}^2(\cY^{|A|},\pi_A(P_{(X_1,...,X_n)}))$. Then by using the bilinearity of inner products we get
	\begin{align*}
	&\la T-P_A(T), h \ra \\
	&= \la T- E[T|X_i:i\in A],g_A(X_i:i\in A)\ra- \sum_{B\subsetneq A} (-1)^{|A|-|B|}\la E(T|X_i:i\in B), g_A(X_i:i\in A) \ra  ,
	\end{align*}
	Now realize that $E[T|X_i:i\in A]$ is the orthogonal projection of $T$ onto the closed linear subspace $\mathcal{H}_A$ of all $\sigma(X_i:i\in A)$-measurable mappings, and therefore by \cite{schilling} Corollary 21.6(ii), we get that $T-E(T|X_i:i\in A)\in \mathcal{H}_A^\perp$, implying that the first term is zero. As to the second term we simply note that for any $B\subsetneq A$              
	\begin{align*}	
	\la E(T|X_i:i\in B), g_A(X_i:i\in A) \ra &= E  E(T|X_i:i\in B) E \lp  g_A(X_i:i\in A) |X_i:i\in B \rp  	\\
	&=0,
	\end{align*}
	by the defining property of random variables in $H_A$. We conclude that $P_A(T)\in H_A$ and that $T-P_A(T)\in H_A^\perp$, proving that $P_A(T)$ is indeed the orthogonal projection of $T$ onto $H_A$. 
	
	As regards  the orthogonal projection $P_A(T)$ when $A=\emptyset$ the formula still holds. In that case we have that the projection onto $H_\emptyset$ is given by $P_{\emptyset}(T)=(-1)^0E(T|\emptyset)=E(T)$, by the above mentioned convention about conditioning on the empty set. Lastly, we can identify $H_\emptyset = L^2(\Omega,\{\Omega,\emptyset\},P|_{\{\Omega,\emptyset\}})$, since  mappings that are measurable with respect to the trivial $\sigma$-algebra are constant and vice versa. Now we know that the orthogonal projection of $T\in L^2(\Omega,\F,P)$ onto $L^2(\Omega,\{\Omega,\emptyset\},P|_{\{\Omega,\emptyset\}})$ is given by the conditional expectation $E(T|\{\Omega,\emptyset\})$ conditioning on the trivial $\sigma$-algebra, which by Theorem 22.4 (xiii) in \cite{schilling} coincides with $E(T)$, proving that the formula holds.  \\ \\
	As regards the two last claims, assume  $T\perp H_B$ for all $B\subset A\subset \{1,...,n\}$ and note that if $|A|=0$ then if $T\perp H_\emptyset$ then $E(TP_{H_{\emptyset}}(T))=E(TE(T))=0\iff E(T)=0$, implying that $E(T|\emptyset)=E(T)=0$, so the assertion holds for $|A|=0$. Here we used that $P_{H_{\emptyset}}(T)=E(T)$ which is easily seen by using the above formula for the projection onto $H_A$ spaces. Now assume that the assertion also holds for any $A\subset \{1,...,n\}$ with $|A|=1,...,k$ hence by induction we are done if it holds for $A$ with $|A|=k+1$. The induction assumption implies that every term with $E(T|X_i:i\in B)=0$ for all $B\subsetneq A$, hence
	\begin{align*}
	P_A(T) &=E(T|X_i:i\in A)+ \sum_{B\subsetneq A} (-1)^{|A|-|B|}E(T|X_i:i\in B) = E(T|X_i:i\in A).
	\end{align*}
	But note that the assumption  $T\perp H_A$ implies that $P_{H_A}(T)=0$. Now using that $P_A(T)=P_{H_A}(T)$, we get  $E(T|X_i:i\in A)=0$, which proves the claim. \\ \\
	The very last claim can be verified by checking that $g(X_i:i\in A)\in \sum_{B\subset A} H_B$ or equivalently
	\begin{align*}
	G_A:&=g(X_i:i\in A)-P_{\sum_{B\subset A} H_B}(g(X_i:i\in A))  \\
	&= g(X_i:i\in A) - \sum_{B\subset A} P_B(g(X_i:i\in A)) =0,
	\end{align*}
	for all	  $g(X_i:i\in A)\in L^2(\Omega, P)$. Fix any such $G_A$ and note that $\{0\}\in H_B$ for all $B\subset A$ implying that $H_C\subset \sum_{B\subset A}H_B$ for all $C\subset A$. Hence $G_A \in \lp \sum_{B\subset A} H_B\rp^\perp \subset H_C^\perp$ for all $C\subset A$, or equivalently $G_A \perp H_B$ for all $B\subset A$, which by the above claim implies that $E(G_A|X_i:i\in A)=0$. But since $G_A$ is $\sigma(X_i:i\in A)$-measurable it follows that $G_A=0$.
\end{p}
\noindent The following theorem called the Hoeffding decomposition theorem gives an explicit representation of any symmetric square-integrable mapping in terms of its projections onto the above mentioned subspaces. This decomposition theorem will later allow us to decompose U-statistics (defined next section) in a beneficial way.
	\begin{theorem}[The Hoeffding decomposition] \label{theorem_Hoeffding_decomposition}
		Let $T:\cX^n \to \R$ be a symmetric measurable mapping such that $T(X_1,...,X_n)\in \mathcal{L}^2(\Omega,\F,P)$. Then we have the following decomposition holds almost surely
		\begin{align*}
		T(X_1,...,X_n) &= \sum_{i=0}^n \sum_{\underset{|A|=i}{A\subset \{1,...,n\}}} P_A(T(X_1,...,X_n))  \\
		&=\sum_{i=0}^n \sum_{\underset{|A|=i}{A\subset \{1,...,n\}}} T_i(X_{A_1},...,X_{A_{i}}),
		\end{align*}
		where $T_i:\cX^i \to \R$ is a symmetric function given by
		\begin{align*}
		T_i(x_1,...,x_i) = \sum_{B\subset \{1,...,i\}} (-1)^{i-|B|} E(T(x_{B_1},...,x_{B_{|B|}},X_{1},...,X_{n-|B|})).
		\end{align*}
			\end{theorem}
			An important thing to note is that for all $A\subset \{1,...,n\}$ with identical cardinality the projections $P_A$ is given by a fixed function with arguments $(X_i:i\in A)$.
\begin{p}
	First note that by the integrability condition we have that $T(X_1,...,X_n)\in \sum_{A\subset \{1,...,n\}}H_A$, by  \cref{lemma_projection_onto_H_A_spaces}(3). Hence $T(X_1,...,X_n)$ is identical to it's projection onto $\sum_{A\subset \{1,...,n\}}H_A$. By the mutual orthogonality of the spaces in $(H_A)_{A\subset \{1,...,n\}}$ that projection is given by the sum of projections onto $H_A$ for $A\subset \{1,...,n\}$. Each of these subspace projections can be expressed by formula in \cref{lemma_projection_onto_H_A_spaces}(i). Thus
	\begin{align*}
	T(X_1,...,X_n) &=  \sum_{A\subset \{1,...,n\}}  P_A(T(X_1,...,X_n)) \\ 
	&= \sum_{i=0}^n \sum_{\underset{|A|=i}{A\subset \{1,...,n\}}}P_A(T(X_1,...,X_n))  \\
	&=  \sum_{i=0}^n \sum_{\underset{|A|=i}{A\subset \{1,...,n\}}}  \sum_{B\subset A} (-1)^{|A|-|B|}E(T(X_1,...,X_n)|X_i:i\in B) \\ 
	&= \sum_{i=0}^n \sum_{\underset{|A|=i}{A\subset \{1,...,n\}}} \sum_{B\subset A} (-1)^{i-|B|}\Psi_{B}(X_{B_1},...,X_{B_{|B|}}),
	\end{align*}
	almost surely, where the mapping $\Psi_B:\cX^{|B|}\to \R$ is a $P_{X}^{|B|}$-almost everywhere unique mapping satisfying that it is a conditional expectation of $T(X_1,...,X_n)$ given $((X_{B_1},...,X_{B_{|B|}})=(x_1,...,x_{|B|}))$. In order words $\Psi_{B}(X_{B_1},...,X_{B_{|B|}})=E(T(X_1,...,X_n)|X_i:i\in B)$ almost surely. \\ \\
	Now fix any $A\subset \{1,...,n\}$ with $|A|=i\in\{0,...,n\}$ and $B\subset A$. By the symmetry of $T$, we have that $ T(X_1,...,X_n) =T(X_{\sigma(1)},...,X_{\sigma(n)})$	for any permutation $\sigma$ of $\{1,...,n\}$, so we may  change the order of the arguments. Hence with $N=\{1,...,n\}$
	\begin{align*}
	T(X_1,...,X_n) = T(X_{B_1},...,X_{B_{|B|}},X_{(N\setminus B)_{1}},...,X_{(N\setminus B)_{|N\setminus B|}}),
	\end{align*}
	which by similar arguments as in Corollary 2.2.4 \cite{beting} implies that for $P_X^{|B|}$-almost all $(x_1,...,x_{|B|})\in \cX^{|B|}$ that
	\begin{align*}
	\Psi_{B}(x_1,...,x_{|B|})&=E(T(X_{B_1},...,X_{B_{|B|}},X_{(N\setminus B)_{1}},...,X_{(N\setminus B)_{|N\setminus B|}})|(X_{B_1},...,X_{B_{|B|}})=(x_1,...,x_{|B|})) \\
	&= E(T(x_1,...,x_{|B|},X_{(N\setminus B)_{1}},...,X_{(N\setminus B)_{|N\setminus B|}})) \\
	&=E(T(x_1,...,x_{|B|},X_{1},...,X_{n-|B|})) \\
	&=: \Psi_{|B|}(x_1,...,x_{|B|}).
	\end{align*}
	Hence every mapping $\Psi_B$ for $B\subset A$ with identical cardinality coincides $P_X^{|B|}$-almost everywhere with $\Psi_{|B|}$. 
	This allows us to change  summation indexes in the following way
	\begin{align*}
	T(X_1,...,X_n)&=\sum_{i=0}^n \sum_{\underset{|A|=i}{A\subset \{1,...,n\}}} \sum_{B\subset A} (-1)^{i-|B|}\Psi_{|B|}(X_{B_1},...,X_{B_{|B|}})  \\
	&=\sum_{i=0}^n \sum_{\underset{|A|=i}{A\subset \{1,...,n\}}} \sum_{B\subset \{1,...,i\}} (-1)^{i-|B|}\Psi_{|B|}(X_{A_{B_1}},...,X_{A_{|B|}}),
	\end{align*}
	almost surely. Now realize that this is the form as stated in the theorem, that is
	\begin{align*}
	T(X_1,...,X_n)=\sum_{i=0}^n \sum_{\underset{|A|=i}{A\subset \{1,...,n\}}} T_i(X_{A_1},....,X_{A_{i}}),
	\end{align*}
	almost surely, where $T_i:\cX^{i}\to \R$ is given by
	\begin{align*}
	T_i(x_1,...,x_i) = \sum_{B\subset \{1,...,i\}} (-1)^{i-|B|} E(T(x_{B_1},...,x_{B_{|B|}},X_{1},...,X_{n-|B|})),
	\end{align*}	
		which by the symmetry of $T$ is itself a symmetric function.
	\end{p}

\newpage
\subsubsection{U-statistics} \label{Appendix_U-statistics} 
 Let $A$ be a index set and consider a family of probability distributions $\mathcal{P}=\{P_\alpha:\alpha\in A\}$ on a measurable space $(\cX,\bF)$ and a functional $\gamma: \mathcal{P}\to \R$. In the terminology of \cite{Theory_of_ustatistics} we assume that $\gamma$ is a regular functional, that is there exists a mapping (refereed to as the kernel) $h:\cX^m \to \R$ that is $P_\alpha^m$-integrable for all $\alpha\in A$, such that
 \begin{align*}
 \gamma(P_\alpha)= \int h(x_1,...,x_m) \, dP_\alpha^m(x_1,...,x_m).
 \end{align*}
 Let $(X_n)_{n\in \N}$ be a sequence of independent and identically distributed random variables each with distribution $P_X\in \mathcal{P}$ and let $X_{1,n}=(X_1,...,X_n)$ for all $n\geq 1$. We want to establish an estimator for $\gamma(P_X)$ and the obvious one is to simply estimate $\gamma(P_X)$ by $h(X_1,...,X_m)$ but in the case that we have $n>m$ observations there is unused samples. Hence we propose the unbiased estimator for $\gamma(P_X)$  given by the arithmetic mean
\begin{align} \label{eq_definition_U_statistics_non_symmetric_kernel}
\tilde{U}_n^m(h,X_{1,n})= \frac{1}{n_{(m)}}\sum_{1\leq i_1 \not = \cdots \not = i_m\leq n} h(X_{i_1},...,X_{i_m}),
\end{align}
where $n_{(m)}=n(n-1)\cdots (n-m+1)=(n-m)!/n!$ and the summation is over all $n_{(m)}$ possible $m$-permutations $(i_1,...,i_m)$ of $(1,...,n)$. We note that one may replace each term with the arithmetic mean of all $m!$ $m$-permutations of $(i_1,...,i_m)$. That is 
\begin{align} \label{eq_temp_1123}
\tilde{U}_n^m(h,X_{1,n}) =\frac{1}{n_{(m)}}\sum_{1\leq i_1 \not = \cdots \not = i_m\leq n} \frac{1}{m!} \sum_{\sigma\in \Pi_m} h(X_{i_{\sigma(1)}},...,X_{i_{\sigma(m)}}),
\end{align}
where $\Pi_m$ is the set of all permutations of $\{1,...,m\}$. 

To see this, fix any $m$-permutation $(i_1,...,i_m)$ of $(1,...,n)$ and define $S(i_1,...,i_m) = \sum_{\sigma \in \Pi_m} h(X_{i_{\sigma(1)}},...,X_{i_{\sigma(m)}})$ and note that this only contains $h(X_{j_1},...,X_{j_m})$ one time for each permutation $(j_1,...,j_m)$ of $(i_1,...,i_m)$ and nothing else. Now consider the expression $\sum_{1\leq i_1 \not = \cdots \not= i_m\leq n}S(i_1,...,i_m)$, the sum of $S(i_1,...,i_m)$ over all possible $m$-permutations $(i_1,...,i_m)$ of $(1,...,n)$. This sum consists solely of terms of the form $h(X_{j_1},...,X_{j_m})$ for $(j_1,...,j_m)$ being a $m$-permutation on $(1,...,n)$. For any fixed $m$-permutation $(j_1,...,j_m)$ of $(1,...,n)$, we shall count how many times $h(X_{j_1},...,X_{j_m})$ occurs in the sum we consider. We note that $h(X_{j_1},...,X_{j_m})$ only occurs in $S(i_1,...,i_m)$ whenever $(i_1,...,i_m)$ is a permutation of $(j_1,...,j_m)$ and in the affirmative it occurs only once. There are $m!$ terms in the sum $\sum_{1\leq i_1 \not = \cdots \not= i_m\leq n}$ such that $(i_1,...,i_m)$ is a permutation of $(j_1,...,j_m)$. Thus for every $m$-permutation $(j_1,...,j_m)$ of $(1,...,n)$ the term $h(X_{j_1},...,X_{j_m})$ appears $m!$ times, hence we conclude that \begin{align*}
\sum_{1\leq i_1 \not = \cdots \not = i_m\leq n} \sum_{\sigma\in \Pi_m} h(X_{i_{\sigma(1)}},...,X_{i_{\sigma(m)}}) = \sum_{1\leq i_1 \not = \cdots \not = i_m \leq n} m! h(X_{i_1},...,X_{i_m}),
\end{align*}
proving that the equality in \cref{eq_temp_1123} is valid. \\ \\
Now define the symmetrized version of $h$ by \begin{align*}
\bar{h}(x_{1},...,x_{m}) = \frac{1}{m!}\sum_{\sigma\in \Pi_m} h(x_{\sigma(1)},...,x_{\sigma(m)}) = \frac{1}{m!} \sum_{1\leq i_1 \not = \cdots \not = i_m \leq m} h(x_{i_1},...,x_{i_m}),
\end{align*}
such that
\begin{align*}
\tilde{U}_n^m(h,X_{1,n}) =\frac{1}{n_{(m)}}\sum_{1\leq i_1 \not = \cdots \not = i_m\leq n} \bar{h}(X_{i_{1}},...,X_{i_m}).
\end{align*}Since $\bar{h}$ is a symmetric mapping, it holds that for each $1\leq i_1 < \cdots <i_m\leq n$, there will be $m!$ terms in the above sum which are permutations of $(i_1,...,i_m)$ contributing the same amount. Hence we may write
\begin{align*}
\tilde{U}_n^m(h,X_{1,n}) &=\frac{m!(n-m)!}{n!}\sum_{1\leq i_1<\cdots <i_n \leq n} \bar{h}(X_{i_{1}},...,X_{i_m})\\
&= {n\choose m}^{-1} \sum_{1\leq i_1 < \cdots <i_m\leq n} \bar{h}(X_{i_{1}},...,X_{i_m}).
\end{align*}
Thus we may without loss of generality restrict the concept of these unbiased estimators to symmetric kernels and define U-statistics as follows. 
\begin{definition}
	Suppose that the kernel $h$ for the regular functional $\gamma$ is symmetric. Then the unbiased estimator $U_n^m(h,X_{1,n})$ for the parameter $\gamma(P_X)$, based on the first $n$-samples $X_{1,n}=(X_1,...,X_n)$ for $n>m$, called the U-statistic with symmetric kernel $h$ of degree $m$, is given by
\begin{align*}
U_n^m(h,X_{1,n})={n\choose m}^{-1} \sum_{1\leq i_1 < \cdots <i_m\leq n} h(X_{i_{1}},...,X_{i_m}).
\end{align*}
If the kernel $h$ is non-symmetric, then we define $\tilde{U}_n^m(h,X_{1,n})$ by \cref{eq_definition_U_statistics_non_symmetric_kernel} and note that
$\tilde{U}_n^m(h,X_{1,n})=U_n^m(\bar{h},X_{1,n})$, where $\bar{h}$ is the symmetrized version of $h$.
\end{definition}
\begin{remark}[Numerical considerations for a given $n$-sample]
	Say we are given an $n$-sample $X_{1,n}=x_{1,n}\in \cX^n$ and want to calculate $U_n^m(h,x_{1,n})$.  If $h$ initially was a symmetric mapping we obviously have that $\bar{h}=h$ and the definition of $U_n^m(h,X_{1,n})$ reduces ($m>1$) the number of computations of $U_n^m(h,x_{1,n})$ compared to the representation \cref{eq_definition_U_statistics_non_symmetric_kernel}, since we only need to deal with ${n\choose m}< \frac{n!}{(n-m)!}$ summands. In the case that $h$ is non-symmetric the representation doesn't matter, both have $\frac{n!}{(n-m)!}$ summands.
\end{remark}
\noindent In the further analysis of U-statistics we need to define the following mappings
\begin{definition}
	For any kernel $h:\cX^m \to \R$ we define the mappings $h_c:\cX^c\to \R$ by
	\begin{align*}
	h_c(x_1,...,x_c)=Eh(x_1,...,x_c,X_{c+1},...,X_{m}),
	\end{align*}
	for all $x_1,...,x_c\in\cX$ and $c=0,...,m$, e.g. $h_0=\gamma(P_X)$ and $h_m=h$. 	
\end{definition}
\noindent Recall that the kernel $h$ is  a $P_X^m$-integrable mappings so the conditional expectations are well-defined and we especially have that $h_c(x_1,...,x_c)$ is a conditional expectation of $h(X_1,...,X_m)$ given $(X_1,...,X_c)=(x_1,...,x_c)$.  That is,
\begin{align*}
h_c(x_1,...,x_c)&=  \int h(x_1,...,x_c,x'_{c+1},...,x'_m) dP_{X}^{m-c}(x'_{c+1},...,x'_m) \\
&=E(h(X_1,...,X_m)|(X_1,...,X_c)=(x_1,...,x_c)) ,
\end{align*}
for $P_X^c$-almost all $(x_1,..,x_c)\in\cX^c$, by similar arguments as in Corollary 2.2.4 of \cite{beting}.\\ \\
Having defined the conditional expectations we will now introduce some rather tedious and (for the moment) unintuitive recursively defined mappings.
\begin{definition} \label{appendix_definition_h^(k)}
	For any symmetric kernel $h:\cX^m \to \R$ we recursively define the mappings $h^{(k)}:\cX^k \to \R$ by
	\begin{align*}
	h^{(k)}(x_1,..,x_k)&=h_k(x_1,..,x_k)- \gamma(P_X)-\sum_{j=1}^{k-1}\sum_{1\leq i_1 < \cdots <i_{j} \leq k } h^{(j)}(x_{i_1},...,x_{i_{j}}),
	\end{align*}
	for $k=1,...,m$. With the convention that $\sum_{j=1}^0=0$ such that $h^{(1)}(x_1)=h_1(x_1) -  \gamma(P_X)$.
\end{definition} 
\noindent We may note that the mappings $h^{(k)}$ defined above, are themselves symmetric kernels for all $1\leq k \leq m$. Furthermore one can show that they possess the following properties.
\begin{theorem} \label{appendix_theorem_properties_of_h^(k)}
	For any symmetric kernel $h:\cX^m \to \R$ it holds that
	\begin{itemize}
		\item[(i)] $(h^{(k)})_c(x_1,...,x_c)=0$ for all $c=1,..,k-1$ and $k=1,...,m$.
		\item[(ii)] $Eh^{(k)}(X_1,...,X_k)=0$ for all $k=1,...,m$.
	\end{itemize}
\end{theorem}
\begin{p}
	See \cite{Lee_U-statistics_theory_and_practice} theorem 2 in section section 1.6
\end{p}
\noindent We will henceforth drop the first parenthesis and apply the convention that the superscript is always read first, i.e. the first equality will now be written as $h^{(k)}_c(x_1,...,x_c)=0$.  These recursively defined mappings $(h^{(k)})_{1\leq k \leq m}$ turns out to exactly be the orthogonal projections of $h(X_1,...,X_m)$ onto the spaces $(H_{\{1,...,k\}})_{1\leq k \leq m}$ defined in \cref{Appendix_Hoeffding_decomposition_section}, whenever the orthogonal projections exists. We stress that the orthogonal projections of $h(X_1,...,X_m)$ exist if $h(X_1,...,X_m)\in L^2(\Omega,\F,P)$, but to define $(h^{(k)})_{1\leq k\leq m}$ it suffices that $h(X_1,...,X_m)\in L^1(\Omega,\F,P)$. We will now show the equality of the recursively defined mappings and orthogonal projections in conjunction with the so-called Hoeffding decomposition of U-statistics.
\begin{theorem}[The Hoeffding decomposition of U-statistics] \label{appendix_theorem_Hoeffding_decomposition_U_stat}
Assuming that $h(X_1,...,X_m)\in \mathcal{L}^2(\Omega,\F,P)$ is a symmetric kernel of degree $m$, then  the orthogonal projection of the U-statistics $U_n^m(h,X_{1,n})$ onto the space $\sum_{B\subset \{1,...,n\}} H_B$, defined in \cref{Appendix_Hoeffding_decomposition_section}, exists. As a consequence we arrive at the following decomposition of the U-statistic into a linear combination of U-statistics with kernels of lower degrees:
\begin{align*}
U_n^m(h,X_{1,n}) -\gamma(P_X) &= \sum_{k=1}^m  {m \choose k} U_n^k(\Psi_k,X_{1,n}) \\
&=\sum_{k=1}^m  {m \choose k} U_n^k(h^{(k)},X_{1,n}),
\end{align*}
almost surely, where the mapping $\Psi_k:\cX^k \to \R$ is a orthogonal projection mapping of $h(X_1,...,X_m)$ onto $H_{\{1,...,k\}}$, that is,  $P_{\{1,...,k\}}(h(X_1,...,X_m))=\Psi_k(X_1,...,X_k)$ and it is given by
\begin{align*}
\Psi_k(x_1,...,x_k) &= \sum_{B\subset \{1,...,k\}} (-1)^{k-|B|} h_{|B|}(x_{B_1},...,x_{B_{|B|}}) \\ \\
&=h^{(k)}(x_1,...,x_k),
\end{align*}
where these last two equalities also hold if $h(X_1,...,X_m)\in \mathcal{L}^1(\Omega,\F,P)$.
\end{theorem}
\begin{p}
	Using the Hoeffding decomposition - \cref{theorem_Hoeffding_decomposition} - on the square integrable (Minkowski's inequality) mapping \begin{align*}
	U_n^m(h,X_{1,n})=T(X_1,...,X_n),
	\end{align*}
	for some measurable mapping $T:\cX^n \to \R$,  we get that $U_n^m(h,X_{1,n})$ can be written as the projection onto the space $\sum_{B\subset \{1,...,n\}}H_B$ in the following way
	\begin{align*}
	U_n^m(h,X_{1,n}) &= P_{\sum_{B\subset \{1,...,n\}}H_B}(U_n^m(h,X_{1,n}))\\
	&=\sum_{i=0}^n   \sum_{\underset{|A|=i}{A\subset \{1,...,n\}}} P_A \lp 
	 {n\choose m}^{-1} \sum_{1\leq t_1< \cdots <t_m \leq n } h(X_{t_1},...,X_{t_m}) \rp \\
	 &= \sum_{i=0}^n    
	 {n\choose m}^{-1} \sum_{\underset{|A|=i}{A\subset \{1,...,n\}}}\sum_{1\leq t_1< \cdots <t_m \leq n }P_A \lp  h(X_{t_1},...,X_{t_m}) \rp .
	\end{align*}
	Note that for any $A\not \subset \{t_1,...,t_m\}$ then $h(X_{t_1},....,X_{t_m})\in \sum_{B \subset \{t_1,...,t_m\}}H_B \perp H_A$ since $H_A\perp H_B$ for all $B\not = A$, implying that $P_A(h(X_{t_1},...,X_{t_m}))=0$, so these terms does not contribute to the above summation. Hence we can for starters remove the terms for $m<i\leq n$ without changing anything. Now for any two $1\leq s_1< \cdots < s_m \leq m$ and $1\leq t_1 < \cdots < t_m \leq n$ with $A \subset \{t_1,...,t_m\}, \{s_1,...,s_m\}$ we have that 
	\begin{align*}
	P_A(h(X_{t_1},...,X_{t_m}))=P_A(h(X_{s_1},...,X_{s_m})) = \Psi_{|A|}(X_i : i \in A).
	\end{align*}
	This is seen by inspecting the representation of $P_A$ in the Hoeffding decomposition (\cref{theorem_Hoeffding_decomposition}) and noting that since  $X_1,...,X_n$ is i.i.d. then $(X_{s_1},...,X_{s_{m-|B|}})\eqd (X_{t_1},...,X_{t_{m-|B|}})$ for any $B\subset A$, implying that the projections are given by identical functions $\Psi_{|A|}$ composed with $(X_i:i\in A)$.  Now for any fixed $A\subset \{1,...,n\}$ with $|A|=i \in \{0,...,m\}$ all $1\leq t_1 < \cdots <t_m \leq n$ that does not contain $A$ yields a zero, but on the other hand as we argued above, every partitioning that does contain $A$ yields the same projection term. When summing over all $m$-partitionings $\sum_{1\leq t_1< \cdots <t_m \leq n }$ we have exactly ${n-i \choose m-i}$ terms which contains $A$, implying that
	\begin{align*}
	U_n^m(h,X_{1,n}) &= \sum_{i=0}^m   \sum_{\underset{|A|=i}{A\subset \{1,...,n\}}} {n\choose m}^{-1}{n-i \choose m-i}  \Psi_{i}(X_i : i \in A) \\
	&=\gamma(P_X)+\sum_{i=1}^m   \sum_{\underset{|A|=i}{A\subset \{1,...,n\}}} {n\choose m}^{-1}{n-i \choose m-i}  \Psi_{i}(X_i : i \in A) ,
	\end{align*}
	where we used that $\Psi_0 =E(h(X_1,...,X_m)|\emptyset)$ as per above mentioned convention (see explicit formula for $P_{\emptyset}$ in \cref{theorem_Hoeffding_decomposition}) is the conditioning on the trivial $\sigma$-algebra $\{\Omega,\emptyset\}$. Conditioning on the trivial $\sigma$-algebra is simply the regular expectation, that is $\Psi_0=Eh(X_1,...,X_m)=\gamma(P_X)$.
	The above binomial coefficient factors can be rewritten as
	\begin{align*}
	{n\choose m}^{-1}  {n-i \choose m-i} &= \frac{m!(n-m)!(n-i)!}{n!(m-i)!(n-m)!} =\frac{m!}{i!(m-i)!} \frac{i!(n-i)!}{n!} = {m \choose i} {n \choose i}^{-1},
	\end{align*}
	and as a consequence we get that \begin{align*}
	U_n^m(h,X_{1,n})-\gamma(P_X) &= \sum_{k=1}^m   {m \choose k} {n \choose k}^{-1} \sum_{1\leq i_1 < \cdots <i_k\leq n}  \Psi_{k}(X_{i_1},...,X_{i_k}) \\
	&=\sum_{k=1}^m  {m \choose k} U_n^k(\Psi_k,X_{1,n}),
	\end{align*}
	where
	\begin{align*}
	U_n^k(\Psi_k,X_{1,n}) :=   {n \choose k}^{-1}\sum_{1\leq t_1 < \cdots < t_k \leq n}   \Psi_{k}(X_{t_1},...,X_{t_k}),
	\end{align*}
	is itself is a $n$-sample U-statistic of degree $k$ with symmetric kernel $\Psi_k$. \\ \\
	Using the explicit formula for $\Psi_k$ found in \cref{theorem_Hoeffding_decomposition} we have that
	\begin{align*}
	\Psi_k(x_1,...,x_k) &= \sum_{B\subset \{1,...,k\}} (-1)^{k-|B|} E(h(x_{B_1},...,x_{B_{|B|}},X_{1},...,X_{n-|B|})) \\
	&=\sum_{B\subset \{1,...,k\}} (-1)^{k-|B|} h_{|B|}(x_{B_1},...,x_{B_{|B|}}) \\
	&= h_k(x_1,...,x_k)+\sum_{j=1}^{k-1} \sum_{\stackrel{B\subset \{1,...,k\}}{|B|=j}}  (-1)^{k-j} h_{j}(x_{B_1},...,x_{B_j})+(-1)^k h_0\\
		&=h_k(x_1,...,x_k)+ \sum_{j=1}^{k-1} \sum_{1\leq i_1 < \cdots < i_j\leq k}  (-1)^{k-j} h_{j}(x_{i_1},...,x_{i_j})+(-1)^{k}\gamma(P_X),
	\end{align*}
	 Now realize that we can re-index (reverse the order of summation) the outer sum over with $d=j-k$ to get
	\begin{align*}
\Psi_k(x_1,...,x_k) &=h_k(x_1,...,x_k)+ \sum_{d=1}^{k-1} \sum_{1\leq i_1 < \cdots < i_{k-d}\leq k}  (-1)^{d} h_{k-d}(x_{i_1},...,x_{i_{k-d}})+(-1)^{k}\gamma(P_X) \\
&=h^{(k)}(x_1,...,x_k),
	\end{align*}
	where we in the last equality used the identity given in equation (11) in section 1.6 of \cite{Lee_U-statistics_theory_and_practice}. This identity is proved by rewriting $h^{(k)}$ is terms of integrals followed by further manipulations. The specific steps are notation heavy and does not provide further insight into the nature of the mappings $h^{(k)}$, hence we refer the reader to \cite{Lee_U-statistics_theory_and_practice} section 1.6 for the proof of the identity. Thus we also have that
	\begin{align*}
	\sum_{k=1}^m  {m \choose k} U_n^k(\Psi_k,X_{1,n})=\sum_{k=1}^m  {m \choose k} U_n^k(h^{(k)},X_{1,n}),
	\end{align*}
	which is what we wanted to show.
\end{p}
\begin{corollary}
	The above decomposition of $U_n^m(h,X_{1,n})$ holds even though $h(X_1,...,X_m)\in \mathcal{L}^1(\Omega,\F,P)$, but the geometric property that $\Psi_k=h^{(k)}$ are projection mappings does no longer hold.
\end{corollary}
\begin{p}
	See page 8-9 in \cite{Ustatistics_in_banach_spaces_Borovskikh} or Lemma $A$ in section 5.1.5 of \cite{serfling2009approximation}.
\end{p}
\newpage
\subsubsection{V-statistics} \label{Appendix_V-statistics}
Consider the exact same set-up as in the above section on U-statistics. That is, we have a sequence of independent and identically distributed random elements $(X_n)_{n\in\N}$ in a measurable space $(\cX,\F)$. We want to estimate 
\begin{align*}
\gamma(P_X) = \int_{\cX^m} h(x_1,...,x_m)dP_X^m(x_1,...,x_m)=Eh(X_1,..,X_m),
\end{align*}
for $m\leq n$ and $h$ a  kernel. \\ \\
We will now construct an estimator for $\gamma(P_X)$ based on the $n$ first samples  $X_{1,n}=\lp X_1,...,X_n\rp.$
\begin{definition}
	The in general biased estimator $V_n^m(h,X_{1,n})$ for the parameter $\gamma(P_X)$, based on the first $n$-samples $X_{1,n}=(X_1,...,X_n)$, called the V-statistic with kernel $h$ of degree $m$, is given by
\begin{align*}
V_n^m(h,X_{1,n}) &= \int_{\cX^m} h(x_1,...,x_m)d\lp P_X^{(n)}\rp^m(x_1,...,x_m) \\
&= \frac{1}{n^m} \sum_{i_1=1}^n \cdots \sum_{i_m=1}^n h(X_{i_1},...,X_{i_m}),
\end{align*}
where $P_X^{(n)}$ is the random empirical measure of $P_X$ based on $X_{1,n}.$
\end{definition}
\newpage
\subsubsection{Various results and definitions} \label{Appendix_various_results_and_defitions_for_U-_and_V-statistics}

\begin{definition}[Degeneracy] \label{Appendix_definition_degenerate_kernel}
	Let $h:\cX^m\to\R$ be a symmetric kernel. The rank of $h$ or the corresponding U- or V-statistic is defined as the smallest integer $r$ such that
	\begin{align*}
	h^{(1)}(X_1)=\cdots=h^{(r-1)}(X_1,...,X_{r-1})= 0,
	\end{align*}
	almost surely and
	\begin{align*}
	h^{(r)}(X_1,...,X_r)\not = 0,
	\end{align*}
	with positive probability.	We say the the kernel is $P_X$-degenerate of order $d=r-1$ and if $r=1$ it is called non-degenerate and in the case that $r=m$ we say that it is completely degenerate.
\end{definition}
\noindent In the literature there is (at least) two non-equivalent ways of defining degeneracy of kernels. The above definition of degeneracy of kernels coincide with that of \cite{Ustatistics_in_banach_spaces_Borovskikh}, \cite{Theory_of_ustatistics} and \cite{SLLN_for_V-statistics}, which allows for stronger results than authors who defines degeneracy as below. For example in \cite{SLLN_for_V-statistics} we have a SLLN for degenerate U-statistics, which has weaker convergence conditions than the square integrability required to define degeneracy, using the definition from \cite{Lee_U-statistics_theory_and_practice}, \cite{serfling2009approximation} and \cite{van2000asymptotic}. They define the degeneracy of kernels that has second moment $Eh(X_1,...,X_m)^2<\i$, in the following way.  If the second moment exists then the following constants are well-defined
\begin{align*}
\sigma_i^2&= \textrm{Var}(h_i(X_1,...,X_i)) \\
&=\textrm{Cov}(h(X_1,...X_i,X_{i+1},...,X_m),h(X_1,...X_i,X'_{i+1},...,X'_m))
\end{align*}
for all $i\geq 1$ and $\sigma_0^2=0$. Then they define the kernel $h$ to be degenerate of order $d$ if
\begin{align*}
0=\sigma_0^2=\cdots =\sigma_{d}^2<\sigma_{d+1}^2<\i.
\end{align*}
For consistency - since we use results from literature defining degeneracy in both ways - we show that definitions are equivalent under the assumption of square integrability of the kernel.
\begin{lemma}
	If $Eh(X_1,...,X_m)^2<\i$ then the above two definitions of degenerate kernels are equivalent.
\end{lemma}
\begin{p}
	Assume that $Eh(X_1,...,X_m)^2<\i$ and that $h$ is $P_X$-degenerate of order $d$, that is
	\begin{align*}
0=h^{(1)}(X_1)=\cdots=h^{(d)}(X_1,...,X_{d}),
\end{align*}
almost surely and $h^{(d+1)}(X_1,...,X_{d+1})\not =0$ with positive probability. Thus by the recursive nature of $(h^{(j)})$ (formally by an induction argument as below) we may realize that $h_i(X_1,...,X_i)=\gamma(P_X)=Eh_i(X_1,...,X_i)$ almost surely  implying that $$\sigma_i^2=\textrm{Var}(h_i(X_1,...,X_i))=0,$$ for all $i\leq d$. By the definition of $h^{(d+1)}$ we see that 
\begin{align*}
0\not =  h^{(d+1)}(X_1,...,X_{d+1}) \iff h_{d+1}(X_1,...,X_{d+1})\not =\gamma(P_X).
\end{align*}
Thus  $h_{d+1}(X_1,...,X_{d+1})$ is with positive probability not equal to its mean, hence we have that $\sigma_{d+1}^2=\textrm{Var}(h_{d+1}(X_1,...,X_{d+1}))>0$. Furthermore
\begin{align*}
Eh_{d+1}(X_1,...,X_{d+1})^2 &= EE(h(X_1,...,X_m)|X_1,...,X_{d+1})^2 \\
& \leq E E(h(X_1,...,X_m)^2 |X_1,...,X_{d+1} ) \\
&= Eh(X_1,...,X_m)^2 \\
&<\i,
\end{align*}
by Jensen's conditional  inequality, proving that $\sigma_{d+1}^2<\i$. \\ \\
Conversely if $0=\sigma_0^2=\cdots \sigma_d^2<\sigma_{d+1}^2<\i$ then we start by inductively showing that $0=h^{(1)}(X_1)=\cdots=h^{(d)}(X_1,...,X_{d})$. We obviously have that $h^{(1)}(X_1)=h_1(X_1)-\gamma(P_X)=0$ almost surely, showing the induction basis. Now for the inductive step assume that 
$h^{(j)}(X_1,...,X_j)=0$ almost surely for all $j\leq d-1$ and note that this also holds for any $j$'element subset $(X_{i_1},...,X_{i_j})$ of $(X_n)_{n\in\N}$. Thus
	\begin{align*}
h^{(j+1)}(X_1,...,X_{j+1})&=h_{j+1}(X_1,..,X_{j+1})- \gamma(P_X)-\sum_{k=1}^{j}\sum_{1\leq i_1 < \cdots <i_{k} \leq j+1 } h^{(k)}(X_{i_1},...,X_{i_{k}}) \\
&\stackrel{a.s.}{=} h_{j+1}(X_1,..,X_{j+1})- \gamma(P_X),
\end{align*}
but since $j\leq d-1$ we have that $\sigma_{j+1}^2=0$ implying that the above difference vanishes, proving that $h^{(j+1)}(X_1,...,X_{j+1})=0$ almost surely. By induction we now have that $0=h^{(1)}(X_1)=\cdots=h^{(d)}(X_1,...,X_{d})$. As argued above these stay almost surely zero for any such independent and identically distributed arguments. Hence the  double sum of   $h^{d+1}(X_1,...,X_{d+1})$ vanishes, such that
\begin{align*}
h^{(d+1)}(X_1,...,X_{d+1})&=h_{d+1}(X_1,..,X_{d+1})- \gamma(P_X),
\end{align*}
almost surely.
Since $0<\sigma_{d+1}^2$ we have that $h_{d+1}(X_1,..,X_{d+1})$ is non-degenerate and hence different from its mean $\gamma(P_X)$ with positive probability. Thus  $h^{(d+1)}(X_1,...,X_{d+1})\not = 0$ with positive probability.
\end{p}
\begin{corollary} \label{appendix_corollary_h^(k)_is_completely_degenerate}
	For any non-zero symmetric kernel $h:\cX^m\to\R$ it holds that $h^{(k)}$ defined in \cref{appendix_definition_h^(k)} are complete degenerate kernels for all $1\leq k \leq m$.
\end{corollary}
\begin{p}
	This is an immediate consequence of \cref{appendix_theorem_properties_of_h^(k)}
\end{p}
\noindent In the thesis we are going to use the following decomposition theorems of V-statistics which are similar to the above proven Hoeffding decomposition of $U$-statistics.
\begin{lemma} \label{appendix_lemma_V-statistic_decomposed_into_V-statistics}
	A centered V-statistic with symmetric kernel $h$ of degree $m$ can be decomposed into a linear combination of V-statistics. That is,
	\begin{align*}
	V_n^m(h,X_{1,n})- \gamma(P_X) = \sum_{c=1}^m {m \choose c} V_n^c(h^{(c)},X_{1,n}).
	\end{align*}
\end{lemma}
\begin{p}
	See section 1.3  in \cite{Ustatistics_in_banach_spaces_Borovskikh}. 
\end{p}
\begin{lemma}  \label{appendix_lemma_V-statistic_decomposed_into_U-statistics}
	A V-statistic with symmetric kernel $h$ of degree $m$ can be decomposed into a linear combination of U-statistics. That is\begin{align*}
	V_n^m(h,X_{1,n}) = \sum_{c=1}^{m} {n\choose c}n^{-m} U_n^c(h_{mc},X_{1,n}),
	\end{align*}
	where $h_{mc}:\cX^c \to \R$ is a symmetric and measurable mapping given by
	\begin{align*}
	h_{mc}(x_1,...,x_c) = \sum_{\stackrel{v_1+\cdots+v_c=m}{v_i\geq 1}} \frac{m!}{v_1!\cdots v_c!}h(x_1^{(v_1)},...,x_c^{(v_c)}),
	\end{align*}
	with $x_i^{(v_j)}=(x_i,...,x_i)\in \cX^{v_j}$.	
\end{lemma}
\begin{p}
	See section 1.3  in \cite{Ustatistics_in_banach_spaces_Borovskikh} or theorem 1 of section 4.2 in \cite{Lee_U-statistics_theory_and_practice}.
\end{p}
\newpage
\subsubsection{Asymptotic results for U- and V-statistics}
In the following let $(\cX,\F)$ be a measurable space and let $(X_i)_{i\in \N}$ be an i.i.d. sequence of random elements in with values in $\cX$ and corresponding distribution $P_X$.
\begin{theorem} [Asymptotic distribution of degenerate U-Statistics] \label{appendix_theorem_asymp_dist_of_2_deg_degenerate_U_stat}
	Let $h:\cX^m \to \R$ be a symmetric and measurable kernel of degree $m\leq n$. If $h$ is $P_X$-degenerate of order 1 and $h(X_1,...,X_m)\in \mathcal{L}^2(\Omega,P)$, then 
	\begin{align*}
	n (U_n^m(h,X_{1,n})-\gamma(P_X)) \convd \frac{m(m-1)}{2} \sum_{i=1}^\i \lambda_i (Z_i^2-1),
	\end{align*}
	as $n$ tends to infinity, where $(Z_i)_{i\in N}$ are independent and identically standard normal distributed, and $(\lambda_i)_{i\in \N}$ are the real eigenvalues (counting algebraic multiplicity) of the operator $A:L^2(\cX,\bF,P_{X})\to L^2(\cX,\bF,P_{X})$ given by
	\begin{align*}
	A(f)(x)&= \int_{\cX} h^{(2)}(x,y)f(y)dP_{X}(y)\\
	&=\int_{\cX} [Eh(x,y,X_3,...,X_m)-\gamma(P_X)]f(y)dP_{X}(y).
	\end{align*}
	That is $(\lambda_i)_{i\in \N}$ consists of every $\lambda$ for which there exists a $f\in L^2(\cX,\bF,P_{X})$ such that
	\begin{align*}
	A(f)-\lambda f = 0,
	\end{align*}
	as a mapping in $L^2(\cX,\F,P_{X})$
	
\end{theorem}
\begin{p}
	See corollary 4.2.2 \cite{Ustatistics_in_banach_spaces_Borovskikh}, theorem 4.3.1 ($m=2$) / corollary 4.4.2 \cite{Theory_of_ustatistics} or theorem 5.5.2 ($\cX=\R$) \cite{serfling2009approximation}. 
	\end{p}
	\begin{theorem}[Strong Law of Large Number for U-statistics] \label{theorem_SLLN_for_U-statistics}
	Assume that $h:\cX^k \to \R$ is a symmetric and measurable function. 
	If
	\begin{align*}
	E|h(X_{1},...,X_{k})| <\i \iff h\in \mathcal{L}^1(\cX^k,P_X^k),
	\end{align*}
	then
	\begin{align*}
	U_n^k(h,P_X) = \frac{1}{{n\choose k}} \sum_{1\leq i_1 < \cdots < i_k \leq n} \cdots \sum_{i_k=1}^n h(X_{i_1},...,X_{i_k}) \convas E(h(X_1,...,X_k)).
	\end{align*}	
	\end{theorem}
	\begin{p}
	See theorem 3.1 \cite{Ustatistics_in_banach_spaces_Borovskikh} or theorem 3.1.2 \cite{Theory_of_ustatistics} or as first proved in \cite{StrongLawOfLargeNumbersForUStatistics_hoeffding1961}.

\end{p}
\begin{theorem}[Strong Law of Large Numbers for Degenerate U-statistics] \label{theorem_SLLN_for_degenerate_U-statistics}
	Let $h:\cX^k\to \R$ be a symmetric and measurable kernel which is $P_X$-degenerate of order $0\leq r-1<k$, and let $k-r/2<s<k$. If $E|h(X_1,...,X_k)|^{r/(s+r-k)}<\i$, then
	\begin{align*}
	n^{k-s}\lp U_n^k(h,X_{1,n})-Eh(X_1,...,X_k) \rp\convas_n 0.
	\end{align*}That is, \begin{align*}
	n^{-s}\sum_{1\leq i_1 < \cdots <i_k \leq n} \lp h(X_{i_1},...,X_{i_k})-Eh(X_1,...,X_k) \rp\convas_n 0.
	\end{align*}
\end{theorem}
\begin{p}
	See theorem 2 \cite{SLLN_for_V-statistics}.
\end{p}
\begin{theorem}[Strong Law of Large Number for V-statistics] \label{theorem_SLLN_for_V-statistics}
	Assume that $h:\cX^k \to \R$ is a symmetric and measurable function. If
	\begin{align*}
	E|h(X_{i_1},...,X_{i_k})|^{\frac{\#\{i_1,...,i_k\}}{k}} <\i, \quad \quad \forall\,\, 1\leq i_1\leq \cdots \leq  i_k\leq k,
	\end{align*}
	then
	\begin{align*}
	V_n = \frac{1}{n^k} \sum_{i_1=1}^n \cdots \sum_{i_k=1}^n h(X_{i_1},...,X_{i_k}) \convas E(h(X_1,...,X_k)).
	\end{align*}	
\end{theorem}
\begin{p}
	see proposition on page 274 in \cite{SLLN_for_V-statistics} or proposition 2.3 in \cite{SLLN_for_V-statistics_2}.
\end{p}
\newpage
\subsection{Integration of Hilbert space valued mappings \label{Integration of Hilbert space valued mappings}} 
This section is an short introduction into Pettis integration and is inspired by the construction method of the Pettis integral seen in \cite{Tensor_products_of_Banach_spaces} and  \cite{TopicsInBanachSpaceIntegration2005}. The goal is to establish some theory  which allows for the  integration of Hilbert space valued mappings $f:\cX \to \cH$, where $(\cX,\F,\mu)$ is a finite measure space $\cH$ is a $\bK$-Hilbert space, where the scalar field $\bK$ is either $\R$ or $\bC$. \\ \\
The weak Pettis integral was introduced by Billy James Pettis in  \cite{pettis1938integration} for Banach space valued mappings but we will restrict ourself to Hilbert space valued mappings, since it suffices for our purpose in this thesis.
There are different notions of integrals with values in Hilbert spaces. Two of these integrals are the \textit{strong} Bochner integral and the above mentioned \textit{weak} Pettis integral.

 The Bochner integral which allows for integration of Banach space valued mappings and \textit{strong} reefers to the fact that the integral is constructed as the limit of integrals of simple mappings which approximates the integrand. That a mapping  has such a sequence of approximating simple mappings is called being strongly measurable. A theorem called the Pettis measurability theorem yields equivalence between being strongly measurable and being both weakly measurable (defined below) and  $\mu$-essentially separably valued (see \cite{Tensor_products_of_Banach_spaces} section 2.3). However the a mapping only needs to be weakly measurable in order to define its Pettis integral. 
 
  This section on Pettis integration was written before restricting the thesis to marginal spaces that are  separable (and as a consequence the Hilbert spaces of attention are separable). This is why, even though we actually have strong measurability of the mappings we intend to integrate, we still use the Pettis integral. As we see in the main part of the thesis, the Pettis integral suffices for our needs (in fact one can show that the Pettis and Bochner integral coincide in our cases). 
 
The name \textit{weak} comes from the fact that we only impose weak conditions on the integral. Unlike the Bochner integral (which is constructed as the limit of simple mappings),  the weak Pettis integral of a sufficiently nice integrand $f:\cX \to \cH$ over $E\in \F$ with respect to $\mu$ is defined as the unique element $\int_E f\, d\mu \in \cH$ which satisfies 
\begin{align*}
h^*\lp \int_E f\, d\mu \rp = \int_E h^* \circ f\, d\mu \in \bK,
\end{align*}
for all $h^* \in \cH^*$ (see below). This \textit{weak} uniquely determining property seems to be a reasonable starting requirement for the integral, since it also holds for the Bochner integral and in some sense also conforms with the interpretation that integrals are related to  sums which exhibit similar properties. In this explanation of the Pettis integral we already claimed uniqueness and existence of the Pettis integral element in $\cH$, and this is essentially what the remaining part of this section sets out to prove. \\ \\
For the rest of the section, assume that $(\cX,\F,\mu)$ is a finite measure space, $\cH$ is a $\R$-Hilbert space and $f:\cX\to \cH$ is some mapping. The arguments for $\bC$-Hilbert spaces follows similarly.  \\ \\
  We also let $\cH^*$ denote its continuous dual space, that is
\begin{align*}
\cH^*= \{h^*:\cH \to \R(=\bK) \, | \, h^* \textit{ is continuous and linear} \},
\end{align*}
is the space of all linear mapping from $\cH$ to $\R$.  Now we define the measurability and integrability conditions our integrands should possess in order for the construction of the Pettis integral to be successful.
\begin{definition}[Weakly measurable mappings]
	A mapping $f:\cX \to \cH$ is said to be weakly $\F$-measurable (or scalarly measurable) if $h^*\circ f$ is $\F/\cB(\R)$-measurable, for all $h^*\in \cH^*$.
\end{definition}
\begin{definition}[Scalar integrable mappings] \label{appendix_defi_scalarly_integrable}
	A weakly $\F$-measurable mapping $f:\cX\to \cH$ is called scalarly $\mu$-integrable if $h^*\circ f\in L^1(\mu)$, for all $h^*\in \cH^*$. \\
\end{definition}
We prove the existence of the Pettis integral for a class of mappings, by proving the existence of the Dunford integral and then establishing a link between them. That is, for a suitable mapping $f:\cX\to \cH$ we prove the existence of the Dunford integral of $f$ over $E\in \F$ with respect to $\mu$, written as $(D)\int_E f\, d\mu$, which is an element in $\cH^{**}=(\cH^*)^*$. Then we introduce the natural injective embedding $ev:\cH\to \cH^{**}$ and if $(D)\int_E f\, d\mu\in ev(\cH)$ we say that $f$ is Pettis integrable over $E$ with respect to $\mu$ and define the Pettis integral as the element in $\cH$ which maps to the Dunford integral via $ev$. Thus the first order of business is to prove the existence of the Dunford integral. \\ \\
Before proceeding to show the existence of the Dunford integral, we have to establish which topology we equip $\cH^*$ with before we even start to talk about its continuous dual. In general we equip every space of continuous linear  mappings between normed vector spaces $X$ and $Y$ with the operator norm. This norm is denoted $\|\cdot \|_{\text{op}}:\text{Lin}(X,Y)\to \overline{\R}$ and equivalently defined by either of the following expressions
\begin{align*}
\|L\|_{\text{op}} &= \sup \{ \|L(x)\|_Y : x\in X \text{ with } \|x\|_X\ \leq 1 \} \\
&= \sup \{ \|L(x)\|_Y/\|x\|_X:x\in X \text{ with } x\not = 0 \} \\
&= \inf \{ c \geq 0 : \|L(x)\|_Y \leq c \|x\|_X \text{ for all } x\in X\}.
\end{align*}
It is well known that this is indeed a norm on the subspace of continuous linear mappings $X^*$, it in fact makes $X^*$ a Banach space (It is important that it is complete for later arguments). Thus the operator norm of a linear mapping is finite (that is, a linear map is bounded)  if and only if the linear map is continuous. To see the only if part simply note that if $\|L\|_{\text{op}}<\i$ there exists a $c\geq 0 $ such that
\begin{align*}
\|L(x+h)-L(x)\|_Y = \|L(h)\|_Y \leq c\|h\|_X \stackrel{h\to 0}{\longrightarrow} 0,
\end{align*}
for any $x\in X$, proving continuity. As a last remark about the operator norm; one should realize that $\|L(x)\|\leq \|L\|_{\text{op}}\|x\|_X$ for any $x\in X$, which trivially follows from the second of the equivalent definitions of $\|\cdot \|_{\text{op}}$. \\ \\ 
Now assume that $f:\cX \to \cH$ is scalarly $\mu$-integrable. For any measurable set $E\in \F$ define the linear mapping $T_E:\cH^* \to L^1(\mu)$ by letting 
\begin{align*}
T_E(h^*) = h^*(1_Ef) =1_E \cdot (h^*\circ f), \quad \quad \quad \forall h^*\in \cH^*.
\end{align*}
Furthermore we define the linear integral operator $I_E:\cH^*\to \R$  by letting
\begin{align*}
I_E(h^*)=\int T_E ( h^* ) \, d\mu = \int_E h^* \circ f \, d\mu, \quad \quad\quad  \forall h^*\in \cH^*.
\end{align*}
Now if $I_E$ is continuous it is an element of $\cH^{**}$ and then we define $(D)\int_E f \, d\mu:= I_E$ and denote it the Dunford integral of $f$ over $E\in \F$ with respect to $\mu$. By linearity of $I_E$ it suffices to show that $I_E$ is continuous in zero. Thus note; for any $h^* \in \cH^*$ that
\begin{align*}
|I_E(h^*)| \leq \int_E |T_E(h^*)|d\mu = \| T_E(h^*)\|_{L^1(\mu)} \leq \|T_E\|_{\text{op}} \|h^*\|_{\text{op}},
\end{align*}
which tends to zero when $h^*$ tends to zero, if $\|T_E\|_{\text{op}}<\i$. That is, if  $T_E$ is a continuous map then $I_E$ is a linear continuous map. $T_E$ is indeed continuous and one can realize this by utilizing the closed graph theorem described below.
\begin{theorem}[Closed graph theorem]
	Assume that $X$ and $Y$ are Banach spaces and $\Lambda:X\to Y$ is a linear map such that, whenever $(x_n)_{n\in\N}\subset X$  and $(\Lambda(x_n))_{n\in\N} \subset Y$ are Cauchy sequences, then $\lim_{n\to \i } \Lambda(x_n)= \Lambda(\lim_{n\to\i} x_n)$. Then $\Lambda$ is a continuous mapping.
\end{theorem}
\begin{p}
	See  theorem 2.15 and its remark in \cite{rudin1991functional} for proof.
\end{p}
Now to use this theorem we simply assume that $(h^*_n)_{n\in\N}\subset \cH^*$ and $(T_E(h^*_n))_{n\in\N}\subset L^1(\mu)$ are Cauchy sequences with limits, say $h^*\in \cH^*$ and $g\in L^1(\mu)$. By the  above Closed graph theorem it suffices to show that $g=T_E(h^*)$ in $L^1(\mu)$ in order to prove that $T_E$ is continuous. Now recall that $L^1(\mu)$-convergence implies convergence in $\mu$-measure which in turn implies that there exists a subsequence $(T_E(h^*_{n_k}))_{k\in\N}$ converging $\mu$-almost everywhere to $g$. On the other hand we note that $\|h^*_n-h^*\|_{op}\to_n 0$ implying that
\begin{align*}
| T_E(h_n^*)(x) - T_E(h^*)(x) | &=  |T_E(h^*_n-h^*)(x)| \leq |(h^*_n-h^*)\circ f (x)| \leq \|h^*_n-h^*\|_{op} \|f(x)\|_\cH,
\end{align*}
tends to zero as $n$ tends to infinity. We have established point-wise convergence of $(T_E(h^*_{n}))_{n\in\N}$ to $T_E(h^*)$ but we also showed the existence of a subsequence hereof converging $\mu$-almost everywhere to $g$. Hence we must have that the limits coincide $\mu$-almost everywhere, that is $T_E(h^*)(x)=1_Eh^*(f(x))=g(x)$ for $\mu$-almost all $x\in \cX$. Identification of mappings up to $\mu$-almost everywhere equality in $L^1(\mu)$ now implies that $T_E(h^*)=g$ in $L^1(\mu)$, so $T_E$ is indeed a linear and continuous map and we may conclude that $\|T_E\|_{\text{op}}<\i$. By the above arguments we therefore have that $I_E$ is linear and continuous implying that $I_E\in \cH^{**}$.
\begin{definition}[Dunford integral]
	Assume that $f:\cX \to \cH$ is scalarly $\mu$-integrable mapping. For every set $E\in \F$ we denote the Dunford integral of $f$ over $E$ with respect to $\mu$ by $(D)\int_E f \, d\mu $ and define it as the unique  functional in $\cH^{**}$ which maps
	\begin{align*}
	\cH^* \ni h^* \mapsto \int_E h^* \circ f \, d\mu.
	\end{align*}
\end{definition}
Now we are close to defining the Pettis integral of $f$ with respect to $\mu$. The above Dunford integral is an element in $\cH^{**}$ but we want an element in $\cH$. Before proceeding we introduce the natural isometric embedding $ev:\cH \to \cH^{**}$, which we are going to use in determining which element of $\cH$ we define as the Pettis integral. In the general theory  for Banach spaces $f$ is Pettis integrable if the image of $ev$ contains the Dunford integral of $f$, and one defines the Pettis integral as the element which is mapped to the Dunford integral by $ev$.  But as we shall see this is always the case whenever the value space of $f$ is a Hilbert space.\\ \\
Define the linear evaluation map $ev:\cH \to \cH^{**}$ by letting 
$ev(h)\in \cH^{**}$ such that
\begin{align*}
ev(h)(h^*)=h^*(h),
\end{align*}
for all $h^*\in \cH^*$. That $ev(h)\in \cH^{**}$ for all $h\in \cH$ follows by noting that $ev(h):\cH^* \to \R$ is clearly linear and  continuity in zero (and by linearity: everywhere) follows from the inequality
$|ev(h)(h^*)|=|h^*(h)|\leq \|h^*\|_{\text{op}} \| h\|_\cH \stackrel{h^* \to 0}{\longrightarrow} 0
$.
\begin{lemma}
	The evaluation mapping $ev:\cH \to \cH^{**}$  is an isometric isomorphism. That is $ev$ is a bijective mapping satisfying 
	\begin{align*}
	\| ev(h_1)-ev(h_2)\|_{\text{op}}  = \|h_1-h_2\|_\cH,
	\end{align*}
	for all $h_1,h_2\in \cH$.
\end{lemma} 
\begin{p}
 First we show that $ev$ is isometric embedding (satisfy the above equation). The last inequality before this lemma yields that
 \begin{align*}
 \|ev(h)\|_{\text{op}}\leq \sup_{h^*\in \cH^*,\|h^*\|_{\text{op}}\leq 1} \|h^*\|_{\text{op}} \| h\|_\cH  =  \|h\|_{\cH}.
 \end{align*}
 Conversely by the corollary to theorem 3.3 (Hahn-Banach theorem) \cite{rudin1991functional} there exists an $a^*\in \cH^*$ with  $|a^*(h)|=\|h\|_\cH$ and $\|a^*\|_{\mathrm{op}}\leq 1$. Thus we get that
 \begin{align*}
 \|ev(h)\|_{\text{op}} = \sup_{h^*\in \cH^*,\|h^*\|_{\text{op}}\leq 1} |h^*(h)| \geq |a^*(h)| = \|h\|_\cH.
 \end{align*}
 Hence $\|ev(h)\|_{\text{op}}=\|h\|_\cH$ for all $h\in \cH$, and by linearity of $ev$ we get that
 \begin{align*}
 \| ev(h_1)-ev(h_2)\|_{\text{op}} = \|ev(h_1-h_2)\|_{\text{op}} = \|h_1-h_2\|_\cH,
 \end{align*}
 for all $h_1,h_2\in \cH$, proving that $ev$  is an isometric (hence also injective and continuous) embedding  of $\cH$ into $\cH^{**}$. 
 It remains to be shown that $ev$ is a surjective mapping. In our scenario, with the value space being a Hilbert space, it is well known that the natural embedding $ev$ into the double continuous dual space is surjective - but in general  - spaces possessing this property are called reflexive. \\
 
 In an effort to keep the thesis self-contained we sketch a proof for the fact that every Hilbert space is reflexive. 
 Note that by Reisz representation theorem (see \cite{schilling}) we get that  $\psi_\cH:\cH \to \cH^*$ given by $\psi_\cH(h_1)(h_2)= \la h_2 , h_1 \ra_\cH $ is a bijection and the inverse  especially fulfils that  
 \begin{align}
 	 h^*(h)=\la h, \psi_\cH^{-1}(h^*) \ra_\cH   \label{Ap_Reisz_prop},
 \end{align}
 for any $h\in \cH$ and $h^*\in \cH^*$. It can easily be verified that $\la\cdot ,\cdot \ra_{\cH^*}:\cH^* \times \cH^* \to \R$ defined by
 \begin{align*}
 \la h^*_1 ,h^*_2\ra_{\cH^*} = \la \psi_\cH^{-1}(h_1^*),\psi_\cH^{-1}(h_2^*) \ra_{\cH},
 \end{align*}
 is indeed a inner product which agrees with the operator norm on $\cH^*$, making it a Hilbert space. Yet again  Reisz representation theorem yields that $\psi_{\cH^*}:\cH^* \to \cH^{**}$ given by $\psi_{\cH^*}(h^*_1)(h^*_2)=\la h^*_2, h^*_1 \ra_{\cH^*}$ is a bijection. The important thing to note is that $\psi_{\cH^*}\circ \psi_{\cH}:\cH \to \cH^{**}$ is surjective and
 \begin{align*}
 \psi_{\cH^*}\circ \psi_{\cH}(h)(h^*) &=\la h^* ,\psi_{\cH}(h)\ra_{\cH^*} \\
 &= \la \psi_\cH^{-1}(h^*),\psi_\cH^{-1}(\psi_{\cH}(h)) \ra_{\cH} \\
 &=\la \psi_\cH^{-1}(h^*),h \ra_{\cH} \\
 &= h^*(h) \\
 &= ev(h)(h^*),
 \end{align*}
 for all $h\in \cH$ and $h^* \in \cH^*$, proving that $ev= \psi_{\cH^*}\circ \psi_{\cH}$ is surjective.
\end{p} 
Hence we know that there always exists a unique $h\in \cH$ such that $ev(h)=(D)\int_E f \, d\mu$ which leads us to the following definition of the Pettis integral.
\begin{definition}[Pettis integral]
	If $f:\cX \to \cH$ is scalarly $\mu$-integrable we say that $f$ is Pettis integrable with respect to $\mu$ and for any $E\in \F$ there exists a unique element $h_E\in \cH$ such that
	$ev(h_E)=(D)\int_E f \, d\mu$. We denote this element $\int_E f \, d\mu$ and call it the Pettis integral of $f$ over $E$ with respect to $\mu$.
\end{definition}
For simplicity lets make en equivalent definition of the Pettis integral which circumvents using the Dunford integral but rather its defining property. 
\begin{theorem} \label{Theorem_definition_Pettis_Integral}
	If $f:\cX \to \cH$ is scalarly $\mu$-integrable and $E\in \F$ then the Pettis integral $\int_E f\, d\mu\in \cH$ is the unique element satisfying 
	\begin{align*} \label{defi_property_of_Pettis_integral}
	h^*\lp  \int_E f \, d\mu \rp = \int _E h^* \circ f \, d\mu,
	\end{align*}
	for all $h\in \cH^*$.
\end{theorem}
\begin{p}
	Recall that the Dunford integral $(D)\int_E f \, d\mu$ is the unique element in $\cH^{**}$ satisfying 
	\begin{align*}
	  \lp  (D)\int_E f \, d\mu \rp  (h^*) = \int_E h^* \circ f \, d\mu ,
	\end{align*}
	for all $h^* \in \cH^*$. But $\int_E f \, d\mu $ is the unique element in $\cH$  satisfying $ev \lp \int_E f \, d\mu\rp =(D)\int_E f \, d\mu$. Combining these two characterisations we get that the Pettis integral $\int_E f\, d\mu$ is the unique element in $\cH$ satisfying
	\begin{align*}
	h^*\lp  \int_E f \, d\mu \rp=ev \lp \int_E f\, d\mu  \rp (h^*)=\int_E h^* \circ f \, d\mu,
	\end{align*}
	for all $h^*\in \cH^*$.
\end{p}
\begin{corollary} \label{appendix_corollary_linearity_of_Pettis_integral}
	The Pettis integral exhibits the same linearity property as the Lebesgue integral. That is, if $f:\cX\to \cH$ and $g:\cX\to \cH$ are both Pettis integrable over $E\in \F$ with respect to $\mu$ then for any $a,b$ in the scalar field of $\cH$ it holds that
	\begin{align*}
  \int a f +bg \, d\mu = a \int f \, d\mu + b\int g \, d\mu.
	\end{align*}
\end{corollary}
\begin{p}
	This is a consequence of the above theorem. Simply note that by the unique defining property of the Pettis integral
	\begin{align*}
	\int a f +bg \, d\mu &= a \int f \, d\mu + b\int g \, d\mu \\&\iff\\ h^* \lp a \int f \, d\mu + b\int g \, d\mu\rp &= \int h^* (af(x)+bg(x)) \, d\mu(x), \quad \forall h^*\in\cH^*.
	\end{align*}
	The latter condition is easily verified by the linearity of $h^*\in \cH^*$. That is,
	\begin{align*}
	h^* \lp a \int f \, d\mu + b\int g \, d\mu\rp &= a \int h^*(f(x)) d\mu(x) +b\int h^*(g(x))\, d\mu(x) \\
	&= \int h^* (af(x)+bg(x)) \, d\mu(x).
	\end{align*}
	
\end{p}
\begin{lemma} \label{lemma_f_is_pettis_integrabel_if_measurable_and_norm_of_f_is_integrable}
	Any $\F/\cB(\cH)$-measurable mapping $f:\cX\to \cH$, if Pettis integrable with respect to $\mu$, if
	\begin{align*}
	\int \| f(x) \|_\cH \, d\mu(x) <\i.
	\end{align*}
\end{lemma}
\begin{p}
		 By \cref{Theorem_definition_Pettis_Integral} is suffices to show that $f$ is $\mu$-scalarly integrable. That is, it suffices to show that $h^*\circ f \in L^1(\mu)$ for all $h^*\in \cH$. \textit{Measurability}: $f$ is an $\F/\cB(\cH)$-measurable mapping and any $h^*\in \cH^*$ is continuous and hence $\cB(\cH)/\cB(\R)$-measurable. We conclude that the composition $h^*\circ \phi$ indeed is $\F/\cB(\R)$-measurable, for any $h^*\in \cH^*$. \textit{Integrability:} Note that for any $h^*\in \cH^*$ we have that $\|h^*\|_{\text{op}}<\i$, and since $|h^*\circ \phi(x)| \leq \|h^*\|_{\text{op}}|\phi(x)|$ for all $x\in \cX$, we get that
		\begin{align*}
		\frac{1}{\|h^*\|_{op}}\int |h^* \circ f(x)| d \mu (x) &\leq \int  \|f(x) \|_\cH \, d\mu (x)  <\i,
		\end{align*}
		if $\int \|f(x)\|_\cH \, d\mu(x)<\i$, proving that $f$ is Pettis integrable with respect to $\mu$.
\end{p}
\begin{lemma} \label{appendix_lemma_linear_map_inside_pettis_integral_also_pettis_integrable}
	If $f:\cX\to\cH_1$ is Pettis integrable with respect to $\mu$ and $L:\cH_1\to\cH_2$ is a continuous linear map, then $L\circ f$ is Pettis integrable with respect to $\mu$ and
	\begin{align*}
	L\lp \int f \, d\mu \rp = \int L\circ f \, d\mu.
	\end{align*}
\end{lemma}
\begin{p}
	First we show that $L\circ f$ is indeed Pettis integrable. Note that if $f$  is scalarly $\mu$-integrable then one simply hote that for any $h^*\in\cH_2^*$ then $h^*\circ L\in \cH_1^*$, hence
	\begin{align*}
	h^*\circ L \circ f \in\cL^{1}(\mu),
	\end{align*}
	proving that also $L\circ f$ is Pettis integrable with respect to $\mu$. Now for the identity note for any $h^*\in\cH_2^*$ we have that $h^*\circ L\in \cH_1^*$, hence
	\begin{align*}
	h^* \lp L \lp \int f \, d\mu \rp \rp = h^* \circ L \lp \int f \, d\mu \rp  = \int h^* \circ L \circ f \, d\mu,
	\end{align*}
	and
	\begin{align*}
	h^* \lp \int L\circ f \, d\mu \rp = \int h^* \circ L \circ f \, d\mu,
	\end{align*}
	for any $h^*\in\cH^*_2$, proving the identity by the unique defining property of the Pettis integral.
	\end{p}
We may also extend the notions of expectations and variance to Hilbert space valued random elements. Let $X:\Omega\to \cH$ be a random Borel element in a Hilbert space $\cH$, defined on some probability space $(\Omega,\F,P)$.
\begin{definition} \label{defi_Hilbert_space_expectation_and_variance}
	If $X$ is Pettis integrable with respect to $P$ (e.g. if $E\|X\|_\cH =  \int \|x\|_\cH \, dX(P)(x) < \i$, by above lemma), then we define the expectation of $X$ as the Pettis integral of $X$ with respect to $P$, that is 
	\begin{align*}
	E(X) = \int X \, dP \in \cH,
	\end{align*}
	and if $E\|X\|_\cH^2 <\i$ then we define the variance of $X$ as the positive real number
	\begin{align*}
	\text{Var}(X) = E\|X-E(X)\|_\cH^2.
	\end{align*}
\end{definition} 
This allows for a different notation of Pettis integrals with respect to probability measures. For any probability space $(\cX,\F,\mu)$, we may construct a random element $X:(\Omega,\bK,P)\to (\cX,\F)$ such that $X(P)=\mu$. Then for any Hilbert space valued mapping $f:\cX\to \cH$, that is Pettis integrable with respect to $\mu$, we have that
\begin{align*}
\int f d\mu =E(f(X)).
\end{align*}
This is easily seen, as $E(f(X))$ fulfils the unique defining property of the Pettis integral. That is 
\begin{align*}
h^*(E(f(X))) = \int h^* \circ f(X(\omega)) \, dP(\omega) = \int h^* \circ f \, dX(P) = \int h^* \circ f \, d\mu,
\end{align*}
for any $h^* \in \cH^*$. 
\newpage
\subsection{Complexification and Realification of a Hilbert space}\label{Appendix_Complexification_and_Realification}
\subsubsection{Complexification of a real Hilbert space} 
Assume that we have a $\R$-Hilbert space $\cH$. We will now associate a complex Hilbert space $\cH^{\bC}$ to $\cH$, which will be useful for our analysis of metric spaces of negative and strong negative type. We follow the general procedure for complexification of real vector spaces presented in \cite{AdvancedLinearAlgebraroman2005_complexification}, and modify it to Hilbert spaces.
\begin{definition} \label{appendix_definition_complexification}
	We define the complexification $\cH^{\bC}$ of $\cH$ as the complex vector space of ordered pairs in $\cH^{\bC}=\cH\times \cH$ with coordinatewise addition
	\begin{align*}
	(x,y)+(x',y')=(x+x',y+y'),
	\end{align*}
	and scalar multiplication 
	\begin{align*}
	(a+ib)(x,y)=(ax-by,ay+bx),
	\end{align*}
	for any $x,x',y,y'\in\cH$ and $a+ib\in\bC$.
	Note that we can use notation alike to the complex numbers by denoting $(x,y)\in \cH^{\bC}$ by $x+iy$, and under this notation we may write $\cH^{\bC}=\{x+iy:x,y\in\cH\}$. With this notation addition and scalar multiplication resembles those of the complex numbers.
\end{definition}
\noindent We can also extend the original inner product $\la \cdot,\cdot  \ra_{\cH}$ on $\cH$ to an inner product $\la \cdot , \cdot\ra_{\cH^{\bC}}$ on $\cH^{\bC}$, in the following way
\begin{theorem}
	The mapping $\la \cdot , \cdot\ra_{\cH^{\bC}}:\cH^{\bC}\times \cH^{\bC}\to \bC$ given by
	\begin{align*}
	\la x+iy,x'+iy'\ra_{\cH^{\bC}} = \la x,x' \ra_{\cH} + \la y,y' \ra_{\cH} + i \lp \la y,x' \ra_{\cH}-\la x,y' \ra_{\cH}\rp,
	\end{align*}
	is an inner product on $\cH^{\bC}$
\end{theorem}
\begin{p}
	Since $\la\cdot , \cdot \ra_{\cH}$ is a $\R$-inner product we obviously have that
	\begin{align*}
	\la x+iy,x'+iy'\ra_{\cH^{\bC}} &= \la x',x \ra_{\cH} + \la y',y \ra_{\cH} - i \lp \la x,y' \ra_{\cH}-\la y,x' \ra_{\cH}\rp \\
	&= \la x',x \ra_{\cH} + \la y',y \ra_{\cH} - i \lp \la y',x \ra_{\cH}-\la x',y \ra_{\cH}\rp \\
	&=\overline{\la x'+iy',x+iy\ra_{\cH^{\bC}}},
	\end{align*}
	i.e. we have conjugate symmetry. Further more note that
	\begin{align*}
	\la (a+ib)(x+iy),x'+iy'\ra_{\cH^{\bC}} =& \la ax-by+i(ay+bx),x'+iy'\ra_{\cH^{\bC}} \\
	=&\la ax-by,x' \ra_{\cH} + \la ay+bx,y' \ra_{\cH} \\
	&+ i \lp \la ay+bx,x' \ra_{\cH}-\la ax-by,y' \ra_{\cH}\rp \\
	=& \la ax,x' \ra_{\cH}-\la by,x' \ra_{\cH} + \la ay,y' \ra_{\cH}+\la bx,y' \ra_{\cH} \\
	&+ i \lp \la ay,x' \ra_{\cH}+ \la bx,x' \ra_{\cH} -\la ax,y' \ra_{\cH}+\la by,y' \ra_{\cH}\rp \\
		=& a \lp  \la x,x' \ra_{\cH}+\la y,y' \ra_{\cH} \rp -b \lp \la y,x' \ra_{\cH} +\la x,y' \ra_{\cH} \rp \\
	&+ i \lp a[\la y,x' \ra_{\cH}-\la x,y' \ra_{\cH} ]+ b[\la x,x' \ra_{\cH} +\la y,y' \ra_{\cH} ]\rp \\
	=& (a+ib) \lp \la x,x' \ra_{\cH} + \la y,y' \ra_{\cH} + i \lp \la y,x' \ra_{\cH}-\la x,y' \ra_{\cH}\rp \rp,	\end{align*}
	and
	\begin{align*}
	\la (x+iy)+(x'+iy'),x''+iy''\ra_{\cH^{\bC}} =& \la x,x'' \ra_{\cH} \la x',x'' \ra_{\cH}+ \la y,y'' \ra_{\cH}+\la y',y'' \ra_{\cH} \\
	&+ i \lp \la y,x'' \ra_{\cH}+ \la y',x'' \ra_{\cH}-\la x,y'' \ra_{\cH}-\la x',y'' \ra_{\cH}\rp \\
	&= \la x+iy,x''+iy''\ra_{\cH^{\bC}}+\la x'+iy',x''+iy''\ra_{\cH^{\bC}},
	\end{align*}
	proving linearity in the first argument. Lastly we need to show finite-definiteness of $\la \cdot , \cdot \ra_{\cH^{\bC}}$, and this follows from the symmetry and finite-definiteness of $\la \cdot , \cdot \ra_{\cH}$. Note that
	\begin{align*}
	\la x+iy , x+iy \ra_{\cH^{\bC}}&= \la x,x \ra_{\cH} + \la y,y \ra_{\cH} + i \lp \la y,x \ra_{\cH}-\la x,y \ra_{\cH}\rp \\
	&= \la x,x \ra_{\cH} + \la y,y \ra_{\cH} + i \lp \la y,x \ra_{\cH}-\la y,x \ra_{\cH}\rp \\
		&= \la x,x \ra_{\cH} + \la y,y \ra_{\cH} \\
		&\geq 0,
	\end{align*}
	proving that $\la \cdot , \cdot \ra_{\cH^{\bC}}$ is an inner product on $\cH^{\bC}$.
\end{p}
\begin{theorem}
	The inner product space $(\cH^{\bC},\la \cdot , \cdot \ra_{\cH^{\bC}})$ is complete and therefore a Hilbert space.
\end{theorem}
\begin{p}
	First note that $\cH$ is complete, meaning that every Cauchy sequence converges. Now consider an arbitrary Cauchy sequence $((x_n+iy_n))_{n\in\N}$ in $\cH^{\bC}$ and note that
	\begin{align*}
	\|x_n+iy_n\|_{\cH^{\bC}}^2 = \|x_n\|_\cH^2+\|y_n\|_\cH^2,
	\end{align*}
	by the equality in the above proof. Using this we see that
	\begin{align*}
	\|(x_n+iy_n)-(x_m+iy_m)\|_{\cH^{\bC}}^2 = \|x_n-x_m\|_\cH^2+\|y_n-y_m\|_\cH^2,
	\end{align*}
	which tends to zero as $n,m\to\infty$ since $((x_n+iy_n))_{n\in\N}$ is Cauchy, and we realize that this happens if and only if both of the addends on the right hand side converge to zero. Thus $(x_n)_{n\in\N}$ and $(y_n)_{n\in\N}$ are both Cauchy sequences in $\cH$. As a consequence they have limits $x$ and $y$ in $\cH$, and we note that
		\begin{align*}
		\|(x_n+iy_n)-(x+iy)\|_{\cH^{\bC}}^2 = \|x_n-x\|_\cH^2+\|y_n-y\|_\cH^2 \to_n 0,
		\end{align*}
		proving that $x_n+iy_n$ converges to $x+iy$ as $n$ tends to infinity. We conclude that every Cauchy sequence in $\cH^{\bC}$ converges in $\cH^{\bC}$, proving that $\cH^{\bC}$ is complete. Hence $(\cH^{\bC},\la \cdot , \cdot \ra_{\cH^{\bC}})$ is a complete inner product space, i.e. a Hilbert space.
\end{p}
\noindent Lastly we introduce a mapping which will be used in the main sections.
\begin{definition} \label{appendix_definition_complexification_map}
	We define $\mathrm{cpx}:\cH\to \cH^{\bC}$ by
	\begin{align*}
	\mathrm{cpx}(x)= x+i0,
	\end{align*}
	and call it the complexification map.
\end{definition}
\noindent We easily see that the complexification map is injective and satisfies the following properties
\begin{align*}
&\mathrm{cpx}(0)=(0+i0)\equiv 0 \\
&\mathrm{cpx}(x+x') = x+x'+i0 = (x+i0)+(x'+i0)= \mathrm{cpx}(x)+\mathrm{cpx}(x') \\
&\mathrm{cpx}(ax)=ax+i0 = a(x+i0) =a\mathrm{cpx}(x), \quad a\in \R,
\end{align*}
so cpx is additive and we may "pull out" real scalars, so cxp is almost a linear map if we disregard the fact that linearity of maps are only defined for maps between vector spaces with the same scalar field.
\begin{theorem} \label{appendix_theorem_complexification_map_is_an_isometric_embedding_and_inner_prod_equality}
	The complexification map $\mathrm{cpx}:(\cH,d_\cH)\to (\cH^{\bC},d_{\cH^{\bC}})$ is an isometric embedding and
	\begin{align*}
	\la \mathrm{cpx}(x),\mathrm{cpx}(x') \ra_{\cH^{\bC}} = \la x,x' \ra_{\cH},
	\end{align*}
	for every $x,x'\in \cH$
\end{theorem}
\begin{p}
	Recall that $\|x+iy\|_{\cH^{\bC}}^2 = \|x\|_\cH^2+\|y\|_\cH^2$ and note that
	\begin{align*}
	d_\cH(x,x')&= \sqrt{\|x-x'\|_{\cH}^2 +\|0-0\|_{\cH}^2} \\
	&= \sqrt{ \| (x-x')+i(0-0) \|_{\cH^{\bC}}^2} \\
	&= \| \mathrm{cpx}(x+(-x'))\|_{\cH^{\bC}} \\
	&= \|\mathrm{cpx}(x)+\mathrm{cpx}(-x')\|_{\cH^{\bC}} \\
		&= \|\mathrm{cpx}(x)-\mathrm{cpx}(x')\|_{\cH^{\bC}} \\
	&= d_{\cH^{\bC}}(\mathrm{cpx}(x),\mathrm{cpx}(x')),
	\end{align*}
	proving that $\mathrm{cpx}$ is an isometric embedding from the real Hilbert space $\cH$ into the corresponding  complexification Hilbert space $\cH^{\bC}$. It further more holds that
	\begin{align*} 
	\la \mathrm{cpx}(x),\mathrm{cpx}(x') \ra_{\cH^{\bC}} &=\la x+i0,x'+i0 \ra_{\cH^{\bC}} \\
	&= \la x,x' \ra_{\cH} + \la 0,0 \ra_{\cH} + i \lp \la 0,x' \ra_{\cH}-\la x,0 \ra_{\cH}\rp \\
	&= \la x,x' \ra_{\cH},
	\end{align*}
	for every $x,x'\in \cH$.
\end{p}
\begin{theorem} \label{appendix_theorem_complexification_is_separable}
	If $\cH$ is a separable Hilbert space, then $\cH^\bC$ is a separable Hilbert space.
\end{theorem}
\begin{p}
	Since $\cH$ is separable, we know that there exists a countable dense subset $D\subset \cH$. Note that $D\times D\subset \cH^{\bC}=\cH\times \cH$ is countable and for any $x+iy\in \cH^\bC$ there exists a sequence $(x_n + iy_n)\subset D\times D$ such that $\|x_n- x\|_\cH\to_n 0$ and $\|y_n-y\|_\cH\to_n 0$, since $D$ is dense in $\cH$. Thus
	\begin{align*}
	\|x_n+iy_n - (x+iy)\|_{\cH^\bC} =\|(x_n-x)+i(y_n-y) \|_{\cH^\bC} = \|x_n-x\|_\cH + \|y_n-y\|_\cH \to_n 0,
	\end{align*}
	proving that $\cH^\bC$ is indeed separable.
	
\end{p}
\subsubsection{Realification of a complex Hilbert space}
Let $(\cH,\la \cdot , \cdot \ra_\cH)$ be a $\bC$-Hilbert space. The realification $\cH^\R$ of $\cH$ is given by $\cH$ where we simply ignore the possibility of scalar multiplying elements with complex scalars and only keep scalar multiplication over $\R$. Let $\mathrm{re}:\cH\to \cH^\R$ denote the identification map between these spaces.
\begin{theorem}
	The mapping $\la \cdot , \cdot \ra_{\cH^\R}:\cH^\R\times \cH^\R \to \R$ given by
	\begin{align*}
	\la x,y \ra_{\cH^\R} = \mathfrak{R}\la x,y\ra_{\cH} = \frac{\la x,y \ra_{\cH} + \la y,x \ra_{\cH}}{2},
	\end{align*}
	is an inner product of $\cH^\R$.
\end{theorem}
\begin{p}
	We obviously have that $\la x,y \ra_{\cH^\R}= \la y,x \ra_{\cH^\R}$ and 
	\begin{align*}
	\la a x,y \ra_{\cH^\R} = \mathfrak{R}a\la x,y\ra_{\cH} =a \mathfrak{R}\la x,y\ra_{\cH}= a\la x,y \ra_{\cH^\R},
	\end{align*}
	for any $a\in\R$, also since $\Re$ preserves addition we also have that 
	\begin{align*}
	\la x+x',y \ra_{\cH^\R} = \Re(\la x,y \ra_{\cH}+\la x',y \ra_{\cH}) = \Re \la x,y \ra_{\cH} + \Re\la x',y \ra_{\cH} =\la x,y \ra_{\cH^\R}+\la x',y \ra_{\cH^\R}.
	\end{align*}
	Lastly since $\la \cdot , \cdot \ra_{\cH}$ has the property of positive-definiteness, so has $\la \cdot,\cdot \ra_{\cH^\R}$. That is
	\begin{align*}
	\la x,x \ra_{\cH^\R} = \frac{\la x,x \ra_{\cH} + \la x,x \ra_{\cH}}{2} \geq 0,
	\end{align*}
	and
	\begin{align*}
	\la x,x \ra_{\cH^\R}=0 \iff \frac{\la x,x \ra_{\cH} + \la x,x \ra_{\cH}}{2}=0 \iff \la x,x \ra_{\cH}=0 \iff x=0.
	\end{align*}
\end{p}
\begin{theorem} \label{appendix_theorem_realification_is_a_separable_hilbert_space}
	$(\cH^\R,\la\cdot, \cdot \ra_{\cH^\R})$ is a $\R$-Hilbert space and if $\cH$ is separable, then $\cH^\R$ is separable.
\end{theorem}
\begin{p}
The inner product space $(\cH^\R,\la\cdot, \cdot \ra_{\cH^\R})$ is actually a $\R$-Hilbert space. To see this note that $(x_n)$  is Cauchy in $\cH$ if and only if $(x_n)$ Cauchy in $\cH^\R$, since $\|x\|_\cH=\|x\|_{\cH^\R}$ for all $x$. Furthermore if $(x_n)$ converges to $x$ in $\cH$ we also have that $(x_n)$ converges to $x$ in $\cH^\R$, proving that $(\cH^\R,\la\cdot, \cdot \ra_{\cH^\R})$ is complete. This furthermore implies that if $\cH$ is separable, then $\cH^\R$ is also separable.
\end{p}
We may note that the identification map $\mathrm{re}:\cH\to \cH^\R$ is an additive isometric since $\|x\|_\cH = \|\mathrm{re}(x)\|_{\cH^\R}$ and $\mathrm{re}(x+x')= \mathrm{re}(x)+\mathrm{re}(x')\in \cH^\R$. Lastly it is bijective, rendering $\mathrm{re}:\cH\to \cH^\R$ an additive isometric isomorphism, hence also a homeomorphism
\newpage
\subsection{Tensor product of Hilbert spaces}  \label{Appendix_Tensor_product_of_Hilbert_spaces}
In this section we are going to construct the tensor product of two Hilbert spaces. This tensor product turns out to be a new Hilbert space, which we can use in the theory of metric spaces of negative and strong negative type. The following construction approach is found in \cite{MethodsOfModernMathematicalPhysics_reed1972functional}, but we try to prove things more carefully. 

We can only construct the tensor product of two Hilbert spaces with the same scalar field, so consider two $\bK$-Hilbert spaces $\cH_1$ and $\cH_2$, both with the same scalar field $\bK=\bC$ or $\R$. Since we in the theory of metric spaces of negative type, may have two Hilbert spaces with different scalar fields on our hands, the approach is to complexify [realify] (see \cref{Appendix_Complexification_and_Realification}) the real [complex] Hilbert spaces and carry on with the following construction.\\ \\ For any $\phi\in \cH_1,\psi \in\cH_2$ we define the  map $\phi\otimes \psi:\cH_1\times \cH_2 \to \bK$, as the tensor product of $\phi$ and $\psi$, given by
\begin{align*}
\phi\otimes \psi (x,y) = \la \phi,x  \ra_{\cH_1} \la \psi,y\ra_{\cH_2}.
\end{align*}
Its easily realized that $h_1\otimes h_2$ is a conjugate bilinear map(anti-linear in both arg.), with the following properties
\begin{align*}
&(\phi+\phi')\otimes \psi= \phi\otimes \psi + \phi'\otimes \psi,\\
&\phi\otimes (\psi+\psi') = \phi\otimes \psi + \phi\otimes \psi', \\
&(a\phi) \otimes \psi=\phi \otimes (a\psi)=a(\phi\otimes \psi) ,\quad a\in \bK, \\
&0\otimes \psi = \phi \otimes 0 \equiv 0 ,
\end{align*}
any we may note that $\phi \otimes \psi$ is the zero mapping if and only if $\phi$ or $\psi$ are zero elements of $\cH_1$ and $\cH_2$ respectively. \\ \\
 Now denote $\mathcal{E}$, the collection of conjugate bilinear mappings from $\cH_1\times \cH_2$ to $\bK$, that can be written as a finite linear combinations of tensor products. That is
\begin{align*}
\mathcal{E}=\left\{\sum_{i=1}^n a_i \phi_{i}\otimes \psi_{i}: n\in\N,\, a\in \mathbb{K}^n,\,  \phi\in\cH_1^n,\, \psi\in \mathcal{H}_2^n\right\},
\end{align*}
where scalar multiplication and summation of maps is defined in the regular fashion $(af+g)(x,y)=af(x,y)+g(x,y)$ for any maps $f,g$ and scalar $a$. \\ \\
For any two tensor products $\phi_1\otimes \psi_1$ and $\phi_2\otimes \psi_2$ we define
\begin{align*}
\la \phi_1 \otimes \psi_1, \phi_2 \otimes \psi_2 \ra = \la \phi_1,\phi_2 \ra_{\cH_1} \la \psi_1,\psi_2 \ra_{\cH_2},
\end{align*}
and extend by linearity in first argument and conjugate linearity in the second argument. That is by letting
\begin{align*}
\left\la \sum_{i=1}^n a_i \phi_{i}\otimes \psi_{i}, \sum_{j=1}^m b_j \gamma_j \otimes \beta_j \right\ra= \sum_{i=1}^n\sum_{j=1}^m a_i \overline{b_j} \la \phi_i \otimes \psi_i, \gamma_j \otimes \beta_j \ra,
\end{align*}
for arbitrary linear combinations of tensor products.
That we need conjugate linearity in the second argument follows from the fact that $\la \phi_1 \otimes \psi_1, b(\phi_2 \otimes \psi_2) \ra = \la \phi_1,b\phi_2 \ra_{\cH_1} \la \psi_1,\psi_2 \ra_{\cH_2} = \overline{b}\la \phi_1 \otimes \psi_1, \phi_2 \otimes \psi_2 \ra $, 
using that  $\la \cdot , \cdot \ra_{\cH_1}$ and $\la\cdot , \cdot \ra_{\cH_2}$ are $\bK$-inner products. Note that we have not said anything about $\la \cdot, \cdot \ra$ being a map on $\cE\times \cE$, this requires that $\la \cdot, \cdot \ra$ is invariant to how one represents the bilinear mappings of $\cE$, which is among other things what we show in the following lemma.
  
\begin{lemma}
	$\cE$ is a $\bK$-vector space and $\la \cdot , \cdot \ra:\cE\times \cE\to \bK$ is a well-defined mapping making $(\cE,\la \cdot , \cdot \ra)$ a $\bK$-inner product space.
\end{lemma} 
\begin{p}
	With the above mentioned definition of scalar multiplication and summation of maps, we see that $\cE$ is closed under summation and scalar multiplication. For any $\alpha\in \bK$ and $x,y\in \cE$ then for some $n_1,n_2\in \N$ and $a_1\in\bK^{n_1}$, $a_2\in\bK^{n_2}$, $\phi_1\in\cH_1^{n_1}$, $\phi_2\in \cH_1^{n_2}$, $\psi_1\in\cH_2^{n_1}$, $\psi_2\in\cH_2^{n_2}$ we have that
	\begin{align*}
	x= \sum_{i=1}^{n_1} a_{1_i} \phi_{1_i}\otimes \psi_{1_i}, \quad \textit{ and } \quad y= \sum_{i=1}^{n_2} a_{2_i} \phi_{2_i}\otimes \psi_{2_i}.
	\end{align*}
	Hence we obviously see that
	\begin{align*}
	\alpha x+y = \sum_{i=1}^{n_1+n_2} a_{_i} \phi_{i}\otimes \psi_{i} \in \cE,
	\end{align*}
	where $a=(\alpha a_1,a_2)\in \bK^{n_1+n_2}$,  $\psi=(\phi_1,\phi_2)\in \cH_1^{n_1+n_2}$ and $\psi=(\psi_1,\psi_2)\in \cH_2^{n_1+n_2}$. Furthermore the zero-element of $\cE$ is the (obviously bilinear) zero mapping $\cH_1\times \cH_2 \ni (x,y) \mapsto 0$.  All other axioms for vector spaces are easily realized, so we conclude that $\cE$ is a $\bK$-vector space. \\ \\
	By the extension of $\la \cdot , \cdot \ra$ to finite linear combinations of tensor products above, we have that
	\begin{align*}
	\la ax+bx',cy+dy'\ra =& \left\la \sum_{i=1}^{n_1+n_2} \theta_i \phi_i \otimes \psi_i , \sum_{i=1}^{m_1+m_2} \lambda_i \gamma_i \otimes \beta_j \right\ra  \\ 
	=& \sum_{i=1}^{n_1+n_2} \sum_{j=1}^{m_1+m_2}\theta_i \overline{\lambda_j} \la \phi_i \otimes \psi_i , \gamma_i\otimes \beta_j \ra \\ 
	=& a\overline{c} \sum_{i=1}^{n_1}\sum_{j=1}^{m_1}\alpha_{1_i}\overline{\delta_{1_j}}\la \phi_{1_i}\otimes \psi_{1_i},\gamma_{1_j}\otimes \beta_{1_j} \ra \\
	&+ a\overline{d} \sum_{i=1}^{n_1}\sum_{j=1}^{m_2}\alpha_{1_i}\overline{\delta_{2_j}}\la \phi_{1_i}\otimes \psi_{1_i},\gamma_{2_j}\otimes \beta_{2_j} \ra \\
	&+b\overline{c}\sum_{i=1}^{n_2}\sum_{j=1}^{m_1}\alpha_{2_i}\overline{\delta_{1_j}}\la \phi_{2_i}\otimes \psi_{2_i},\gamma_{1_j}\otimes \beta_{1_j} \ra\\\
	&+b\overline{d}\sum_{i=1}^{n_2}\sum_{j=1}^{m_2}\alpha_{2_i}\overline{\delta_{2_j}}\la \phi_{2_i}\otimes \psi_{2_i},\gamma_{2_j}\otimes \beta_{2_j} \ra \\
	=& a\overline{c}\la x,y \ra +a\overline{d}\la x,y'\ra+b\overline{c}\la x',y \ra+b\overline{d}\la x',y' \ra,
	\end{align*}
	proving that $\la \cdot , \cdot\ra$ is sesquilinear (sesqui meaning one-and-a-half, the same as conjugate linear), where $a,b,c,d\in \bK$ and
	\begin{align*}
	x &= \sum_{i=1}^{n_1} \alpha_{1_i} \psi_{1_i} \otimes \psi_{1_i}, \quad x'= \sum_{i=1}^{n_2} \alpha_{2_i} \psi_{2_i} \otimes \psi_{2_i}, \\
	y&= \sum_{j=1}^{m_1} \delta_{1_j} \gamma_{1_j} \otimes \beta_{1_j}, \quad \, \, y'= \sum_{j=1}^{m_2} \delta_{2_j} \gamma_{2_j} \otimes \beta_{2_j},
	\end{align*}
	with $\alpha=(\alpha_1,\alpha_2)\in\bK^{n_1+n_2}$, $\delta=(\delta_1,\delta_2)\in \bK^{m_1+m_2}$, $\theta=(a\alpha_1,b\alpha_2)\in \bK^{n_1+n_2}$, $\lambda=(c\delta_1,d\delta_2)\in \bK^{m_1+m_2}$, $\phi=(\phi_1,\phi_2)\in \cH_1^{n_1+n_2}$, $\psi=(\psi_1,\psi_2)\in \cH_2^{n_1+n_2}$, $\gamma=(\gamma_1,\gamma_2)\in \cH_1^{m_1+m_2}$, $\beta=(\beta_1,\beta_2)\in \cH_2^{m_1+m_2}$ and  $n_1,n_2,m_1,m_2\in \N$. \\ \\
	In order for $\la \cdot , \cdot \ra$ to be a well-defined mapping on $\cE\times \cE$, we must show that no matter what finite linear combination we express $x\in \cE$ and $y\in \cE$ in, then $\la x,y \ra$ stay the same. Let $x_1,x_2$ be two finite linear combination representations of $x\in \cE$ and let $y_1,y_2$ be two finite linear combination representations of $y\in \cE$. Then we note that $x_1-x_2\equiv 0 = 0(x_1-x_2)$ and $y_1-y_2\equiv 0= 0(y_1-y_2)$, hence
	\begin{align*}
	\la x_1 , y_1 \ra - \la x_2, y_2\ra &=  \la x_1,y_1 \ra - \la x_2,y_1 \ra + \la x_2,y_1 \ra - \la x_2, y_2 \ra \\
	&= \la x_1-x_2 , y_1 \ra + \la x_2,y_1-y_2 \ra \\
	&= \la 0(x_1-x_2) , y_1 \ra + \la x_2,0(y_1-y_2)\ra \\
	&= 0\la x_1-x_2 , y_1 \ra + 0\la x_2,y_1-y_2\ra \\
	&=0,
	\end{align*}
 proving that $\la\cdot , \cdot \ra$ is a well-defined mapping on $\cE \times \cE$. Thus we have proved that $\la \cdot , \cdot \ra: \cE\times \cE \to \bK$ is a well-defined mapping of sesquilinear form, and it only remains to be shown that $\la\cdot , \cdot \ra$ possess the positive-definiteness property, i.e. $\la x,x \ra \geq 0$ with $\la x,x \ra=0 \iff x\equiv 0$, in order to confirm that it is a inner product on $\cE$. \\  \\
 Let $f\in\cE$, and assume that a finite linear combination representation is given by
 \begin{align*}
 f= \sum_{i=1}^n \alpha_i \phi_i \otimes \psi_i,
 \end{align*}
 for some $\alpha\in \bK^n$, $\phi\in \cH_1^n$, $\psi\in \cH_2^n$ and $n\in \N$. Consider the following subspaces
 \begin{align*}
 \text{span}(\phi_1,...,\phi_n) &=  \{a_1 \phi_1+ \cdots + a_n \phi_n:a\in \bK^n\}\subset \cH_1,\\
 \text{span}(\psi_1,...,\psi_n) &=  \{a_1 \psi_1+ \cdots + a_n \psi_n:a\in \bK^n\} \subset \cH_2,
 \end{align*}
 and note that there exists two orthonormal bases $\{\gamma_1,...,\gamma_{m_1}\}\subset \cH_1$ and $\{\beta_1,...,\beta_{m_2}\}$ with $m_1,m_2\leq n$ for the $\text{span}(\phi_1,...,\phi_n)$ and $  \text{span}(\psi_1,...,\psi_n)$ respectively (reduce to independent vectors and utilize Gram-Schmidt procedure). Now note that for any $(x,y)\in \cH_1\times \cH_2$, we can represent $\alpha_i\phi_i$ and $\psi_i$ in terms of linear combinations of the basis elements for every $1\leq i \leq n$, yielding
 \begin{align*}
 f(x,y) &= \sum_{i=1}^n \alpha_i \phi_i \otimes \psi_i (x,y) = \sum_{i=1}^n \la x, \alpha_i \phi_i \ra_{\cH_1} \la y, \psi_i \ra_{\cH_2} \\
 &= \sum_{i=1}^n \Big\la x, \sum_{j_1=1}^{m_1}  a_{{i_j}_1} \gamma_{j_1}\Big\ra_{\cH_1} \Big\la y,  \sum_{j_2=1}^{m_2}b_{{i_j}_2} \beta_{j_2} \Big\ra_{\cH_2} \\
 &= \sum_{i=1}^n \sum_{j_1=1}^{m_1}  \sum_{j_2=1}^{m_2} \overline{a_{{i_j}_1} b_{{i_j}_2}}\la x,  \gamma_{j_1} \ra_{\cH_1} \la y,  \beta_{j_2} \ra_{\cH_2} \\
&= \sum_{i=1}^n \sum_{j_1=1}^{m_1}  \sum_{j_2=1}^{m_2} \overline{a_{{i_j}_1} b_{{i_j}_2}} \gamma_{j_1}\otimes\beta_{j_2}   (x,y) =  \sum_{j_1=1}^{m_1}  \sum_{j_2=1}^{m_2} c_{j_1,j_2} \gamma_{j_1}\otimes\beta_{j_2}   (x,y).
 \end{align*}
 where $
 c_{j_1,j_2} = \sum_{i=1}^n \overline{a_{{i_j}_1} b_{{i_j}_2}}
$. Hence
 \begin{align*}
  \la f, f\ra  &= \left\la  \sum_{j_1=1}^{m_1}  \sum_{j_2=1}^{m_2} c_{j_1,j_2} \gamma_{j_1}\otimes\beta_{j_2}, \sum_{j_1=1}^{m_1}  \sum_{j_2=1}^{m_2} c_{j'_1,j'_2} \gamma_{j_1}\otimes\beta_{j_2} \right\ra \\
  &= \sum_{j_1=1}^{m_1}  \sum_{j_2=1}^{m_2}  \sum_{j_1'=1}^{m_1}  \sum_{j_2'=1}^{m_2} c_{j_1,j_2} \overline{c_{j'_1,j'_2}} \la \gamma_{j_1},\gamma_{j'_1} \ra_{\cH_2} \la \beta_{j_2} , \beta_{j'_2} \ra_{\cH_2},
 \end{align*} 
 and note by the orthogonality that each term is zero if $j_1\not = j'_1$ or $j_2\not = j'_2$ and for those terms where $j_1 = j'_1$ and $j_2 = j'_2$ the inner products become $\la \gamma_{j_1},\gamma_{j'_1} \ra_{\cH_1}=\la \beta_{j_2},\beta_{j'_2} \ra_{\cH_2}=1  $. Hence
 \begin{align*}
 \la f , f \ra =  \sum_{j_1=1}^{m_1}  \sum_{j_2=1}^{m_2}    c_{j_1,j_2} \overline{c_{j_1,j_2}} =\sum_{j_1=1}^{m_1}  \sum_{j_2=1}^{m_2}    |c_{j_1,j_2}|^2 .
\end{align*}
 This proves that $\la f, f\ra \geq 0$ and we have that
\begin{align*}
\la f, f \ra =0 \iff \forall (j_1,j_2)\in\{1,...,m_1\}\times \{1,...,m_2\}: c_{j_1,j_2}=0,
 \end{align*}
 but note that $f=\sum_{j_1=1}^{m_1}  \sum_{j_2=1}^{m_2} c_{j_1,j_2} \gamma_{j_1}\otimes\beta_{j_2}$ and since $\{\gamma_1,...,\gamma_{m_1}\}\subset \cH_1$ and $\{\beta_1,...,\beta_{m_2}\}$ are orthonormal bases then $\gamma_{j_1}\otimes \beta_{j_2}\not \equiv 0$, proving that $f\equiv 0$ if and only if every factor $c_{j_1,j_2}$ is zero. We conclude that
\begin{align*}
\la f, f \ra =0 \iff f\equiv 0,
\end{align*}
proving that $\la\cdot , \cdot \ra$ is an inner product on $\cE$.	 
\end{p}
\noindent First some notes about how the completion of $\cE$ and hence also the tensor product is constructed: Let $\| \cdot \|$ be the inner product induces norm on $\cE$, and let $E$ be the $\bK$-vector space of all Cauchy sequences
\begin{align*}
E= \{(x_n): (x_n) \textit{ is Cauchy in } (\cE,\|\cdot\|)\},
\end{align*}
with scalar multiplication and addition defined by $\lambda(x_n)=(\lambda x_n)$ and $(x_n)+(y_n)=(x_n+y_n)$ and zero element $(0)$.
Define the semi-norm $s:E\to \R$ by
\begin{align*}
s((x_n))= \lim_{n\to\i} \|x_n\|,
\end{align*}
and note that the limit always exists and is an element of $[0,\i)$ since $(\|x_n\|)$ is Cauchy in $\R$ which in complete. This is a semi-norm since there may be Cauchy sequences $(x_n)\not = (0)$ with $s((x_n))=0$. But if we define $N=\{(x_n)\in E:s((x_n))=0\}$ then $s$ is a  norm on the quotient $\bK$-vector space $\tilde{\cE}=E/N$ given by the space of equivalence classes $[(x_n)]=(x_n)+N=\{(x_n)+(y_n):(y_n)\in N\}$, which identifies elements $(x_n)\sim (y_n)$ if $(x_n)-(y_n)\in N$. This $\bK$-vector space has scalar multiplication and addition defined by $\lambda[(x_n)]=[\lambda(x_n)]$ and $[(x_n)]+[(y_n)]=[(x_n)+(y_n)]$. That is, we have the following normed $\bK$-vector space
\begin{align*}
(\tilde{\cE},\|\cdot\|_\sim) = (\{[(x_n)]: (x_n) \in E\},\|\cdot\|_\sim),
\end{align*}
where $\|[(x_n)]\|_\sim=\lim_{n\to\i}\|x_n\|$. Note that $\|\cdot\|_\sim$ is well-defined on $\tilde{\cE}$ by this definition, since for any two $[(x_n)]=[(y_n)]$ then $(x_n)-(y_n)=(x_n+y_n)\in N$ and hence $\lim_{n\to\i}| \, \|x_n\| - \|y_n\| \, | \leq \lim_{n\to\i} \|x_n-y_n\| =0$, proving that $\|[(x_n)]\|_\sim=\|[(y_n)]\|_\sim$. \\ \\
We can now define the linear (by the scalar multiplication and addition defined on the space $E$ and $\tilde{\cE}$) operator $\iota: \cE \to \tilde{\cE}$ by
\begin{align*}
\iota(x) = [(x)],
\end{align*}
where $(x)$ is the constant Cauchy sequence $(x,x,x,...)\in E$. By the linearity of $\iota$ we see that
\begin{align*}
\| x-y\| =  \lim_{n\to \i} \| x-y\| = \|[(x-y)]\|_\sim =\| \iota(x-y)\|_\sim=\|\iota(x)-\iota(y)\|_\sim,
\end{align*}
 for any $x,y\in\cE$, proving that $\iota$ is an isometry between $(\cE,\|\cdot\|)$ and $(\tilde{\cE},\|\cdot\|_\sim)$. In other words we have that $(\cE,\|\cdot\|)$ is isometrically embeddable via $\iota$ into its completion $(\tilde{\cE},\|\cdot\|_\sim)$. 
\begin{definition}
	We define the tensor product of $\cH_1$ and $\cH_2$ by $(\cH_1\otimes \cH_2,\|\cdot\|_{\cH_1\otimes \cH_2})=(\tilde{\cE},\|\cdot\|_\sim)$, i.e. the $\bK$-Banach space given by the completion of $\cE$ with respect to the induced norm. Whenever a simple tensor $\phi\otimes \psi$ is presented it will henceforth, without mention,  represent the corresponding element in the completion $\iota(\phi\otimes \psi)\in \cH_1\otimes \cH_2$.
\end{definition}

\noindent By proposition 1.9 \cite{conway1990ACourseInFunctionalAnalysis_reference_for_completion_inner_product_space} there exists an inner product $\la \cdot , \cdot \ra_{\cH_1\otimes \cH_2}$ on $\cH_1\otimes \cH_2$ with induced norm coinciding with $\|\cdot \|_{\cH_1\otimes \cH_2}$ rendering $(\cH_1\otimes \cH_2,\la\cdot, \cdot \ra_{\cH_1\otimes \cH_2})$ a $\bK$-Hilbert space. This inner product furthermore satisfies $\la x,y \ra = \la \iota(x),\iota(y) \ra_{\cH_1\otimes \cH_2}$ for any $x,y\in \cE$. This especially entails that
\begin{align*}
\la \phi_1\otimes \psi_1, \phi_2 \otimes \psi_2 \ra_{\cH_1\otimes \cH_2} = \la \phi_1,\phi_2\ra_{\cH_1}\la \psi_1,\psi_2 \ra_{\cH_2},
\end{align*}
for any $\phi_1,\phi_2\in\cH_1$ and $\psi_1,\psi_2\in\cH_2$. \\ \\
The following proof is only true for separable Hilbert spaces $\cH_1$ and $\cH_2$.
\begin{lemma}\label{appendix_lemma_tensor_map_is_measurable}
	If $\cH_1$ and $\cH_2$ are separable Hilbert spaces, then the map $\cH_1\times \cH_2 \ni (\psi,\phi) \mapsto \phi \otimes \psi \in \cH_1\otimes \cH_2$ is a  $\cB(\cH_1)\otimes \cB(\cH_2)/\cB(\cH_1\otimes \cH_2)$-measurable.
\end{lemma}
\begin{p}
	Fist note that the map in question if the composition $ \iota\circ f: \cH_1\times \cH_2 \to \cH_1\otimes \cH_2$, where $f:\cH_1\times \cH_2 \to \cE$ is given by the simple tensor, that is $f(\phi,\psi) = \phi\otimes \psi \in \cE$ for any $(\phi,\psi)\in\cH_1\times \cH_2$. The isometric embedding into the completion $\iota:(\cE,\|\cdot\|)\to (\cH_1\otimes \cH_2,\la\cdot , \cdot \ra_{\cH_1\otimes \cH_2})$ is continuous and therefore measurable with respect to the Borel $\sigma$-algebra induced by the respective norms, i.e.  $\cB(\cE)/\cB(\cH_1\otimes \cH_2)$-measurable. Thus it suffices to show that $f$ is $\cB(\cH_1)\otimes \cB(\cH_2)/\cB(\cE)$-measurable. Since $\cH_1$ and $\cH_2$ are separable spaces we note that $\cB(\cH_1)\otimes \cB(\cH_2)= \cB(\cH_1\times \cH_2)$ and therefore it suffices to show that $f$ is a continuous mapping. Fix any point $(\phi_0,\psi_0)\in \cH_1\times \cH_2$ and note that
	\begin{align*}
	\| f(\phi,\psi) \| = \|\phi \otimes \psi\| = \la \phi,\phi\ra_{\cH_1}\la \psi,\psi \ra_{\cH_2} = \|\phi\|_{\cH_2}\|\psi\|_{\cH_2}.
	\end{align*}
	Thus for any $\epsilon >0$, let $\delta=\min\{1, \epsilon/(1+\|\phi_0\|_{\cH_1}+\|\psi_0\|_{\cH_2})\}$ and notice that for any $(\phi,\psi)\in B((\phi_0,\psi_0),\delta)=B(\phi_0,\delta)\times B(\psi_0,\delta)$ (see \cref{Appendix_Product_spaces_section} regarding that the maximum metric is a metric generating the product topology on $\cH_1\times \cH_2$) we have that
	\begin{align*}
	\| f(\phi,\psi)-f(\phi_0,\psi_0)\| &= \| \phi\otimes \psi -\phi_0 \otimes \psi_0\| \\
	&= \| (\phi-\phi_0)\otimes (\psi-\psi_0) +\phi\otimes \psi_0 +\phi_0\otimes \psi- 2\phi_0\otimes \psi_0 \| \\
	&= \| (\phi-\phi_0)\otimes (\psi-\psi_0) + (\phi-\phi_0)\otimes \psi_0 + \phi_0\otimes (\psi-\psi_0)\| \\
	&\leq \|\phi-\phi_0\|_{\cH_1}\| \psi-\psi_0  \|_{\cH_2} +\|\phi-\phi_0\|_{\cH_1}\|\psi_0\|_{\cH_2}+ \|\phi_0\|_{\cH_1}\|\psi-\psi_0\|_{\cH_2} \\
	&< \delta^2 +\delta \|\psi_0\|_{\cH_2}+ \|\phi_0\|_{\cH_1}\delta \\
	&\leq \delta (1+ \|\psi_0\|_{\cH_2}+\|\phi_0\|_{\cH_1}) \\
	&\leq \epsilon,
	\end{align*}
	proving continuity of $f$ in every point $(\phi_0,\psi_0)\in \cH_1\times \cH_2$, which concludes the proof.
\end{p}
\begin{theorem} \label{appendix_theorem_orthonormal_basis_tensor_product}
	$\{e_{1,i}\otimes e_{2,j} : i\in I,j\in J\}$ is an orthonormal basis for $\cH_1\otimes \cH_2$ if $\{e_{1,i}:i\in I\}$  and $\{e_{2,j}:j\in J\}$ are orthonormal bases for $\cH_1$ and $\cH_2$ respectively. Furthermore since $\cH_1$ and $\cH_2$ are both separable Hilbert spaces we know that $I$ and $J$ are either finite or countably infinite index sets.
\end{theorem}
\begin{p}
	That any orthonormal basis for a separable Hilbert space is at most countably infinite, follows from standard Hilbert space theory, the arguments can also be found in \cref{remark_eigenvalues_of_S_and_limit_dist}. Assume without loss of generality that both Hilbert spaces $\cH_1$ and $\cH_2$ are infinite dimensional. We know that every Hilbert space has an orthonormal basis, so fix any two arbitrary orthonormal bases $\{e_{1,j}:j\in \N\}$  and $\{e_{2,j}:j\in \N\}$ of $\cH_1$ and $\cH_2$ respectively. Denote $B=\{e_{1,i}\otimes e_{2,j} : i,j\in\N\}$ and note that it suffices to show that $\iota(\cE)\subset \overline{\mathrm{span}(B)}$. To see this note that in the affirmative, then $\cH_1\otimes \cH_2 = \overline{\iota(\cE)} \subset \overline{\mathrm{span}(B)}$ and $\mathrm{span}(B)\subset \cH_1\otimes \cH_2$, where the overline denotes the closure with respect to the topology on $\cH_1\otimes \cH_2$. The fact that $\cH_1\otimes \cH_2 = \overline{\iota(\cE)}$, which follows by showing that $\iota(\cE)$ is dense in $\cH_1\otimes \cH_2$ (omitted, trivial $\epsilon/\delta$-proof) and the fact that the closure of a dense set is equal to the space itself. 
	
	Since $i(\cE)$ is the space of finite linear combinations of elements $\phi\otimes \psi$ for $\phi\in\cH_1$, $\psi\in\cH_2$, it suffices to show that $\psi\otimes \phi \in \overline{\mathrm{span}(B)}$ for any $\phi\in\cH_1$ and $\psi\in\cH_2$, where
	\begin{align*}
	\mathrm{span}(B) = \Big\{\sum_{i=1}^n\sum_{j=1}^m \lambda_{i,j} e_{i}\otimes e_{j} : n,m\in\N,\lambda_{i,j}\in\bK,e_{i}\otimes e_{j}\in B\Big\}.
	\end{align*} 
	Assume without loss of generality that both $\cH_1$ and $\cH_2$ are infinite dimensional Hilbert spaces, rendering the orthonormal bases $\{e_{1,i}:i\in I\}$  and $\{e_{2,j}:j\in J\}$ countably infinite, so we may enumerate them by the natural numbers. Thus fix $\phi\in\cH_1$ and $\psi\in\cH_2$ and note that since $\{e_{1,i}:i\in \N\}$  and $\{e_{2,j}:j\in \N\}$ are orthonormal bases for $\cH_1$ and $\cH_2$ respectively, we get that
	\begin{align*}
	\phi = \sum_{i=1}^\i  \la e_{1,i} , \phi \ra e_{1,j} \quad \quad \text{and} \quad \quad \psi = \sum_{j=1}^\i \la e_{2,j} , \psi \ra e_{2,j},
	\end{align*}
	where we equalities are understood as convergence in the norms on $\cH_1$ and $\cH_2$. Now realize that $\sum_{i=1}^\i \sum_{j=1}^\i \la e_{1,i} , \phi \ra_{\cH_1} \la e_{2,j} , \psi \ra_{\cH_2} e_{1,j}\otimes  e_{2,j}\in \overline{\mathrm{span}(B)}$, since it converges in the topology of $\cH_1\otimes \cH_2$. This is seen by noting that such a series converge if and only if 
	\begin{align*}
	\sum_{i=1}^\i \sum_{j=1}^\i \|\la e_{1,i} , \phi \ra_{\cH_1} \la e_{2,j},\psi \ra_{\cH_2} e_{1,i}\otimes e_{2,j}\|_{\cH_1\otimes \cH_2}^2 ,
	\end{align*}
	converges, but we note that this equals 
	\begin{align*}
	\sum_{i=1}^\i \sum_{j=1}^\i |\la e_{1,i} , \phi \ra_{\cH_1}|^2 \, | \la e_{2,j} , \psi \ra_{\cH_2}|^2 \|e_{1,i}\|_{\cH_1}^2 \|e_{2,j}\|_{\cH_2}^2 \leq \|\phi\|_{\cH_1}^2 \|\psi\|_{\cH_2}^2 <\i,
	\end{align*}
	by the orthonormality of the bases and the inequality $\sum_{i=1}^n | \la \phi,e_{1,i}\ra_{\cH_1}|^2 \leq \|\phi\|_{\cH_1}^2$ known as Bessel's inequality. With $\phi_n =\sum_{i=1}^n \la e_{1,i} , \phi \ra_{\cH_1}  e_{1,j}$ and $\psi_m = \sum_{j=1}^m 
	\la e_{2,j} , \psi \ra_{\cH_2} e_{2,j}$, linearity and the triangle inequality yields that
	\begin{align*}
	&\Big\| \phi\otimes \psi -\sum_{i=1}^n \sum_{j=1}^m \la e_{1,i} , \phi \ra_{\cH_1} \la e_{2,j} , \psi \ra_{\cH_2} e_{1,j}\otimes  e_{2,j} \Big\|_{\cH_1\otimes \cH_2} \\
	=&\Big\| \phi\otimes \psi - \phi_n \otimes  \psi_m  \Big\|_{\cH_1\otimes \cH_2} \\
	 \leq & \| \phi \otimes \psi - \phi \otimes \psi_m \|_{\cH_1\otimes \cH_2} + \| \phi\otimes \psi_m - \phi_n \otimes \psi_m \|_{\cH_1\otimes \cH_2} \\
	 =& \|\phi \otimes (\psi-\psi_m)\|_{\cH_1\otimes \cH_2} + \| (\phi-\phi_n)\otimes \psi_m) \|_{\cH_1\otimes \cH_2}\\
	=& \|\phi\|_{\cH_1}\|\psi-\psi_m\|_{\cH_2} + \| \phi-\phi_m\|_{\cH_1} \|\psi_m\|_{\cH_2},
	\end{align*}
	which converge to zero when $m$ and $n$ tends to infinity, by the basis representation of $\phi$ and $\psi$ above and the fact that $\|\phi_n\|_{\cH_1}\to_n \|\phi\|_{\cH_1}<\i$ (reverse triangle inequality). Since $\phi\otimes \psi$ can be written as a limit of elements in $\mathrm{span}(B)$ it must lie in the closure $\overline{\mathrm{span}(B)}$.		
\end{p} \newpage
\subsection{Characteristic functions of random elements in $\R^n$} \label{Appendix_Characteristic_fucntions}
We will briefly  introduce the theory of characteristic functions of random elements in $n$-dimensional Euclidean spaces. To this end, let $\inner{\cdot}{\cdot}$ denote the regular inner product on $\R^n$, and let \pspace be some probability space. We refer to \cite{folland1999real} for the theory concerning integration of complex valued functions and to \cite{schilling} for measure and integration theory. \\ \\
Let $X:\lp \Omega, \F, P\rp \to \lp \R^n , \cB(\R^n)\rp$ be a measurable mapping (random vector) and let $P_X$ denote the push-forward measure on $\R^n$, i.e. $P_X :=X(P)= P\circ X^{-1}$. 

\begin{definition} \label{defi_char}
	The characteristic function $\phi_X: \R^n \to \bC$ of $X$ (or the law $P_X$ on $\R^n$) is defined as
	\begin{align*}
	\phi_X (t)&= \int_{\R^n} e^{i\inner{t}{x}}dP_X(x)= \int_{\R^n} e^{it^\t x}dP_X(x) \\
	&=\int_{\R^n} \cos \lp t^\t x\rp dP_X(x)+i\int_{\R^n} \sin \lp t^\t x\rp dP_X(x),
	\end{align*}
	for any $t\in \R^n$. For a Borel probability measures $\mu$ on $\R^n$ we denote the corresponding characteristic function by $\hat{\mu}:\R^n \to \bC$
\end{definition}
We note that the characteristic function is always well-defined by realizing that $|e^{iy}|=1$. \\ \\
Now lets prove some properties of characteristic functions of random vectors. Note these properties equally apply to probability measures on $\R^n$, in the sense that we can think of the identity mapping $Id:\R^n \to \R^n$ as the random vector with sample space $\Omega=\R^n$. .
\begin{theorem}[Properties of the characteristic function]\label{appendix_theorem_properties_char}
	Let $X$ and $Y$ be random elements in $\R^n$ and $\R^m$ respectively, then the following holds
	\begin{enumerate}[label={(\arabic*)}., ref=(\arabic*)]
		\item\label{char_properties_simple} $\phi_X(0)=1$, $|\phi_X(t)|\leq 1$ and $\bar{\phi}_X(t)=\phi_X(-t)$  for all $t\in \R^n$.
		\item\label{char_prooerties_cont} The mapping $t \mapsto \phi_X(t)$ is uniformly continuous.
		\item\label{char_prooerties_unique} If $n=m$: $\phi_X(t)=\phi_Y(t)$ for all $t\in\R^n$ if and only if $X\eqd Y$.
		\item\label{char_prooerties_indep} $\phi_{X,Y}(t,s)=\phi_X(t) \phi_Y(s)$ for all $(t,s)\in \R^n \times \R^m$ if and only if $X\independent Y$.
		\item\label{char_prooerties_sum} If $n=m$ then $\phi_{X+Y}(t)=\phi_{X,Y}(t,t)$ for all $t\in\R^n$.
		\item\label{char_prooerties_relax} The conditions of \cref{char_prooerties_unique} and \cref{char_prooerties_indep} are equivalent to requiring that they hold almost everywhere with respect to the Lebesgue measure on $\R^n$ and $\R^{n+m}$ respectively.
		\item\label{char_prooerties_transformation} For any $p\geq1$, let $\alpha\in \R^p$ and $B:\R^n \to \R^p$ be a linear mapping then $\alpha+BX$  has characteristic function given by $\phi_{\alpha+BX}(t)=e^{i\inner{t}{\alpha}}\phi_X(B^Tt)$ for all $t\in \R^p$.
		\item\label{char_prooerties_sym} $X$ has a symmetric distribution around $a\in\R^n$ if and only if $\Im \phi_{X-a}=0$.
	\end{enumerate}
\end{theorem}
\begin{p}
	(1): It is trivial that $\phi_X(0)=1$ since $\exp(0)=1$. Furthermore let $t\in\R^n$ and note that by the triangle inequality for complex valued integrals we have that
	\begin{align*}
	|\phi_X(t)| \leq \int_{\R^n} \lv e^{i\inner{t}{x}}\rv dP_X(x) =P_X(\R^n)=1,
	\end{align*}
	and as regards the last claim simply note that
	\begin{align*}
	\bar{\phi}_X(t) &= \int_{\R^n} \cos \inner{t}{x}dP_X(x)-i\int_{\R^n} \sin \inner{t}{x}dP_X(x) =\phi_X(-t),
	\end{align*}
	since cosine is even, sinus is odd and $t\mapsto \langle t,x \rangle$ is linear. \\ \\
	(2): The proof proceeds analogously to the one-dimensional case: By the continuity of the map $h\mapsto e^{i\inner{h}{x}}$ for any $x\in\R^n$, we get that for any sequence $(h_k)_{k\in\N}\subset \R^n$ with $\lim_{h\to\i} h_k = 0$
	\begin{align*}
	\lim_{k\to \i } \int_{\R^n} \lv e^{i \inner{h_k}{x}}-1 \rv dP_X(x) =0,
	\end{align*}
	as $k$ tends to infinity by Lebesgue's dominated convergence theorem. Hence for any $\ep>0$ there exists a $\delta>0$ such that 
	\begin{align*}
	\lv \phi_X(t)-\phi_X(s)\rv &= \lv \int_{\R^n}e^{i \inner{t}{x}}-e^{i \inner{s}{x}} dP_X(x)\rv  \\
	& \leq  \int_{\R^n} \lv e^{i \inner{t}{x}}-e^{i \inner{s}{x}} \rv dP_X(x) \\
	& =  \int_{\R^n} \lv e^{i \inner{s}{x}} \rv \lv e^{i \inner{t-s}{x}}-1 \rv dP_X(x) \\
	& =  \int_{\R^n} \lv e^{i \inner{t-s}{x}}-1 \rv dP_X(x) < \ep,
	\end{align*}
	for all $t,s\in\R^n$ with $|t-s|<\delta$, proving uniform continuity of the characteristic function. \\ \\
	(3): see \cite{dudley2002real} theorem 9.5.1. \\ \\
	(4): \textit{If}  follows rather trivially by noting that $X \independent Y \iff (X,Y)(P)=X(P)\otimes Y(P)$ on  $(\R^{n+m},\cB(\R^{n+m}))=(\R^n \times \R^m , \cB(\R^n)\otimes \cB(\R^m))$ allowing us to utilize Fubini's theorem to conclude that
	\begin{align*}
	\phi_{X,Y}(t,s) &= \int_{\R^n} \int_{\R^m} e^{i\inner{(t,s)}{(x,y)}} dP_Y(y) dP_X(x) \\
	&= \int_{\R^n} \int_{\R^m} e^{i\inner{t}{x}} e^{i\inner{s}{y}} dP_Y(y) dP_X(x) \\
	&= \phi_X(t)\phi_Y(s),
	\end{align*}
	for any $(t,s)\in\R^n\times \R^m$. \textit{Only if} follows after realizing that by the above we have that the distribution $X(P)\otimes Y(P)$  on $(\R^n \times \R^m , \cB(\R^n)\otimes \cB(\R^m))$ has characteristic function $\phi_X(t)\phi_Y(s)$ for $(t,s)\in\R^n \times \R^m$ coinciding with the characteristic function of $(X,Y)$. Hence by (2) we have that $(X,Y)(P)=X(P)\otimes Y(P)$ proving independence. \\ \\ 
	(5): Simply note that for any $t\in\R^n$ \begin{align*}
	\phi_{X+Y}(t)=E\lp e^{it^\t(X+Y)}\rp=E\lp e^{i (t^\t X+t^\t Y)}\rp=E\lp e^{i(t,t)^\t(X,Y)}\rp=\phi_{X,Y}(t,t).
	\end{align*}
	(6): The result follows from an easy application of \cref{char_properties_simple}. Assume for contradiction that $\phi_X=\phi_Y$ $\lambda^n$-almost everywhere and that there exists a $t\in \R^n$ such that $\phi_X(t)\not = \phi_Y(t)$. Let $\ep:=|\phi_X(t)-\phi_Y(t)|>0$ and note that by the continuity of $s\mapsto |\phi_X(s)-\phi_Y(s)|$ there exists a $\delta >0$ such that
	$
	|\phi_X(s)-\phi_Y(s)| > \ep/2 \text{ whenever } s\in B(t,\delta)
	$, where $B(t,\delta)$ denotes the open $\delta$-ball of $t$. However we now have that $\lambda^n(\{s\in \R^n : \phi_Y(s)\not = \phi_Y(s)\})\geq \lambda^n(B(t,\delta))=(\sqrt{\pi}\delta )^n/\Gamma(n/2+1)>0$ \Lightning. The equivalence for \cref{char_prooerties_indep} follows exactly similarly. \\ \\
	(7): Let $p\geq 1$ and simply note that $
	\phi_{\alpha+BX}(t) = E\lp e^{i \inner{t}{\alpha+BX}}\rp = e^{i \inner{t}{\alpha}}E\lp e^{i t^TBX}\rp =  e^{i \inner{t}{\alpha}}\phi_X(B^Tt)
	$
	for any $t\in \R^p$. \\ \\
	(8): By definition $X$ is symmetric around $a\in\R^n$ if $X-a \eqd a-X$, and if this is the case then \cref{char_prooerties_unique} and \cref{char_prooerties_transformation} yields that
	\begin{align*}
	\phi_{X-a}(t) &= \phi_{(-I_n)(-I_n)(X-a)}(t) = \phi_{(-I_n)(X-a)}((-I_n)^\t t) \\
	&=\phi_{a-X}(-t) = \phi_{X-a}(-t) = \bar{\phi}_{X-a}(t),		\end{align*}
	for any $t\in\R^n$, where $I_n:\R^n \to \R^n$ is the identity mapping, proving that $\Im \phi_{X-a}=0$. For the converse realize that if $\Im \phi_{X-a}=0$ then by the same arguments as above we get that \begin{align*}
	\phi_{X-a}(t)= \bar{\phi}_{X-a}(t)= \phi_{X-a}(-t) = \phi_{a-X}(t),
	\end{align*}
	for any $t\in\R^n$, which by \cref{char_prooerties_unique} proves that $X-a\eqd a-X$.
	
\end{p}
\newpage
\subsection{Miscellaneous}
\begin{lemma} \label{lemma_help_injectivity_proof}
	If $\theta\in M(\cX\times \cY)$ and $f\in \cL^1_\R(\pi_1(|\theta|))$ then the set function $\theta_f:\cB(\cY)\to \R$ given by \begin{align*}
	\theta_f(B) = \int f(x)1_B(y) \, d\theta(x,y),
	\end{align*}
	is a finite signed measure. If furthermore $ (x,y) \mapsto f(x)g(y)\in\cL^1(\theta)$ and $g$ is measurable then $g\in \cL^{1}(\theta_f)$ and $$\int g(y) \,d\theta_f(y) =\int f(x)g(y) \, \theta(x,y).$$
\end{lemma}
\begin{p}
	We obviously have that $\theta_f(\emptyset)=0$ and for any disjoint sequence of sets $(A_i)\subset \cB(\cY)$ we have that\begin{align*}
	\theta_f \lp  \bigcup_{i=1}^\i A_i\rp &= \int f(x) \sum_{i=1}^\i 1_{A_i}(y) \, d\theta(x,y) \\
	&= \lim_{n\to\i} \sum_{i=1}^n\int f(x)1_{A_i}(y) \, d\theta(x,y) \\
	&= \sum_{i=1}^\i \theta_f(A_i),
	\end{align*}
	by the dominated convergence theorem for signed measures, since $|f(x)\sum_{i=1}^n1_{A_i}(y)|\leq |f(x)|\in \cL^1(\theta)$. Hence $\theta_f$ is a signed measure, but it furthermore holds that $|\theta_f(B)|\leq \int |f(x)|\, d|\theta|(x,y) <\i$ for any $B\in\cB(\cY)$, proving that $\theta_f:\cB(\cY)\to \R$ is a finite signed measure. Now note that
	\begin{align*}
	\theta_f(B) =& \int f^+(x)1_B(y)\, d\theta^+(x,y) +\int f^-(x)1_B(y)\, d\theta^-(x,y) \\
	&- \int f^+(x)1_B(y)\, d\theta^-(x,y) - \int f^-(x)1_B(y)\, d\theta^+(x,y),
	\end{align*}
	and as a consequence 
	\begin{align*}
	|\theta_f|(B) \leq& \int f^+(x)1_B(y)\, d\theta^+(x,y) +\int f^-(x)1_B(y)\, d\theta^-(x,y) \\
	&+ \int f^+(x)1_B(y)\, d\theta^-(x,y) + \int f^-(x)1_B(y)\, d\theta^+(x,y) \\
	=& \int |f(x)|1_B(y) \, d|\theta|(x,y) \\
	=& (|\theta|_{|f|})(B),
	\end{align*}
	by the same reasoning as in \cref{lemma_M(X)_and_M^1(X)_are_vector_spaces}. So $g\in\cL^1(\theta_f)$ if $g\in \cL^1(|\theta|_{|f|})$, and as we shall see this is indeed the case if $(x,y)\mapsto f(x)g(y)\in \cL^1(\theta)$. If $g$ is measurable there exists a sequence of positive simple mappings $(|g|_n)$ such that $|g|_n\to_n |g|$ point-wise with $|g|_n\leq |g|_{n+1}$ for all $n\in\N$. By the monotone convergence theorem we  have that
	\begin{align*}
	\int |f(x)g(y)|\, d|\theta|(x,y) &= \lim_{n\to\i} \int |f(x)||g|_n(y) \, d|\theta|(x,y) \\
	&=\lim_{n\to\i} \int |f(x)|\lp \sum_{i=1}^{k(n)}c_{n,i}1_{A_{n,i}}(y)\rp \, d|\theta|(x,y) \\ 
	&=\lim_{n\to\i} \sum_{i=1}^{k(n)} c_{n,i} \int |f(x)|1_{A_{n,i}}(y) \, d|\theta|(x,y) \\
	&= \lim_{n\to\i} \sum_{i=1}^{k(n)} c_{n,i} |\theta|_{|f|}(A_{n,i}) \\
	&=\lim_{n\to\i} \int |g|_n(y) \, d|\theta|_{|f|}(y) \\
	&= \int |g(y)| \, d|\theta|_{|f|}(y) \\
	&\geq \int |g(y)| \, d|\theta_f|(y),
	\end{align*}
	proving that $g\in\cL^1(\theta_f)$ if $(x,y)\mapsto f(x)g(y)\in \cL^1(\theta)$.  Hence if $g$ is measurable and $(x,y)\mapsto f(x)g(y)\in \cL^1(\theta)$ then there exists a sequence of simple mappings $(g_n)$ such that $g_n\to_n g$ point-wise with $|g_n|\leq |g|$ for all $n\in\N$. By the dominated convergence theorem for signed measures we have that
	\begin{align*}
	\int g(y) \, d\theta_f(y) &= \lim_{n\to \i} \int g_n(y) \, d\theta_f(y)\\
	&= \lim_{n\to\i} \int  f(x) g_n(y) \, d\theta(x,y) \\
	&=  \int  f(x) g(y) \, d\theta(x,y),
	\end{align*}
	where we used that $|f(x)g_n(y)|\leq |f(x)g(y)|\in\cL^1(\theta)$.
\end{p}
\begin{lemma}\label{lemma_inequality_f}
	It holds that
		\begin{align*}
		\frac{|f_i(z_1,z_2,z_3,z_4)|}{2}\leq \left\{
		\begin{array}{lll}
		d_i(z_1,z_4), & d_i(z_2,z_3), & d_i(z_1,z_2) \lor d_i(z_1,z_3), \\
		d_i(z_1,z_2) \lor d_i(z_1,z_4), & d_i(z_1,z_2) \lor d_i(z_2,z_3), & d_i(z_1,z_2) \lor d_i(z_2,z_4), \\
		d_i(z_1,z_4)  \lor d_i(z_1,z_3), & d_i(z_1,z_4)  \lor d_i(z_2,z_3), & d_i(z_1,z_4)  \lor d_i(z_2,z_4), \\
		d_i(z_2,z_3) \lor d_i(z_1,z_3), & d_i(z_2,z_3) \lor d_i(z_2,z_4),  & d_i(z_3,z_4) \lor d_i(z_1,z_3), \\ 
		d_i(z_3,z_4) \lor d_i(z_1,z_4), & d_i(z_3,z_4) \lor d_i(z_2,z_3), &d_i(z_3,z_4) \lor d_i(z_2,z_4). \\ 	\end{array}
		\right.
		\end{align*}
\end{lemma}
\begin{p}
	
		Initially recall that for any metric
		\begin{equation} \label{temp_inequality}
		d(a,b)-d(a,c)\leq d(b,c),
		\end{equation}
		by the triangle inequality. \\ \\
	Fix $z_1,z_2,z_3,z_4$ in $\cX$ or $\cY$ and note that if $f=f(z_1,z_2,z_3,z_4)>0$ then \begin{align*}
	|f|=d(z_1,z_2)-d(z_1,z_3)-d(z_2,z_4)+d(z_3,z_4).
	\end{align*}
	Now we create upper bounds for $|f|$ by using the triangle inequality on the positive terms
	
	\begin{itemize}
		\item Expanding first term with
		\begin{itemize}
			\item $z_4$ as an intermediate point
			\begin{align*}
			|f| & \leq d(z_1,z_4)+d(z_2,z_4)-d(z_1,z_3)-d(z_2,z_4)+d(z_3,z_4) \\
			&= d(z_1,z_4)-d(z_1,z_3)+d(z_3,z_4) \\
			&\leq \left\{ \begin{array}{c}
			2d(z_3,z_4) \quad \quad \textit{\cref{temp_inequality} on term 1 and 2}\\
			2d(z_1,z_4) \quad \quad \textit{\cref{temp_inequality} on term 2 and 3}
			\end{array}\right.
			\end{align*}
			\item $z_3$ as an intermediate point
			\begin{align*}
			|f| & \leq d(z_1,z_3)+d(z_2,z_3)-d(z_1,z_3)-d(z_2,z_4)+d(z_3,z_4) \\
			&= d(z_2,z_3)-d(z_2,z_4)+d(z_3,z_4) \\
			&\leq \left\{ \begin{array}{c}
			2d(z_3,z_4) \quad \quad \textit{\cref{temp_inequality} on term 1 and 2}\\
			2d(z_2,z_3) \quad \quad \textit{\cref{temp_inequality} on term 2 and 3}
			\end{array}\right.
			\end{align*}
		\end{itemize}
		\item Expanding fourth term with 
				\begin{itemize}
					\item $z_1$ as an intermediate point
					\begin{align*}
					|f| & \leq d(z_1,z_2)-d(z_1,z_3)-d(z_2,z_4)+d(z_1,z_3)+d(z_1,z_4) \\
					&= d(z_1,z_2)-d(z_2,z_4)+d(z_1,z_4)  \\
					&\leq \left\{ \begin{array}{c}
					2d(z_1,z_4) \quad \quad \textit{\cref{temp_inequality} on term 1 and 2}\\
					2d(z_1,z_2) \quad \quad \textit{\cref{temp_inequality} on term 2 and 3}
					\end{array}\right.
					\end{align*}
					\item $z_2$ as an intermediate point
					\begin{align*}
					|f| & \leq d(z_1,z_2)-d(z_1,z_3)-d(z_2,z_4)+d(z_2,z_3)+d(z_2,z_4) \\
					&= d(z_1,z_2)-d(z_1,z_3)+d(z_2,z_3) \\
					&\leq \left\{ \begin{array}{c}
					2d(z_2,z_3) \quad \quad \textit{\cref{temp_inequality} on term 1 and 2}\\
					2d(z_1,z_2) \quad \quad \textit{\cref{temp_inequality} on term 2 and 3}
					\end{array}\right.
					\end{align*}
				\end{itemize}
	\end{itemize}
	Thus if  $f>0$ then
	\begin{align*}
	|f| \leq \left\{  \begin{array}{c}
	2d(z_1,z_2) \\
	2d(z_1,z_4) \\
	2d(z_2,z_3) \\
	2d(z_3,z_4)
	\end{array}\right.
	\end{align*}
	If $f<0$ then by using the above inequalities we get
	\begin{align*}
	|f(z_1,z_2,z_3,z_4)| &= d(z_1,z_3)-d(z_1,z_2)-d(z_3,z_4)+d(z_2,z_4) \\
	&= f(x_1,x_3,x_2,x_4) \leq \left\{  \begin{array}{c}
	2d(z_1,z_3) \\
	2d(z_1,z_4) \\
	2d(z_2,z_3) \\
	2d(z_2,z_4)
	\end{array}\right.
	\end{align*}
	proving that in general
		\begin{align*}
		\frac{|f_i(z_1,z_2,z_3,z_4)|}{2}\leq \left\{
		\begin{array}{lll}
		d_i(z_1,z_4), & d_i(z_2,z_3), & d_i(z_1,z_2) \lor d_i(z_1,z_3), \\
		d_i(z_1,z_2) \lor d_i(z_1,z_4), & d_i(z_1,z_2) \lor d_i(z_2,z_3), & d_i(z_1,z_2) \lor d_i(z_2,z_4), \\
		d_i(z_1,z_4)  \lor d_i(z_1,z_3), & d_i(z_1,z_4)  \lor d_i(z_2,z_3), & d_i(z_1,z_4)  \lor d_i(z_2,z_4), \\
		d_i(z_2,z_3) \lor d_i(z_1,z_3), & d_i(z_2,z_3) \lor d_i(z_2,z_4),  & d_i(z_3,z_4) \lor d_i(z_1,z_3), \\ 
		d_i(z_3,z_4) \lor d_i(z_1,z_4), & d_i(z_3,z_4) \lor d_i(z_2,z_3), &d_i(z_3,z_4) \lor d_i(z_2,z_4). \\ 	\end{array}
		\right. ,
		\end{align*}
	for $i=\cX$ or $i=\cY$.	
\end{p}

\bibliographystyle{plain}
\bibliography{bibfile}
\end{document}